\documentclass{article}
\usepackage[a4paper, total={6in, 8in}]{geometry}
\usepackage[utf8]{inputenc}
\usepackage{amsmath, amsfonts,amsthm}
\usepackage{float}
\usepackage{cite}
\usepackage{graphicx}
\usepackage{textcomp}
\usepackage{xcolor}
\usepackage{tikz}
\usepackage{pgfplots}
\usepackage[export]{adjustbox}
\usetikzlibrary{spy}
\usepackage{graphicx}
\usepackage{tabularx}
\usepackage{xcolor}
\usepackage{amssymb,amsmath,amsfonts, amsthm}
\usepackage{tikz}
\usepackage{multirow}
\usepackage{booktabs}
\usepackage{blindtext}
\usepackage{multicol}
\usepackage{tikzscale}
\usepackage{numprint}
\npdecimalsign{.}
\nprounddigits{3} 
\usepackage{caption}
\usepackage{subcaption}
\usepackage{parskip}
\usepackage{enumitem}
\usepackage{float}
\usepackage{pict2e}
\usepackage{siunitx}
\usepackage[english]{babel}
\usepackage{ulem}
\usepackage{authblk}
\newcommand{\C}{\mathbb C}

\newcommand{\R}{\mathbb R}

\newcommand{\CC}{\mathbf{c}}                       
\newcommand{\Sd}{\mathbf{S}}

\newcommand{\XX}{\mathbf{x}} 
\newcommand{\PP}{\mathbf{p}}
\newcommand{\QQ}{\mathbf{q}}
                      
\newcommand{\ZZ}{\mathbf{z}}

\newcommand{\YY}{\mathbf{y}} 
\newcommand{\UU}{\mathbf{u}} 
    
\newcommand{\VV}{\mathbf{v}}

\newcommand{\Ad}{\mathbf A}                        
\newcommand{\Wd}{\mathbf W} 
\newcommand{\Fd}{\mathbf F} 
                       
\newcommand{\Id}{\mathbf I}                        
\newcommand{\Cd}{\mathbf C}                        
\newcommand{\Ed}{\mathbf E}                        

\newcommand{\Td}{\mathbf T}

\newcommand{\herm}{{\scriptstyle \boldsymbol{\mathsf{H}}}}
\newcommand{\trans}{{\scriptstyle \boldsymbol{\mathsf{T}}}}

\newcommand{\LLambda}{\boldsymbol{\Lambda}} 

\usepackage{mathtools}

\usepackage{mathtools}
\usepackage[percent]{overpic}
\usepackage{algorithm,algcompatible,amsmath}

\DeclareMathOperator{\diag}{diag}

\algnewcommand\INPUT{\item[\textbf{Input:}]}%
\algnewcommand\PARAMETER{\item[\textbf{Parameters:}]}%
\algnewcommand\OUTPUT{\item[\textbf{Output:}]}%
\usepackage[
 colorlinks=false, backref=page
]{hyperref}
\renewcommand*{\backref}[1]{}
\renewcommand*{\backrefalt}[4]{%
    \ifcase #1 (Not cited).%
    \or        (Cited on page~#2).%
    \else      (Cited on pages~#2).%
    \fi}

\colorlet{lred}{red!80}
\colorlet{cmix}{blue!80!red} 
\colorlet{lgreen}{green!80}
\colorlet{lblue}{blue!80}

\newtheorem{theorem}{Theorem}

\newtheorem{assumption}[theorem]{Assumption}

\newtheorem{remark}[theorem]{Remark}
\newtheorem{proposition}[theorem]{Proposition}
\newtheorem{corollary}[theorem]{Corollary}
\newcommand{\stkout}[1]{\ifmmode\text{\sout{\ensuremath{#1}}}\else\sout{#1}\fi}
\usepackage{algorithmicx}
\usepackage{comment}

\usepackage{array}

\newcolumntype{L}[1]{>{\raggedright\let\newline\\\arraybackslash\hspace{0pt}}m{#1}}
\newcolumntype{C}[1]{>{\centering\let\newline\\\arraybackslash\hspace{0pt}}m{#1}}
\newcolumntype{R}[1]{>{\raggedleft\let\newline\\\arraybackslash\hspace{0pt}}m{#1}}

\title{Learning Regularization Parameter-Maps for 
Variational
Image Reconstruction using Deep Neural Networks and Algorithm Unrolling}

\author[1]{Andreas Kofler}
\author[2,3]{Fabian Altekrüger}
\author[3]{Fatima Antarou Ba}
\author[1]{Christoph Kolbitsch}
\author[4,5]{Evangelos Papoutsellis}
\author[1]{David Schote}
\author[6]{Clemens Sirotenko}
\author[1]{Felix Frederik Zimmermann}
\author[7]{Kostas Papafitsoros}
\affil[1]{Physikalisch-Technische Bundesanstalt (PTB), Braunschweig and Berlin, Germany~~~~~~~~~~}
\affil[2]{Humboldt-Universit\"at zu Berlin, Department of Mathematics, Berlin, Germany~~~~~~~~~~~~~~~}
\affil[3]{Technische Universit\"at Berlin, Institute of Mathematics, Berlin, Germany~~~~~~~~~~~~~~~~~~~~~}
\affil[4]{Finden Ltd, Rutherford Appleton Laboratory,  Harwell Campus, Didcot, United Kingdom~~~~}
\affil[5]{Science and Technology Facilities Council, Harwell Campus, Didcot, United Kingdom~~~~~~~}
\affil[6]{Weierstrass Institute for Applied Analysis and Stochastics, Berlin, Germany~~~~~~~~~~~~~~~~~~}
\affil[7]{School of Mathematical Sciences, Queen Mary University of London, United Kingdom~~~~~\vspace{1em}}
\affil[ ]{\textit {\{andreas.kofler, felix.zimmermann, david.schote, christoph.kolbitsch\}@ptb.de}}
\affil[ ]{\textit {fabian.altekrueger@hu-berlin.de, fatimaba@math.tu-berlin.de, epapoutsellis@gmail.com,   }}
\affil[ ]{\textit {sirotenko@wias-berlin.de, k.papafitsoros@qmul.ac.uk}}
\begin{document}

\maketitle

\begin{abstract}
We introduce a method for fast estimation of data-adapted, spatio-temporally dependent regularization parameter-maps for variational image reconstruction, focusing on total variation (TV)-minimization. Our approach is inspired by recent developments in algorithm unrolling using deep neural networks (NNs), and relies on two distinct sub-networks. The first sub-network estimates the regularization parameter-map from the  input data. The second sub-network unrolls $T$ iterations of an iterative algorithm which approximately solves the corresponding TV-minimization problem incorporating the previously estimated regularization parameter-map. The overall network is trained end-to-end in a supervised learning fashion using pairs of clean-corrupted data but crucially without the need of having access to labels for the optimal regularization parameter-maps. We prove consistency of the unrolled scheme by showing that the unrolled energy functional used for the supervised learning $\Gamma$-converges as $T$ tends to infinity, to the corresponding functional that incorporates the exact solution map of the TV-minimization problem. We apply and evaluate our method on a variety of large scale and dynamic imaging problems in which the automatic computation of such parameters has been so far challenging:  2D dynamic cardiac MRI reconstruction, quantitative brain MRI reconstruction,  low-dose CT and dynamic image denoising. The proposed method consistently improves the TV-reconstructions using scalar parameters and the obtained parameter-maps adapt well to each imaging problem and data by leading to the preservation of detailed features. Although the choice of the regularization parameter-maps is data-driven and based on NNs, the proposed algorithm is entirely interpretable since it inherits the properties of the respective iterative reconstruction method from which the network is implicitly defined.
\end{abstract}

\section{Introduction}

Inverse imaging problems can often be described as 
\begin{equation}\label{eq:forward_problem}
    \ZZ = \Ad \XX_{\mathrm{true}} + \mathbf{e}
\end{equation}
where $\XX_{\mathrm{true}}\in V^{n}$ with $V \in \{\R, \C\}$ is the object to be imaged, $\Ad: V^{n}\to V^{m}$ is  a linear operator which models the data-acquisition process, $\mathbf{e}\in V^{m}$ denotes some random noise component and $\ZZ\in V^{m}$ represents the measured data. The goal is to reconstruct $\XX_{\mathrm{true}}$ or at least a good enough approximation of it given the data $\ZZ$. In practice, problem \eqref{eq:forward_problem} is however ill-posed for various reasons. For example, in Magnetic Resonance Imaging (MRI) which is known to suffer from long acquisition times, the measurement process is often accelerated by undersampling in the raw-data domain, the so-called $k$-space, leading to an underdetermined systems. In low-dose CT, where one reduces the radiation exposure of the patient by reducing the energy of the photons emitted from the X-ray source, the measured data is noisy. Further, different inherent properties of the operator $\Ad$ also often determine how well-posed the problem is.
Therefore, the reconstruction procedure requires the use of regularization methods to be able to obtain high quality images and particularly in medical imaging, images with diagnostic accuracy. A prominent approach is to formulate the reconstruction as a minimization problem
\begin{equation}\label{eq:reg_problem}
    \underset{\XX }{\min}  \, d(\Ad \XX, \ZZ) + \mathcal{R}(\XX),
\end{equation}
where $d(\,\cdot \,, \, \cdot \,)$ denotes a data-discrepancy measure and $\mathcal{R}(\,\cdot\,)$ a regularization term. Typical choices for $\mathcal{R}$ vary from $\mathcal{R}(\,\cdot\,) = \| \,\cdot\,\|_2^2$ for the well-known Tikhonov regularization \cite{Tikhonov} or  $\mathcal{R}(\,\cdot\,) = \| \Td \,\cdot\,\|_1$ for methods enforcing sparsity in some basis \cite{Donoho_1994, Chang_2000}.
One of the most widely applied methods is the so-called Total Variation (TV) regularization \cite{ROF92, ChambolleLions}. 
Remaining in the finite dimensional  setting, and choosing the square of the $\ell_{2}$ norm as  data discrepancy (appropriate for Gaussian noise), the reconstruction problem is formulated as 
\begin{equation}\label{eq:tv_min_problem}
    \underset{\XX }{\min}\; \frac{1}{2} \| \Ad \XX - \ZZ\|_2^2 + \lambda \| \nabla \XX \|_1.
\end{equation}
Here $\nabla$ denotes a finite-differences operator and $\lambda>0$ is a scalar regularization parameter that balances the effect of the two terms. This means that the regularization imposed on the sought image is given by sparsity in the gradient domain of the image measured with respect to the $\ell_1$-norm. One reason for the great success of this method lies in its simple intuition, interpretability as well as its interesting mathematical properties. As a result, in the last decades, it has driven both theoretical as well as applied research fields such as biomedical engineering, inverse problems, optimization and geometric measure theory among others \cite{chambollepock2016, Burger2013, scherzer2009variational}, with the  complete list of publications in which the approach is investigated for different reconstruction problems in different imaging modalities being quite extensive. In addition, there exist nowadays numerous algorithms with proven convergence guarantees \cite{chambolle2011first, wang2008new, stadler, hint_kun, chambollepock2016, Tai_augmented} as well as extensions to overcome inherent limitations of the structural properties of the solutions of the problem \eqref{eq:tv_min_problem}, e.g.\  the total generalized variation (TGV) \cite{TGV}, for solving for the well-known TV staircasing artefacts (blocky-like, piecewise constant structures).\\
A crucial aspect which impacts the quality and the usefulness of the images which can be reconstructed by solving problem \eqref{eq:tv_min_problem} is the careful choice of the parameter $\lambda$. Underestimating $\lambda$ yields poor regularization, while overestimating it yields too smooth images with artificial ``cartoon-like" appearance. Particularly in medical imaging applications, where images are at the basis of diagnostic decisions and therapy planning, a proper choice of any regularization parameter is crucial. There exist quite a few methods regarding the automatic choice of a scalar parameter $\lambda$, placed either in the TV term or in the data-discrepancy term, based on the discrepancy principle, $L$-curve methods and others, e.g. \cite{Chuan_discrepancy, Golub1979, CALVETTI2000423, Montagner2014}.\\
However, employing  one single scalar parameter $\lambda$ which globally dictates the strength of the regularization over the entire image might seem sub-optimal for various obvious reasons. Depending on the application, it might be desirable to maintain locally higher data-fidelity instead of enforcing visually appealing but rather wrong image features. In that case, one can replace the parameter $\lambda\in \R_{+}$ in \eqref{eq:tv_min_problem} with a spatially varying and pixel/voxel dependent one,  denoted now by $\boldsymbol{\Lambda}\in \R_{+}^{qn}$,  
with $q$ denoting the number of directions for which the partial derivatives are computed. Implementation-wise that translates to a stack of diagonal operators which contain a regularization parameter for each single pixel/voxel in the respective gradient domain of the image, resulting in a problem of the form 
\begin{equation}\label{eq:tv_min_problem_spatial}
\underset{\XX }{\min}\; \frac{1}{2} \| \Ad \XX - \ZZ\|_2^2 +\|\boldsymbol{\Lambda} \nabla \XX\|_{1}.
\end{equation}
An automatic choice for such spatially varying regularization parameter is rather challenging, as the number of its components drastically increases. Towards that task, bilevel optimization techniques have been employed during the last years, which have the following general formulation:
\begin{equation}\label{general_bilevel}
\left \{
\begin{aligned}
&\min_{\boldsymbol{\Lambda}}\; \sum_{i=1}^{M}l(\XX^{i}(\boldsymbol{\Lambda}),\XX_{\mathrm{true}}^{i})\\
&\text{subject to }\;\;\XX^{i}(\boldsymbol{\Lambda})=\underset{\XX}{\operatorname{argmin}}\;\frac{1}{2} \| \Ad \XX - \ZZ_{i}\|_2^2 +  \| \boldsymbol{\Lambda}\nabla \XX \|_1, \quad i=1,\ldots,M.
\end{aligned}\right.
\end{equation}
Here, $(\ZZ_{i}, \XX_{\mathrm{true}}^{i})_{i=1}^{M}$ are $M$ pairs of measured data  and corresponding ground truth, and $l$ is a suitable upper level objective. For instance, in the case where $l(x_{1},x_{2})=l_{\mathrm{PSNR}}(x_{1},x_{2}):=\|x_{1}-x_{2}\|_{2}^{2}$, the bilevel problem \eqref{general_bilevel} aims to compute the parameters $\boldsymbol{\Lambda}$ which are ``on the average the best ones" (i.e.\ PSNR-maximizing), for the given $M$ data pairs. The idea is that, given some new data $\ZZ_{\mathrm{test}}$ 
that has been measured in a similar way as $(\ZZ_{i})_{i=1}^{M}$, solving \eqref{eq:tv_min_problem_spatial} (in the ``online phase") with the offline-computed $\boldsymbol{\Lambda}$ will yield a good reconstruction.
 This scheme has been extensively studied both for scalar and spatially varying regularization parameters. However, in practice it has mainly been applied  for image denoising (i.e.\ $\Ad=\mathbf{I}_n$) and for scalar or coarse patch-based parameters \cite{pockbilevel, bilevellearning, noiselearning, Chung_Reyes_Schoenlieb, DelosReyes2021}. An extension for learning the optimal sampling pattern in MRI \cite{Sherry_2020}, as well as extensions to non-local and higher order regularizers \cite{D_elia_2021, TGV_learning2} have been considered as well. Further, unsupervised approaches, employing upper level energies that do not depend on the ground truth $\XX_{\mathrm{true}}$, i.e., $l=l(\XX(\boldsymbol{\Lambda}))$ and $M=1$, have also been  considered in a series of works \cite{bilevelconvex, bilevelTGV, bilevel_handbook, HiRa17, hintermuellerPartII}. The upper level energy considered there aims to constrain localized versions of the image residuals  $\Ad \XX - \ZZ$ within a certain tight corridor around the variance of the (Gaussian) noise $\mathbf{e}$, which is assumed to be known. Even though these bilevel optimization methods are typically accompanied by elegant mathematical theories, there exist limitations  on the computational time they require in order to give satisfactory results. For instance, employing these methods for 2D or even 3D dynamic imaging requires a vast computation effort and as a result,  these limitations pose a challenge for the application in modern medical imaging modalities and hence in the clinical routine.\\
Recently, methods that are based on neural networks (NNs) have been proposed for the task of the estimation of such regularization parameter-maps. In \cite{Afkham_2021}, the authors employ a classical supervised learning approach in order to learn the map from the data $\ZZ$ to the optimal scalar regularization parameter $\lambda$. The pipeline consists again of an offline and an online phase. More precisely, given again $M$ pairs of measured data and corresponding ground truth images $(\ZZ_{i}, \XX_{\mathrm{true}}^{i})_{i=1}^{M}$, during the first part of the offline phase, a corresponding family of optimal regularization parameters $(\lambda_{i})_{i=1}^{M}$ is computed, e.g.\ by employing a scheme like \eqref{general_bilevel} separately for each $i$. Then, in the second part of the offline phase,  using the training data $\mathcal{D}=\{(\lambda_{i}, \ZZ_{i})_{i=1}^{M}\}$, the parameters $\Theta$  of a NN $\mathcal{N}_{\Theta}$ are learned by minimizing
\begin{equation}\label{sup_learn_param_off1}
\min_{\Theta} \mathcal{L}(\Theta) := \frac{1}{M} \sum_{i=1}^{M} l( \mathcal{N}_{\Theta}(\ZZ_{i}),\lambda_{i}),
\end{equation}
for a suitable loss function $l$. Once an estimate of the optimal parameters $\Theta$ has been learned, one passes to the online phase, and given some new data $\ZZ_{\mathrm{test}}$, the regularization parameter is simply calculated by applying the learned network to $\ZZ_{\mathrm{test}}$, i.e., $\lambda_{\Theta}=\mathcal{N}_{\Theta}(\ZZ_{\mathrm{test}})$ and the classical image reconstruction problem 
\begin{equation}\label{sup_learn_param_off2}
    \underset{\XX }{\min}\; \frac{1}{2} \| \Ad \XX - \ZZ_{\mathrm{test}}\|_2^2 +  \lambda_{\Theta}\|  \nabla \XX \|_1, 
   \end{equation}
is solved by an appropriate algorithm. The idea is that, due to the good generalizability and adaptability of NNs on unseen data, the computed regularization parameter $\lambda_{\Theta}=\mathcal{N}_{\Theta}(\ZZ_{\mathrm{test}})$ will be better adapted to $\ZZ_{\mathrm{test}}$ than the ``average" $\lambda$ of the (scalar parameter version of) bilevel optimization approach \eqref{general_bilevel}.
The authors in \cite{Afkham_2021} apply this pipeline to learn scalar regularization parameters for computerized tomography reconstruction and image deblurring. Nevertheless, the computational burden for computing offline the training data as well as solving \eqref{sup_learn_param_off2} for high dimensional (3D and dynamic) problems still remains. A similar approach, where the supervised learning problem \eqref{sup_learn_param_off1} is performed at the level of small image-patches, was performed in \cite{Nekhili_2022} for the image denoising problem.\\
In this work, inspired by the recent success of unrolled NNs  \cite{Monga_2021}, and targeting a variety of inverse problems including dynamic ones,  we apply a different strategy for the construction of the regularization parameter-maps. We construct an unrolled NN which corresponds to an implementation of an iterative scheme of finite length to approach the solution of problem \eqref{eq:tv_min_problem} assuming a \textit{fixed} regularization parameter-map. Within the unrolled NN, the regularization parameter-map is estimated from the input data and is used throughout the whole reconstruction scheme.  To be more precise, given some initial estimate $\XX_{0}$ we work with an iterative scheme
\begin{equation}\label{iteraive_scheme}
\XX_{T}=S^{T}(\XX_{0},\ZZ, \boldsymbol{\Lambda},\Ad), \quad T=0,1,2,\ldots
\end{equation}
for which we know that $\XX_{T}\to S^*(\ZZ,\boldsymbol{\Lambda},\Ad)$ as $T\to\infty$ where $S^*(\ZZ,\boldsymbol{\Lambda}, \Ad)$ is a solution of \eqref{eq:tv_min_problem_spatial}. 
We note that sometimes,  we will drop the dependence of $S^{T}$ on $\XX_{0}$, $\ZZ$, $\Ad$, for notational convenience. Then, for some fixed number of iterations $T\in \mathbb{N}$, our unrolled NN reads as follows:
\begin{equation}\label{unrolled_intro}
\left \{
\begin{aligned}
\boldsymbol{\Lambda}_{\Theta}&=\mathrm{NET}_{\Theta}(\XX_0),\\
\XX_{1}&=S^{1}(\XX_{0},\ZZ, \boldsymbol{\Lambda}_{\Theta},\Ad),\\
\XX_{2}&=S^{2}(\XX_{0},\ZZ, \boldsymbol{\Lambda}_{\Theta},\Ad),\\
\vdots\\
\XX_{T}&=S^{T}(\XX_{0},\ZZ, \boldsymbol{\Lambda}_{\Theta},\Ad).
\end{aligned}
\right.
\end{equation}
Here, $\mathrm{NET}_{\Theta}$ denotes some convolutional NN (CNN) with learnable parameters $\Theta$. We denote by $\mathcal{N}_{\Theta}^{T}$ the overall resulting network, i.e. 
\[\mathcal{N}_{\Theta}^{T}(\XX_0)=S^{T}(\XX_{0},\ZZ, \boldsymbol{\Lambda}_{\Theta},\Ad)
=S^{T}(\XX_{0},\ZZ, \mathrm{NET}_{\Theta}(\XX_0),\Ad).\] The unrolled NN can then be end-to-end trained in a supervised manner on a set of input-target image-pairs.
This resulting network can be identified as a pipeline that combines in a sequential way 
\begin{itemize}
    \item the estimation of the regularization parameter-map which is adapted to the  data $\ZZ$ (and hence in medical imaging to the new patient) and 
    \item the iterative scheme that solves the  image reconstruction problem.
\end{itemize}
In particular, given a new unseen input data $\ZZ_{\mathrm{test}}$, the regularization  parameter-map $\boldsymbol{\Lambda}_{\Theta}$ is estimated and  stays fixed. Then, the reconstruction problem is solved by unrolling an appropriate algorithm. As such, the resulting method is entirely interpretable and naturally inherits all convergence properties of the initial reconstruction algorithm since the data-driven component  merely lies in the choice of the parameter-map.\\
 Our approach can be considered to belong to the family of  recently developed image reconstruction methods that combine elements both from \textit{model-based} and \textit{data-driven} regularization approaches. This is a modern and active field of research where interpretability and convergence guarantees from the traditional variational image reconstruction approaches are combined with the flexibility and adaptability of the deep-learning based methods. These combined approaches which aim to bring together the best of both worlds can result from instance by learning the regularization functional $\mathcal{R}$ from data and embed it into a scheme like\eqref{eq:reg_problem}, see \cite{Lunz_2018, Li_2020}, by enforcing the reconstruction to be close to an output of a network via $\mathcal{R}$, e.g.\  $\mathcal{R}(\,\cdot\,) = \|\,\cdot\, - u_{\Theta}(\mathbf{s})\|_2^2$ for some network $u_{\Theta}$ with trainable parameters $\Theta$ and input $\mathbf{s}$ \cite{schlemper2017deep, hyun2018deep, Kofler_2020, Duff_2021}, by substituting proximal operators in classical iterative schemes by learned NN denoisers (in a ``plug-and-play'' fashion) \cite{Meihardt_2017, Romano_2017}, or by using learned iterative schemes \cite{Adler_2017, Adler_2018, Hammernik_2018, Kofler_2021}, see also the review papers \cite{McCann_2017, arridge_maass_oktem_schonlieb_2019, Monga_2021}. Since one of our choices for the iterative scheme \eqref{iteraive_scheme} will be  the Primal-Dual Hybrid Gradient method (PDHG) of Chambolle and Pock \cite{chambolle2011first}, our approach is  related to the Learned Primal-Dual method \cite{Adler_2018}, where the proximal operators in the primal and dual step of PDHG are fully substituted by learnable networks. Here, we decrease the complexity but increase the interpretability by keeping the iterations of the iterative scheme  untouched and put all the power of NNs in the estimation of the input-dependent regularization parameter-map, given in the  first line of \eqref{unrolled_intro}. As a result, our approach can be regarded as an intermediate approach between \cite{Afkham_2021} and  \cite{Adler_2018}.
 One of the main reasons we follow this approach is because, apart from the increased interpretability, we also target dynamic imaging applications and we are particularly interested in the interplay between the learned temporal and spatial regularization. As far as we are aware and in contrast to static imaging problems,  there are no existing works on automatically computing regularization parameters for dynamic problems that are  both spatially and temporally varying.
Furthermore, because the ``black-box'' nature of CNNs in entirely put on the estimation of the regularization parameter-maps, the probability to possibly observe   instabilities of the method in the sense of \cite{Colbrook_2022} is rather small. From a theoretical point of view, at least for denoising, it can be shown that for smooth regularization parameter-maps, no artefacts (i.e.\ new discontinuities) can appear in the reconstructions and for rougher weights, any creation of new discontinuities can be controlled  \cite{caselles2007discontinuity, mine_spatial, jalalzai2014discontinuities}. Moreover, even the worst-case of locally very large produced regularization weights  will only result in  a locally flat area in the image with controlled values.  
Further, from a practical point of view, in preliminary experiments, we have observed that even fully random regularization parameter-maps yield reconstructions whose artefacts can be at worst similar to the ones which would result from a locally too low or too strong TV-regularization.

We evaluate the proposed approach on a variety of reconstruction problems such as accelerated cardiac cine MRI, quantitative MRI, dynamic image denoising and low-dose computerized tomography (CT). We show that the proposed approach significantly improves the reconstruction results which can be obtained by the respective methods using  only scalar regularization values, and better preserves the fine scale details by adapting the regularization strength to the given data. We finally stress that even though here we focus on TV regularization only, the proposed framework can be in principle adapted to more sophisticated regularization methods 

The rest of the paper is organized as follows. In Section \ref{sec:spat_temp_vrm} we review spatio-temporal TV-based regularization and introduce  notation. In Section \ref{sec:proposed_method}, we present our proposed approach for obtaining a data/patient-adaptive spatial or spatio-temporal regularization parameter-map. We investigate theoretical aspects of the proposed approach in Section~\ref{sec:consistency_ana_scheme}, focusing on the consistency of the unrolled scheme. In Section \ref{sec:experiments}, we conduct experiments to evaluate the proposed method on  different imaging problems. We conclude the work in Section \ref{sec:conclusion} by discussing some aspects of the proposed approach, its limitations and possible future research directions.

\section{Spatio-Temporal Variational  Regularization Models}\label{sec:spat_temp_vrm}
In this section we introduce in more detail the different spatio-temporal regularization functionals and we review relevant works from the literature. In parallel, we also fix the different notations  for the several regularization parameters (scalar and spatially/spatio-temporally varying) in relationship to the way these are computed, e.g.\ supervised, unsupervised, ground truth-based, NNs-based.

\subsection{Spatio-Temporal Total Variation}
Setting $V^{n}:=V^{n_{x}\times n_{y} \times n_{t}}$, $n_{x}, n_{y}, n_{t}\in \mathbb{N}$, we define the discrete spatio-temporal gradient operator $\nabla : V^{n}\to (V\times V\times V)^{n}$ as 
\begin{equation}\label{def_spatio_temp_grad}
\nabla \XX(z)= [\nabla_x \XX(z), \nabla_y \XX(z), \nabla_t \XX(z)]^\trans, \quad \XX \in V^{n},
\end{equation}
where $\nabla_x, \nabla_y, \nabla_t$ are finite difference operators along the corresponding dimension. Note that in the case $V=\mathbb{C}$,  $\nabla_x \XX(z):=[\nabla_x \mathrm{Re}(\XX(z)), \nabla_x \mathrm{Im}(\XX(z))]$, with $\mathrm{Re}(\XX(z))$ and $\mathrm{Im}(\XX(z))$ denoting the real and the imaginary part of $\XX(z)\in \mathbb{C}$, respectively and similarly for $\nabla_y \XX(z)$ and $\nabla_t \XX(z)$. Here, $z\in [1, \ldots, n_{x}]\times [1, \ldots, n_{y}]\times [1, \ldots, n_{t}]:=\mathcal{I}$ denotes the set of indices. The (anisotropic) total variation of $\XX\in V^{n}$ is the $\ell_{1,1}$-norm of $\nabla \XX$
\begin{equation}\label{standard_TV}
\mathrm{TV}(\XX)=\|\nabla \XX\|_{1}=\sum_{z\in \mathcal{I}}|\nabla \XX(z)|_{1}:= \sum_{z\in \mathcal{I}} |\nabla_x \XX(z)| + |\nabla_y \XX(z)| +  |\nabla_t \XX(z)|,
\end{equation}
and the isotropic version is defined analogously using the $\ell_{2,1}$-norm, i.e., with the Euclidean norm $|\cdot|_{2}$ being used instead of $|\cdot|_{1}$ in \eqref{standard_TV}. In this work we employ the anisotropic version of TV and we mention that for $V=\mathbb{C}$ we also set $|\nabla_x \XX(z)|=|\nabla_x \mathrm{Re}(\XX(z))| + |\nabla_x \mathrm{Im}(\XX(z))|$, and similarly for the $y$- and $t$-direction.

\subsection{Notations on the Different Regularization Weights and Corresponding  Spatio-Temporal Total Variation Functionals}
In general, we will denote scalar and spatially (and/or temporally) varying regularization parameters with $\lambda$ and $\boldsymbol{\Lambda}$, respectively. Whenever such a parameter is the output of  a NN (with weights $\Theta$), the subindex $\Theta$ will be used, i.e., $\lambda_{\Theta}$ or $\boldsymbol{\Lambda}_{\Theta}$. If such a parameter produces the best corresponding TV-reconstruction with respect to the PSNR  for some given data $\ZZ$, it will be denoted by $\lambda_{\mathrm{P}}$ or $\boldsymbol{\Lambda}_{\mathrm{P}}$. For instance, these would be the optimal parameters that are solutions to the bilevel scheme \eqref{general_bilevel} when the upper level energy $l_{\mathrm{PSNR}}$ is used and $M=1$. In that case, the training and the test image coincide. If the ``best" is understood as ``on average" based on some training data, i.e.\ bilevel scheme \eqref{general_bilevel}  with $l_{\mathrm{PSNR}}$ and $M>1$, we denote these parameters as $\lambda_{\tilde{\mathrm{P}}}$ or $\boldsymbol{\Lambda}_{\tilde{\mathrm{P}}}$. In that case the test image is not part of the training data.

On the other hand, we use superindices to index whether the different components of the regularization parameters that correspond to the different dimensions are the same or not. For instance,
\begin{align}
\lambda^{x,y,t}&=(\lambda^{x}, \lambda^{y}, \lambda^{t}) \in \R_{+}^{3}\label{def_lambda_x_y_t},\\
\lambda^{xy,t}&=(\lambda^{xy}, \lambda^{xy}, \lambda^{t}) \in \R_{+}^{3}\label{def_lambda_xy_t},\\
\lambda^{xyt}&=(\lambda^{xyt}, \lambda^{xyt}, \lambda^{xyt}) \in \R_{+}^{3},\label{def_lambda_xyt}
\end{align}
denote  parameters that weight all the components differently, weight only the spatial components equally, and weight all the components equally respectively.  For instance,  the use of \eqref{def_lambda_x_y_t}  leads to the following version of weighted TV:

\begin{equation}\label{scalar_lambda_x_y_t_TV}
\mathrm{TV}_{\lambda^{x,y,t}}(\XX):=\|\lambda^{x,y,t}\nabla \XX\|_{1}=\sum_{z\in \mathcal{I}}\lambda^{x,y,t}|\nabla \XX(z)|_{1}:= \sum_{z\in \mathcal{I}} \lambda^{x}|\nabla_x \XX(z)| +\lambda^{y} |\nabla_y \XX(z)| +  \lambda^{t}|\nabla_t \XX(z)|.
\end{equation}
In contrast, $\lambda^{xy,t}=(\lambda^{xy}, \lambda^{xy}, \lambda^{t})$ denotes a parameter where the spatial components $x$ and $y$ are weighted equally. Analogously, we define the spatio-temporally varying versions, generally denoted by $\LLambda\in \mathbb{R}_{+}^{qn}$. In particular, we define
\begin{align}
\boldsymbol{\Lambda}^{x,y,t}&=(\boldsymbol{\Lambda}^{x}, \boldsymbol{\Lambda}^{y},\boldsymbol{\Lambda}^{t}) \in  (\R_{+}^{n})^{3},\label{def_Lambda_x_y_t}\\
\boldsymbol{\Lambda}^{xy,t}&=(\boldsymbol{\Lambda}^{xy}, \boldsymbol{\Lambda}^{xy}, \boldsymbol{\Lambda}^{t}) \in (\R_{+}^{n})^{3},\label{def_Lambda_xy_t}\\
\boldsymbol{\Lambda}^{xyt}&=(\boldsymbol{\Lambda}^{xyt}, \boldsymbol{\Lambda}^{xyt}, \boldsymbol{\Lambda}^{xyt}) \in (\R_{+}^{n})^{3},\label{def_Lambda_xyt}
\end{align}
with \eqref{def_Lambda_xy_t}, for instance, leading to the following version of weighted TV 
\begin{align}\label{spatial_Lambda_x_y_t_TV}
\mathrm{TV}_{\boldsymbol{\Lambda}^{xy,t} }(\XX)&:=\| \boldsymbol{\Lambda}^{xy,t}\nabla \XX\|_{1}=\sum_{z\in \mathcal{I}}|\boldsymbol{\Lambda}^{xy,t}(z)\nabla \XX(z)|_{1} \\
&:=\sum_{z\in \mathcal{I}} \boldsymbol{\Lambda}^{xy}(z)|\nabla_x \XX(z)| + \boldsymbol{\Lambda}^{xy}(z)|\nabla_y \XX(z)| +  \boldsymbol{\Lambda}^{t}(z)|\nabla_t \XX(z)|.
\end{align}
Here the multiplication of $\boldsymbol{\Lambda}^{xy,t}$ and $\nabla\XX$ is considered component-wise. The full notations of the type, e.g.\  $\boldsymbol{\Lambda}_{\Theta}^{xy,t}$,  $\boldsymbol{\Lambda}_{\mathrm{P}}^{xy,t}$,  $\lambda_{\mathrm{P}}^{xy,t}$  have the obvious meaning. 
By definition, we have the following, easily checked, inequalities for the PSNRs of the corresponding reconstructed images:
\begin{align}
&\lambda_{\tilde{\mathrm{P}}}^{xy,t}\preceq  \lambda_{\mathrm{P}}^{xy,t}\preceq \boldsymbol{\Lambda}_{\mathrm{P}}^{xy,t},\\
&\lambda_{\tilde{\mathrm{P}}}^{x,y,t}\preceq  \lambda_{\mathrm{P}}^{x,y,t}\preceq \boldsymbol{\Lambda}_{\mathrm{P}}^{x,y,t},\\
&\lambda_{\mathrm{P}}^{xy,t}\preceq \lambda_{\mathrm{P}}^{x,y,t},\\
&\boldsymbol{\Lambda}_{\mathrm{P}}^{xy,t}\preceq \boldsymbol{\Lambda}_{\mathrm{P}}^{x,y,t}.
\end{align}
Of particular interest will be the comparison of the above to the reconstructions that correspond to the  parameters $\boldsymbol{\Lambda}_{\Theta}^{x,y,t}$/$\boldsymbol{\Lambda}_{\Theta}^{xy,t}$/$\boldsymbol{\Lambda}_{\Theta}^{xyt}$ that we learn through our unrolled scheme. We also note that the  quantities correspond to spatial regularization only, i.e., static imaging tasks are defined straightforwardly by omitting the temporal component $t$.

\subsection{Related Literature on Spatio-Temporal Total Variation-type Functionals}
There have been quite a few related works in the dynamic inverse problems literature that employ regularization functionals of the type \eqref{scalar_lambda_x_y_t_TV}, or higher order extensions. Even though, we will also use our approach for static tasks, we briefly review these works since the literature on computing spatio-temporal regularization parameters for dynamic problems is essentially void. 
In \cite{TGV_video}, the authors use a regularization functional defined as an infimal convolution of functionals of the type \eqref{scalar_lambda_x_y_t_TV} for video denoising and decompression, an approach which splits the image sequence into two components with little change in space and time respectively. A bilevel approach for dynamic denoising is considered in \cite{BilevelDynamicDenoising}. A higher extension of the approach, applied to dynamic MRI was investigated in \cite{Schloegl_2017} and for dynamic PET in \cite{bergounioux:hal-01694064}. Regularization of the type \eqref{scalar_lambda_x_y_t_TV} has been also considered for dynamic tomographic imaging \cite{Papoutsellis_2021} and dynamic cardiac MRI \cite{wang2013compressed}. We stress however that in all these works  all regularization parameters are scalar and they are manually selected. 
Here we allow for better flexibility in the regularization by automatically computing  regularization parameters that are both spatially and temporally dependent, and we let these parameters guide the decoupling to static and moving parts in the image sequence.

\section{Proposed Unrolled Network Structure}\label{sec:proposed_method}

As described in \eqref{unrolled_intro}, our network architecture $\mathcal{N}_{\Theta}^T$  consists of two parts. 
The first part of the network is concerned with the determination of the regularization parameter-map $\boldsymbol{\Lambda}_{\Theta}$ which is subsequently fed into the second part. We describe this procedure in more detail in Section \ref{subsubsec:obtaining_lambda}. Assuming the regularization parameter-map $\LLambda_{\Theta}$ is fixed in the second module of the network,  $\LLambda_{\Theta}$ is fed into an unrolled iterative scheme of length $T$ which, if run until convergence, exactly solves 
\begin{equation}\label{eq:tv_min_problem_parameter_map}
    \underset{\XX }{\min}\; \frac{1}{2} \|\Ad \XX - \ZZ \|_2^2 + \|\LLambda_{\Theta} \nabla\XX\|_1.
\end{equation}
For the latter we choose $T$ iterations of the PDHG algorithm \cite{chambolle2011first}, which we briefly recall here.
\begin{figure}[!t]
    \centering
    \includegraphics[width=0.9\linewidth]{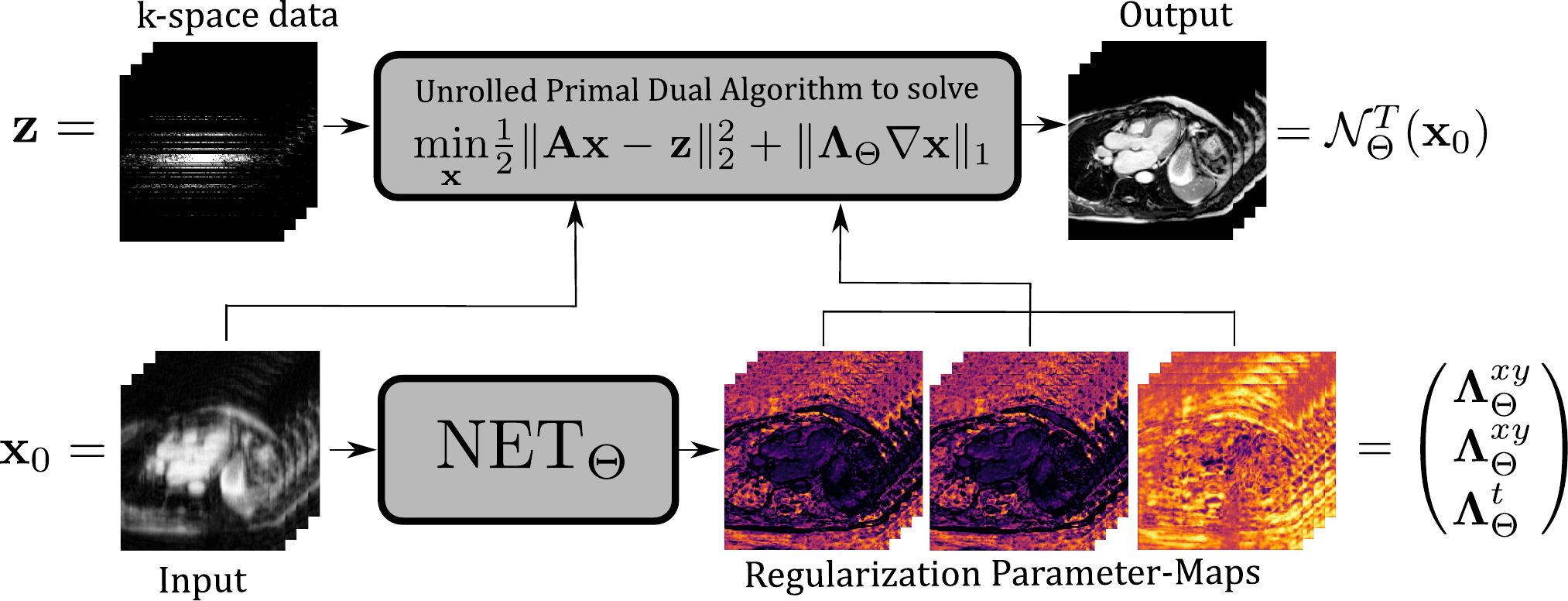}
    \caption{\small{Illustration of the proposed network architecture for a dynamic cardiac MR image reconstruction problem. The network consists of a sub-network which estimates the regularization parameter-maps and a sub-network which reconstructs the image using the PDHG algorithm described in Algorithm \ref{algo:pdhg_algo}. First, a spatio-temporal parameter-map $ \boldsymbol{\Lambda_{\Theta}}$, here as in \eqref{def_Lambda_xy_t},  is estimated by applying $\mathrm{NET}_{\Theta}$ to an input image $\XX_0$. The regularization parameter-map is then used within the reconstruction network which assumes the parameter-map to be fixed. The regularization parameter-map is trained such that the output of the PDHG algorithm is close to a reference image. } }
    \label{fig:tv_primal_dual_nn.pdf}
\end{figure}

By denoting $X=V^{n}$, $Z=V^{m}, Q=V^{qn}$, the image reconstruction problem \eqref{eq:tv_min_problem_spatial} can be equivalently formulated as 
 \begin{equation}\label{eq:fKx_gx}
\underset{\XX\in X }{\min}\; f(\mathbf{K}\XX) + g(\XX),
\end{equation} 
with $f:Y=Z\times Q\to \mathbb{R}$ where
\begin{align}
    f(\YY)=f(\PP, \QQ) := f_1(\PP) + f_2(\QQ) = \frac{1}{2} \| \PP - \ZZ \|_2^2 + \| \boldsymbol{\Lambda}_{\Theta} \, \QQ \|_1, \quad
    \mathbf{K} := \begin{bmatrix}
        \Ad\\
        \nabla
        \end{bmatrix},\quad
    g(\XX) := \mathbf{0} \label{eq:identification3}.
\end{align}
Here, the variables $\PP, \QQ$ belong to the  finite dimensional Euclidean spaces $Z$ and $Q$ that correspond to the specificities, e.g.\ dimensions, of each problem, and $\mathbf{K}:X\to Y$.

The PDHG algorithm for solving  problems of the general form  \eqref{eq:fKx_gx}, i.e., with $f$ and $g$ convex as well as $\mathbf{K}$ bounded and linear is described in Algorithm \ref{algo:pdhg_algo}. Recall that, for a convex function $h$ and scalar $\sigma>0$, the proximal operator $\mathrm{prox}_{\sigma h}$  is defined as
\begin{equation}\label{eq:prox}
    \mathrm{prox}_{\sigma h}(\overline{\YY}):= \underset{\YY}{\operatorname{argmin}} \, \frac{1}{2}\| \YY - \overline{\YY} \|_2^2  + \sigma h(\YY),
\end{equation}
while the convex conjugate of $h$ is defined as
\begin{equation}\label{eq:conjugate}
    h^\ast(\overline{\YY}):= \underset{\YY }{\max} \,  \langle \YY , \overline{\YY} \rangle  - h(\YY).
\end{equation}

In order to be consistent with our purposes, we have stated the Algorithm \ref{algo:pdhg_algo} such that it terminates in $T$ iterations with an output $\XX_{T}$. However, from standard convergence analysis it holds that $\XX_{T}\to \XX^{\ast}$ as $T\to\infty$, where $\XX^{\ast}$ solves \eqref{eq:tv_min_problem_spatial}.

\begin{algorithm}[h]
        \caption{Unrolled PDHG algorithm \cite{chambolle2011first}}\label{algo:pdhg_algo}
  \begin{algorithmic}[1]
  \INPUT $L = \| \mathbf{K}\|$,\; $\tau\sigma\le  1/L^2$, \;  $\theta = 1$, \; $\text{initial guess}~\XX_0$
  \PARAMETER number of iterations $T>0$ 
   \OUTPUT reconstructed image $\XX_{T}$ 
  \STATE $\bar{\XX_0} = \XX_0$
  \STATE $\YY_0 = \mathbf{0}$
    \FOR {$k < T$  }
    \STATE $\YY_{k+1} = \mathrm{prox}_{\sigma f^\ast}( \YY_k + \sigma \mathbf{K}  \bar{\XX}_k)$ 
     \STATE $\XX_{k+1} = \mathrm{prox}_{\tau g}( \XX_k - \tau \mathbf{K}^\trans  \YY_{k+1})$
     \STATE $\bar{\XX}_{k+1} = \XX_{k+1} + \theta(\XX_{k+1} - \XX_k)$ 
    \ENDFOR
  \end{algorithmic}
\end{algorithm}

\begin{remark}\label{rem:Lambda_in_algo}
Later we recall the precise form of Algorithm \eqref{algo:pdhg_algo} for the problem \eqref{eq:tv_min_problem_spatial}. Here, we only mention that the $\mathrm{prox}_{\sigma f^{\ast}}$ for $f$ as in \eqref{eq:identification3}, decouples  to $\mathrm{prox}_{\sigma f_{1}^{\ast}}$ and  $\mathrm{prox}_{\sigma f_{2}^{\ast}}$, with the latter acting as a pointwise projection (``clipping'') onto the bilateral set $[-\boldsymbol{\Lambda}_{i}, \boldsymbol{\Lambda}_{i}]$, $i=1,\ldots, qn$. In particular, the map $\LLambda \mapsto \mathrm{prox}_{\sigma f_{2}^{\ast}}(\QQ)$ is Lipschitz with constant one, for every $\QQ$. We remark that this is the only place where the parameter $\LLambda$ appears in the  version of Algorithm \eqref{algo:pdhg_algo} for the problem \eqref{eq:tv_min_problem_spatial}.
\end{remark}

\begin{remark}\label{rem:PD3O}
We mention that for the low-dose CT application which we consider in Section \ref{sec:CT}, we will be using a generalization of PDHG, namely the PD3O algorithm \cite{Yan2018} which is better adapted for the Kullback-Leibler divergence fidelity term used there. We give more details  later in that section.
\end{remark}

\subsection{Obtaining the Regularization Parameter-Map  Via a CNN}\label{subsubsec:obtaining_lambda}
In our set-up, $\boldsymbol{\Lambda}_{\Theta}$ is the output of a CNN with parameters $\Theta$, denoted by $\mathrm{NET}_{\Theta}$, which takes as an input an initial image $\XX_0$, i.e., $\boldsymbol{\Lambda}_{\Theta} = \mathrm{NET}_{\Theta}(\XX_0)$.
Depending on the structure of the considered imaging problem, we can explore different possibilities for the construction of the latter. For instance, for a dynamic imaging problem, i.e., 2D + time, we might prefer to attribute equal importance to the $x$- and $y$-direction, but use a different parameter-map for the temporal component resulting in

 \begin{equation}\label{eq:spatio_temporal_lambda_map}
    \boldsymbol{\Lambda}_{\Theta}=(\boldsymbol{\Lambda}_{\Theta}^{xy}, \boldsymbol{\Lambda}_{\Theta}^{xy},\boldsymbol{\Lambda}_{\Theta}^{t}).
\end{equation}

This choice is motivated by the later shown cardiac cine MRI reconstruction problem. There, the temporal dimension is the one which on the one hand exhibits the largest correlation to be exploited by the TV-method, but on the other hand also the one which contains the diagnostic information and therefore requires special care to ensure that important features are preserved.\\
For a 3D imaging problem, one could for example attribute equal importance to all spatial-directions or opt for a construction as in \eqref{eq:spatio_temporal_lambda_map}, if for example the $z$-direction has a different resolution than the $x$- and $y$-directions. Moreover, for complex-valued images, it seems intuitive to share the same regularization map across the real and the imaginary parts of the images.\\
The core of the overall network, denoted by  $u_{\Theta}$, consists of  a (sub-)CNN with high expressive capabilities such as the U-Net \cite{Ronneberger2015}. 
To constrain the regularization parameter-maps to be strictly positive, we then apply a softplus activation function $\phi$ and, as last operation, we  multiply the output by a positive parameter $t>0$. 
Empirically, we have experienced that the network's training benefits in terms of faster convergence if the order of the scale of the output is properly set depending on the application. This can be achieved either by accordingly initializing the weights of the network $u_{\Theta}$, or in a simpler way, as we do here by scaling the output of the CNN. Summarizing, given an input image $\XX_0$, we estimate the corresponding regularization parameter-map by 
\begin{equation}\label{eq:lambda_reg_map}
    \boldsymbol{\Lambda}_{\Theta}=\mathrm{NET}_{\Theta}(\XX_0) = t\,  \phi( u_{\Theta}(\XX_0)).
\end{equation}
We finally recall that the overall network has the form
\begin{equation}\label{eq:final_net}
\mathcal{N}_{\Theta}^{T}(\XX_0)=
S^{T}(\XX_{0},\ZZ, \mathrm{NET}_{\Theta}(\XX_0),\Ad).
\end{equation}

\begin{remark}
 Note that in our set-up we use the same quantity $\XX_{0}$ as the input for the CNN-block $\mathrm{NET}_{\Theta}(\XX_{0})$ as well as the initialization for the unrolled PDHG $S^{T}(\XX_{0},\ldots)$. According to our experience, this produces satisfactory results, see also the discussion in Section \ref{sec:initialization}. However this is not a hard constraint of the method and one could also further experiment with having different values for these variables. 
\end{remark}

\begin{figure}[!t]
    \centering
    \includegraphics[width=0.8\linewidth]{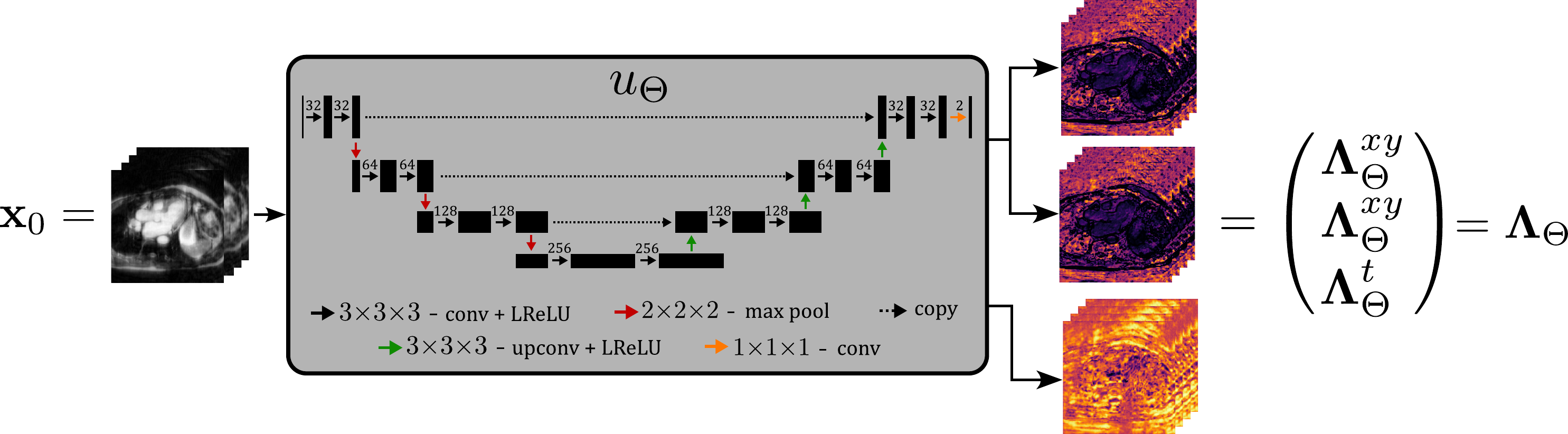}
    \caption{\small{
The CNN-block  $u_{\Theta}$ is a U-Net which we hyper-parametrize by the number of encoding stages, the number of convolutional layers and the initial number of filters applied to the input-image.  The filter-sizes as well as the number of output channels of $u_{\Theta}$ depend on the considered application. The network shown here is a 3D U-Net with four encoding stages, two convolutional layers per stage and 32 initially applied filters. Here, the CNN is constructed to yield two different components of the regularization parameter-map, i.e.\ $\LLambda_{\Theta}^{xy}$ and $\LLambda_{\Theta}^t$, which are then used to construct the final spatio-temporal regularization parameter-map $\LLambda_{\Theta}^{xy,t}$ according to \eqref{def_Lambda_xy_t}, i.e.\ by $\LLambda_{\Theta}^{xy,t} = (\LLambda_{\Theta}^{xy}, \LLambda_{\Theta}^{xy}, \LLambda_{\Theta}^{t})$.}}
    \label{fig:u_Theta.pdf}
\end{figure}

\subsection{Network Training}
Training the network $\mathcal{N}_{\Theta}^T$ refers to minimizing a chosen energy-function $\mathcal{L}$ over a set of input-target training pairs $\mathcal{D}=\{ (\XX_0^i, \XX_{\mathrm{true}}^i)_{i=1}^M : \XX_0^i = \Ad^\ddagger \ZZ_i, \,\ZZ_i:= \Ad \XX_{\mathrm{true}}^i + \mathbf{e}_i\}$. Here $\Ad^\ddagger$ denotes some reconstruction operator, e.g.  the pseudo-inverse of $\Ad$, which provides  the inputs $\XX_{0}^{i}$ and  $\mathbf{e}_i$ are noisy terms. Using an appropriate loss function $l$  and potentially some regularization function $r$ for the weights $\Theta$, we end up with the following energy-function:
\begin{equation}\label{eq:loss_fct}
\mathcal{L}(\Theta) = \frac{1}{M} \sum_{i=1}^{M} l \big(\mathcal{N}_{\Theta}^T( \XX_0^{i})\, , \, \XX_{\mathrm{true}}^{i}\big)+ r(\Theta).
\end{equation}

Note that $r$ can also encode some constraints on $\Theta$ by being an indicator of some set.
Note that the network is trained end-to-end from the initial reconstruction to its estimate. Therefore, the set $\Theta$ is adjusted such that the estimated parameter-map $\boldsymbol{\Lambda}_{\Theta}$ is appropriate for a subsequent reconstruction using a suitable reconstruction algorithm as for example, the primal-dual method in Algorithm \ref{algo:pdhg_algo}. This conceptually highly differs from approaches as in \cite{Afkham_2021}, in which a network is trained to estimate the best scalar regularization parameter which is previously obtained by a time-consuming grid-search. First of all, in \cite{Afkham_2021} the learning procedure is entirely decoupled from the employed reconstruction algorithm. Second, opposed to our approach, the method  requires access to a target regularization parameter, meaning that a generalization of \cite{Afkham_2021} to regularization parameter-maps would require access to entire target regularization parameter-maps which can typically only obtained by even more time-consuming approaches. Our approach, in contrast, allows to implicitly learn the regularization parameter-maps by unrolling the reconstruction algorithm and thus only requires access to ground truth target images.

\section{Consistency Analysis of the Unrolled Scheme}\label{sec:consistency_ana_scheme}

In addition to the practical advantages of the proposed method which will be highlighted in Section \ref{sec:experiments}, we want to discuss some of the emerging theoretical questions and in particular some consistency results when we let the number of unrolled iterations $T\to\infty$. We note that there are papers that study hyperparameter search with bilevel optimization and unrolled optimization methods, see e.g. \cite{liu2020generic, liu2022general, ochs2015bilevel}. Although some of the latter articles provide consistency analysis in different contexts, we think that none of the techniques presented there can be applied to our problem.

We will  be using the space and operator notation $X, Y, Z$ and  $\mathbf{K}$ as these were defined in the previous section and we will also set $\mathcal{V}:=X\times Y$. For simplicity here we work with the real-valued case, i.e.\ $V=\mathbb{R}$.
Recall that the solution of the convex variational problem \eqref{eq:tv_min_problem_spatial}
and the corresponding $T$-th iterate of the unrolled algorithm are
 denoted by $\XX^{\ast}=S^*(\ZZ,\boldsymbol{\Lambda})$ and  $\XX_{T}=S^T(\XX_{0},\ZZ,\boldsymbol{\Lambda})$. Recall also that for the ease of notation we sometimes suppress the dependence of $S^{T}$ on the initialization $\XX_0$ of Algorithm \ref{algo:pdhg_algo}, as well as the one of the dual variable $\YY_{0}$. We then have $S^T(\XX_0,\ZZ,\boldsymbol{\Lambda}) \to S^*(\ZZ,\boldsymbol{\Lambda})$ as $T\to\infty$.
Furthermore, for this section we consider a more general fidelity term $d$, such that $f_{1}(\cdot):=d(\cdot,\ZZ)$.

Let us now consider the learning framework as presented in Section \ref{sec:proposed_method} 
\begin{equation}
\min_{\Theta\in \mathbb{R}^{\ell}} \mathcal{L}^T(\Theta) := \frac{1}{M} \sum_{i=1}^M l (\mathcal{N}_\Theta^T(\XX_0^i)\,,\, \XX_{\text{true}}^i) + r(\Theta), \label{eq:empirical_risk_T} 
\end{equation}
as well as the corresponding training scheme where no unrolling is taking place, i.e.,
\begin{equation}
 \min_{\Theta \in \mathbb{R}^{\ell}} \mathcal{L}^*(\Theta) := \frac{1}{M} \sum_{i=1}^M l (\mathcal{N}_\Theta^*(\XX_0^i)\,,\, \XX_{\mathrm{true}}^i) + r(\Theta), \label{eq:empirical_risk}
 \end{equation}
where we used analogously the notation $\mathcal{N}_\Theta^*(\XX_{0}) := S^*(\ZZ,\LLambda_\Theta(\XX_{0}))$. Our target will be to  show  convergence of ($\epsilon$)-minimizers
of \eqref{eq:empirical_risk_T} to minimizers of \eqref{eq:empirical_risk} as $T\to\infty$ under appropriate conditions via a $\Gamma$-convergence argument. 

Naturally, in order to guarantee  existence of minimizers for the problems \eqref{eq:empirical_risk_T} and  \eqref{eq:empirical_risk},   the functionals $\mathcal{L}^{T}$ and $\mathcal{L^{\ast}}$ must be coercive, in addition to the standard lower semicontinuity assumptions. However, it is not so clear if this can be achieved without imposing coercivity via the regularization function $r$,  which can be the case when e.g.\  $r$ is some norm in $\mathbb{R}^{\ell}$ or an indicator function of a bounded set. Even though strictly speaking, it is not needed for our main consistency result Corollary \ref{cor:consistency}, we will assume that  the minimization problems \eqref{eq:empirical_risk_T} and  \eqref{eq:empirical_risk} indeed admit solutions.
Of course in deep learning practice, one does not compute minimizers for these problems,  but rather it is  aimed that the energy is decreased up to some degree based on a validation set, in order to guarantee generalizability. However, the analysis presented here can serve as a starting point to further show consistency in the level of stationary points and/or energy decrease using validation sets.

Below we summarize  a series of assumptions which we will need next:

\begin{assumption}\label{main_assumptions}
We assume that the following hold:
\begin{enumerate}
\item[(i)] The operator $\Ad:X\to Z$ is injective. 
\item [(ii)] The fidelity term $d(\cdot,\ZZ)$ is $\mu_{\ZZ}$-strongly convex and Lipschitz continuously differentiable for every $\ZZ\in Z$. 
We denote by $L_{\ZZ}$ the corresponding Lipschitz constant.
\item[(iii)] The parameters $\sigma, \tau>0$ in Algorithm \ref{algo:pdhg_algo} are small enough such that the matrix
\begin{equation}\label{M_st}
\mathbf{M} = 
\begin{pmatrix}
\frac{1}{\tau} \Id & -\mathbf{K}^{\trans} \\ 
-\mathbf{K} & \frac{1}{\sigma} \Id 
\end{pmatrix}
\end{equation}
is symmetric, positive definite and thus defines a norm in $\mathcal{V}$.
Then there exist $c,C>0$ such that
\begin{equation}
c \| \VV \|_{2} \leq \|\VV \|_{\mathbf{M}} := \sqrt{\langle \mathbf{M} \VV , \VV \rangle} \leq C  \| \VV \|_{2} \quad \text{for every $\VV \in \mathcal{V}$.}
\end{equation}
\item[(iv)] The regularization function $r:\mathbb{R}^{\ell}\to\overline{\mathbb{R}}:= \mathbb{R} \cup \lbrace + \infty \rbrace$ is proper and lower semicontinuous.
\item[(v)] The loss function $l: X \times X \to \mathbb{R}$ is continuous.
\item[(vi)] The activation functions in the U-Net $u_{\Theta}$ are continuous. 
\end{enumerate}
\end{assumption}

\begin{remark}\label{rem:injectivity}
The injectivity of the operator $\Ad$,
together with the strong convexity of $d$, is used  in order to ensure that $\XX\mapsto d(\Ad \XX, \ZZ)$ is strongly convex.
This indeed guarantees uniqueness of the solution for the variational problem, and in particular the map $\boldsymbol{\Lambda} \mapsto S^*(\ZZ,\boldsymbol{\Lambda})$ is well-defined and single-valued. For the applications we will consider in Section \ref{sec:experiments}, i.e., denoising, MRI with multiple receiver coils and CT with enough angular views and detectors, this injectivity assumption is satisfied. We note however it might be possible to drop this injectivity assumption following \cite{vaiter2015model}, \cite{vaiter2012robust},  or \cite{zhang2016one}.
\end{remark}

We start with the following Proposition \ref{prop:unroll_properties} regarding Lipschitz continuity and equicontinuity of the iterates $S^T(\ZZ,\boldsymbol{\Lambda}) $ with respect to $\boldsymbol{\Lambda}$. Note that the convergence $\boldsymbol{\Lambda}_{T}\to \LLambda$ as $T\to\infty$, in (iii) below, is merely part of a technical condition and it is not associated to the structure of our unrolled scheme where, as we have pointed out, the CNN-output $\boldsymbol{\Lambda}_{\Theta}=\mathrm{NET}_{\Theta}(\XX_0)$ remains unchanged.

\begin{proposition}\label{prop:unroll_properties}
Assuming $(i)$-$(iii)$ of Assumption \ref{main_assumptions}, the following statements hold:
\begin{enumerate}[label=(\roman*)]
\item The solution map $\boldsymbol{\Lambda} \mapsto S^*(\ZZ,\boldsymbol{\Lambda})$ is Lipschitz continuous for every $\ZZ\in Z$. In particular the following bound holds for every $\LLambda_{1}, \LLambda_{2}\in \mathbb{R}_{+}^{qn}$,
\begin{equation}
\| S^*(\LLambda_1,\ZZ) - S^*(\LLambda_2,\ZZ) \|_{2} \leq \frac{2 \| \nabla \|}{\lambda_{\min}(\Ad^\trans \Ad) \mu_\ZZ} \|\LLambda_1 - \LLambda_2\|_{2}.
\end{equation}
\item The map $\boldsymbol{\Lambda} \mapsto S^T(\XX_0,\ZZ,\boldsymbol{\Lambda})$ is Lipschitz continuous for every $\ZZ\in Z$, $\XX_0 \in X$ and  $T\in \mathbb{N}$. 
\item For $\|\LLambda\|_{2} \leq \overline{\LLambda}$ we obtain the following sub-linear rate, for $\VV_{0}:=(\XX_{0},\YY_{0})$ being the initial iterates of Algorithm \ref{algo:pdhg_algo}
\begin{equation}
\|S^{T}(\XX_{0}, \ZZ,\LLambda)-S^{\ast}(\ZZ,\LLambda)\|_{2} \leq \frac{3 C_{\ZZ,\Ad}}{T^{1/4}}  \left( 1 + 
\| \VV_0 - \VV^*(\LLambda,\ZZ) \|_\mathbf{M} \right),
\end{equation} 
where 
\begin{equation}\label{def_C_z_A}
C_{\ZZ,\Ad} := \frac{\max \left( C L_\ZZ \| \Ad \| \, , \, 4C \overline{\LLambda} \, ,\, 2\, ,\,  \lambda_{\min}(\Ad^\trans \Ad)\mu_\ZZ \right)}{\lambda_{\min}(\Ad^\trans \Ad) \mu_\ZZ},
\end{equation}
with $\lambda_{\min}(\Ad^\trans \Ad)$ denoting the smallest eigenvalue of $\Ad^\trans \Ad$.\\
\item Whenever $\boldsymbol{\Lambda}_{T}\to \boldsymbol{\Lambda}$ as $T\to\infty$, it holds $S^T(\XX_0,\ZZ,\boldsymbol{\Lambda}_{T})\to S^*(\ZZ,\boldsymbol{\Lambda})$ for every $\XX_{0} \in X$, $\ZZ \in Z$.
\end{enumerate}
\end{proposition}

\begin{proof}
(i) This statement is proved similarly to  e.g.\ in \cite[Theorem 4.1]{Reyes_Villacis} and it is strongly based on the $\mu_{\ZZ}$-strong convexity of the map  $d(\cdot,\ZZ)$ and the injectivity of $\Ad$. The statement (ii) can also be seen easily since the only dependence of $\boldsymbol{\Lambda}$ in the unrolled PDHG scheme is via the pointwise projection onto $[-\boldsymbol{\Lambda}_{i},\boldsymbol{\Lambda}_{i}]$ which is a Lipschitz map, recall Remark \ref{rem:Lambda_in_algo}. As a result, the map $\boldsymbol{\Lambda} \mapsto S^T(\XX_{0}, \ZZ,\boldsymbol{\Lambda})$ is Lipschitz,  as a composition of Lipschitz functions. 

The proof of (iii) is more involved. We fix a  ball of radius $\overline{\LLambda}$ centered at the origin, denoted by $B_{\bar{\LLambda}}\subset \mathbb{R}_{+}^{qn}$  and let $\LLambda\in B_{\overline{\LLambda}}$ be arbitrary. In what follows, we initially  suppress the dependence of all variables on $\LLambda$. 
Define the primal-dual gap
\begin{equation}
L(\XX,\YY) := \langle \mathbf{K} \XX, \YY \rangle - f^*(\YY) + g(\XX),
\end{equation}
and denote by $\VV_{T}:=(\XX_{T},\YY_{T})$ the iterates of the Algorithm \ref{algo:pdhg_algo} and by $\VV^{\ast}=(\XX^{\ast},\YY^{\ast})$ the corresponding limits. Then, the following estimate holds, see \cite[Corollary 1]{lu2022infimal} for a proof,
\begin{align}\label{lu2022infimal_inequality}
L(\XX_T,\YY) - L(\XX,\YY_T) \leq \frac{1}{\sqrt{T}} \left( \|\VV_0 - \VV^{*} \|_{\mathbf{M}}^{2}  +  
\|
\VV_0 - \VV^{*} \|_{\mathbf{M}} \|\VV - \VV^{*} \|_{\mathbf{M}} \right), 
\end{align}
where $\VV=(\XX,\YY)$ is arbitrary.
We can thus take the supremum over $\YY\in \partial f(K\XX_{T})$ in both sides in \eqref{lu2022infimal_inequality} and estimate the left hand side as follows
\begin{equation}\label{ineq_1}
\sup_{\YY\in \partial f(\mathbf{K}\XX_{T})}L(\XX_T,\YY) - L(\XX,\YY_T)
\ge \sup_{\YY\in \partial f(\mathbf{K}\XX_{T})} \langle \mathbf{K}\XX_{T}, \YY \rangle - f^{\ast}(\YY) - \langle \mathbf{K}\XX,\YY_{T} \rangle+f^{\ast}(\YY_{T})\ge f(\mathbf{K}\XX_{T})-f(\mathbf{K}\XX),
\end{equation}
where we used the fact that $ \langle \mathbf{K}\XX_{T}, \YY \rangle - f^{\ast}(\YY)=f(\mathbf{K}\XX_{T})$ if and only if $\YY\in \partial f(\mathbf{K}\XX_{T})$. By setting $\XX=\XX^{\ast}$, using the $\mu_{\ZZ}$-strong convexity of $f_{1}(\cdot)=d(\cdot,\ZZ)$, the convexity of $f_{2}$, together with $\YY^{\ast}\in \partial f (\mathbf{K}\XX^{\ast})$ we 
deduce
\begin{equation}\label{ineq_2}
f(\mathbf{K}\XX_{T}) - f(\mathbf{K}\XX^{\ast})\ge \langle \mathbf{K}^{\ast}\YY^{\ast}, \XX_{T}-\XX^{\ast} \rangle + \frac{\mu_{\ZZ}}{2} \| \Ad\XX_{T}-\Ad\XX^{\ast}\|_{2}^{2}.
\end{equation}
Taking into account that $\mathbf{K}^{\ast} \YY^{\ast}=0$ (taking limits at line 6 of Algorithm \ref{algo:pdhg_algo}, using the fact that $\mathrm{prox}_{\tau g}=Id$), using the injectivity of $\Ad$, we infer from \eqref{ineq_1} and \eqref{ineq_2}
\begin{equation}\label{lu2022infimal_inequality_v2}
\|\XX_{T}-\XX^{\ast}\|_{2}^{2}\le \frac{2}{\lambda_{\min}(\Ad^{\trans} \Ad)\mu_{\ZZ}\sqrt{T}} \sup_{\YY\in \partial f(\mathbf{K}\XX_{T})}  \left( \|\VV_0 - \VV^{*} \|_{\mathbf{M}}^{2}  +  
\|
\VV_0 - \VV^{*} \|_{\mathbf{M}} \|\VV - \VV^{*} \|_{\mathbf{M}} \right). 
\end{equation}

We proceed by estimating the last term in \eqref{lu2022infimal_inequality_v2} again making the dependence of $\LLambda$ explicit. Thus, recalling that $\VV=(\XX^{\ast}(\LLambda), \YY)$ with $(\PP,\QQ)=:\YY\in  \partial f(\mathbf{K}\XX_{T}(\LLambda))$ arbitrary, we have 
\begin{align}
\|\VV-\VV^{\ast}(\LLambda)\|_{\mathbf{M}}
&\le C \sqrt{ \|\XX^{\ast}(\LLambda) - \XX^{\ast} (\LLambda)\|_{2}^{2} + \|\PP-\PP^{\ast}(\LLambda)\|_{2}^{2}  + \|\QQ-\QQ^{\ast}(\LLambda)\|_{2}^{2}}\nonumber\\
&= C \sqrt{ \|\PP-\PP^{\ast}(\LLambda)\|_{2}^{2}  + \|\QQ-\QQ^{\ast}(\LLambda)\|_{2}^{2}}\nonumber\\
& \le C \sqrt{ \|\nabla f_{1}(\Ad\XX_{T}(\LLambda)) - \nabla f_{1} (\Ad \XX^{\ast} (\LLambda)) \|_{2}^{2} + 4\overline{\LLambda}^{2}}\nonumber\\
& \le C\left ( L_{\ZZ}\|\Ad\| \|\XX_{T}(\LLambda)-\XX^{\ast}(\LLambda)\|_{2} + 2 \overline{\LLambda}\right ),\label{last_bound}
\end{align}

where the last inequality used the fact that $\sqrt{a^{2}+b^{2}}\le a+ b$ for $a,b\ge 0$.
We also used the relationship $\PP^{\ast}(\LLambda)=\nabla f_{1}(\Ad \XX^{\ast}(\LLambda))$,   the Lipschitz continuity of $\nabla f_{1}$, as well as  $\QQ^{\ast}(\LLambda)\in \partial f_{2} (\nabla \XX^{\ast} (\LLambda))$ which implies that $\QQ^{\ast}(\LLambda)\in B_{\overline{\LLambda}}$.
By combining \eqref{lu2022infimal_inequality_v2}, \eqref{first_bound} and \eqref{last_bound}, and by defining
\[r_{0}(\LLambda):=\|\VV_{0}-\VV^{\ast}(\LLambda)\|_{\mathbf{M}},\]
we end up to 
\begin{equation*}\label{combined_bounds}
\|\XX_{T}(\LLambda)-\XX^{\ast}(\LLambda)\|_{2}^{2}\le \frac{2}{\lambda_{\min}(\Ad^\trans \Ad) \mu_\ZZ \sqrt{T}} \left( r_0(\LLambda)^2 + r_0(\LLambda) C L_\ZZ \| \Ad \| \|\XX_{T}(\LLambda)-\XX^{\ast}(\LLambda)\|_{2}  + 2 C r_0(\LLambda)\overline{\LLambda} \right).
\end{equation*} 
By setting $r_T(\LLambda) := \|\XX_{T}(\LLambda)-\XX^{\ast}(\LLambda)\|_{2}$ and
\begin{equation}
C_1 := \frac{C L_\ZZ \| \Ad \|}{\lambda_{\min}(\Ad^\trans \Ad) \mu_\ZZ} 
\quad
C_2 := \frac{4 C \overline{\LLambda}}{\lambda_{\min}(\Ad^\trans \Ad) \mu_\ZZ}
\quad
C_3 := \frac{2}{\lambda_{\min}(\Ad^\trans \Ad) \mu_\ZZ},
\end{equation}
we infer
\begin{align*}
r_T(\LLambda)^2 - \frac{2 C_1}{\sqrt{T}} r_T(\LLambda)r_0(\LLambda)   
+ \frac{C_1^2 r_0(\LLambda)^2}{T}
\le 
\frac{2}{\sqrt{T}} \left( C_3 r_0(\LLambda)^2 + C_2r_0(\LLambda) \right) + \frac{C_1^2 r_0(\LLambda)^2}{T}.
\end{align*}
After applying the binomial formula, this yields
\begin{align}
r_T(\LLambda) 
&\le
\frac{C_1 r_0(\LLambda)}{\sqrt{T}}
+
\sqrt{\frac{2}{\sqrt{T}} \left( C_3 r_0(\LLambda)^2 + C_2 r_0(\LLambda) \right) + \frac{C_1^2 r_0(\LLambda)^2}{T}} \\
&\le
\frac{C_{\ZZ ,\Ad} r_0(\LLambda)}{\sqrt{T}} + 
\sqrt{\frac{C_{\ZZ ,\Ad}}{\sqrt{T}} \left(  r_0(\LLambda)^2 +  r_0(\LLambda) \right) + \frac{C_{\ZZ ,\Ad}^2 r_0(\LLambda)^2}{T}}\\
&\leq
\frac{3 C_{\ZZ,\Ad}}{T^{1/4}} (r_0(\LLambda) + 1),
\end{align}
where the last inequality uses basic estimates, like $\sqrt{T} \leq T$, again $\sqrt{a^2 + b^2} \leq a + b$ for $a,b \geq0$ and the fact that $C_{\ZZ,\Ad} \geq 1$ by its definition \eqref{def_C_z_A}. This proves (iii). 

To show (iv) let $\LLambda_{T}\to \LLambda$ and fix 
$\overline{\LLambda}:=\sup_{T\in\mathbb{N}}\|\LLambda_{T}\|_{2}$. By (iii) we have that
\begin{equation}
\|\XX_{T}(\LLambda_T)-\XX^{\ast}(\LLambda_T)\|_{2} \leq \frac{3 C_{\ZZ,\Ad}}{T^{1/4}}  \left(1 + 
\| \VV_0 - \VV^*(\LLambda_T) \|_\mathbf{M} \right), \label{eq:bound_lambda_T}
\end{equation}
where we can further estimate by norm-equivalence
\begin{equation}\label{first_term}
\|\VV_{0}-\VV^{\ast}(\LLambda_T)\|_{\mathbf{M}}
\le C \sqrt{ \|\XX_{0}-\XX^{\ast}(\LLambda_T)\|_{2}^{2} + \|\PP_{0}-\PP^{\ast}(\LLambda_T)\|_{2}^{2} + \|\QQ_{0}-\QQ^{\ast}(\LLambda_T)\|_{2}^{2} }.
\end{equation}

Using again the boundedness of  $(\LLambda_{T})_{T\in\mathbb{N}}$, the continuity of $\nabla f_{1}$ and the relationships $\PP^{\ast}(\LLambda)=\nabla f_{1}(\Ad \XX^{\ast}(\LLambda))$ and  $\QQ^{\ast}(\LLambda)\in \partial f_{2} (\nabla \XX^{\ast} (\LLambda))$, we conclude that there exists a constant $\hat{C}>0$ independent of $T$, such that 
\begin{equation}\label{first_bound}
\|\VV_{0}-\VV^{\ast}(\LLambda_{T})\|_{\mathbf{M}}\le \hat{C}, \quad \text{ for every $T\in \mathbb{N}$}.
\end{equation}
Thus we deduce
\[\|\XX_{T}(\LLambda_{T})-\XX^{\ast}(\LLambda_{T})\|_{2}\to 0, \quad \text{as  }T\to\infty.\]
We finally use the triangle inequality to obtain
\[\|\XX_{T}(\LLambda_{T})-\XX^{\ast}(\LLambda)\|_{2} \le \|\XX_{T}(\LLambda_{T})-\XX^{\ast}(\LLambda_{T})\|_{2}+\|\XX^{\ast}(\LLambda_{T})-\XX^{\ast}(\LLambda)\|_{2}\to 0,\]
as $T\to\infty$, where we have also used $(i)$.
\end{proof}

 We can now proceed with our main  result.

\begin{theorem}
Let the Assumption \ref{main_assumptions} hold, let the training set $\mathcal{D}$ be fixed and consider the sequence of functionals $\mathcal{L}^{T}:\mathbb{R}^{\ell} \to \overline{\mathbb{R}}$, $T\in\mathbb{N}$, as well as $\mathcal{L}^{\ast}:\mathbb{R}^{\ell} \to \overline{\mathbb{R}}$ defined as in \eqref{eq:empirical_risk_T} and \eqref{eq:empirical_risk}. Then we have that $\mathcal{L}^{T}$ $\Gamma$-converges to $\mathcal{L}^{\ast}$ as $T\to\infty$.

\end{theorem}
\begin{proof}
It suffices to check the conditions in the definition of $\Gamma$-convergence \cite{dalmasogamma}, i.e.,:
\begin{enumerate}
\item[$(i)$] For all $\Theta_{T}\to \Theta$,  it holds  $\mathcal{L}^{\ast}(\Theta)\le \liminf_{T\to\infty} \mathcal{L}^{T}(\Theta_{T})$.
\item[$(ii)$] For all $\Theta\in \mathbb{R}^{\ell}$, there exists  $\Theta_{T}\to \Theta \text{ such that } \limsup_{T\to\infty}  \mathcal{L}^{T}(\Theta_{T})\le \mathcal{L}^{\ast}(\Theta)$.
\end{enumerate}
The first condition holds due to the lower semicontinuity of $r$, the continuity of the map $\Theta\mapsto \LLambda_{\Theta}$, Proposition \ref{prop:unroll_properties} $(iv)$ and the continuity of the loss function $l$.
The fact that the map $\Theta\mapsto \LLambda_{\Theta}$ is continuous, follows by the continuity of all constituent functions, in particular, from the continuity of the  activation functions of the U-Net $u_{\Theta}$. The second condition follows similarly, setting $\Theta_{T}:=\Theta$ for all $T\in\mathbb{N}$ and using the convergence of the iterative scheme, i.e. $S^{T}(\XX_{0}, \ZZ,\LLambda)\to S^{\ast}(\ZZ,\LLambda)$ as $T\to\infty$, as well as the continuity of the other involved functions.
\end{proof}
The following consistency result follows directly from the $\Gamma$-convergence, see \cite[Corollary 7.20]{dalmasogamma} for a proof. 
\begin{corollary}[Consistency of the unrolled scheme]\label{cor:consistency}
Let the Assumption \ref{main_assumptions} hold, let the training set $\mathcal{D}$ be fixed  and let $\epsilon_{T}\to 0$. Suppose that $\Theta_{T}$ is an $\epsilon_{T}$-minimizer of $\mathcal{L}^{T}$ i.e.\ $\mathcal{L}^{T}(\Theta_{T})\le \inf_{\Theta\in\mathbb{R}^{\ell}}\mathcal{L}^{T}(\Theta)+\epsilon_{T}$. Then, if $\Theta$ is an accumulation point  of $(\Theta_{T})_{T\in\mathbb{N}}$ it is a minimizer of $\mathcal{L}^{\ast}$ and $\mathcal{L}^{\ast}(\Theta)=\limsup_{T\to\infty}\mathcal{L}^{T}(\Theta_{T})$.
\end{corollary}

\section{Applications}\label{sec:experiments}

In the following, we apply our proposed method to several different imaging problems to demonstrate its versatility. The considered imaging problems differ in terms of the operator $\Ad$ and, more importantly, on the number of dimensions, e.g.\ 2D, 3D or  2D+time as well as on the specific role the dynamic component plays in the respective problem.
All images were evaluated in terms of PSNR, normalized root mean-squared error (NRMSE), structural similarity index measure\cite{wang2004image} (SSIM) and blur effect \cite{crete2007blur}.\\
\texttt{Python} code is available at  \href{https://www.github.com/koflera/LearningRegularizationParameterMaps}{\url{github.com/koflera/LearningRegularizationParameterMaps}}.

\subsection{Initialization for the Unrolled PDHG}\label{sec:initialization}
In general, an initial image for the PDHG can be directly reconstructed from the measured data by applying the adjoint of the forward operator, i.e., $\XX_0:=\Ad^\herm \ZZ$.  Often, the  set-up for realistic imaging problems is that $\Ad$ is given by a tall operator, i.e., $m > n$. Therefore, to obtain a better estimate of the unknown image from which one can estimate $\LLambda_{\Theta}$ by applying the CNN $\mathrm{NET}_{\Theta}$, one can consider the normal equation 
 \begin{equation}\label{eq:normal_eqs}
 \Ad^\herm \Ad \XX = \Ad^\herm \ZZ,
 \end{equation}
 and approximately solve it. As $\Ad$ is typically constructed such that the normal operator $\Ad^\herm \Ad$ is invertible, in the absence of noise, i.e., $\ZZ \in \mathrm{range}(\Ad)$,  solving \eqref{eq:normal_eqs} using an iterative scheme to approximate $\XX^{\dagger}:= (\Ad^\herm \Ad)^{-1}\Ad^\herm \ZZ$ would allow for a perfect reconstruction of the ground truth image.  However, in the presence of noise, early stopping is required to avoid a noise amplification during the iterations. An approximate solution of  \eqref{eq:normal_eqs} can then be used as a better initial estimate for the PDHG method as well as the image from which the CNN $\mathrm{NET}_{\Theta}$ estimates the different components of the regularization map $\boldsymbol{\Lambda}_{\Theta}$. 
For the case that the considered imaging problem is not overdetermined, i.e., $m\leq n$, e.g.\ for image denoising, one simply uses $\XX_0:=\Ad^\herm \ZZ$ as the input of $\mathrm{NET}_{\Theta}$.\\

\subsection{Dynamic Cardiac MR Image Reconstruction}\label{subsec:dyn_mri}
Here, we apply the proposed NN to a dynamic cardiac MR image reconstruction problem.  The problem consists of a set of independent 2D problems from which static images of the heart can be reconstructed. By stacking the different images along time, one can obtain a sequence of images which cover the entire cardiac cycle, also referred to as cardiac cine MRI. In clinical practice, cardiac cine MRI can be used to assess the cardiac function, see e.g.\ \cite{von2017representation}. Due to the structure of the problem, the temporal dimension is the one which offers the greatest potential to exploit the sparsity of the image in its gradient domain. However, a careful choice of the regularization parameter-map is required to ensure that the cardiac motion as well as smaller diagnostic image details are well-preserved after the reconstruction.
\subsubsection{Problem Formulation}
 For a complex-valued dynamic 2D MR image with vector representation $\XX \in \C^N$ with $n=n_x\cdot n_y\cdot n_t$, the  forward operator in \eqref{eq:forward_problem} is given as 
\begin{equation}\label{eq:operator_mri}
    \Ad:= (\Id_{n_c} \otimes \Ed) \Cd,
\end{equation}
where $\Id_{n_c}$ denotes the $n_c\times n_c$-sized identity-operator with $n_c$ being the number of receiver coils used for the data acquisition. The operator $\Cd$ is a tall operator which contains the coil-sensitivity maps, i.e. $\Cd = [\Cd_1, \ldots, \Cd_{n_c}]^\trans$ with $\Cd_k = \diag(\CC_k)$ and $\CC_k \in \C^n$, $k=1,\ldots,n_c$. Let $\Ed_{I_t}:=\Sd_{I_t} \Fd$  be an operator which acquires the $k$-space data of a static 2D MR image $\XX_t$ at time-point $t$ by sampling the $k$-space coefficients indexed by the set $I_t \subset J$, where $J=\{1,\ldots, n_{xy}\}$ with $n_{xy}:=n_x \cdot n_y$. 
Thereby, the mask $\Sd_{I_t}\in \{0,1\}^{m_t \times n_{xy}}$ with $m_t< n_{xy}$ for all $t=1,\ldots,n_t$ models the undersampling process. Undersampling the Fourier-space data is employed in order to accelerate the data acquisition process which usually takes place during a breathhold of the patient. Finally, the encoding operator $\Ed$ is given by 
\begin{equation}
\Ed:=
    \begin{pmatrix}
    \Ed_{I_1} & \mathbf{0}_{m_1\times n_{xy}} & \mathbf{0}_{m_1\times n_{xy}} & \cdots & \mathbf{0}_{m_1\times n_{xy}} \\
    \mathbf{0}_{m_2\times n_{xy}} & \Ed_{I_2}& \mathbf{0}_{m_2\times n_{xy}} & \cdots & \mathbf{0}_{m_2\times n_{xy}}  \\
    \vdots & \vdots & \vdots &   & \vdots &  \\
    \mathbf{0}_{m_{n_t}\times n_{xy}} & \mathbf{0}_{m_{n_t}\times n_{xy}} & \mathbf{0}_{m_{n_t}\times n_{xy}} & \cdots & \Ed_{I_{n_t}}
    \end{pmatrix},
\end{equation}
where $\mathbf{0}_{m_t \times n_{xy} } \in \{0\}^{m_t \times n_{xy}}$ denotes a $m_t \times n_{xy}$-sized zero-matrix. Typically, the number of receiver coils $n_c$ is chosen to ensure that $m:=n_c\cdot (m_1 + \ldots + m_t)>n$, i.e.\ problem \eqref{eq:forward_problem} is overdetermined when \eqref{eq:operator_mri} is the forward model.\\

\subsubsection{PDHG for Dynamic Multi-Coil MRI}

For the sake of completeness, we briefly summarize the PDHG-algorithm based on the identification mentioned in \eqref{eq:identification3}.
Recall the definition of $f_2$ from   \eqref{eq:identification3}. Since 
\begin{equation}\label{eq:clipping_op}
    \big(\mathrm{prox}_{\tau f_2^\ast}(\QQ)\big)_i =
    \left \{
    \begin{aligned}
    - & (\LLambda_{\Theta})_i, \quad  &\QQ_i &\in \big(-\infty, -(\LLambda_{\Theta})_i \big) \\
     & \QQ_i,  \quad &\QQ_i &\in \big[-(\LLambda_{\Theta})_i , (\LLambda_{\Theta})_i ) \big] \\
     & (\LLambda_{\Theta})_i, &\QQ_i &\in \big( (\LLambda_{\Theta})_i, \infty \big)
    \end{aligned}
    \right. ,
\end{equation}
the proximal operator $\mathrm{prox}_{\tau f_2^\ast}$ acts by ``clipping" each entry in the vector $\QQ$ if its magnitude exceeds the corresponding entry in $\boldsymbol{\Lambda}_{\Theta}$ and we therefore abbreviate it as $\mathrm{prox}_{\tau f_2^\ast}:=\mathrm{clip}_{\boldsymbol{\Lambda}_{\Theta}}$ to emphasize its dependence on the regularization parameter-map $\LLambda_{\Theta}$. The algorithm is summarized in Algorithm \ref{algo:tv_reco_algo_mri}.

\begin{algorithm}
\caption{Unrolled PDHG algorithm  for general linear operator $\Ad$ with $d(\,\cdot \, , \, \cdot \,) = \frac{1}{2}\| \cdot - \cdot\|_2^2$ and \textit{fixed} regularization parameter-map $\boldsymbol{\Lambda}_{\Theta}$ (adapted from \cite{sidky2012convex})}\label{algo:tv_reco_algo_mri}
  \begin{algorithmic}[1]
  \INPUT $L =\| [\Ad, \nabla]^\trans \|$, \; $\tau = 1/L$, \; $\sigma = 1/L$,  \; $\theta = 1$, \;$\text{initial guess}~\XX_0$
   \OUTPUT reconstructed image $\XX_{\mathrm{TV}}$ 
  \STATE $\bar{\XX_0} = \XX_0$
  \STATE $\PP_0 = \mathbf{0}$
  \STATE $\QQ_0 = \mathbf{0}$
    \FOR {$k < T$ }
    \STATE $\PP_{k+1} = ( \PP_k + \sigma ( \Ad \bar{\XX}_k - \YY) / ( 1 + \sigma)$ 
     \STATE $\QQ_{k+1}  = \mathrm{clip}_{\boldsymbol{\Lambda}_{\Theta}}(\QQ_k + \sigma \nabla \bar{\XX}_k)$  
     \STATE $\XX_{k+1} = \XX_k - \tau \Ad^\herm \PP_{k+1} - \tau \nabla^\trans \QQ_{k+1} $
     \STATE $\bar{\XX}_{k+1} = \XX_{k+1} + \theta ( \XX_{k+1} - \XX_k)$ 
    \ENDFOR
   \STATE $\XX_{\mathrm{TV}} = \XX_T$
  \end{algorithmic}
\end{algorithm}

\subsubsection{Experimental Set-Up}
\label{subsubsec:experiments:cardiac}
We used a set of 216 cardiac cine MR images of the study \cite{kolbitsch2014cardiac} which we split in portion of 144/36/36 for training, validation and testing. The images have shape $n_x\times n_y\times n_t = 160 \times 160 \times 30$ and a resolution of $2 \times 2$ mm$^2$ with a slice thickness of 8 mm$^2$. The number of receiver coils is $n_c=12$. We retrospectively simulated $k$-space data according to \eqref{eq:forward_problem} using the model in  \eqref{eq:operator_mri} as the forward operator simulating acceleration factors of $R=4,6,8$ with complex-valued Gaussian noise with standard deviation $\sigma = 0.15, 0.30, 0.45$.\\
As described in Section \ref{subsubsec:obtaining_lambda}, we constructed $u_{\Theta}$ such that it yields two different parameter-maps. One for the spatial $x$- and $y$-directions and one for the temporal direction, i.e., $\boldsymbol{\Lambda}_{\Theta}:=(\boldsymbol{\Lambda}_{\Theta}^{xy}, \boldsymbol{\Lambda}_{\Theta}^{xy},\boldsymbol{\Lambda}_{\Theta}^{t})$.
The CNN $u_{\Theta}$ here corresponds to a 3D U-Net with two input-channels (for the real and the imaginary part of the image, respectively), three encoding stages, two convolutional layers per stage and an initial number of eight filters which are applied to the input image. As in Figure \ref{fig:u_Theta.pdf}, the last layer consists of a $1\times1\times 1$ convolution with two output channels (the first for the parameter-map for the $x$- and $y$-directions, the second for the parameter-map for the $t$-direction) and the softplus activation function $\phi$. Note that the gradients of the real and the imaginary parts of the images share the same regularization parameter-map. The scaling factor $t$ in \eqref{eq:lambda_reg_map} was set to $t=0.1$. The overall number of trainable parameters of $u_{\Theta}$ is 97\,290.
To reduce training times, the network was trained on patches of shape $n_x^\prime \times n_y^\prime \times n_t^\prime = 160 \times 160 \times 16$. The network's number of overall iterations was set to $T=256$ during training, while at test time, we used $T=4096$ iterations. The reason for the different number of iterations at training and test time is discussed later in Subsection \ref{subsec:choosing_T}. The parameters $\sigma, \tau$ and $\theta$ were trained as well and constrained to be in the intervals $(0,1/L), (0,1/L)$ and $(0,1)$, respectively, by using a sigmoid activation-function. Despite of the training, we mention that no noteworthy changes were visible after training, i.e. $\sigma \approx \tau \approx 1/L$. Not training $\sigma, \tau$ and $\theta$ also led to similar results as the ones shown later. As training routine, we used the Adam optimizer \cite{kingma2014adam} with initial learning rate of $10^{-4}$ to minimize the mean squared error (MSE) between the reconstructed image and the target image. We trained all networks for 200 epochs while evaluating the network 25 times over the entire training and validation datasets. We then used the model configuration for which the MSE on the validation set was the lowest. 

\subsubsection{Results}
\begin{figure}
\centering
\begin{minipage}{\linewidth}
    \begin{minipage}{\linewidth}
\hspace{0.8cm} PDHG $\lambda^{xyt}_{\tilde{\mathrm{P}}}$ \hspace{0.4cm} PDHG $\lambda^{xyt}_{\mathrm{P}}$ \hspace{0.4cm}  PDHG $\lambda^{xy,t}_{\tilde{\mathrm{P}}}$ \hspace{0.4cm} PDHG $\lambda^{xy,t}_{\mathrm{P}}$ \hspace{0.4cm} PDHG $\LLambda^{xy,t}_\Theta$ \hspace{0.4cm} Target/ZF
    \end{minipage} 
    \begin{minipage}{\linewidth}
    \rotatebox{90}{
        \begin{minipage}{0.4\linewidth}
        \hspace{-4.6cm }$R=8$ \hspace{2.9cm} $R=6$ \hspace{3.cm} $R=4$
        \end{minipage}
    }
     \begin{minipage}{\linewidth}
    \begin{minipage}{\linewidth}
    \resizebox{\linewidth}{!}{
    \includegraphics[height=3cm]{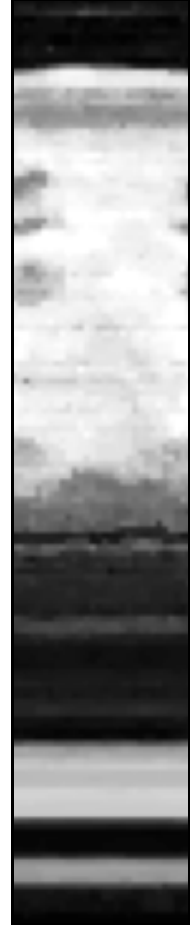}\hspace{-0.1cm}
    \begin{tikzpicture}[spy using outlines={rectangle, white, magnification=2, size=1.cm, connect spies}]
    \node[anchor=south west,inner sep=0]  at (0,0) {\includegraphics[height=3cm]{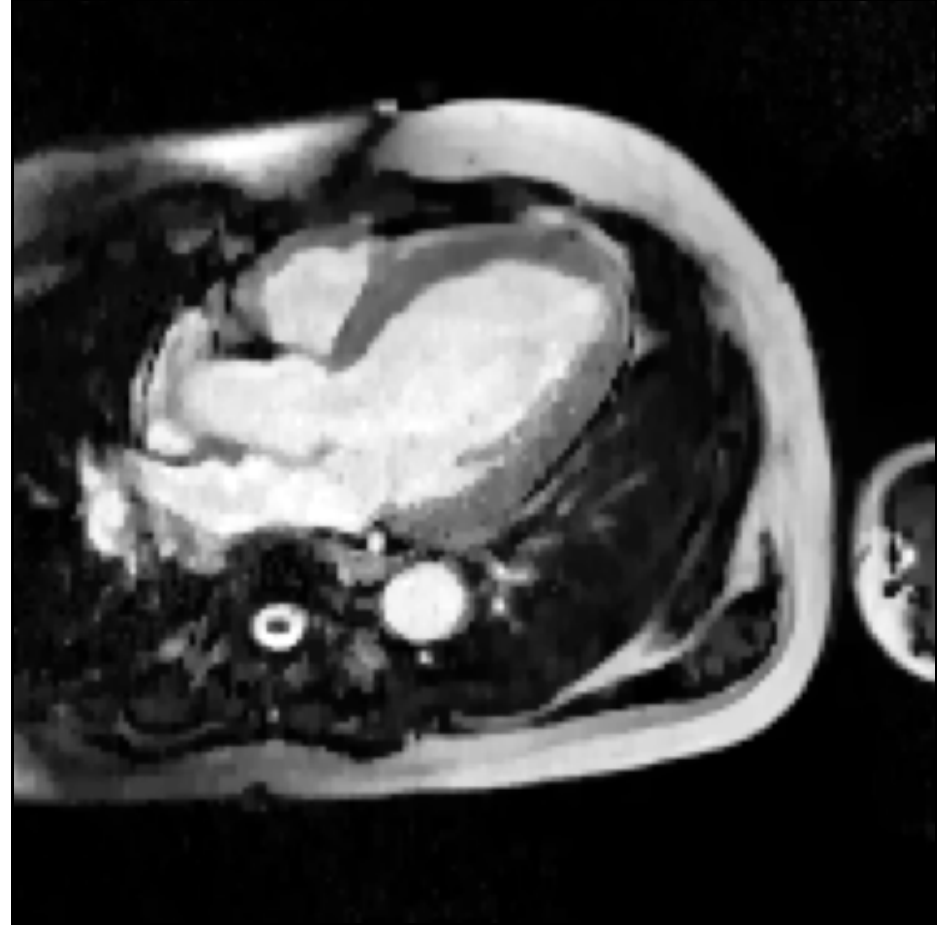}};
    \spy on (0.8, 1.5) in node [left] at (3.0, 2.5);
    \end{tikzpicture}
    \includegraphics[height=3cm]{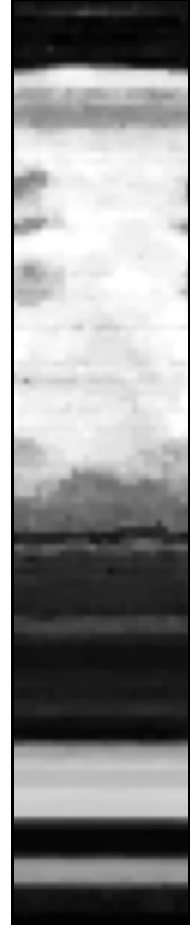}\hspace{-0.1cm}
    \begin{tikzpicture}[spy using outlines={rectangle, white, magnification=2, size=1.cm, connect spies}]
    \node[anchor=south west,inner sep=0]  at (0,0) {\includegraphics[height=3cm]{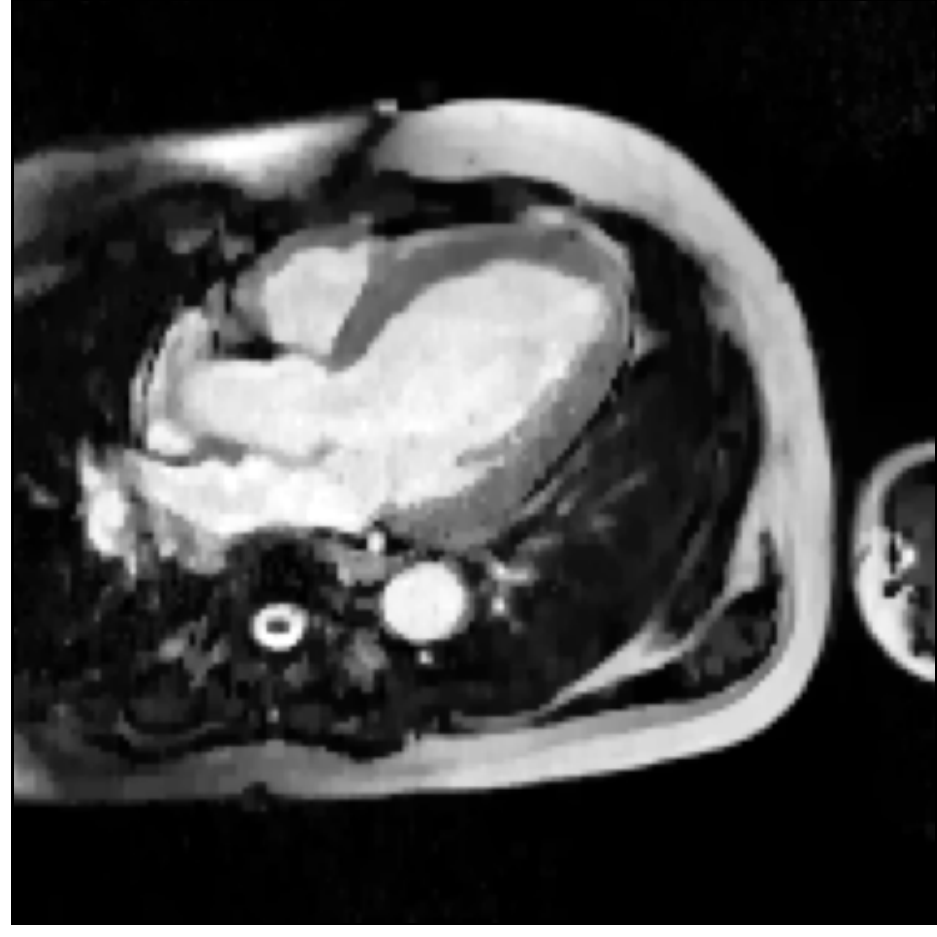}};
    \spy on (0.8, 1.5) in node [left] at (3.0, 2.5);
    \end{tikzpicture}
    \includegraphics[height=3cm]{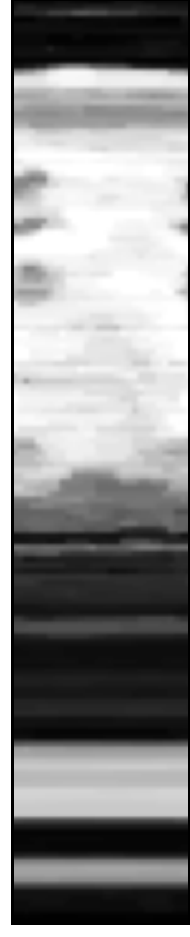}\hspace{-0.1cm}
    \begin{tikzpicture}[spy using outlines={rectangle, white, magnification=2, size=1.cm, connect spies}]
    \node[anchor=south west,inner sep=0]  at (0,0) {\includegraphics[height=3cm]{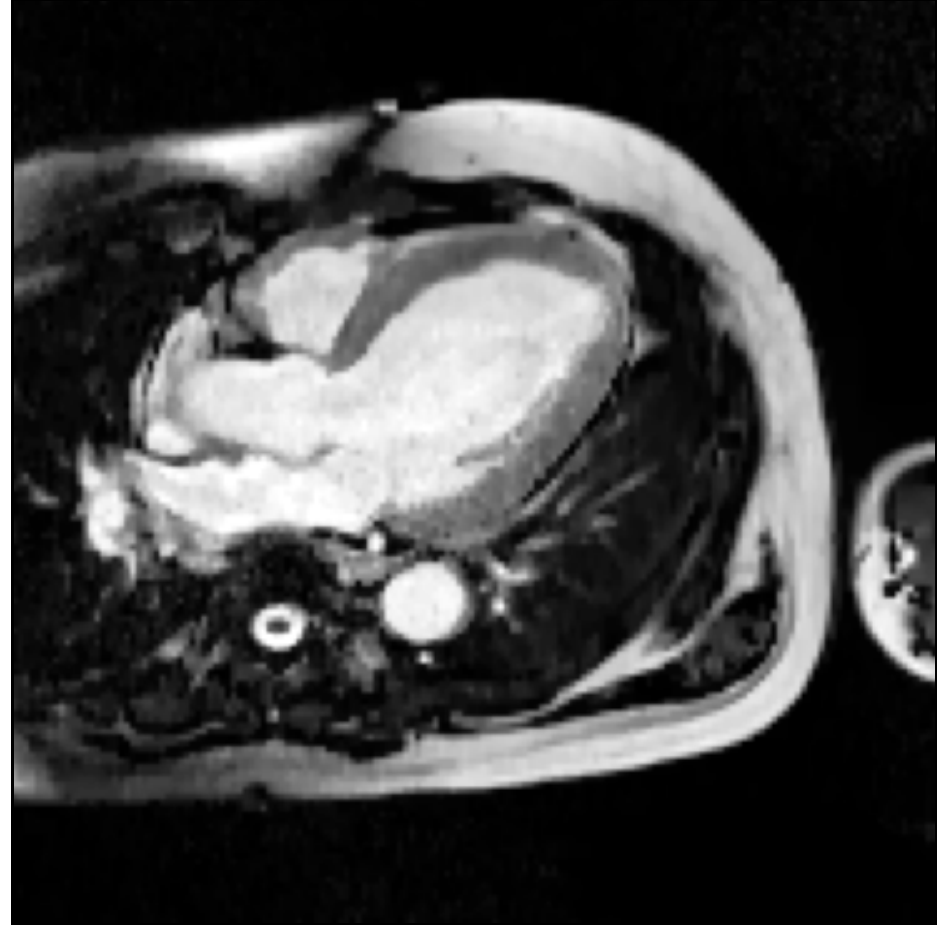}};
    \spy on (0.8, 1.5) in node [left] at (3.0, 2.5);
    \end{tikzpicture}
    \includegraphics[height=3cm]{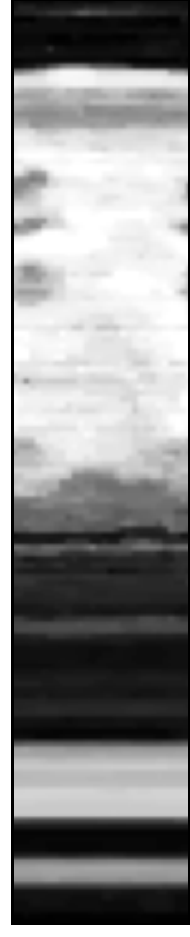}\hspace{-0.1cm}
    \begin{tikzpicture}[spy using outlines={rectangle, white, magnification=2, size=1.cm, connect spies}]
    \node[anchor=south west,inner sep=0]  at (0,0) {\includegraphics[height=3cm]{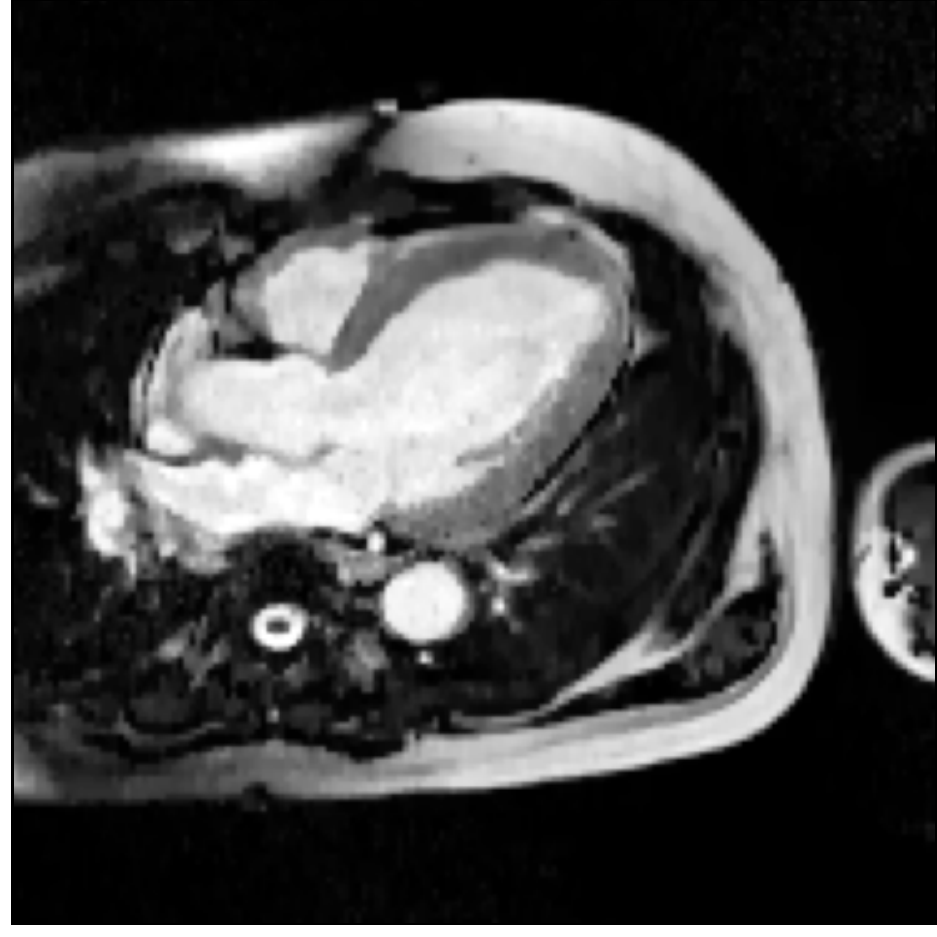}};
    \spy on (0.8, 1.5) in node [left] at (3.0, 2.5);
    \end{tikzpicture}
    \includegraphics[height=3cm]{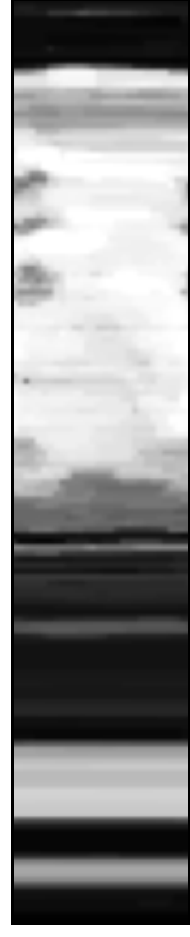}\hspace{-0.1cm}
    \begin{tikzpicture}[spy using outlines={rectangle, white, magnification=2, size=1.cm, connect spies}]
    \node[anchor=south west,inner sep=0]  at (0,0) {\includegraphics[height=3cm]{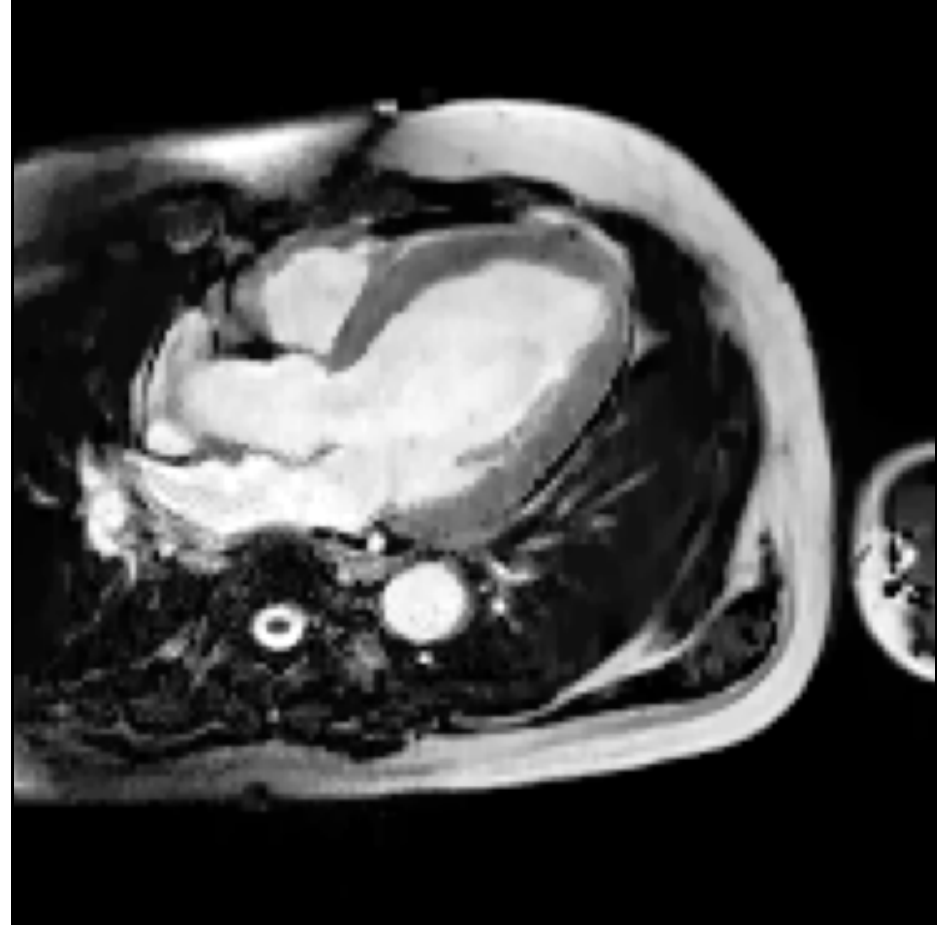}};
    \spy on (0.8, 1.5) in node [left] at (3.0, 2.5);
    \end{tikzpicture}
    \includegraphics[height=3cm]{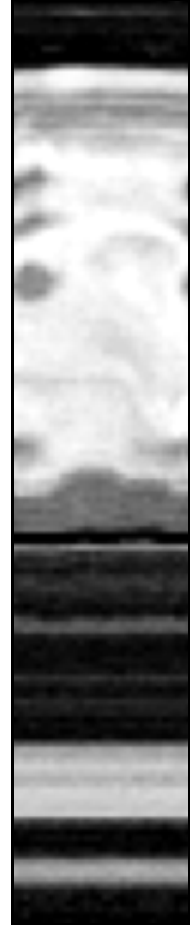}\hspace{-0.1cm}
    \begin{tikzpicture}[spy using outlines={rectangle, white, magnification=2, size=1.cm, connect spies}]
    \node[anchor=south west,inner sep=0]  at (0,0) {\includegraphics[height=3cm]{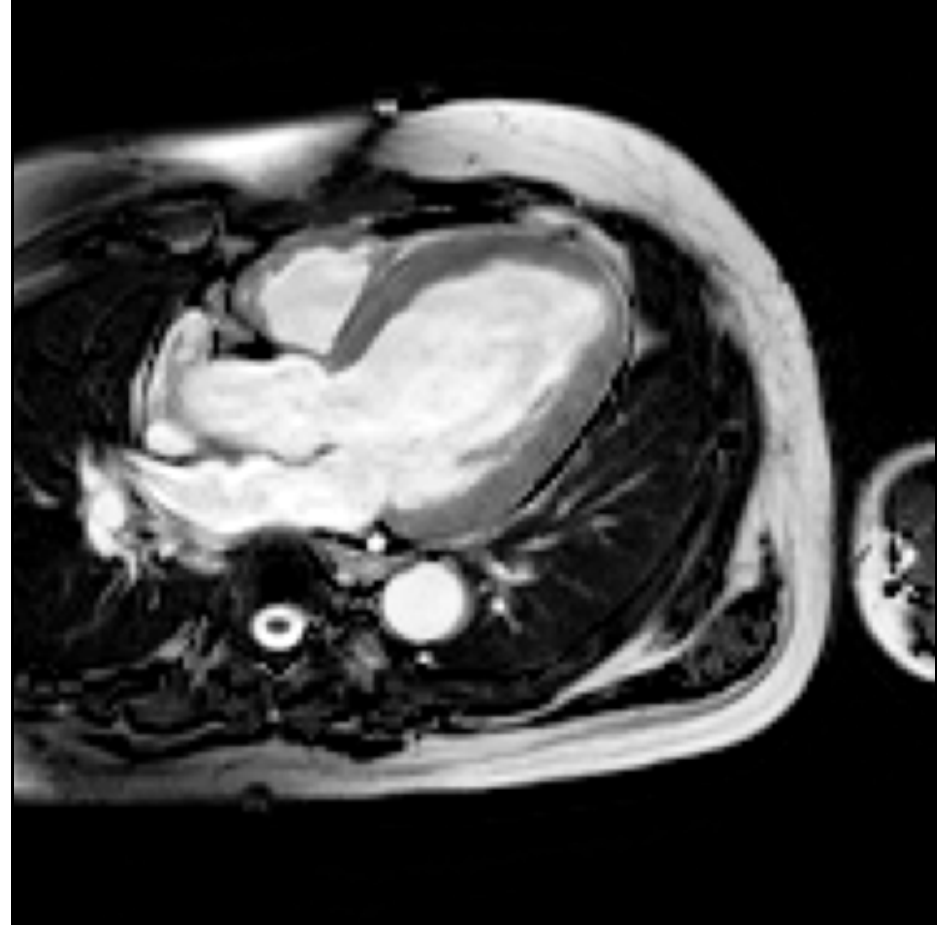}};
    \spy on (0.8, 1.5) in node [left] at (3.0, 2.5);
    \end{tikzpicture}
    }
    \resizebox{\linewidth}{!}{
    \includegraphics[height=3cm]{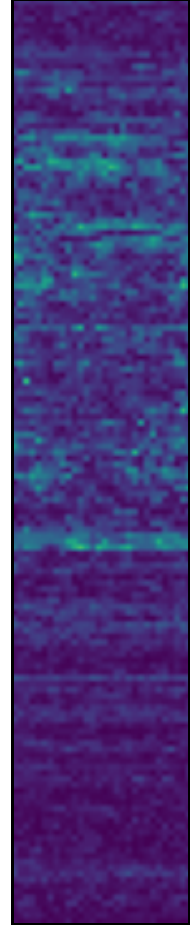}\hspace{-0.1cm}
    \begin{tikzpicture}[spy using outlines={rectangle, white, magnification=2, size=1.cm, connect spies}]
    \node[anchor=south west,inner sep=0]  at (0,0) {\includegraphics[height=3cm]{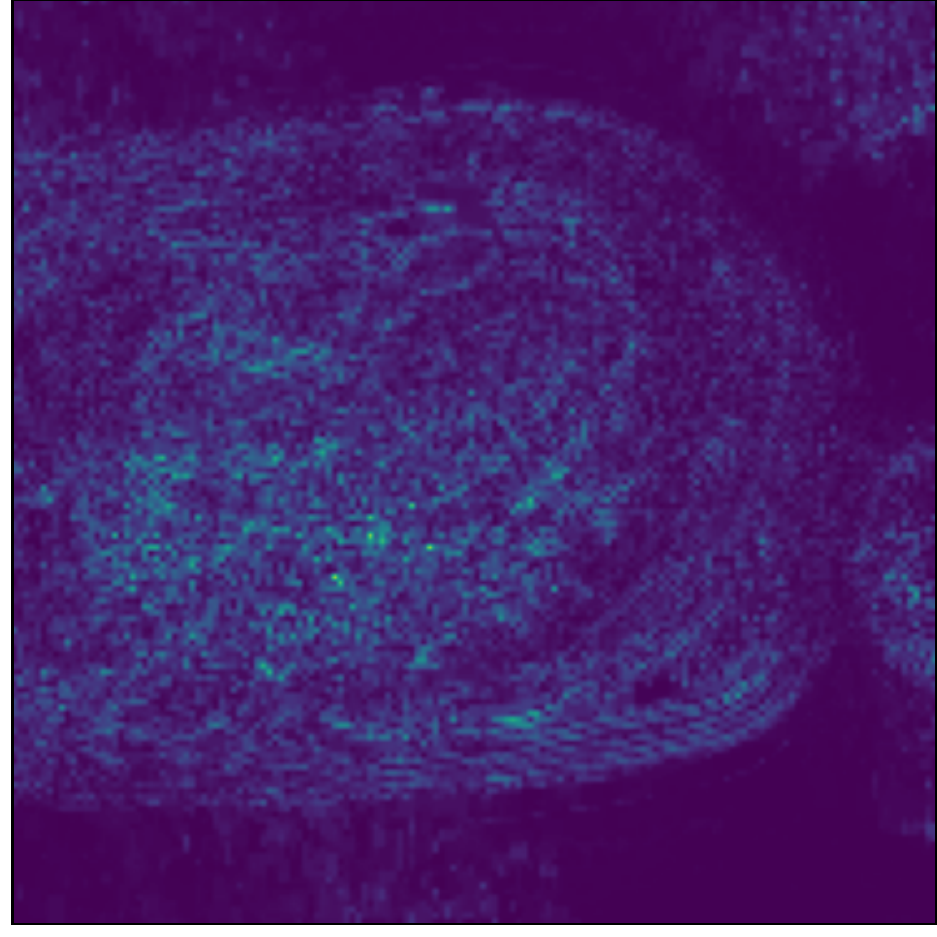}};
    \spy on (0.8, 1.5) in node [left] at (3.0, 2.5);
    \end{tikzpicture}
    \includegraphics[height=3cm]{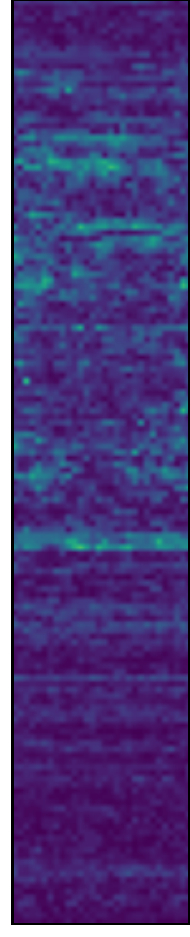}\hspace{-0.1cm}
    \begin{tikzpicture}[spy using outlines={rectangle, white, magnification=2, size=1.cm, connect spies}]
    \node[anchor=south west,inner sep=0]  at (0,0) {\includegraphics[height=3cm]{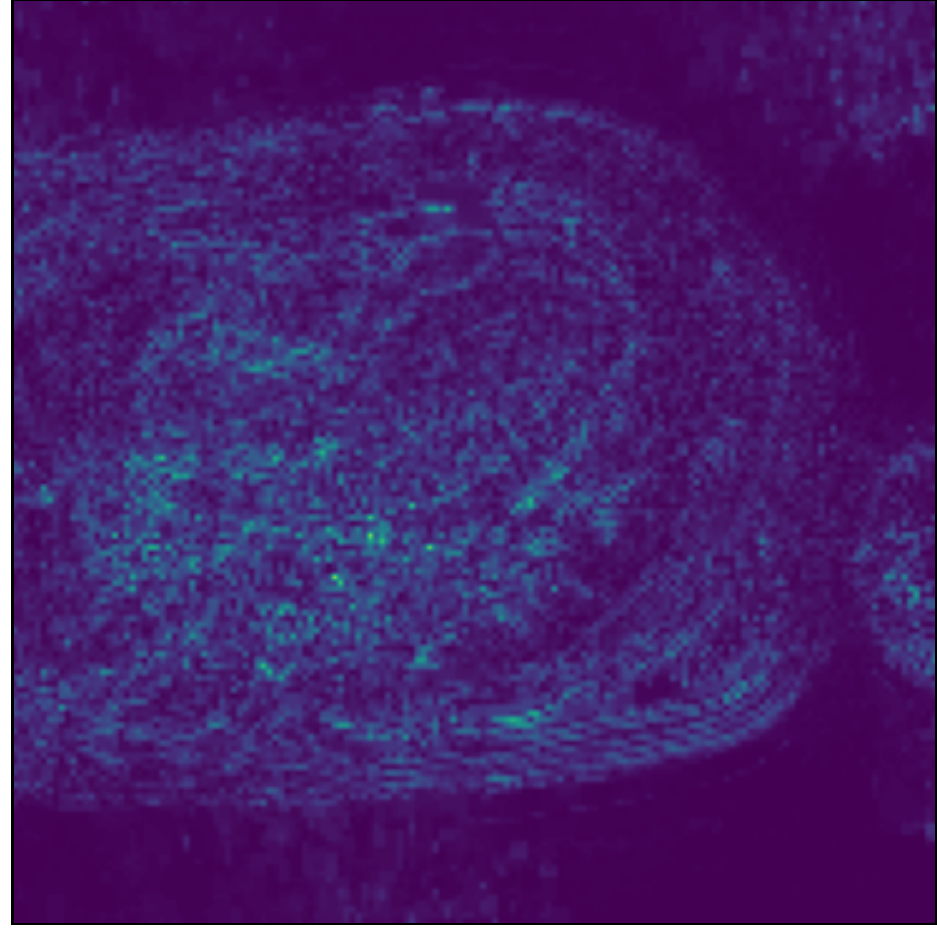}};
    \spy on (0.8, 1.5) in node [left] at (3.0, 2.5);
    \end{tikzpicture}
    \includegraphics[height=3cm]{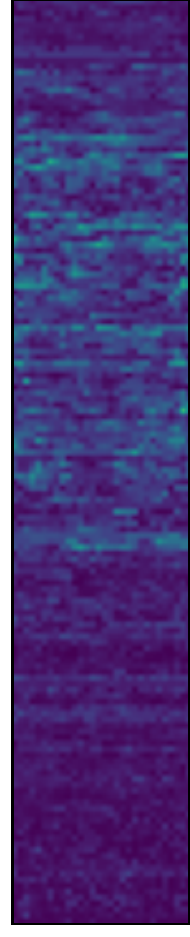}\hspace{-0.1cm}
    \begin{tikzpicture}[spy using outlines={rectangle, white, magnification=2, size=1.cm, connect spies}]
    \node[anchor=south west,inner sep=0]  at (0,0) {\includegraphics[height=3cm]{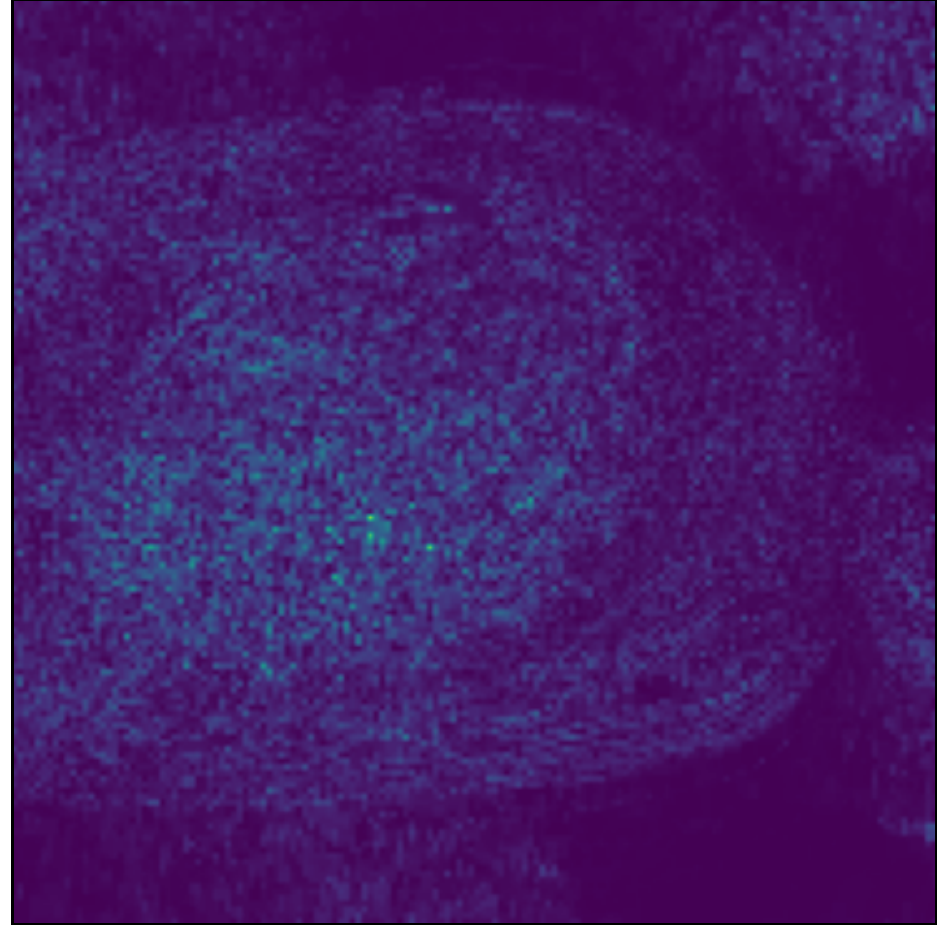}};
    \spy on (0.8, 1.5) in node [left] at (3.0, 2.5);
    \end{tikzpicture}
    \includegraphics[height=3cm]{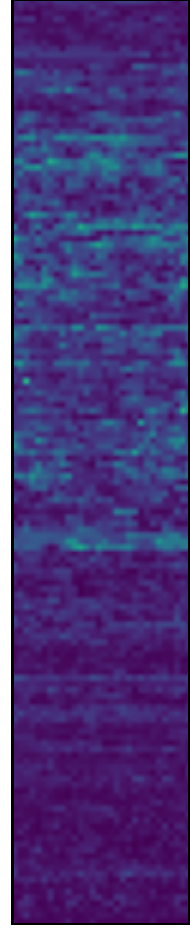}\hspace{-0.1cm}
    \begin{tikzpicture}[spy using outlines={rectangle, white, magnification=2, size=1.cm, connect spies}]
    \node[anchor=south west,inner sep=0]  at (0,0) {\includegraphics[height=3cm]{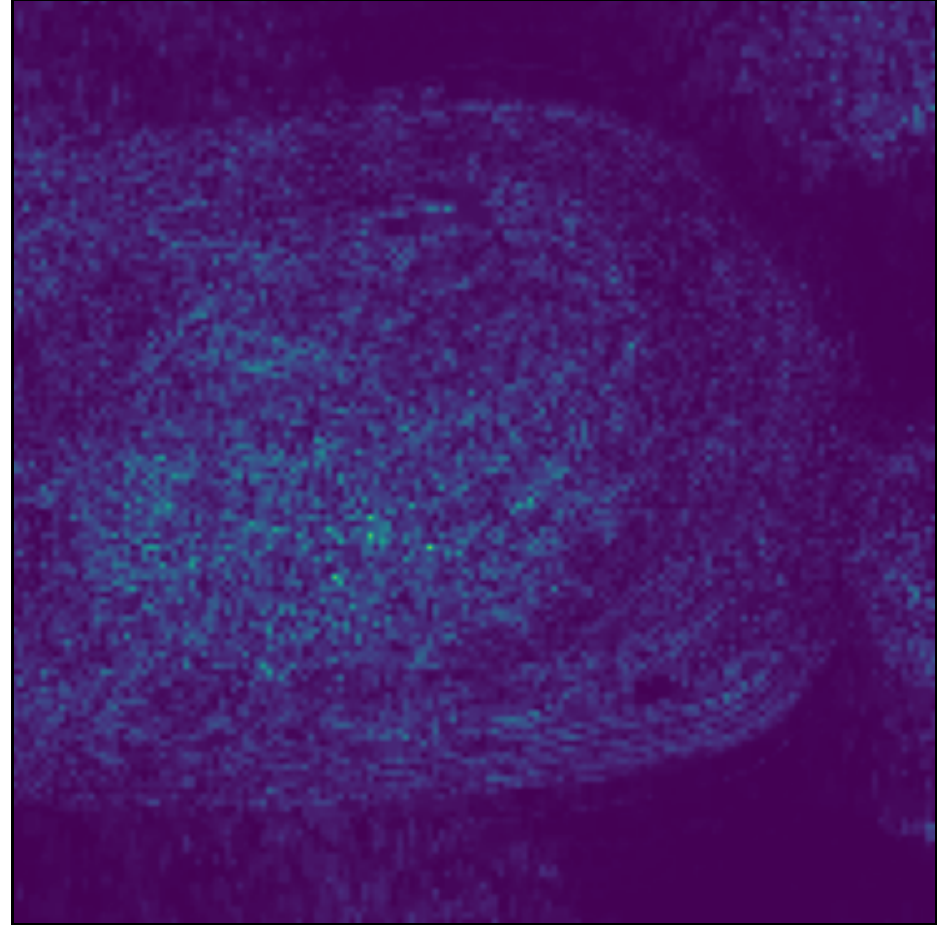}};
    \spy on (0.8, 1.5) in node [left] at (3.0, 2.5);
    \end{tikzpicture}
    \includegraphics[height=3cm]{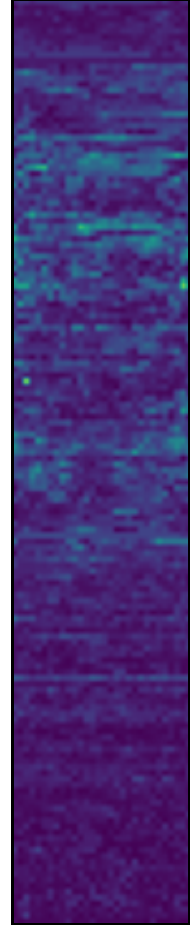}\hspace{-0.1cm}
    \begin{tikzpicture}[spy using outlines={rectangle, white, magnification=2, size=1.cm, connect spies}]
    \node[anchor=south west,inner sep=0]  at (0,0) {\includegraphics[height=3cm]{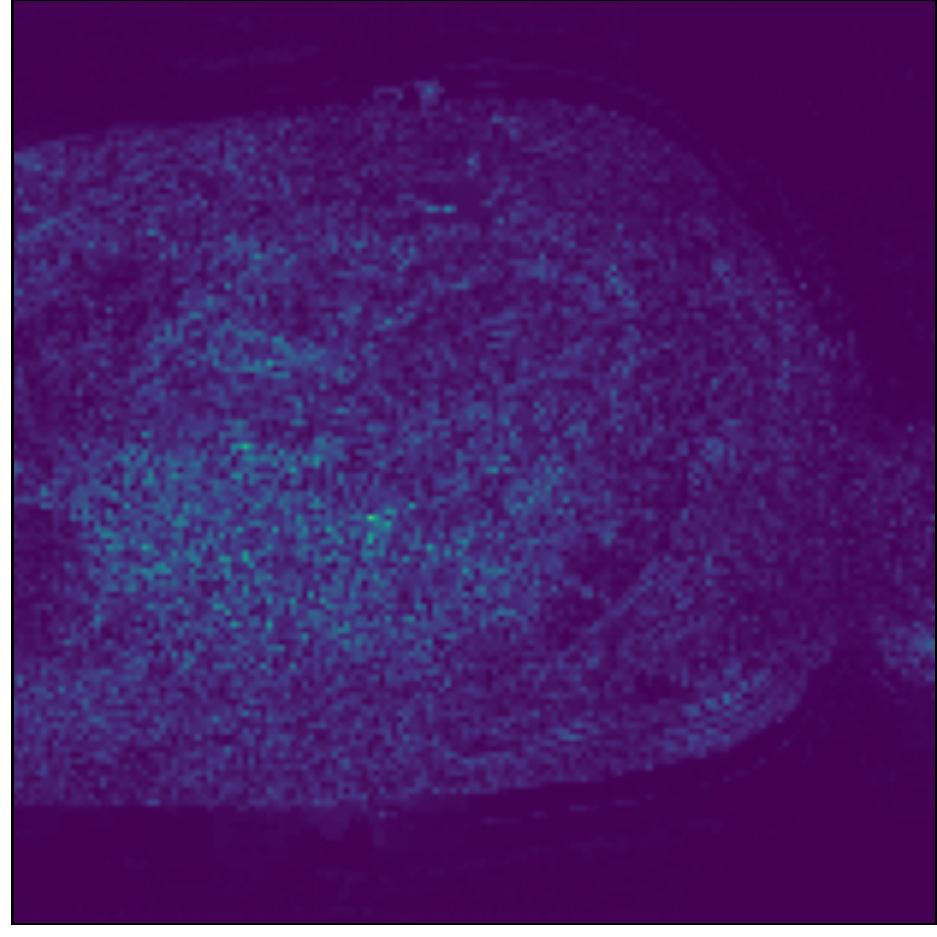}};
    \spy on (0.8, 1.5) in node [left] at (3.0, 2.5);
    \end{tikzpicture}
    \includegraphics[height=3cm]{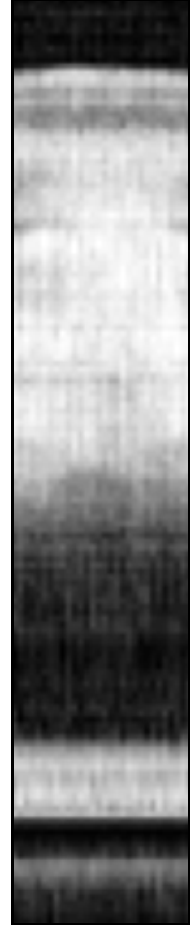}\hspace{-0.1cm}
    \begin{tikzpicture}[spy using outlines={rectangle, white, magnification=2, size=1.cm, connect spies}]
    \node[anchor=south west,inner sep=0]  at (0,0) {\includegraphics[height=3cm]{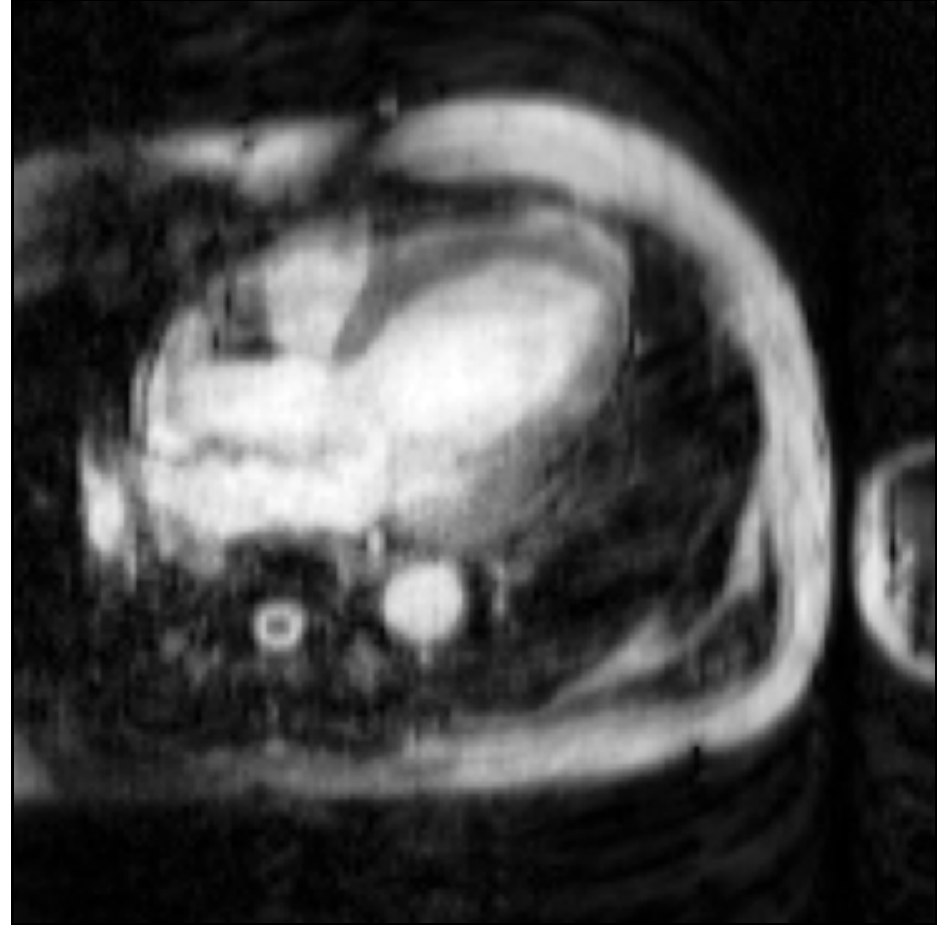}};
    \spy on (0.8, 1.5) in node [left] at (3.0, 2.5);
    \end{tikzpicture}
    }
    \end{minipage}

    \begin{minipage}{\linewidth}
    \resizebox{\linewidth}{!}{
    \includegraphics[height=3cm]{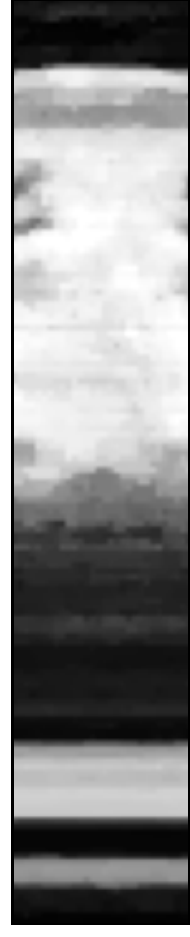}\hspace{-0.1cm}
    \begin{tikzpicture}[spy using outlines={rectangle, white, magnification=2, size=1.cm, connect spies}]
    \node[anchor=south west,inner sep=0]  at (0,0) {\includegraphics[height=3cm]{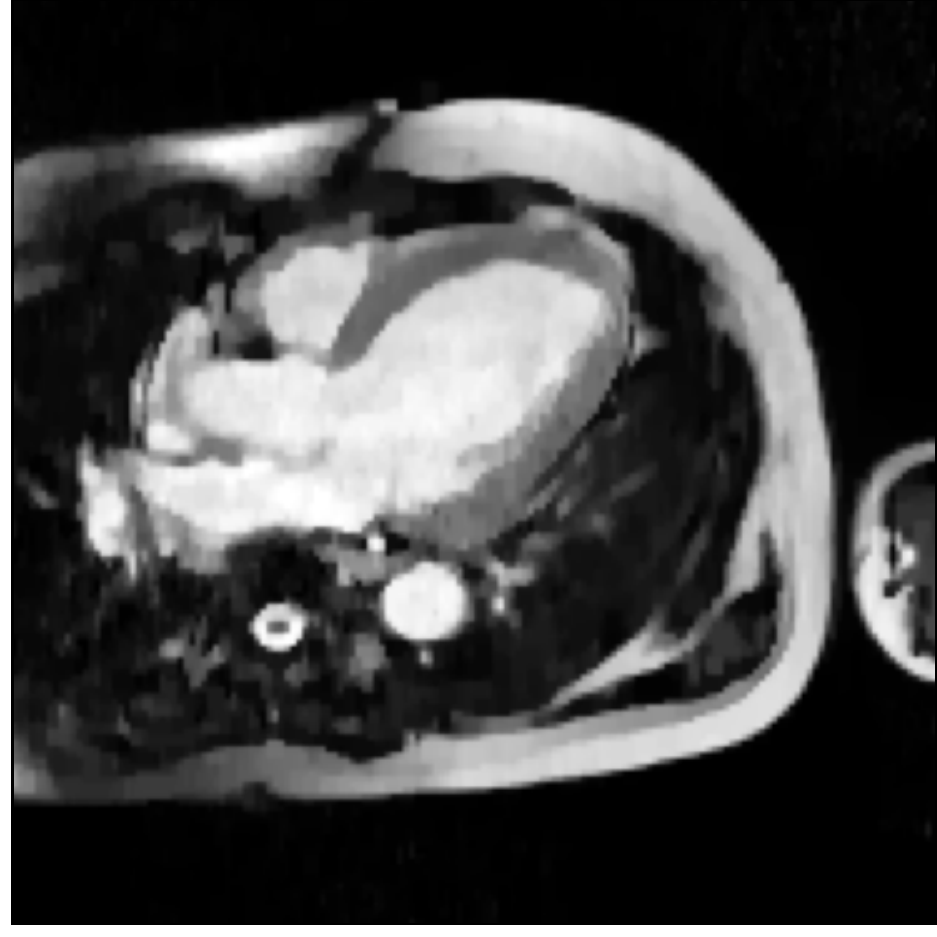}};
    \spy on (0.8, 1.5) in node [left] at (3.0, 2.5);
    \end{tikzpicture}
    \includegraphics[height=3cm]{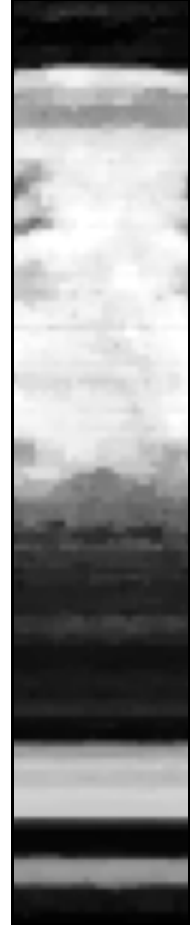}\hspace{-0.1cm}
    \begin{tikzpicture}[spy using outlines={rectangle, white, magnification=2, size=1.cm, connect spies}]
    \node[anchor=south west,inner sep=0]  at (0,0) {\includegraphics[height=3cm]{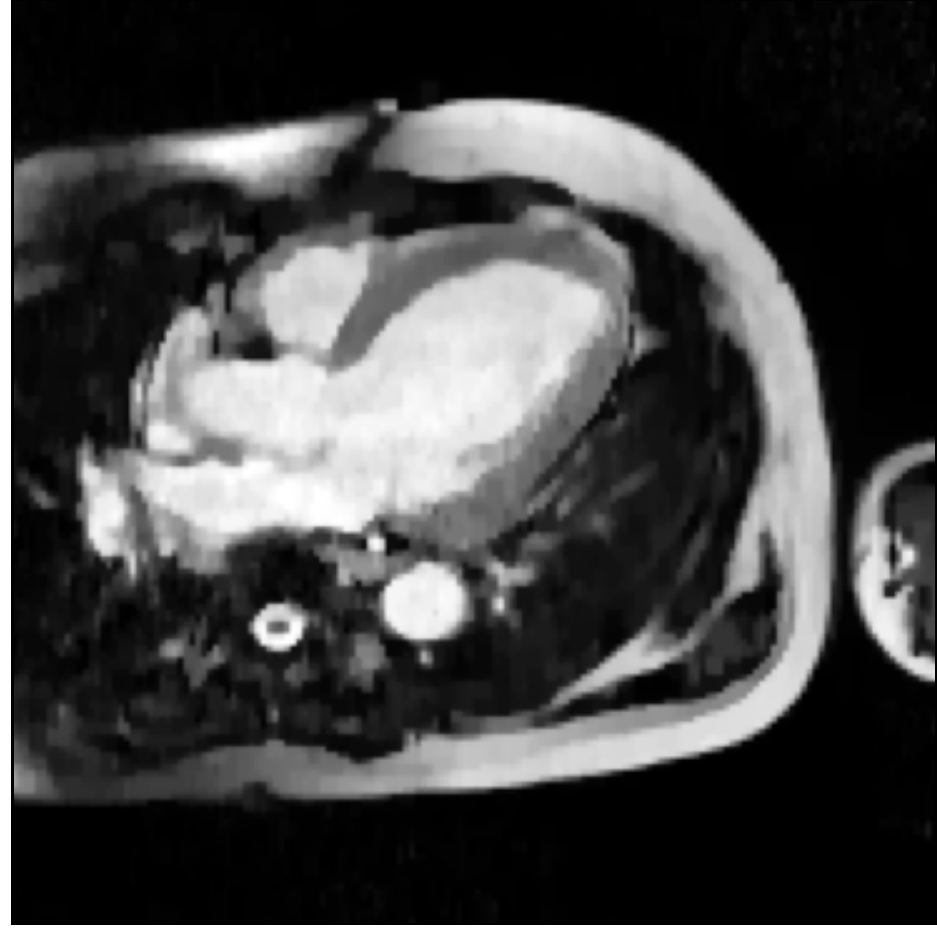}};
    \spy on (0.8, 1.5) in node [left] at (3.0, 2.5);
    \end{tikzpicture}
    \includegraphics[height=3cm]{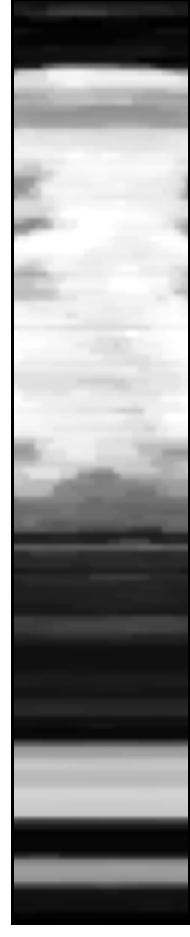}\hspace{-0.1cm}
    \begin{tikzpicture}[spy using outlines={rectangle, white, magnification=2, size=1.cm, connect spies}]
    \node[anchor=south west,inner sep=0]  at (0,0) {\includegraphics[height=3cm]{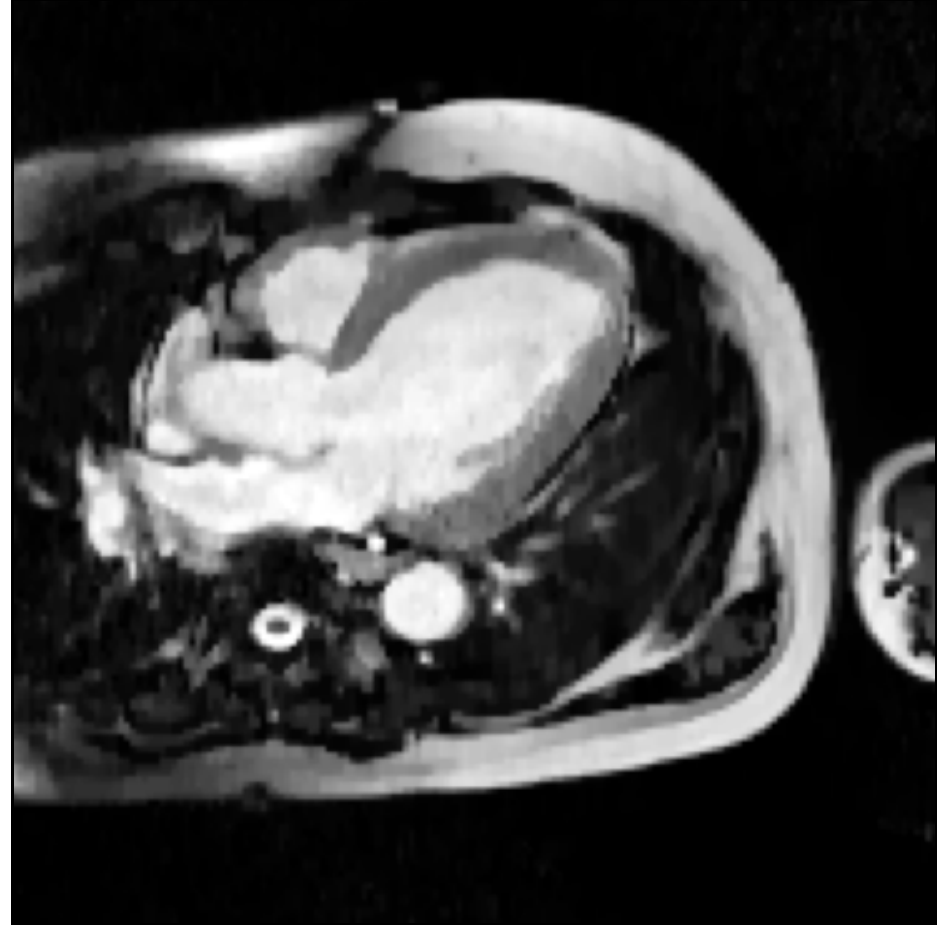}};
    \spy on (0.8, 1.5) in node [left] at (3.0, 2.5);
    \end{tikzpicture}
    \includegraphics[height=3cm]{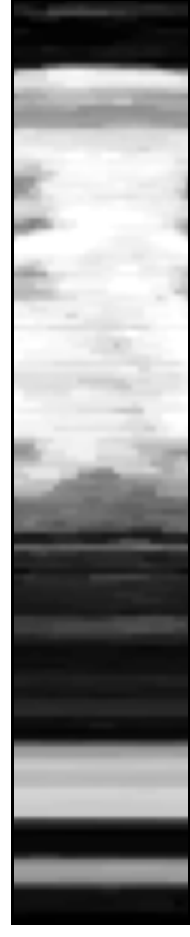}\hspace{-0.1cm}
    \begin{tikzpicture}[spy using outlines={rectangle, white, magnification=2, size=1.cm, connect spies}]
    \node[anchor=south west,inner sep=0]  at (0,0) {\includegraphics[height=3cm]{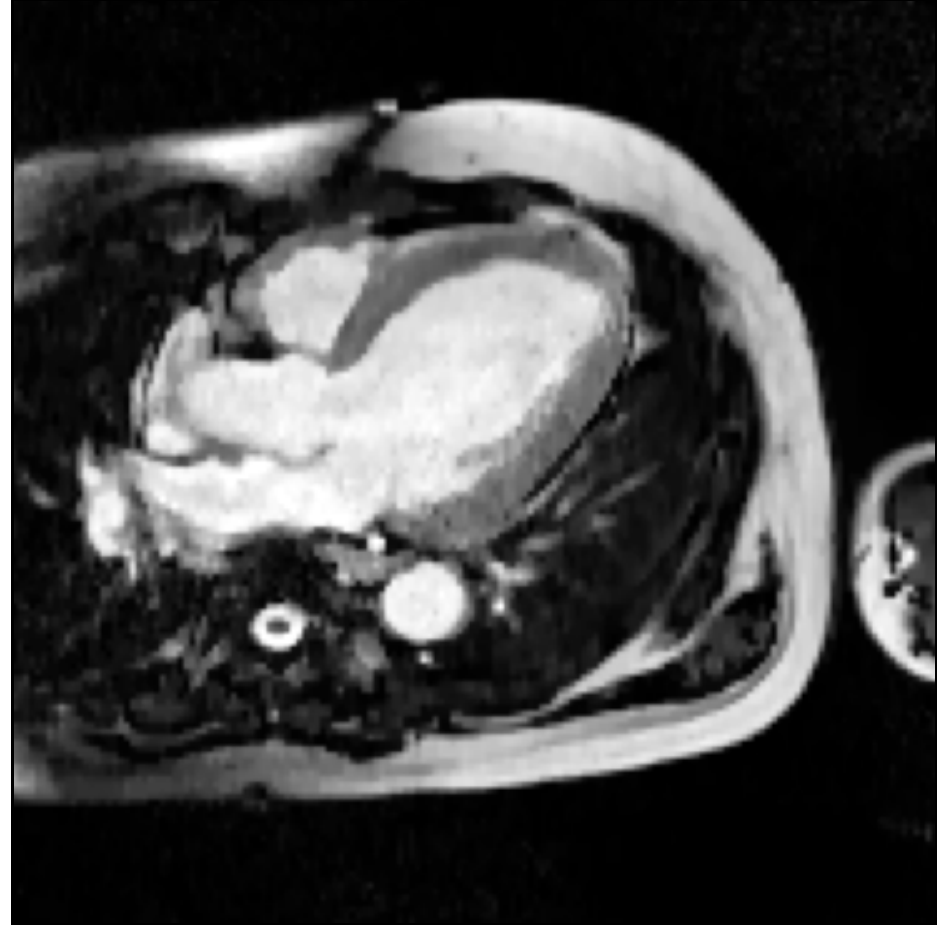}};
    \spy on (0.8, 1.5) in node [left] at (3.0, 2.5);
    \end{tikzpicture}
    \includegraphics[height=3cm]{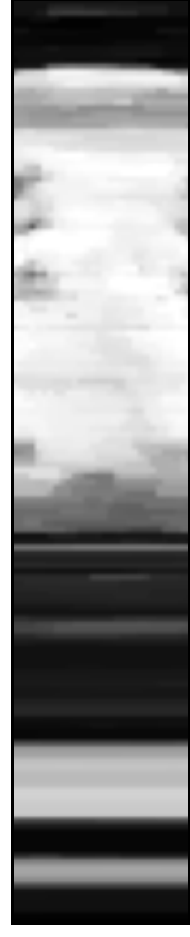}\hspace{-0.1cm}
    \begin{tikzpicture}[spy using outlines={rectangle, white, magnification=2, size=1.cm, connect spies}]
    \node[anchor=south west,inner sep=0]  at (0,0) {\includegraphics[height=3cm]{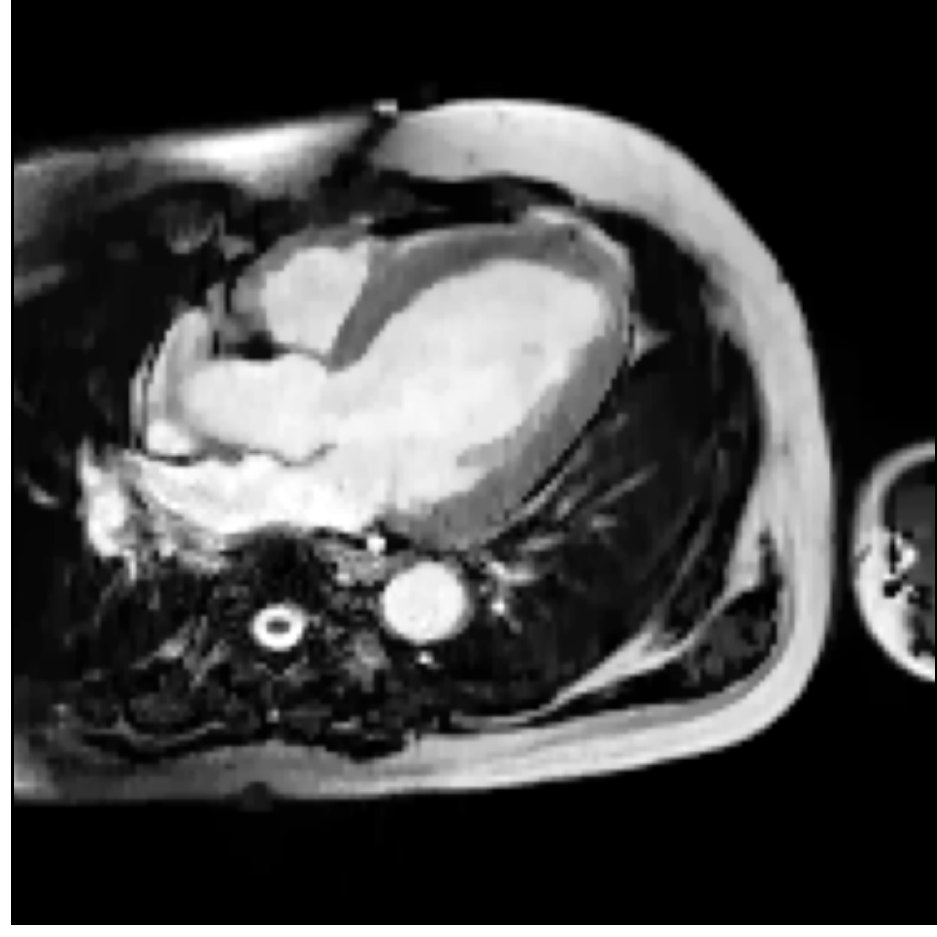}};
    \spy on (0.8, 1.5) in node [left] at (3.0, 2.5);
    \end{tikzpicture}
    \includegraphics[height=3cm]{figures/MRI/mri_results/xf_yt.pdf}\hspace{-0.1cm}
    \begin{tikzpicture}[spy using outlines={rectangle, white, magnification=2, size=1.cm, connect spies}]
    \node[anchor=south west,inner sep=0]  at (0,0) {\includegraphics[height=3cm]{figures/MRI/mri_results/xf_xy.pdf}};
    \spy on (0.8, 1.5) in node [left] at (3.0, 2.5);
    \end{tikzpicture}
    }
    \resizebox{\linewidth}{!}{
    \includegraphics[height=3cm]{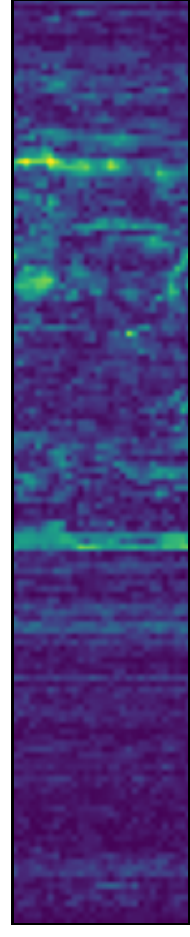}\hspace{-0.1cm}
    \begin{tikzpicture}[spy using outlines={rectangle, white, magnification=2, size=1.cm, connect spies}]
    \node[anchor=south west,inner sep=0]  at (0,0) {\includegraphics[height=3cm]{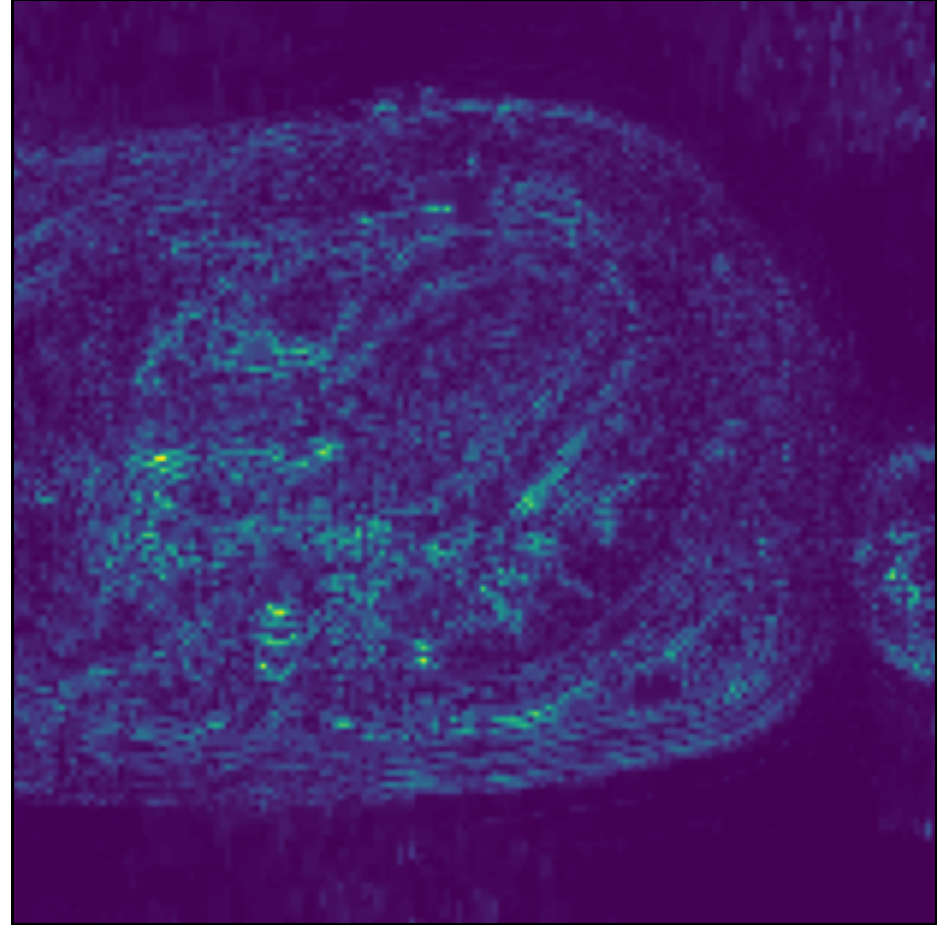}};
    \spy on (0.8, 1.5) in node [left] at (3.0, 2.5);
    \end{tikzpicture}
    \includegraphics[height=3cm]{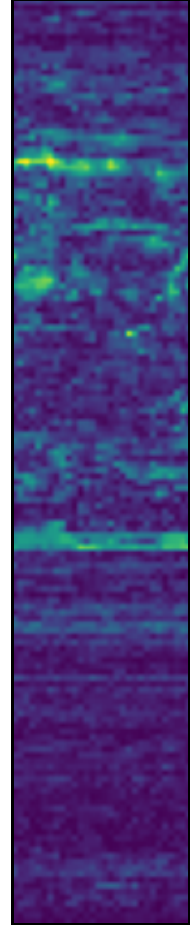}\hspace{-0.1cm}
    \begin{tikzpicture}[spy using outlines={rectangle, white, magnification=2, size=1.cm, connect spies}]
    \node[anchor=south west,inner sep=0]  at (0,0) {\includegraphics[height=3cm]{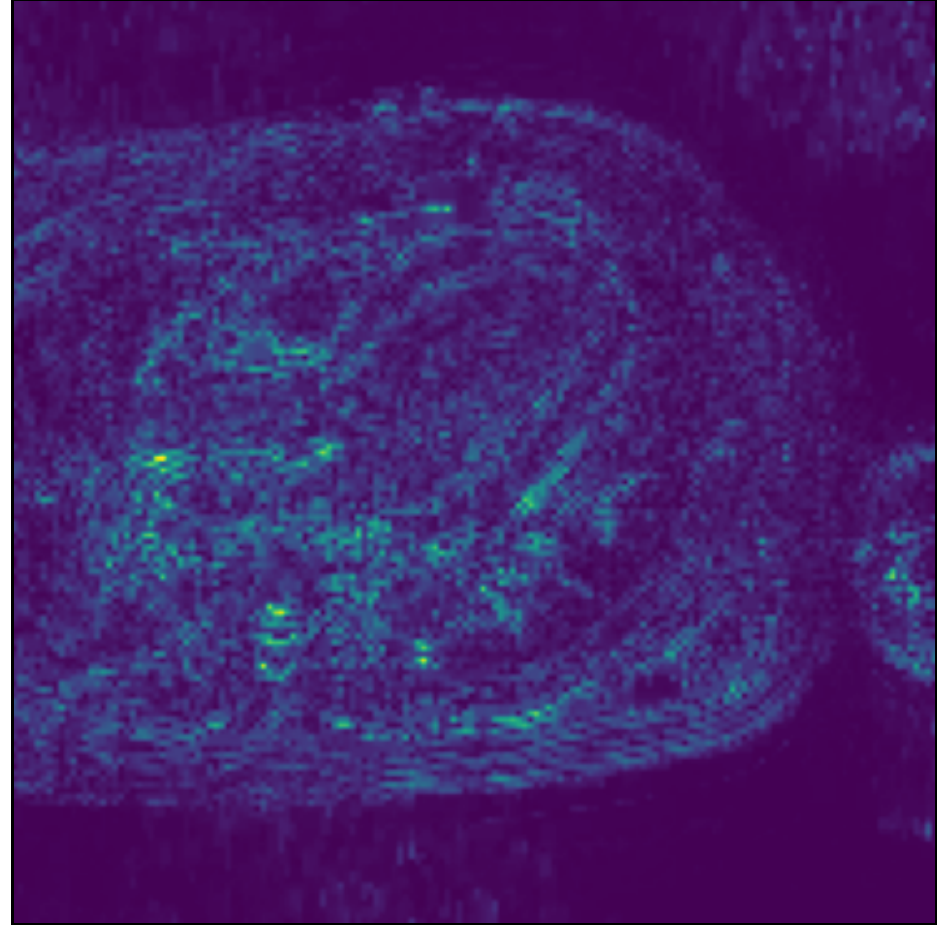}};
    \spy on (0.8, 1.5) in node [left] at (3.0, 2.5);
    \end{tikzpicture}
    \includegraphics[height=3cm]{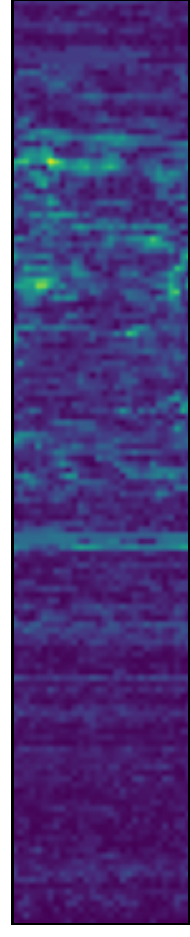}\hspace{-0.1cm}
    \begin{tikzpicture}[spy using outlines={rectangle, white, magnification=2, size=1.cm, connect spies}]
    \node[anchor=south west,inner sep=0]  at (0,0) {\includegraphics[height=3cm]{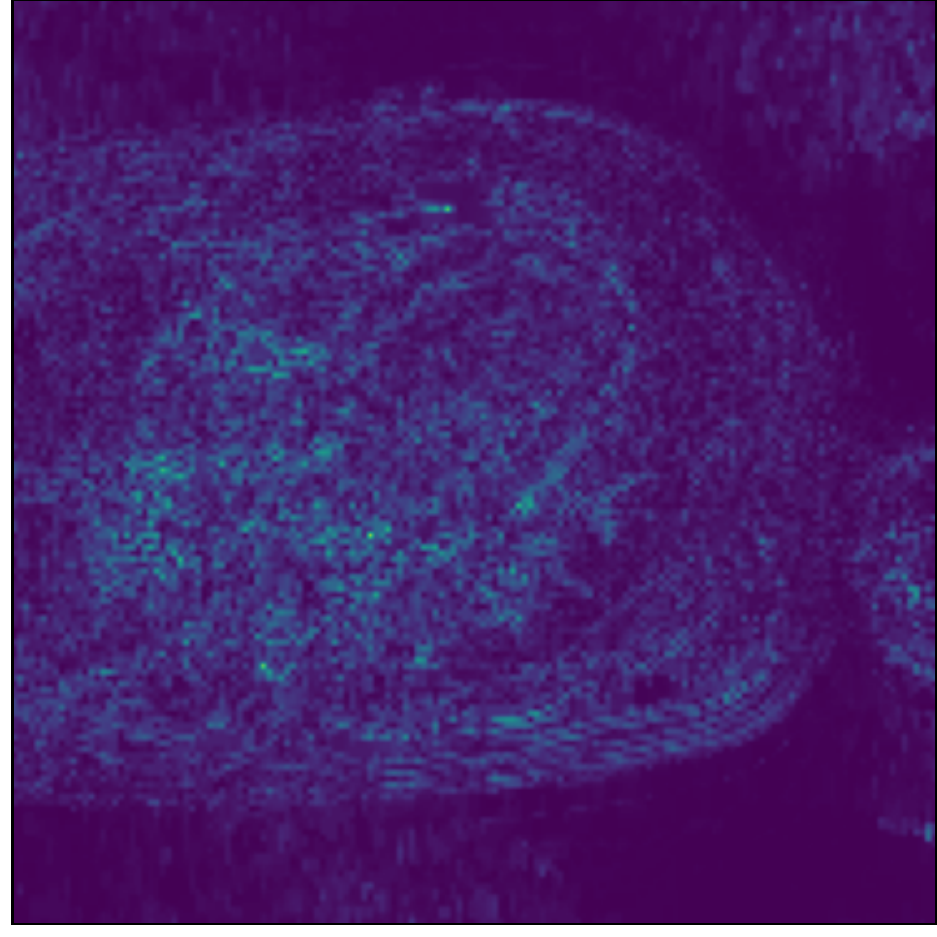}};
    \spy on (0.8, 1.5) in node [left] at (3.0, 2.5);
    \end{tikzpicture}
    \includegraphics[height=3cm]{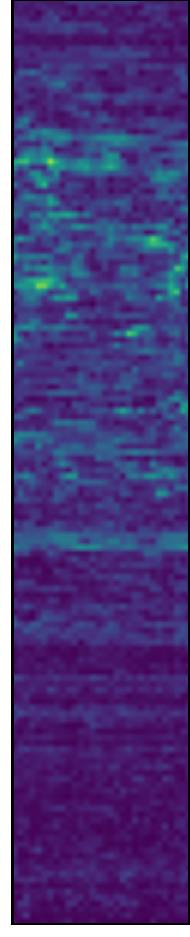}\hspace{-0.1cm}
    \begin{tikzpicture}[spy using outlines={rectangle, white, magnification=2, size=1.cm, connect spies}]
    \node[anchor=south west,inner sep=0]  at (0,0) {\includegraphics[height=3cm]{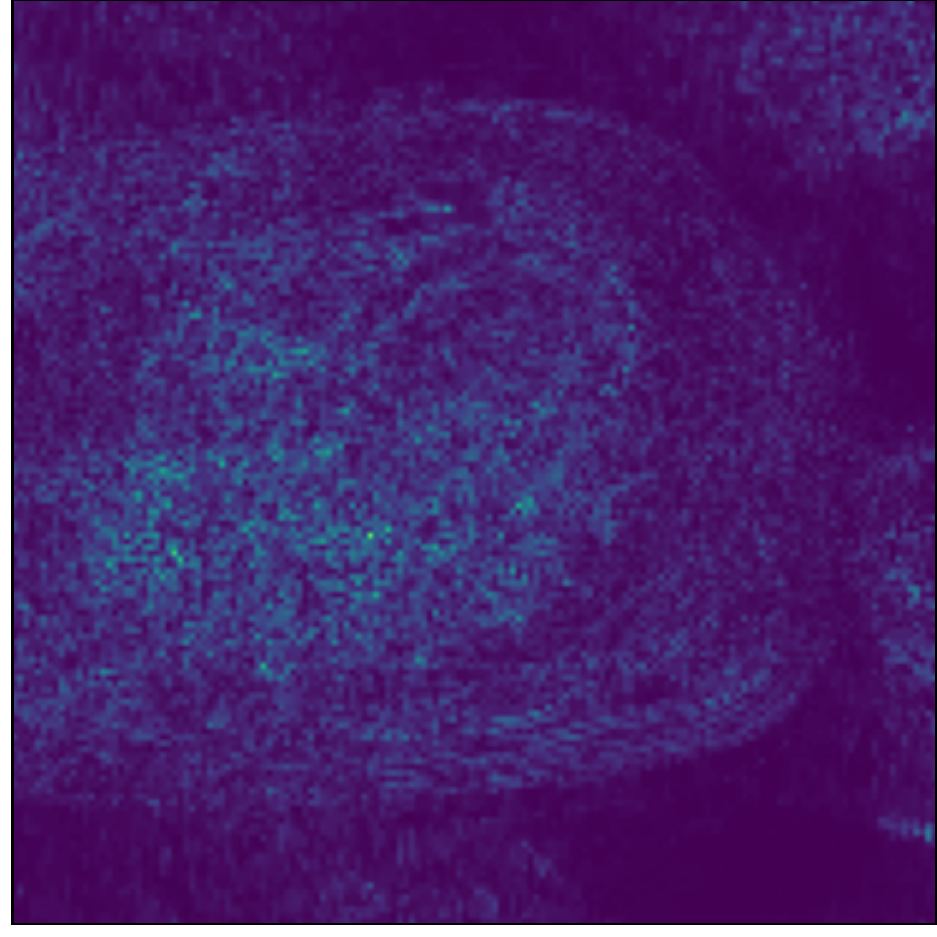}};
    \spy on (0.8, 1.5) in node [left] at (3.0, 2.5);
    \end{tikzpicture}
    \includegraphics[height=3cm]{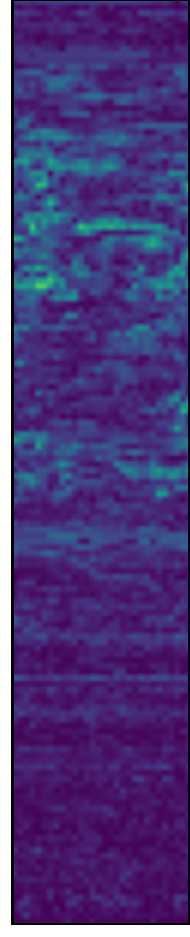}\hspace{-0.1cm}
    \begin{tikzpicture}[spy using outlines={rectangle, white, magnification=2, size=1.cm, connect spies}]
    \node[anchor=south west,inner sep=0]  at (0,0) {\includegraphics[height=3cm]{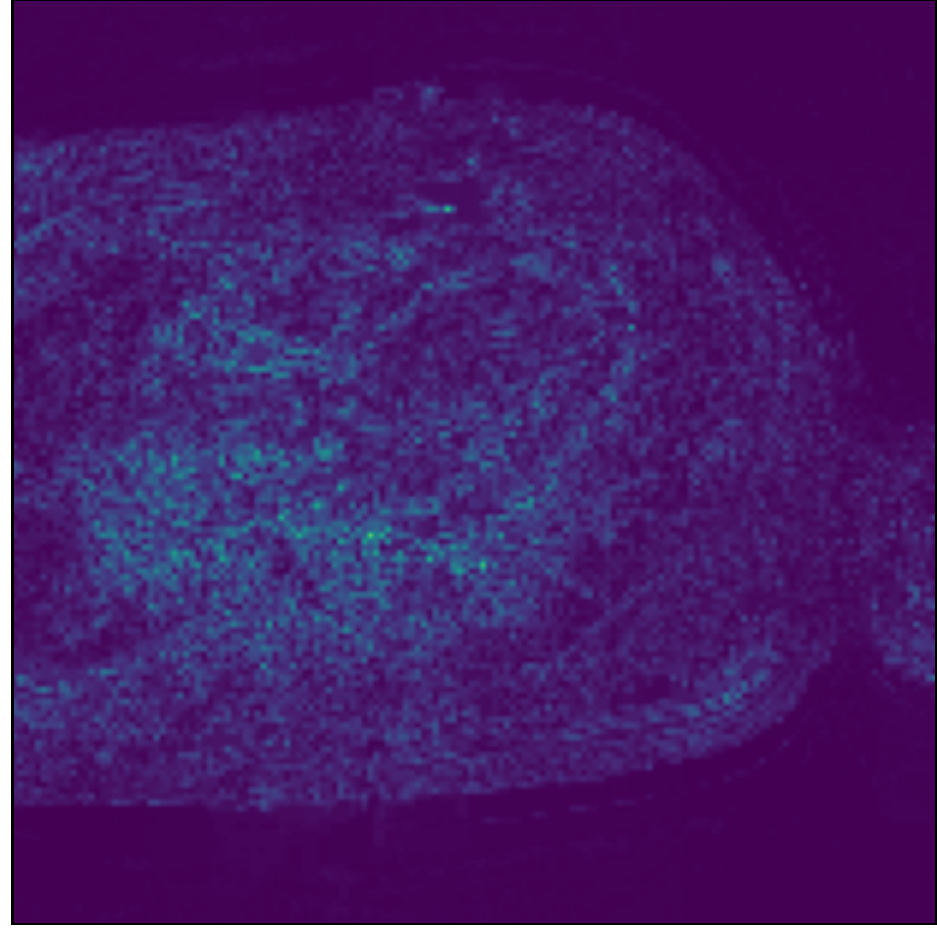}};
    \spy on (0.8, 1.5) in node [left] at (3.0, 2.5);
    \end{tikzpicture}
    \includegraphics[height=3cm]{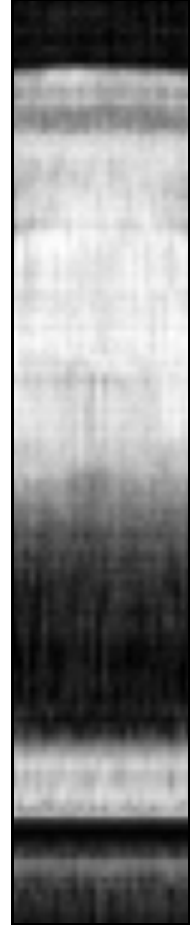}\hspace{-0.1cm}
    \begin{tikzpicture}[spy using outlines={rectangle, white, magnification=2, size=1.cm, connect spies}]
    \node[anchor=south west,inner sep=0]  at (0,0) {\includegraphics[height=3cm]{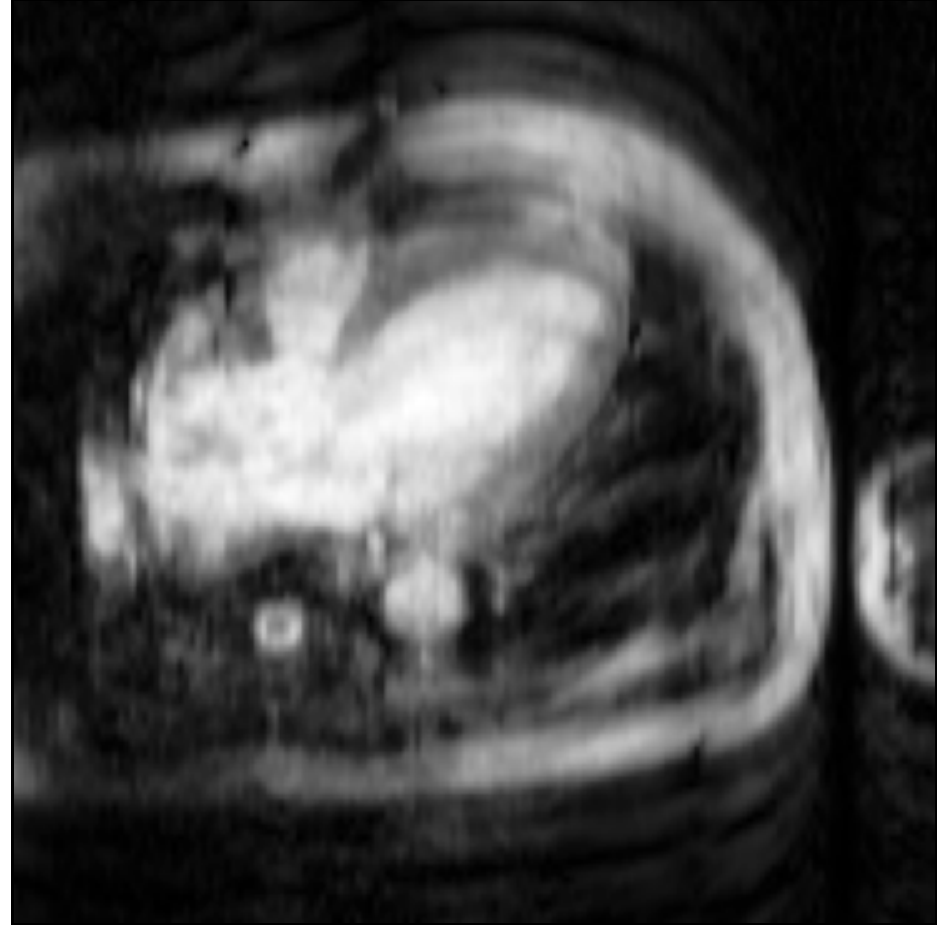}};
    \spy on (0.8, 1.5) in node [left] at (3.0, 2.5);
    \end{tikzpicture}
    }
    \end{minipage}

    \begin{minipage}{\linewidth}
    \resizebox{\linewidth}{!}{
    \includegraphics[height=3cm]{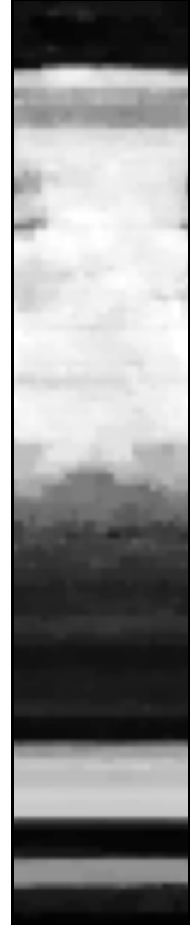}\hspace{-0.1cm}
    \begin{tikzpicture}[spy using outlines={rectangle, white, magnification=2, size=1.cm, connect spies}]
    \node[anchor=south west,inner sep=0]  at (0,0) {\includegraphics[height=3cm]{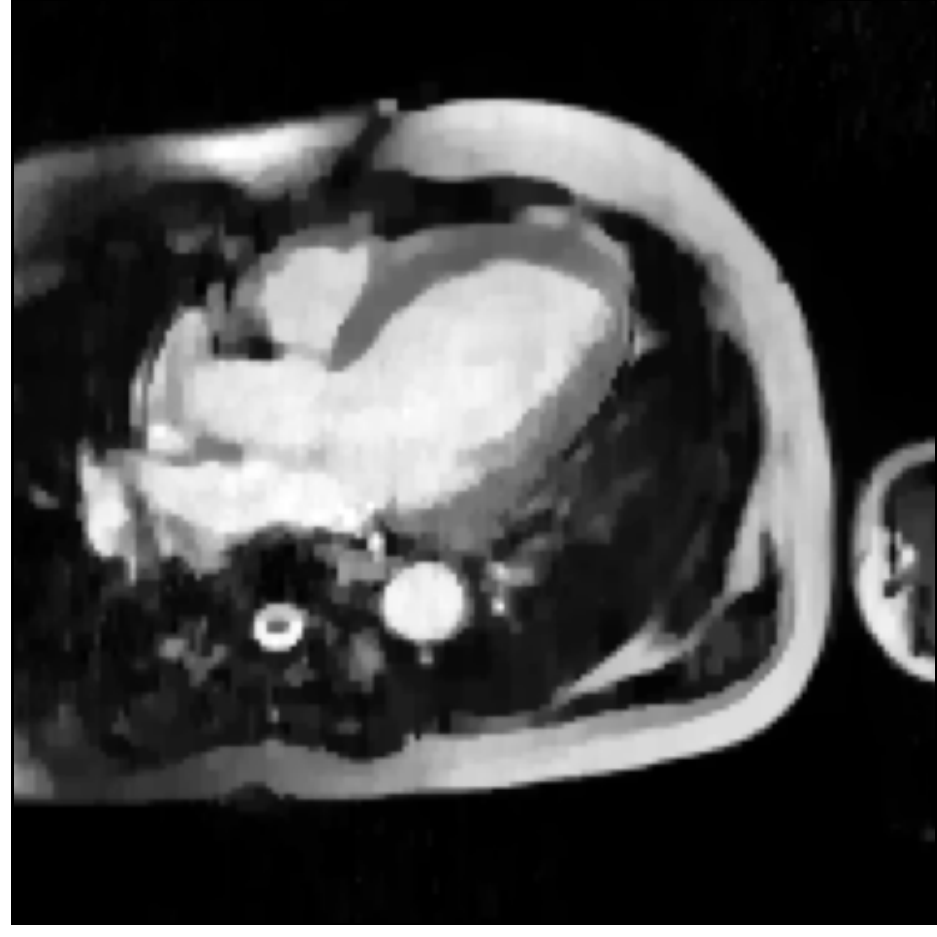}};
    \spy on (0.8, 1.5) in node [left] at (3.0, 2.5);
    \end{tikzpicture}
    \includegraphics[height=3cm]{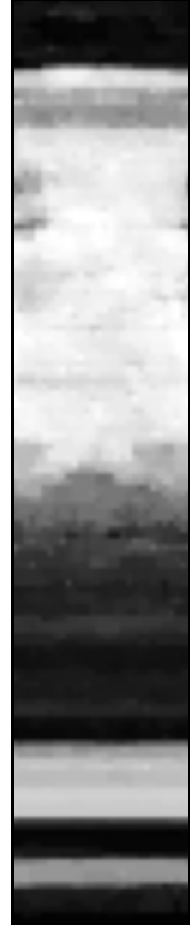}\hspace{-0.1cm}
    \begin{tikzpicture}[spy using outlines={rectangle, white, magnification=2, size=1.cm, connect spies}]
    \node[anchor=south west,inner sep=0]  at (0,0) {\includegraphics[height=3cm]{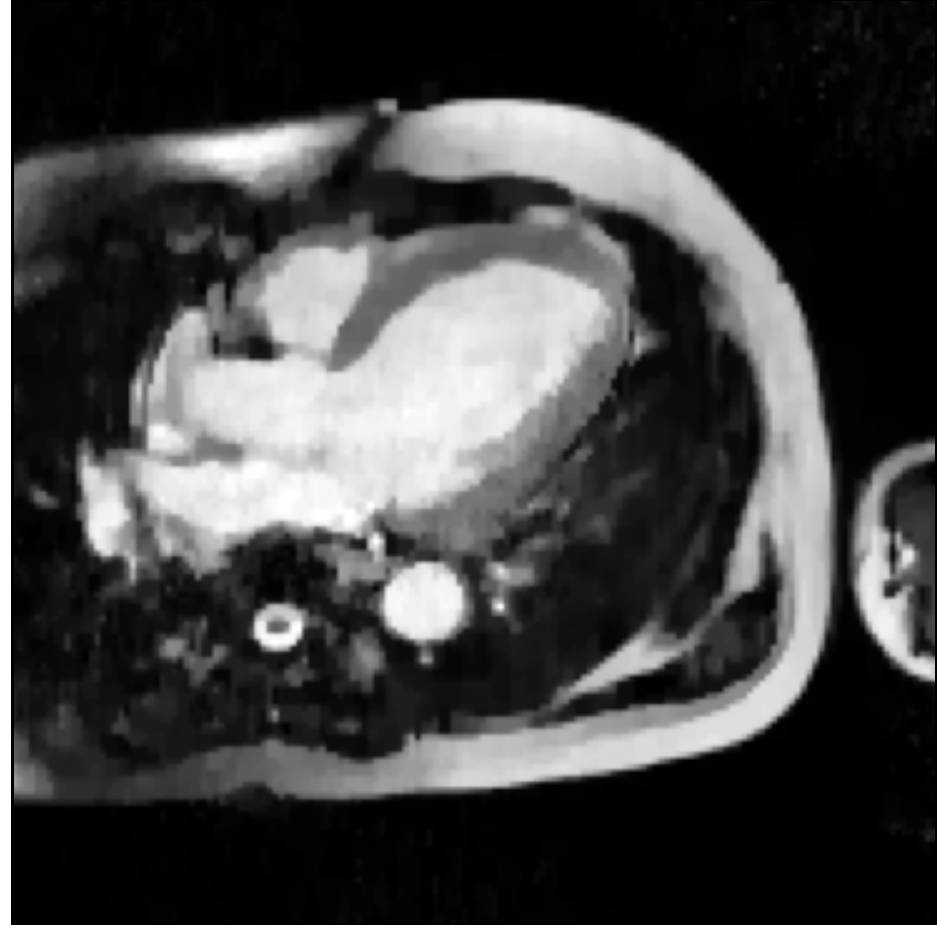}};
    \spy on (0.8, 1.5) in node [left] at (3.0, 2.5);
    \end{tikzpicture}
    \includegraphics[height=3cm]{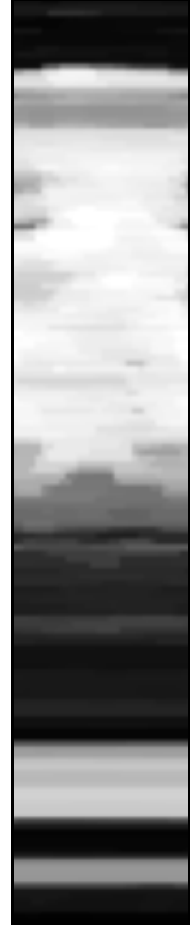}\hspace{-0.1cm}
    \begin{tikzpicture}[spy using outlines={rectangle, white, magnification=2, size=1.cm, connect spies}]
    \node[anchor=south west,inner sep=0]  at (0,0) {\includegraphics[height=3cm]{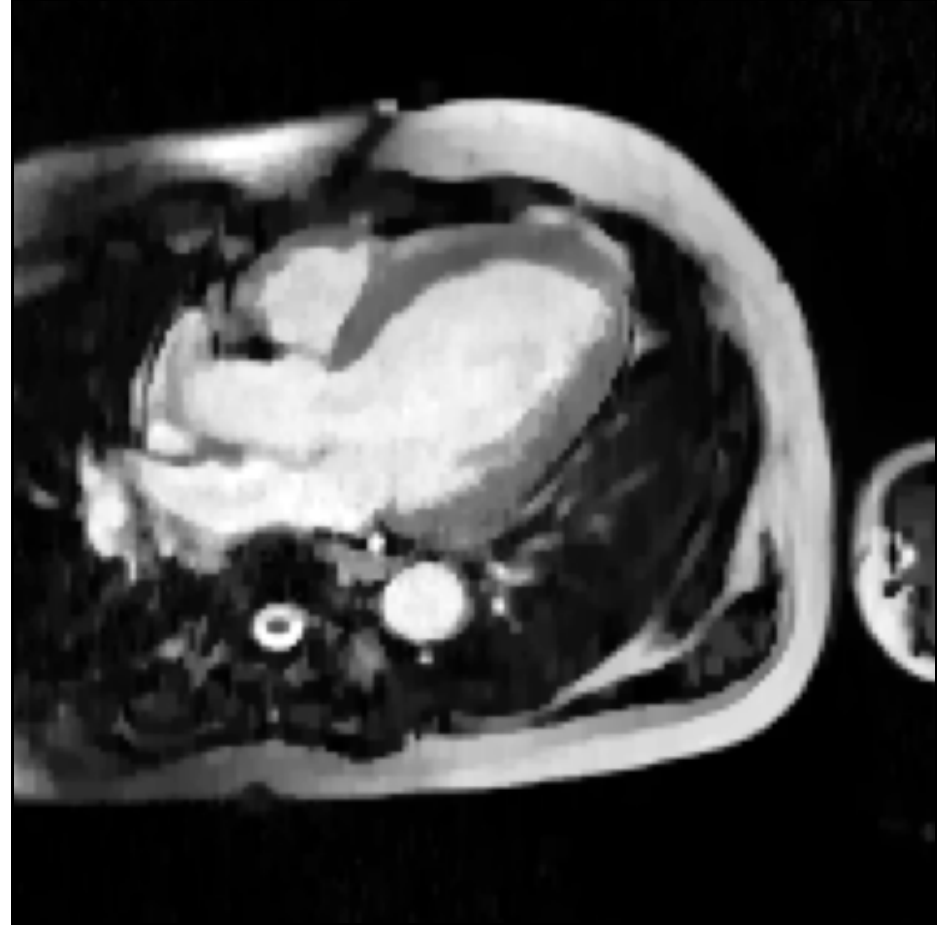}};
    \spy on (0.8, 1.5) in node [left] at (3.0, 2.5);
    \end{tikzpicture}
    \includegraphics[height=3cm]{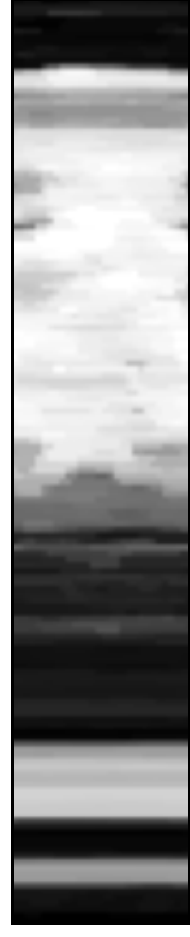}\hspace{-0.1cm}
    \begin{tikzpicture}[spy using outlines={rectangle, white, magnification=2, size=1.cm, connect spies}]
    \node[anchor=south west,inner sep=0]  at (0,0) {\includegraphics[height=3cm]{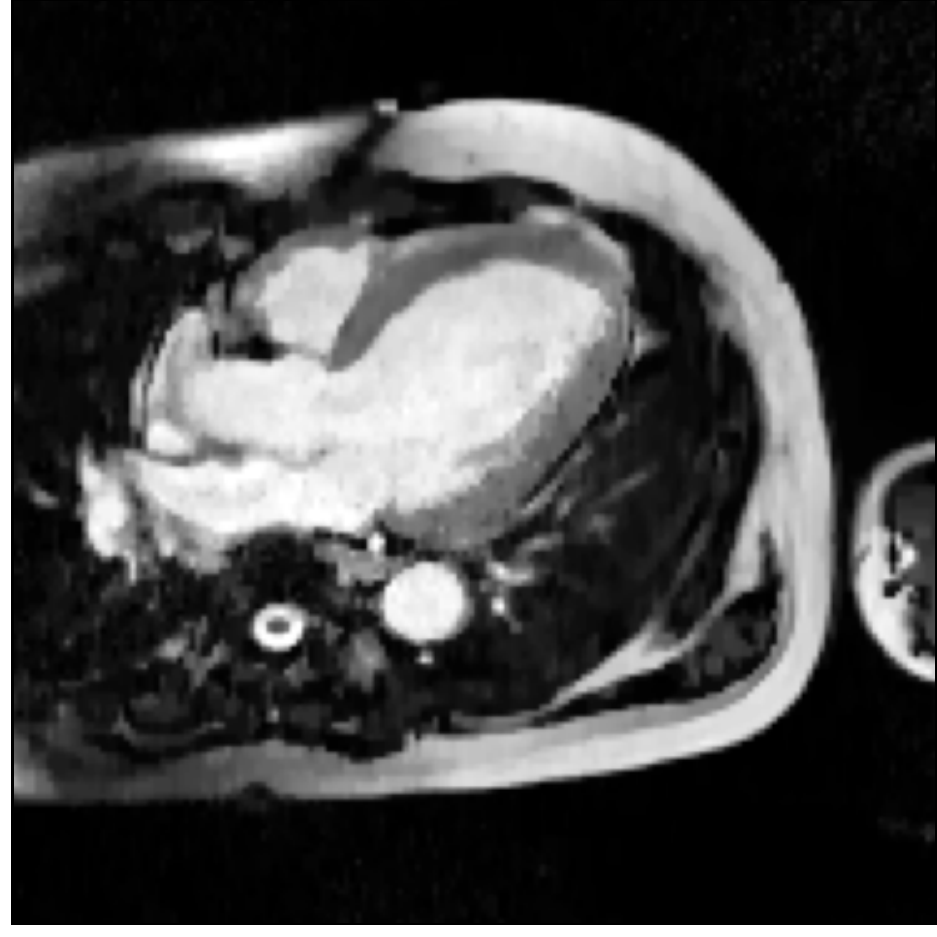}};
    \spy on (0.8, 1.5) in node [left] at (3.0, 2.5);
    \end{tikzpicture}
    \includegraphics[height=3cm]{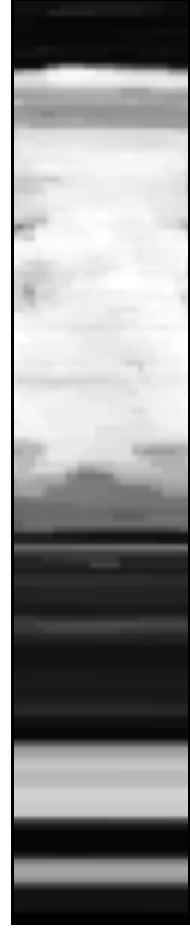}\hspace{-0.1cm}
    \begin{tikzpicture}[spy using outlines={rectangle, white, magnification=2, size=1.cm, connect spies}]
    \node[anchor=south west,inner sep=0]  at (0,0) {\includegraphics[height=3cm]{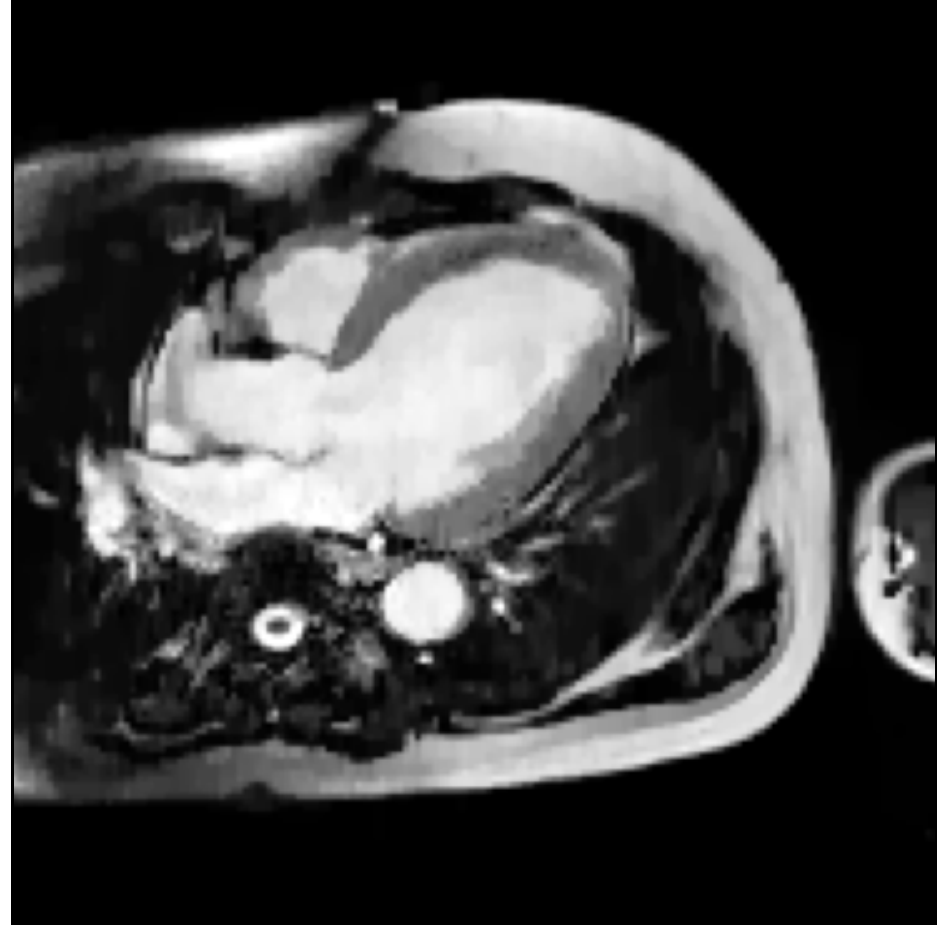}};
    \spy on (0.8, 1.5) in node [left] at (3.0, 2.5);
    \end{tikzpicture}
    \includegraphics[height=3cm]{figures/MRI/mri_results/xf_yt.pdf}\hspace{-0.1cm}
    \begin{tikzpicture}[spy using outlines={rectangle, white, magnification=2, size=1.cm, connect spies}]
    \node[anchor=south west,inner sep=0]  at (0,0) {\includegraphics[height=3cm]{figures/MRI/mri_results/xf_xy.pdf}};
    \spy on (0.8, 1.5) in node [left] at (3.0, 2.5);
    \end{tikzpicture}
    }
    \resizebox{\linewidth}{!}{
    \includegraphics[height=3cm]{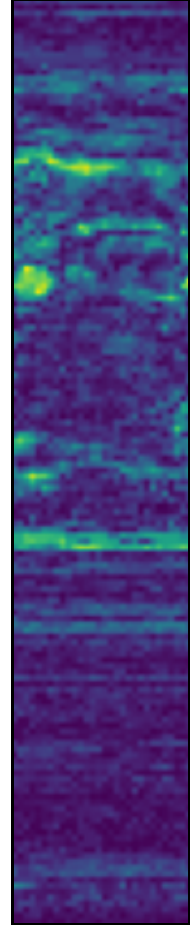}\hspace{-0.1cm}
    \begin{tikzpicture}[spy using outlines={rectangle, white, magnification=2, size=1.cm, connect spies}]
    \node[anchor=south west,inner sep=0]  at (0,0) {\includegraphics[height=3cm]{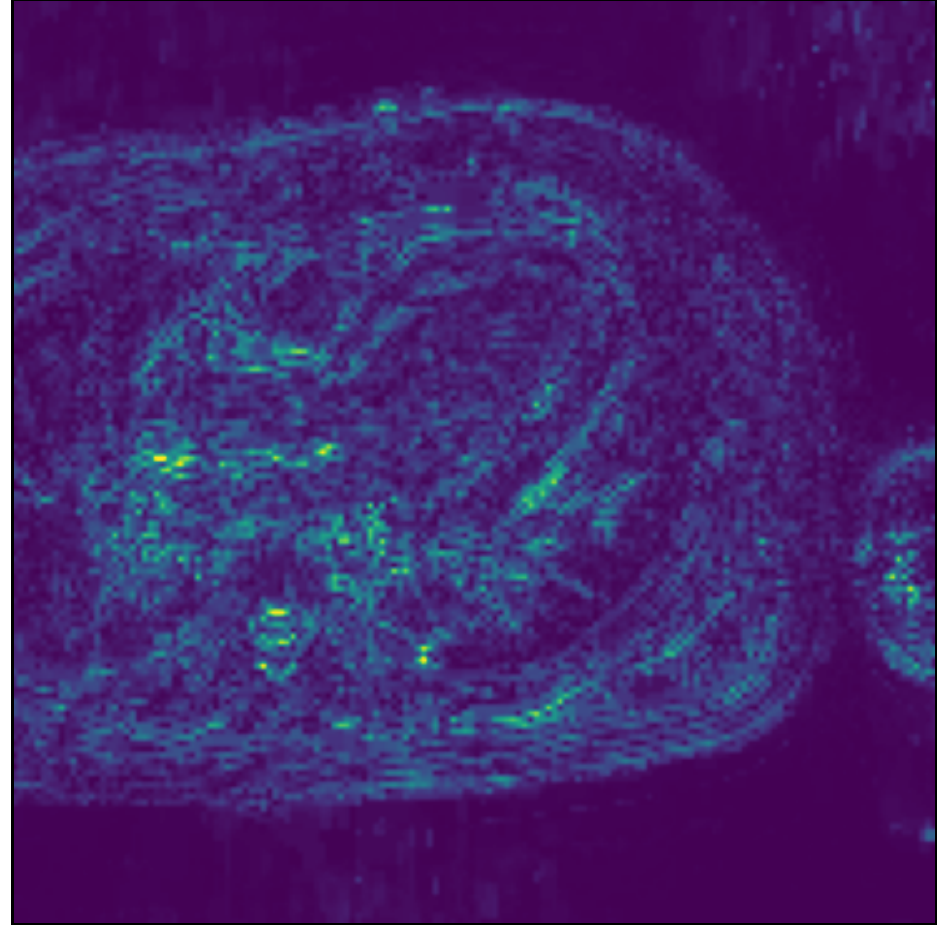}};
    \spy on (0.8, 1.5) in node [left] at (3.0, 2.5);
    \end{tikzpicture}
    \includegraphics[height=3cm]{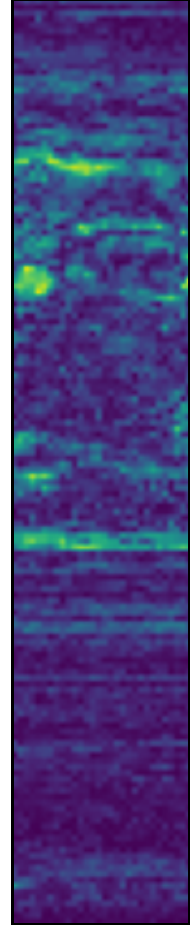}\hspace{-0.1cm}
    \begin{tikzpicture}[spy using outlines={rectangle, white, magnification=2, size=1.cm, connect spies}]
    \node[anchor=south west,inner sep=0]  at (0,0) {\includegraphics[height=3cm]{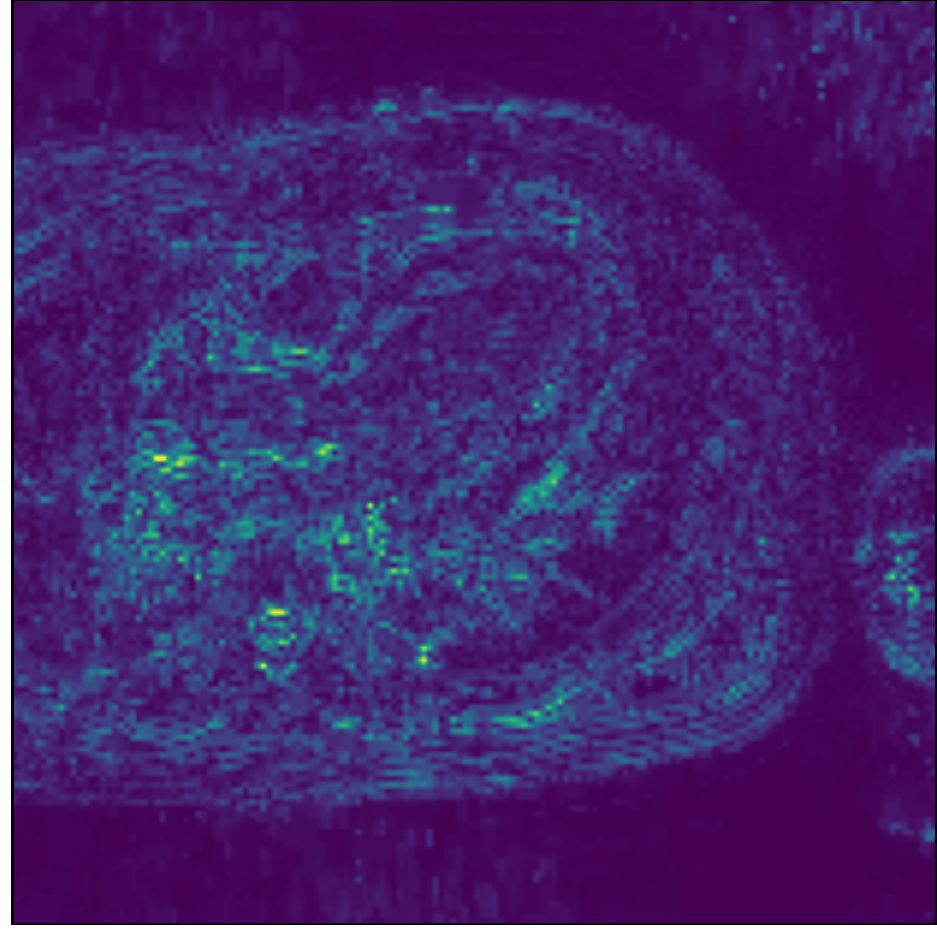}};
    \spy on (0.8, 1.5) in node [left] at (3.0, 2.5);
    \end{tikzpicture}
    \includegraphics[height=3cm]{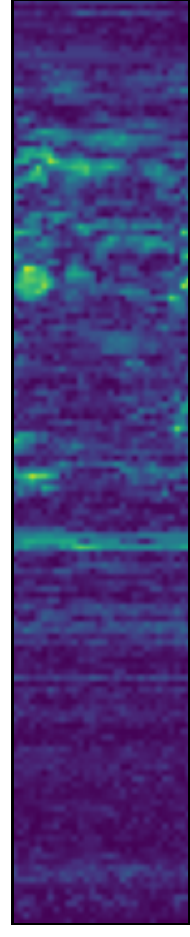}\hspace{-0.1cm}
    \begin{tikzpicture}[spy using outlines={rectangle, white, magnification=2, size=1.cm, connect spies}]
    \node[anchor=south west,inner sep=0]  at (0,0) {\includegraphics[height=3cm]{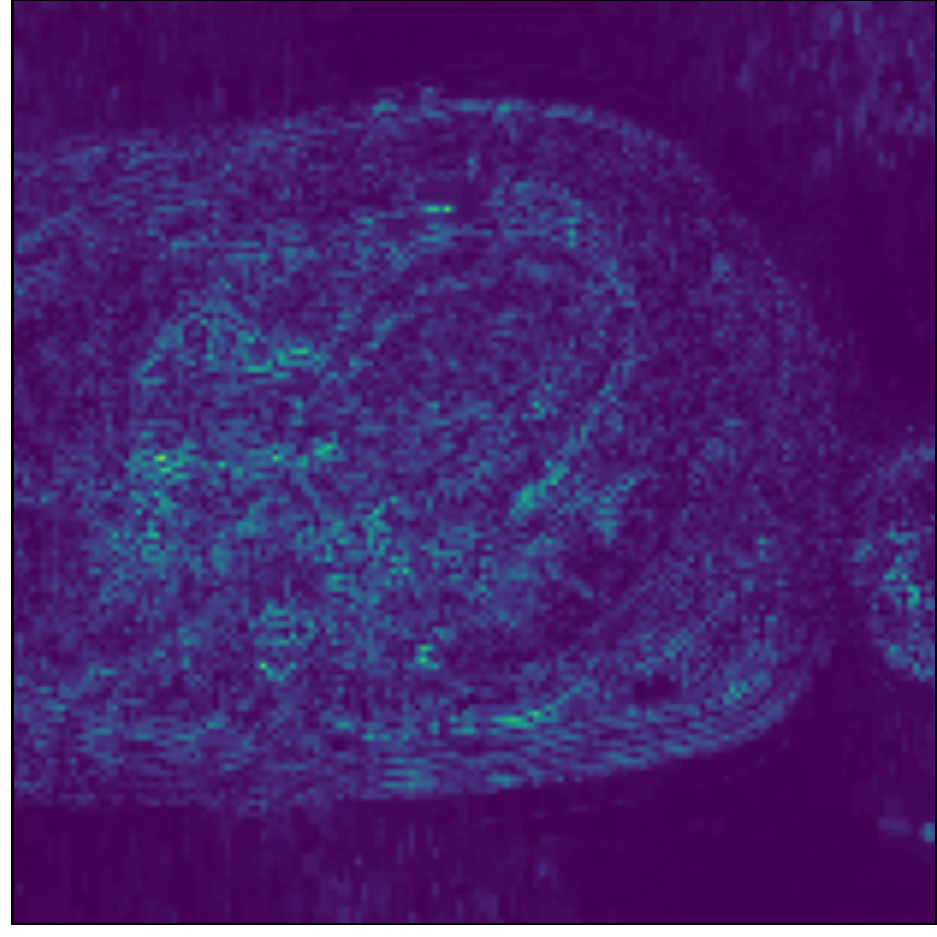}};
    \spy on (0.8, 1.5) in node [left] at (3.0, 2.5);
    \end{tikzpicture}
    \includegraphics[height=3cm]{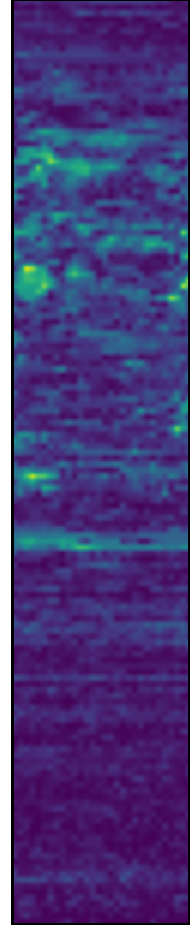}\hspace{-0.1cm}
    \begin{tikzpicture}[spy using outlines={rectangle, white, magnification=2, size=1.cm, connect spies}]
    \node[anchor=south west,inner sep=0]  at (0,0) {\includegraphics[height=3cm]{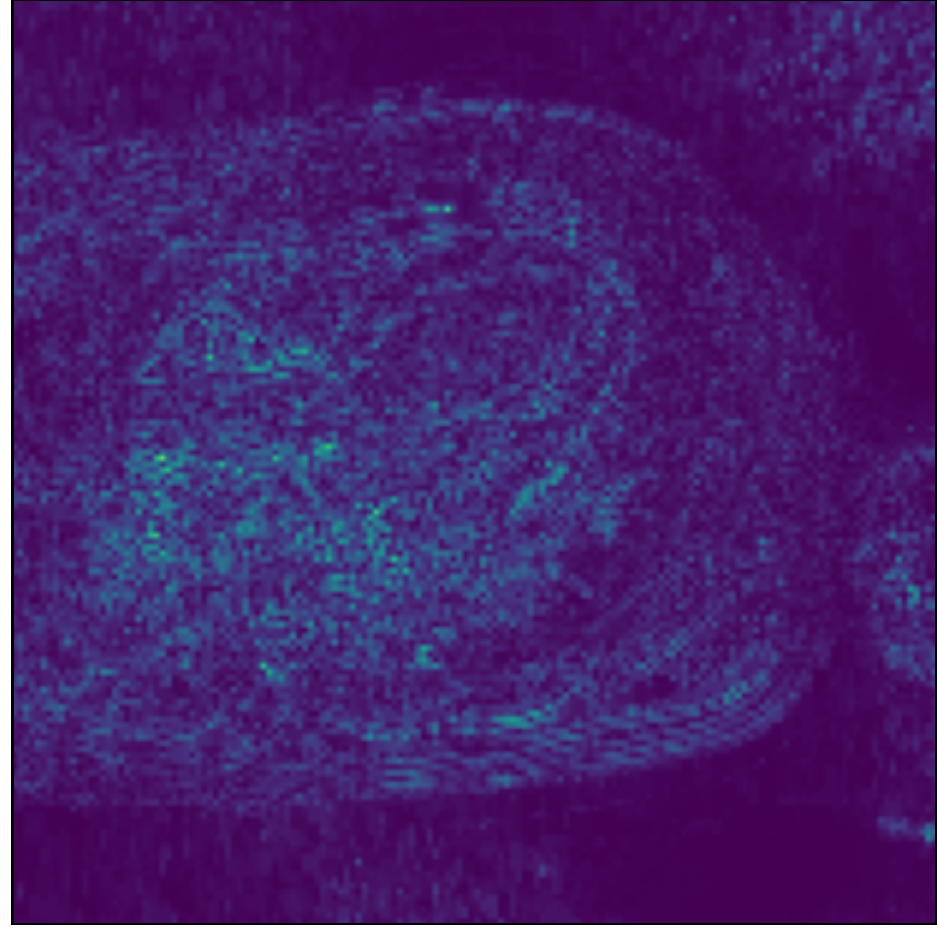}};
    \spy on (0.8, 1.5) in node [left] at (3.0, 2.5);
    \end{tikzpicture}
    \includegraphics[height=3cm]{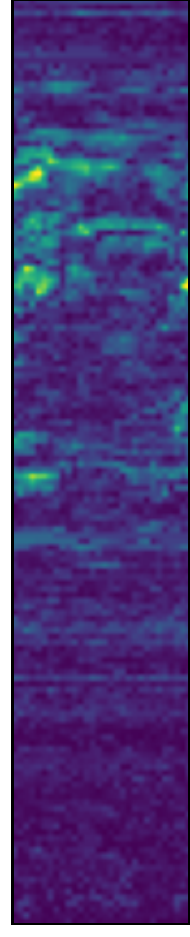}\hspace{-0.1cm}
    \begin{tikzpicture}[spy using outlines={rectangle, white, magnification=2, size=1.cm, connect spies}]
    \node[anchor=south west,inner sep=0]  at (0,0) {\includegraphics[height=3cm]{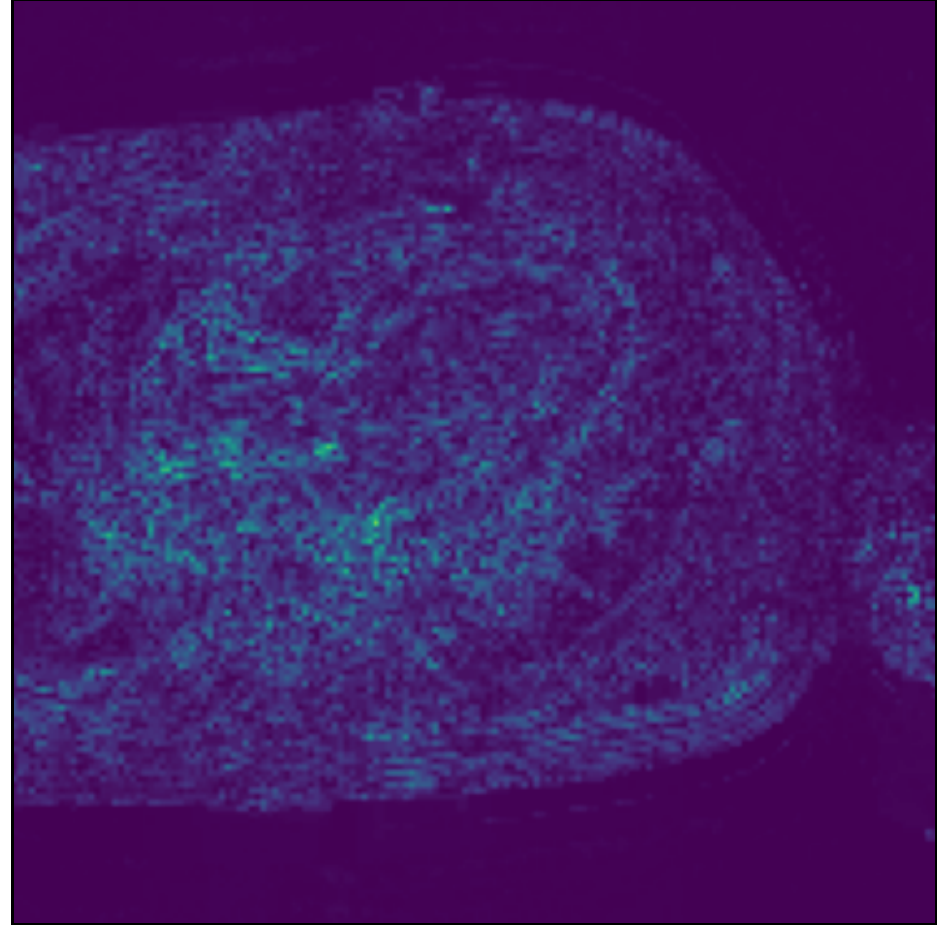}};
    \spy on (0.8, 1.5) in node [left] at (3.0, 2.5);
    \end{tikzpicture}
    \includegraphics[height=3cm]{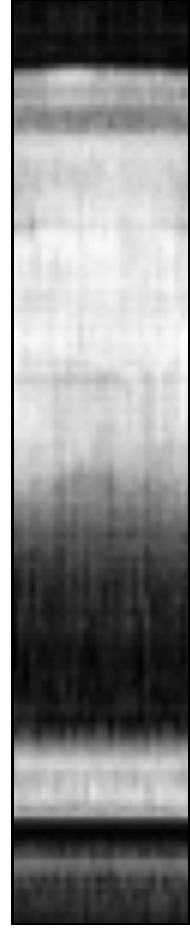}\hspace{-0.1cm}
    \begin{tikzpicture}[spy using outlines={rectangle, white, magnification=2, size=1.cm, connect spies}]
    \node[anchor=south west,inner sep=0]  at (0,0) {\includegraphics[height=3cm]{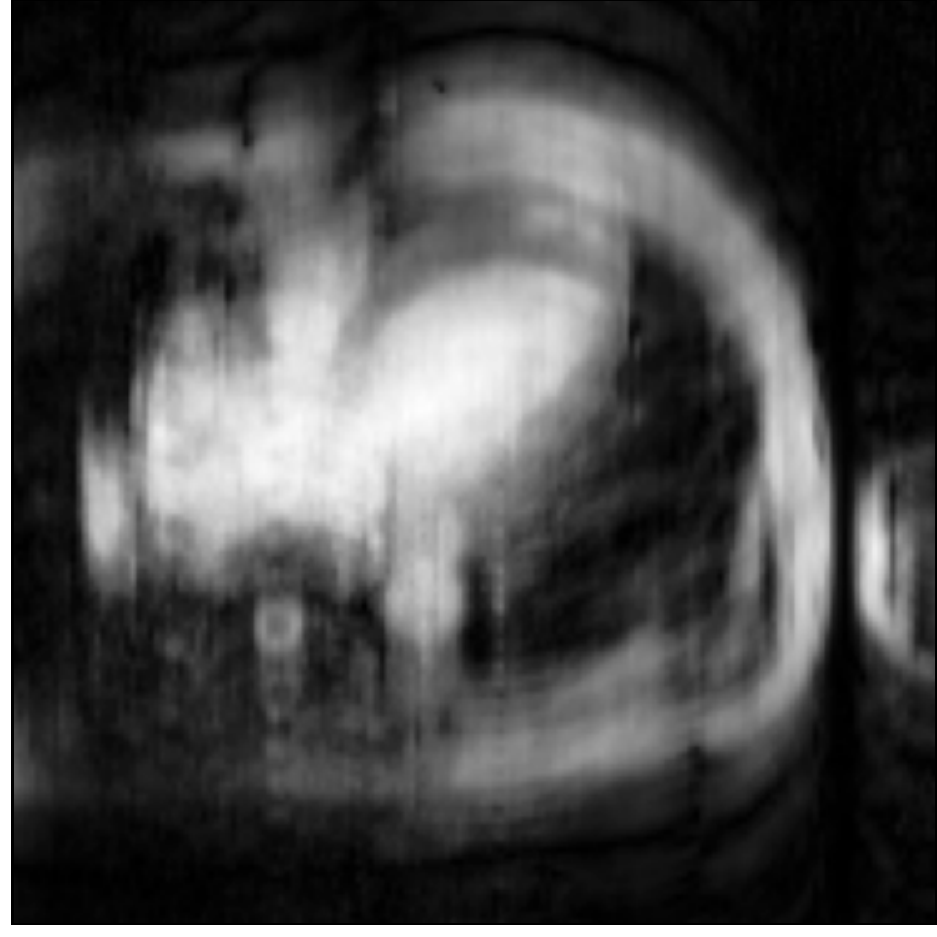}};
    \spy on (0.8, 1.5) in node [left] at (3.0, 2.5);
    \end{tikzpicture}
    }
    \end{minipage}
    \end{minipage}
    \end{minipage}
    \end{minipage} 
    \caption{An example of images reconstructed with the primal-dual scheme in Algorithm \ref{algo:tv_reco_algo_mri} for different choices of regularization parameters and acceleration factors $R=4,6,8$. 
    Single scalar regularization parameter $\lambda_{\mathrm{P}}^{xyt}$ and $\lambda_{\tilde{\mathrm{P}}}^{xyt}$, two scalar regularization parameters for differently weighted spatial and temporal components, $\lambda_{\mathrm{P}}^{xy,t}$ and $\lambda_{\tilde{\mathrm{P}}}^{xy,t}$, the proposed spatially and temporal dependent parameter-map $\boldsymbol{\Lambda}_{\Theta}^{xy,t}$ obtained with the network $\mathcal{N}_{\Theta}^T$. The last column shows the target image and the zero-filled reconstruction. Note again that the results for $\lambda_{\mathrm{P}}^{xyt}$ and $\lambda_{\mathrm{P}}^{xy,t}$ were obtained performing a grid-search for $\lambda^{xyt}>0$ and $\lambda^{xy},\lambda^t>0$, assuming the ground truth image to be known. Therefore, the results for $\lambda_{\mathrm{P}}^{xyt}$ and $\lambda_{\mathrm{P}}^{xy,t}$ cannot be obtained in practice and merely serve for illustrating that the proposed $\boldsymbol{\Lambda}_{\Theta}^{xy,t}$ yields competitive results. }
    \label{fig:mri_results}
    \label{fig:my_label}
\end{figure}

Figure \ref{fig:mri_results} shows an example of a single  frame of the reconstructed MR image sequences  for an acceleration factor of $R=6$ using several approaches. We show the reconstructions that correspond to the single scalar parameter $\lambda_{\mathrm{P}}^{xyt}$ as well as  to the scalar parameter pair (one spatial and one temporal)  $\lambda_{\mathrm{P}}^{xy,t}=(\lambda_{\mathrm{P}}^{xy}, \lambda_{\mathrm{P}}^{xy}, \lambda_{\mathrm{P}}^{t})$ which are the parameters that maximize the PSNR of entire cine MR image and are obtained via a grid search by making use of the corresponding ground truth image. We also show the results that correspond to the parameters $\lambda_{\tilde{\mathrm{P}}}^{xyt}$ and $\lambda_{\tilde{\mathrm{P}}}^{xy,t}$ which are respectively the single and the pair of scalar parameters that on average maximize the PSNR over the training set. These were obtained by treating the scalar regularization parameters as trainable parameters and training them by minimizing \eqref{eq:loss_fct}. We finally show the results  for our  estimated parameter-map $\boldsymbol{\Lambda}_{\Theta}^{xy,t}$ with the proposed method. As observed, for all choices of the regularization parameters, the error with respect to the target image was significantly reduced compared to the initial zero-filled reconstruction. Further, we can see how the use of the estimated parameter-map yields the most accurate reconstruction and the best preservation of image details.

\begin{figure}[!h]
    \centering
    \includegraphics[width=0.5\linewidth]{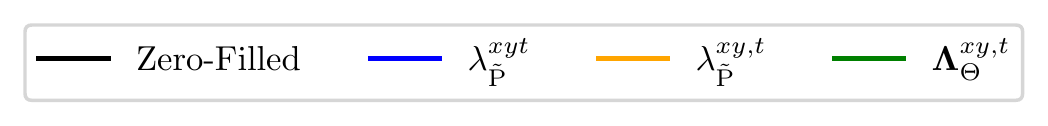}\\
    \includegraphics[width=0.325\linewidth]{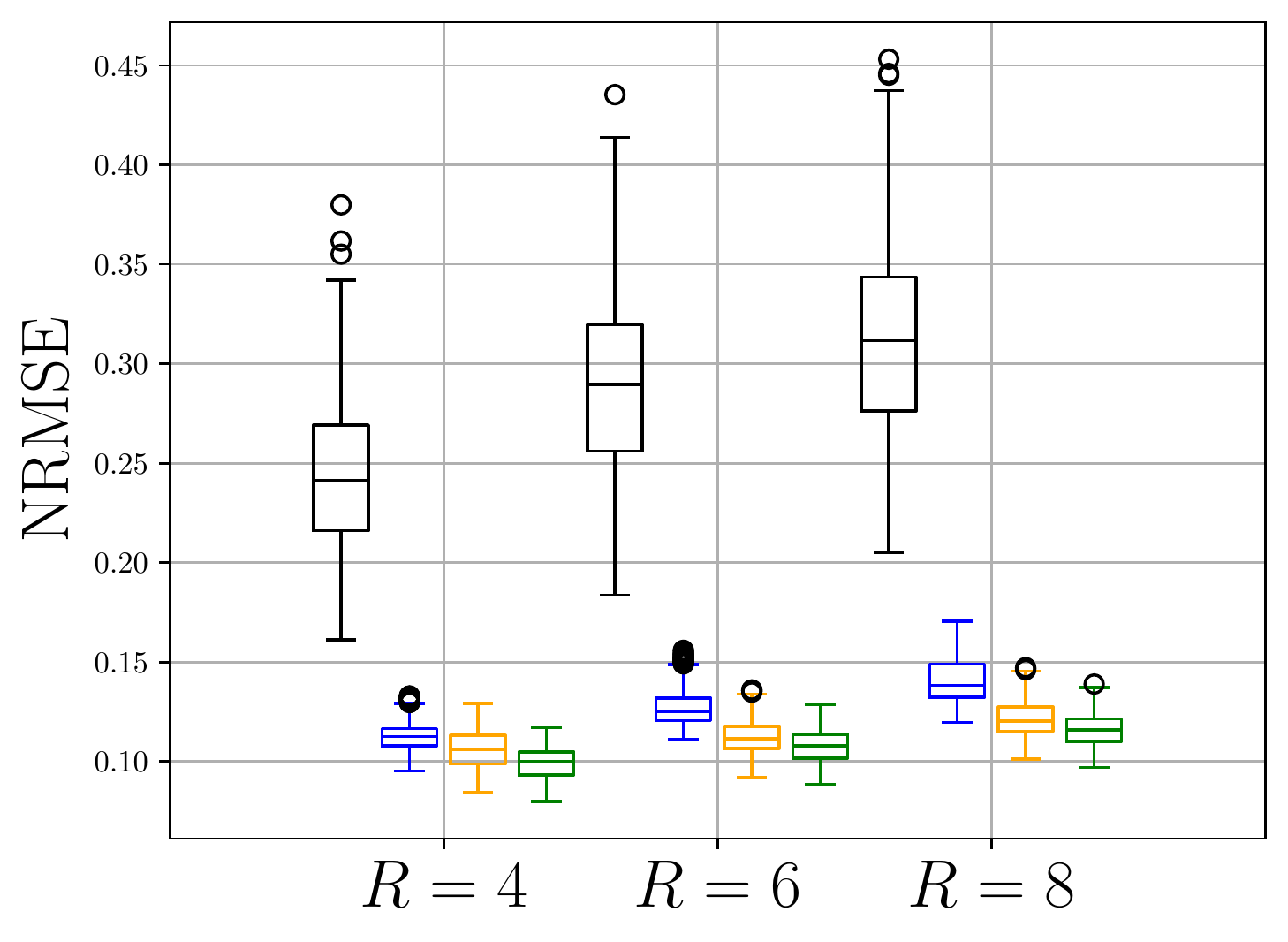}
    \includegraphics[width=0.325\linewidth]{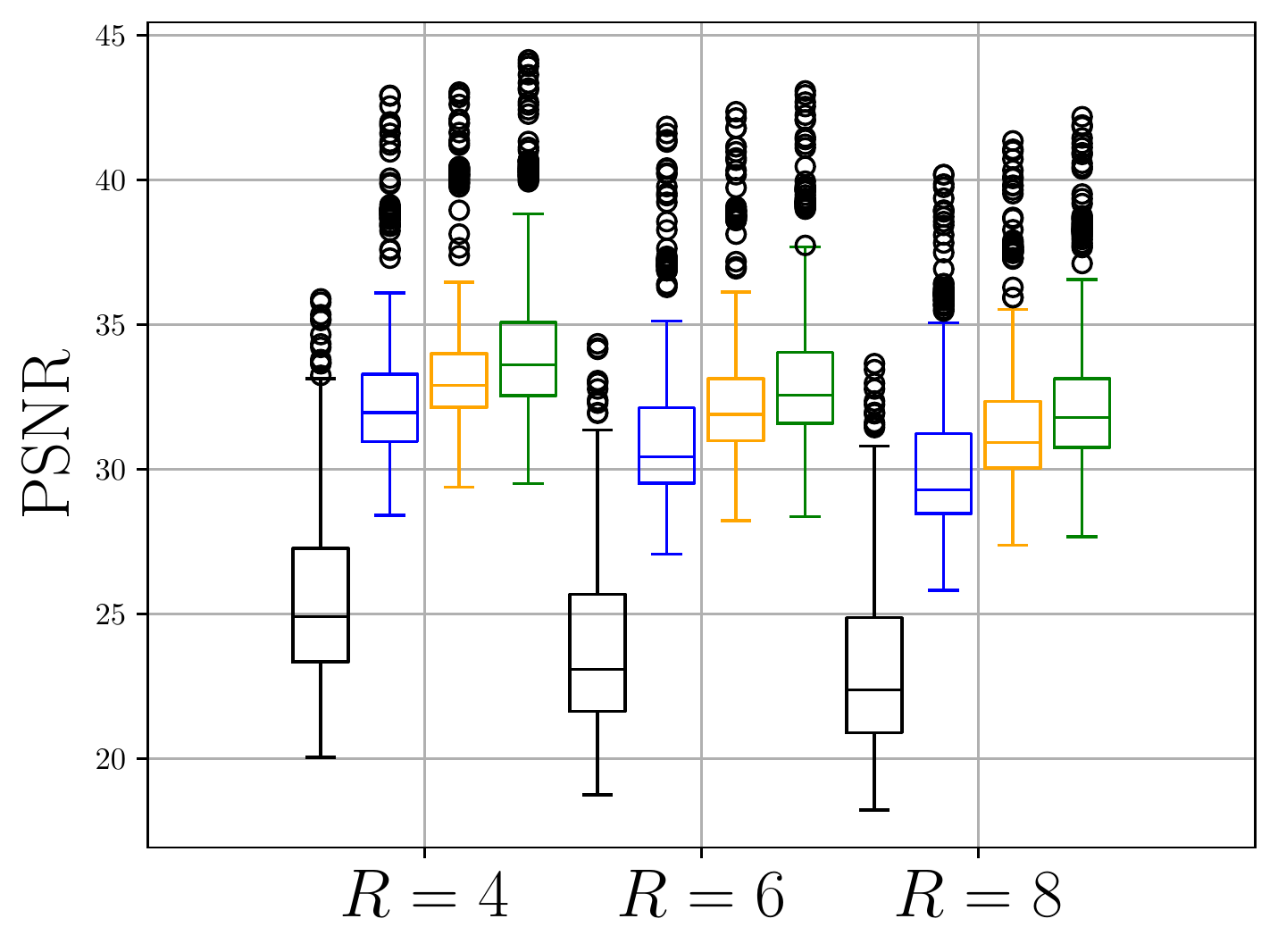}
    \includegraphics[width=0.325\linewidth]{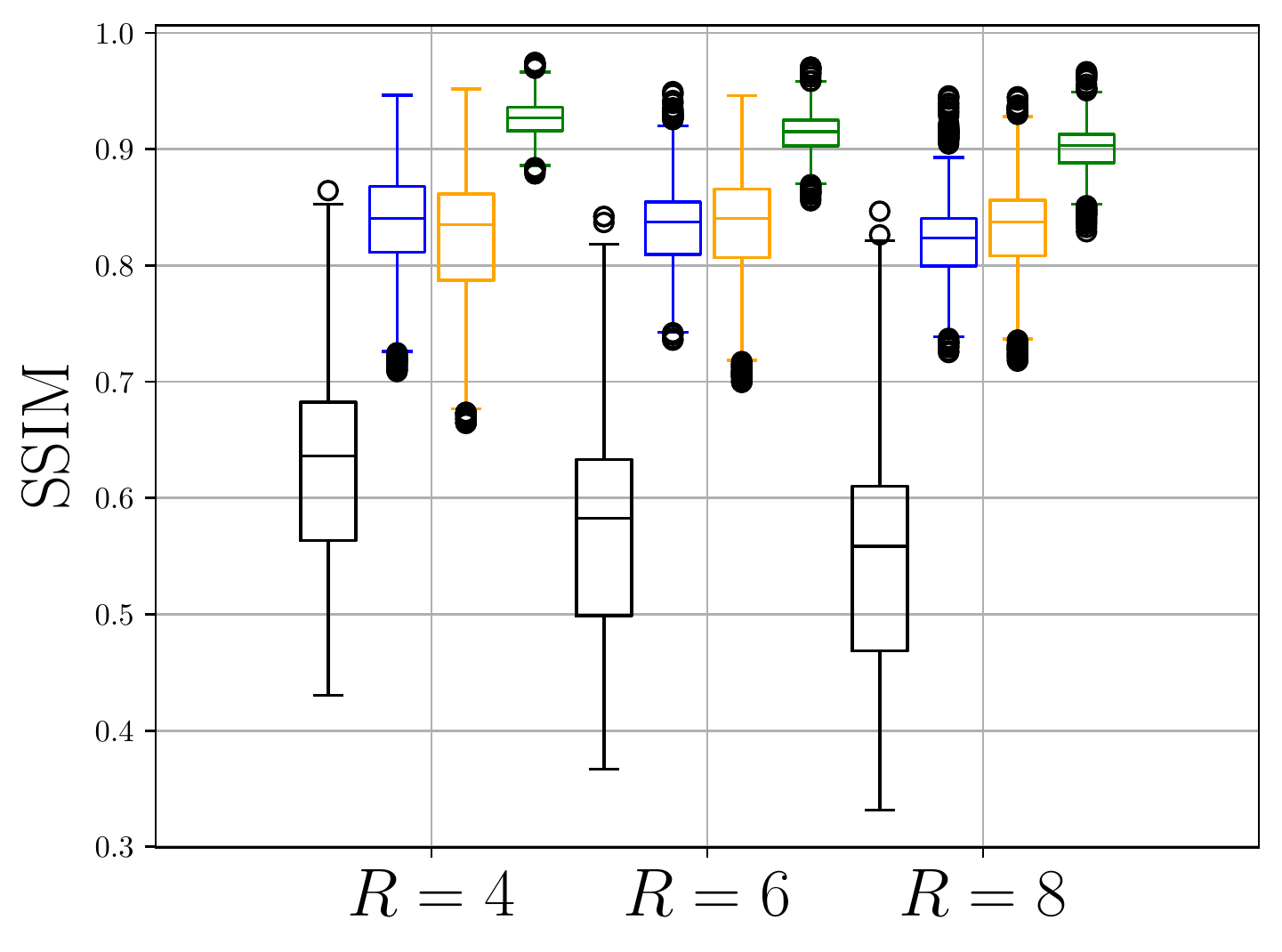}
    \caption{Box-plots summarizing the reconstruction results in terms of PSNR, NRMSE and SSIM obtained with the PDHG for a dynamic cardiac MR image reconstruction problem for different choices of the regularization parameter. Zero-filled reconstruction (black), single scalar regularization parameter  ($\lambda_{\tilde{\mathrm{P}}}^{xyt}$, blue), two scalar regularization parameters; one for the spatial $x$- and $y$-direction, one for the temporal $t$-direction ($\lambda_{\tilde{\mathrm{P}}}^{xy,t}$, blue) and the proposed spatially and temporal dependent parameter-map $\boldsymbol{\Lambda}_{\Theta}^{xy,t}$ obtained with a CNN ($\mathrm{NET}_{\Theta}$, green).}
    \label{fig:box_plots_mri}
\end{figure}

Figure \ref{fig:box_plots_mri} summarizes the results obtained over the test set with the help of box-plots. Compared to the initial zero-filled reconstruction, an improvement is clearly visible for all choices of the regularization parameter with respect to all reported measures and for all acceleration factors. In addition we see how allowing the temporal direction to be differently regularized than the two spatial dimensions positively influences the results compared to having one global parameter $\lambda$ (orange vs blue). Last, we see how using the proposed method to estimate an entire spatio-temporal parameter-map $\boldsymbol{\Lambda}_{\Theta}$ further surpasses the scalar regularization parameter-maps (green vs orange and blue), especially in terms of SSIM.
Table \ref{tab:dyn_mri_results} lists the mean and the standard deviation of all TV-reconstructions. The results are consistent with the ones from the box-plots.

\begin{figure}
    \centering
        \begin{overpic}[width=0.96\linewidth]{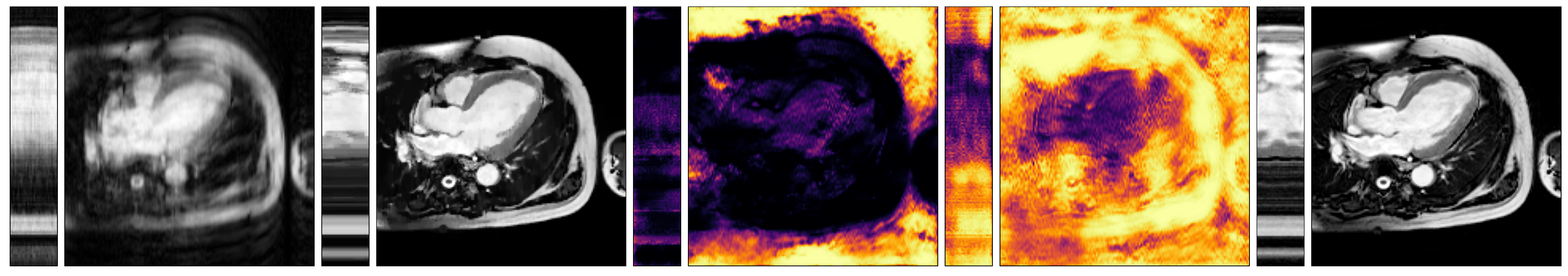}
        \put(4,18){ {\textcolor{black}{\normalsize{Zero-Filled}}}} 
        \put(24,18){ {\textcolor{black}{\normalsize{PDHG $\boldsymbol{\Lambda}_{\Theta}^{xy,t}$}}}} 
        \put(49,18){ {\textcolor{black}{\normalsize{$\boldsymbol{\Lambda}_{\Theta}^{xy}$}}}} 
        \put(69,18){ {\textcolor{black}{\normalsize{$\boldsymbol{\Lambda}_{\Theta}^{t}$}}}} 
        \put(86,18){ {\textcolor{black}{\normalsize{Target}}}} 
        \end{overpic} 
        \begin{overpic}[width=0.96\linewidth]{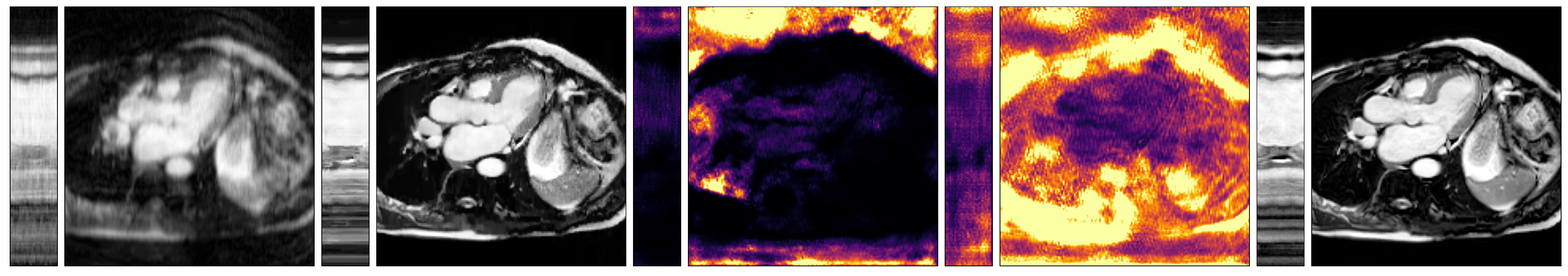}
        \end{overpic}
        \begin{overpic}[width=0.96\linewidth]{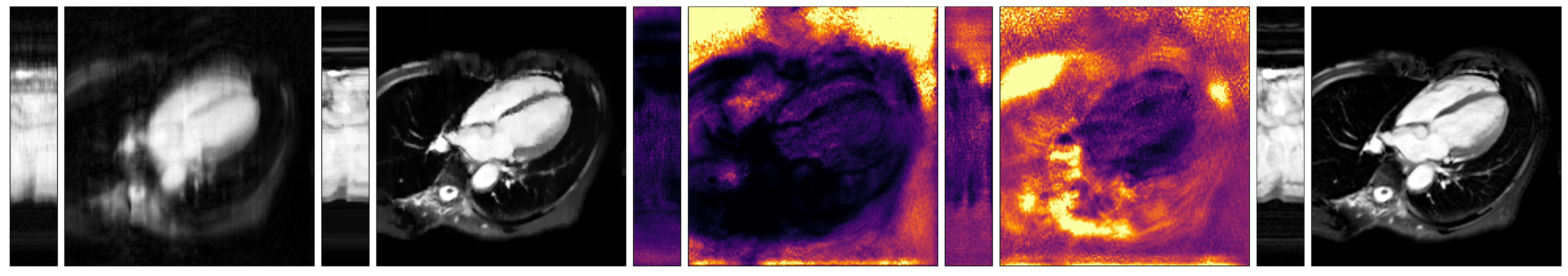}
        \end{overpic}
        \includegraphics[width=0.8\linewidth]{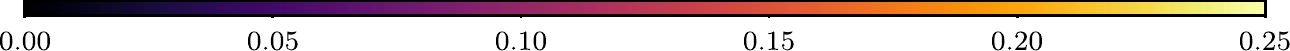}
    \caption{Different examples of reconstruction results and regularization parameter-maps for an acceleration factor of $R=6$. From left to right for each row: zero-filled reconstruction, PDHG-reconstruction with $T=4096$ obtained with the proposed CNN-based regularization spatio-temporal parameter-map $\boldsymbol{\Lambda}_{\Theta}^{xy,t}=(\boldsymbol{\Lambda}_{\Theta}^{xy}, \boldsymbol{\Lambda}_{\Theta}^{xy}, \boldsymbol{\Lambda}_{\Theta}^{t})$, spatial parameter-map $\boldsymbol{\Lambda}_{\Theta}^{xy}$, temporal parameter-map  $\boldsymbol{\Lambda}_{\Theta}^{t}$ and  target ground truth image. Both parameter-maps $\boldsymbol{\Lambda}_{\Theta}^{xy}$ and $\boldsymbol{\Lambda}_{\Theta}^t$ are displayed on the scale $[0,0.25]$.}
    \label{fig:lambda_maps_mri}
\end{figure}

Figure \ref{fig:lambda_maps_mri} shows an example of a spatio-temporal regularization parameter-map which was estimated using the proposed approach for an acceleration factor of $R=6$. The network $u_{\Theta}$ estimates the regularization parameter-map to be pointwise relatively consistenly higher than the spatially required regularization. This result is in fact expected as the temporal dimension is the one for which the gradients of the images are the sparsest because of the high temporal correlation. Further, we see how the network consistently predicts both the spatial regularization  as well as the temporal regularization to be less strong in the area where most of the movement is expected, i.e.\ in the cardiac region.

\begin{table}[h]
\centering
\footnotesize\rm
\begin{tabular}{r|r|rcl|rcl|rcl}

&    & \multicolumn{3}{c}{\textbf{PDHG - $\lambda_{\tilde{\mathrm{P}}}^{xyt}$}}  & \multicolumn{3}{c}{\textbf{PDHG - $\lambda_{\tilde{\mathrm{P}}}^{xy,t}$}} & \multicolumn{3}{c}{\textbf{PDHG - $\boldsymbol{\Lambda}_{\Theta}^{xy,t}$}}  \\[5pt]
\hline
        &  \textbf{SSIM}    & $0.836$ &  $\pm$ & $0.048$    &  $0.824$ &  $\pm$ & $0.059$     & \textbf{0.927} &  $\pm$ & $0.016$   \\
$R=4$   &  \textbf{PSNR}    & $32.35$ &  $\pm$ & $2.21$   &  $33.19$ &  $\pm$ & $2.09$    & \textbf{33.91} & $ \pm$ & $ 2.23$    \\      
        &  \textbf{NRMSE}   & $0.113$ &  $\pm$ & $0.007$    &  $0.106$ &  $\pm$ & $0.010$     & \textbf{0.099} & $ \pm $ & $ 0.008$   \\   
        &  \textbf{Blur}   & $0.369$ &  $\pm$ & $0.018$    &  \textbf{0.353} & $ \pm $ & $0.017$     & $0.359$ &  $\pm$ & $0.019$   \\   
\hline
        &  \textbf{SSIM}    & $0.833$ &  $\pm$ & $0.038$    &  $0.834$ &  $\pm$ & $0.048$     & \textbf{0.915} & $ \pm $ & $ 0.018$ \\
$R=6$   &  \textbf{PSNR}    & $30.98$ &  $\pm$ & $2.28$   &  $32.28$ &  $\pm$ & $2.10$    & \textbf{32.94} & $ \pm $ & $2.21$   \\  
        &  \textbf{NRMSE}   & $0.127$ &  $\pm$ & $0.009$    &  0.112 & $ \pm $& $0.008$     & \textbf{0.108}& $ \pm $& $0.008$   \\  
        &  \textbf{Blur}   & $0.391$ &  $\pm$ & $0.021$    &  $0.369$ &  $\pm$ & $0.018$     & \textbf{0.365} & $ \pm $ & $ 0.019$   \\   
\hline
        &  \textbf{SSIM}    & $0.822$ &  $\pm$ & $0.035$    & $0.832$ &  $\pm$ & $0.042$      & \textbf{0.902}& $ \pm $&$0.021$   \\
$R=8$   &  \textbf{PSNR}    & $29.95$ &  $\pm$ & $2.32$   & $31.37$ &  $\pm$ & $2.14$     & \textbf{32.10} & $ \pm $ & $ 2.18$   \\  
        &  \textbf{NRMSE}   & $0.141$ &  $\pm$ & $0.011$    & 0.122& $ \pm $& $0.009$      & \textbf{0.117} & $ \pm $& $0.008$   \\
        &  \textbf{Blur}   & $0.408$ &  $\pm$ & $0.025$    &  $0.382$ &  $\pm$ & $0.020$     & \textbf{0.371}& $ \pm $& $0.020$   \\   
\hline 
\end{tabular}
\caption{Mean and standard deviation of the measures SSIM, PSNR and NRMSE and Blur obtained over the test set. The TV-reconstruction using the proposed spatio-temporal parameter-maps $\boldsymbol{\Lambda}_{\Theta}^{xy,t}$ improves the results especially in terms of SSIM and PSNR.}\label{tab:dyn_mri_results}
\end{table}

\begin{remark}
    From Algorithm \ref{algo:tv_reco_algo_mri}, we see that the considered PDHG algorithm for solving problem \eqref{eq:tv_min_problem_parameter_map} involves the repeated \textit{separate} application of the forward and the adjoint operators $\Ad$ and $\Ad^\herm$. Depending on the considered problem - more precisely, on the operator of the data-acquisition - this aspect can be problematic.
For example, for non-Cartesian sampling trajectories in MRI, the separate application of $\Ad^\herm$ and $\Ad$ is considerably slower than the composition of $\Ad^\herm\Ad$ because the latter can be efficiently approximated by the Toeplitz-kernel trick \cite{feichtinger1995efficient}, see e.g.\ \cite{mani2015fast,tamir2017t2} for applications. Thereby, the composition of the forward and the adjoint can be approximated by $\Ad^\herm  \Ad \approx \Fd^\herm \Wd \Fd$, where $\Wd$ are Toeplitz-kernels which can be estimated depending on the sampling trajectories and $\Fd^\herm$ and $\Fd$ denote efficient implementations of the FFT.
Therefore, choosing a different reconstruction algorithm that requires the application of  $\Ad^\herm \Ad$ rather than $\Ad^\herm$ and $\Ad$ separately, e.g.\ \cite{wang2008new, goldstein2009split}, may be a viable option for non-Cartesian MRI.
\end{remark}
\subsection{Quantitative MRI Reconstruction}\label{subsec:qmri}
Here, we apply the proposed method to estimate voxel-wise regularization parameter-maps to be used in a quantitative brain MRI reconstruction problem. Similar to the previous case study, the problem consists of different decoupled 2D problems. However, the third temporal dimension contains information about the changing magnetization and thus over time, the contrast of the images changes. Moreover, the speed at which the contrast changes is voxel-depending. This suggests that, different from the previously shown dynamic MRI example, the dynamic component of the estimated regularization parameter-maps should also change over time and thus regularize each time point of the images differently.

\subsubsection{Problem Formulation}

Formally, the data-acquisition process for quantitative MRI reconstruction problems is given by 

\begin{equation}\label{eq:qmri_forward_problem}
    \ZZ = \Ad\, q(\mathbf{u}) + \mathbf{e},
\end{equation}
where the operator $\Ad$ takes the exact form as in Subsection \ref{subsec:dyn_mri}. However, instead of acquiring the $k$-space data of a sequence of qualitative 2D images $\XX = [\XX_1, \ldots, \XX_{N_t}]^\trans$ with similar image contrast, the operator $\Ad$ collects the $k$-space data of the qualitative images defined by 
\begin{equation}\label{eq:q_forward_model}
    \XX_t = q_t(\mathbf{u}),
\end{equation}
where $q_t:\mathbb{R}^{un} \rightarrow \mathbb{R}^n$ combines the vector containing the $u$ quantitative parameters $\mathbf{u}=[\mathbf{u}_1, \ldots, \mathbf{u}_u]^\trans$ to a qualitative image by a non-linear signal-model $q_t$. In the following, we will consider the inversion recovery signal model for $T_1$-mapping given by 
\begin{align}\label{eq:signal_model_qmri}
    q_t: \mathbb{R}^{3 n}  &\rightarrow \mathbb{C}^n \\
     [\mathbf{T}_1, \operatorname{Re}(\mathbf{M}_0), \operatorname{Im}(\mathbf{M}_0) ]^\trans &\mapsto q_t\left(\mathbf{T}_1, \operatorname{Re}(\mathbf{M}_0), \operatorname{Im}(\mathbf{M}_0 )\right)= \mathbf{M}_0 (1-2 e^{-t/\mathbf{T_1}}) ,
\end{align}
where the vector $\mathbf{T}_1$ denotes the longitudinal relaxation times for all pixels and $\operatorname{Re}(\mathbf{M}_0)$ and $\operatorname{Im}(\mathbf{M}_0)$ denote real and imaginary parts of the steady-state magnetization, respectively \cite{haase1990}.

Note that in quantitative MR imaging, one is ultimately interested in the quantities contained in the vector $\mathbf{u}$. However, often, qualitative images are first reconstructed (using some regularization method) as an intermediate step, from which then the vector $\mathbf{u}$ is estimated in a second step using non-linear regression methods, see for example \cite{tamir2020computational}. We can formulate the image reconstruction problem by 
\begin{equation}\label{eq:qmri_recon_problem}
\left \{
\begin{aligned}
\underset{\XX}{\min} 
\;\frac{1}{2} \|\Ad \XX - \ZZ\|_2^2 +    \| \boldsymbol{\Lambda}_{\Theta} \nabla \XX  \|_1,\\
\text{ s.t. } \XX_t = q_{t_i}(\mathbf{u})\, \quad 1\leq i\leq N_t.
\end{aligned}
\right.
\end{equation}
First, we train the proposed NN to estimate appropriate pixel-dependent regularization parameter-maps $\boldsymbol{\Lambda}_{\Theta}^{xy,t}$ to solve the TV-minimization problem and in a second step, perform a pixel-wise regression to obtain the vector $\mathbf{u}$.

\subsubsection{Experimental Set-Up}
We used the BrainWeb \cite{brainweb} dataset of 20 segmented healthy human heads as a basis to generate a quantitative MRI dataset with known ground truth. The subjects were split 17/1/2 for training, validation and testing. We considered axial slices, rescaled to $192\times192$ pixels. In each axial slice, we sampled for each tissue class from uniform distributions around anatomically plausible values the complex magnetization and the longitudinal relaxation rate $R_1=1/T_1$. The phase of the magnetization was further modulated by low amplitude 2D polynomials, approximating residual phases in the acquisition model. Following the signal model \eqref{eq:signal_model_qmri}, we generated images for the inversion times 0.05\,s, 0.1\,s, 0.2\,s, 0.35\,s, 0.5\,s, 1.0\,s, 1.5\,s, 2\,s, 3\,s, 4\,s and transformed them into (undersampled) Cartesian k-space. The number of simulated receiver coils was 8. The acceleration factor $R$ was chosen from 4, 6, and 8 for comparisons. In each case, complex Gaussian noise with $\sigma$ randomly chosen from [0.04, 0.4] was added in k-space. The proposed unrolled network was used to reconstruct the (qualitative) images at different inversion time points.

Similar to Section \ref{subsec:dyn_mri}, we choose a simple 3D U-Net with two downsampling steps, two 3D convolution layers for each encoder and decoder block with LeakyReLu activation, and 8 initial filters, resulting in only 97402 parameters. We used a scaled softplus activation, $\beta \phi(x / \beta)$ with $\beta=5$, for the final activation to keep the predicted regularization strength positive. We initialized the bias of the final convolution layer with -1 (empirically chosen) to stabilize training by starting at a low regularization.  We trained the network with AdamW \cite{adamw} (weight decay $10^{-4}$), cosine annealing learning rate schedule with linear warmup over one epoch with a maximum learning rate of $10^{-2}$, and a batch size of 4. The number of iterations of the unrolled PDHG is set to $T=32$ during warmup and $T=128$ for the rest of the training. Again, $\sigma$, $\tau$ and $\theta$ were trainable parameters.
Optimization was done by minimizing the MSE between the ground truth images and the obtained images after masking out non-brain regions. To find $\lambda_{\mathrm{P}}$ and $\lambda_{\tilde{\mathrm{P}}}$, grid searches, similar to those in the dynamic MRI case, were performed with a fixed number of $T=256$ iterations.
For evaluation, we used $T=256$ iterations of PDHG. We calculated the PSNR and SSIM of the reconstructed images. As a comparison method, we also performed standard iterative MRI reconstruction (CG with early stopping) without any TV regularization \cite{pruessmann2001advances}. We determined the optimal number of iterations based on the MSE to the ground truth images.
Finally, we performed a pixelwise regression on the reconstructed images $\XX$ using the Broyden–Fletcher–Goldfarb–Shanno (BFGS) algorithm, minimizing $\| |q_{t_i}(\mathbf{u})|-|\XX| \|_2^2$, to obtain $T_1$-maps and calculated the RMSE.

\subsubsection{Results}
\begin{figure}
\hspace{0.05\textwidth}CG-SENSE\hspace{0.04\textwidth}PDHG $\lambda^{xy,t}_{\tilde{\mathrm{P}}}$\hspace{0.035\textwidth}PDHG $\lambda^{xy,t}_{\mathrm{P}}$\hspace{0.03\textwidth}PDHG $\LLambda^{xy,t}_\Theta$ \hspace{0.02\textwidth}Target/ZF \hspace{0.05\textwidth}$\LLambda^{xy}_\Theta$ / $\LLambda^{t}_\Theta$ 

\centering
	\rotatebox[origin=c]{90}{Example 1}\,\,\begin{subfigure}[p]{\textwidth}
		\includegraphics[width=0.15\textwidth]{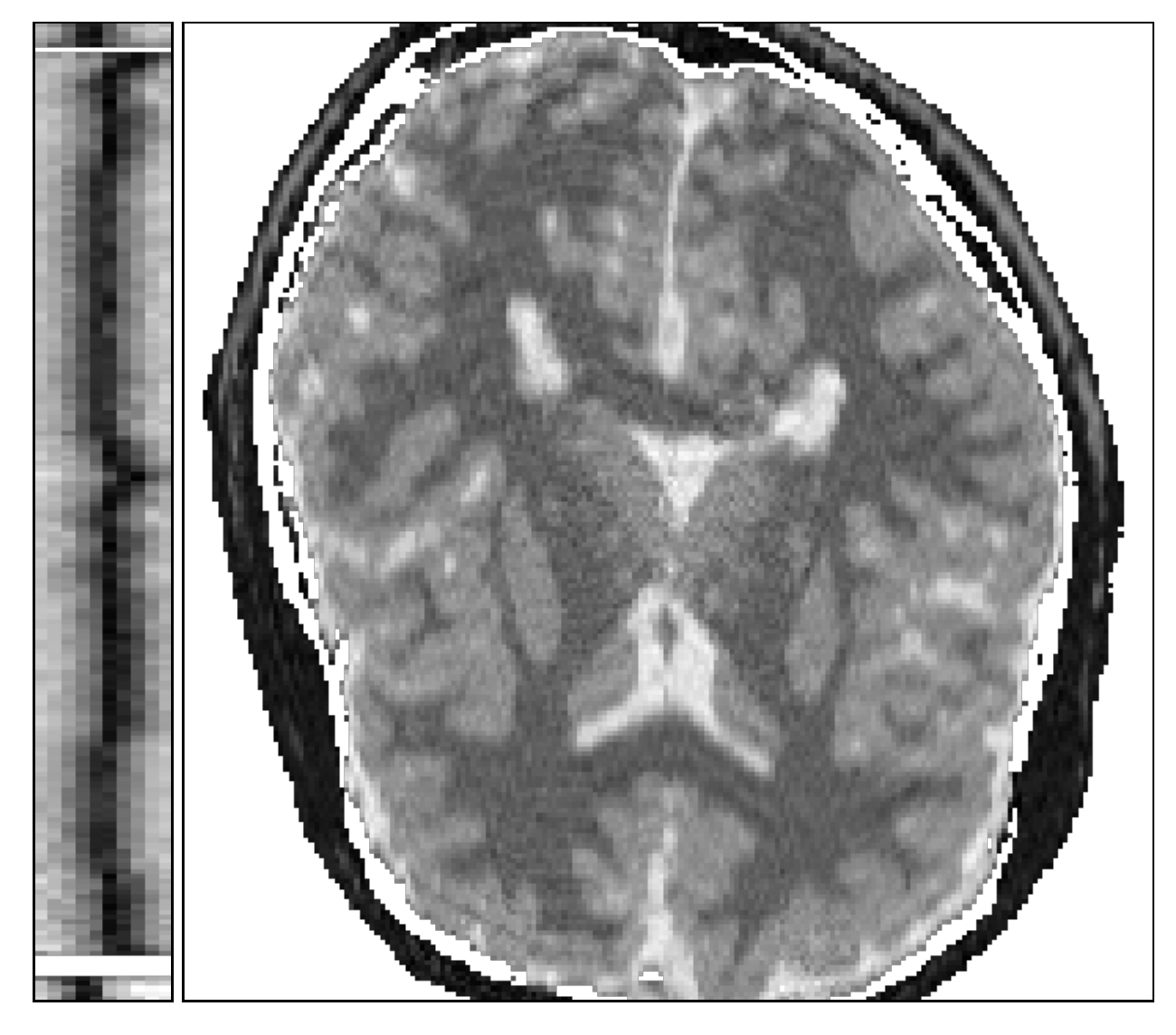}
		\includegraphics[width=0.15\textwidth]{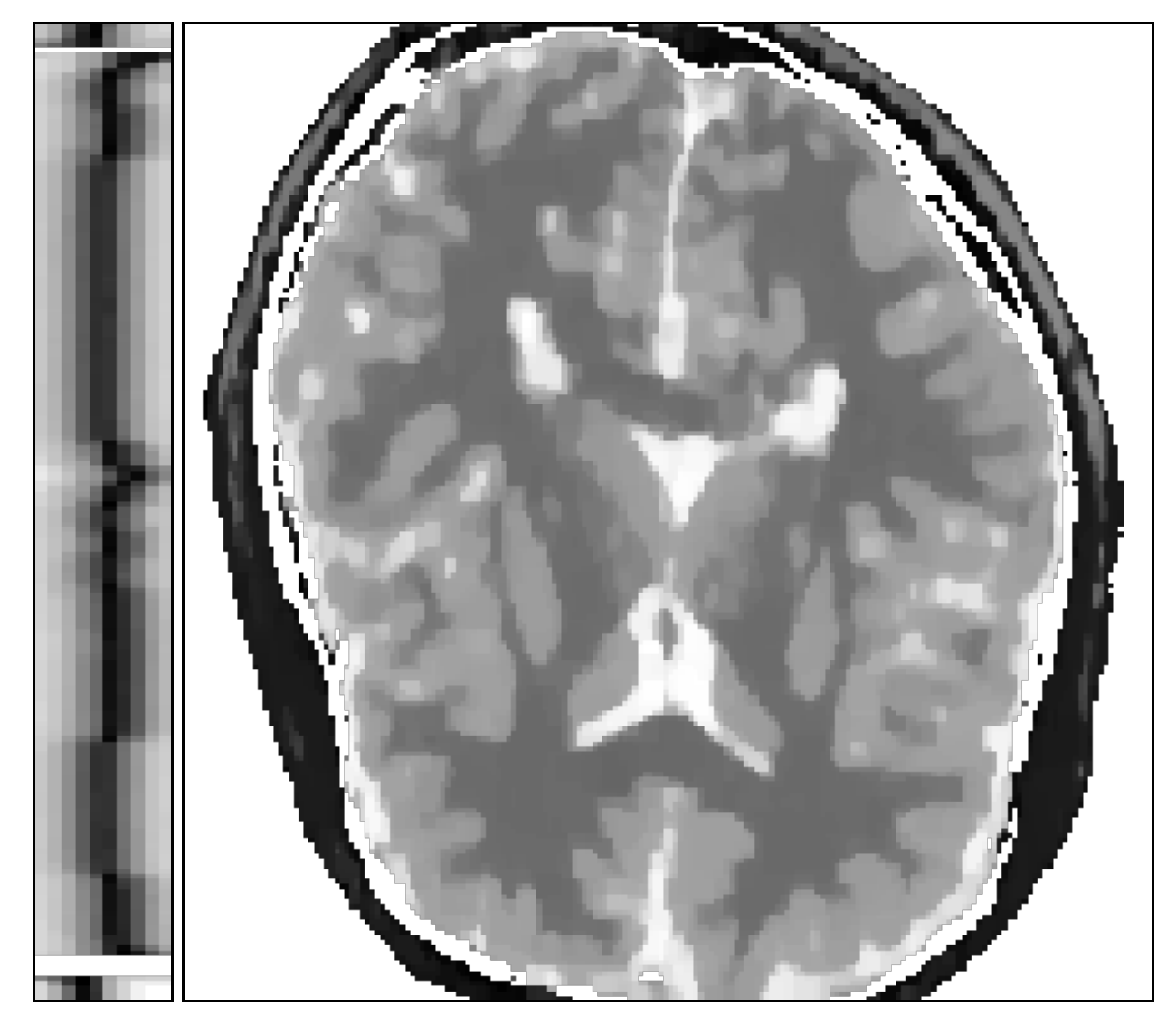}
		\includegraphics[width=0.15\textwidth]{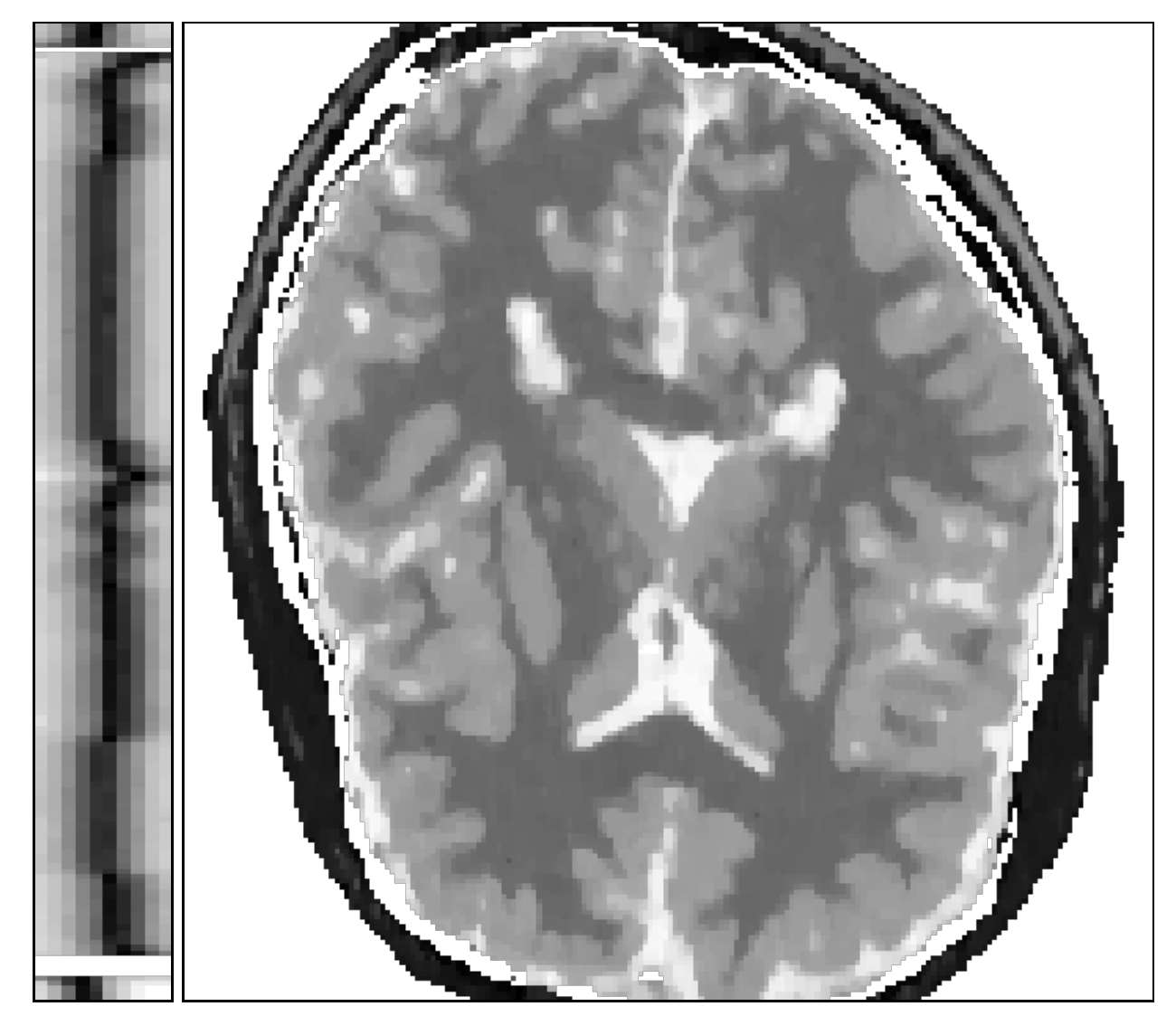}
		\includegraphics[width=0.15\textwidth]{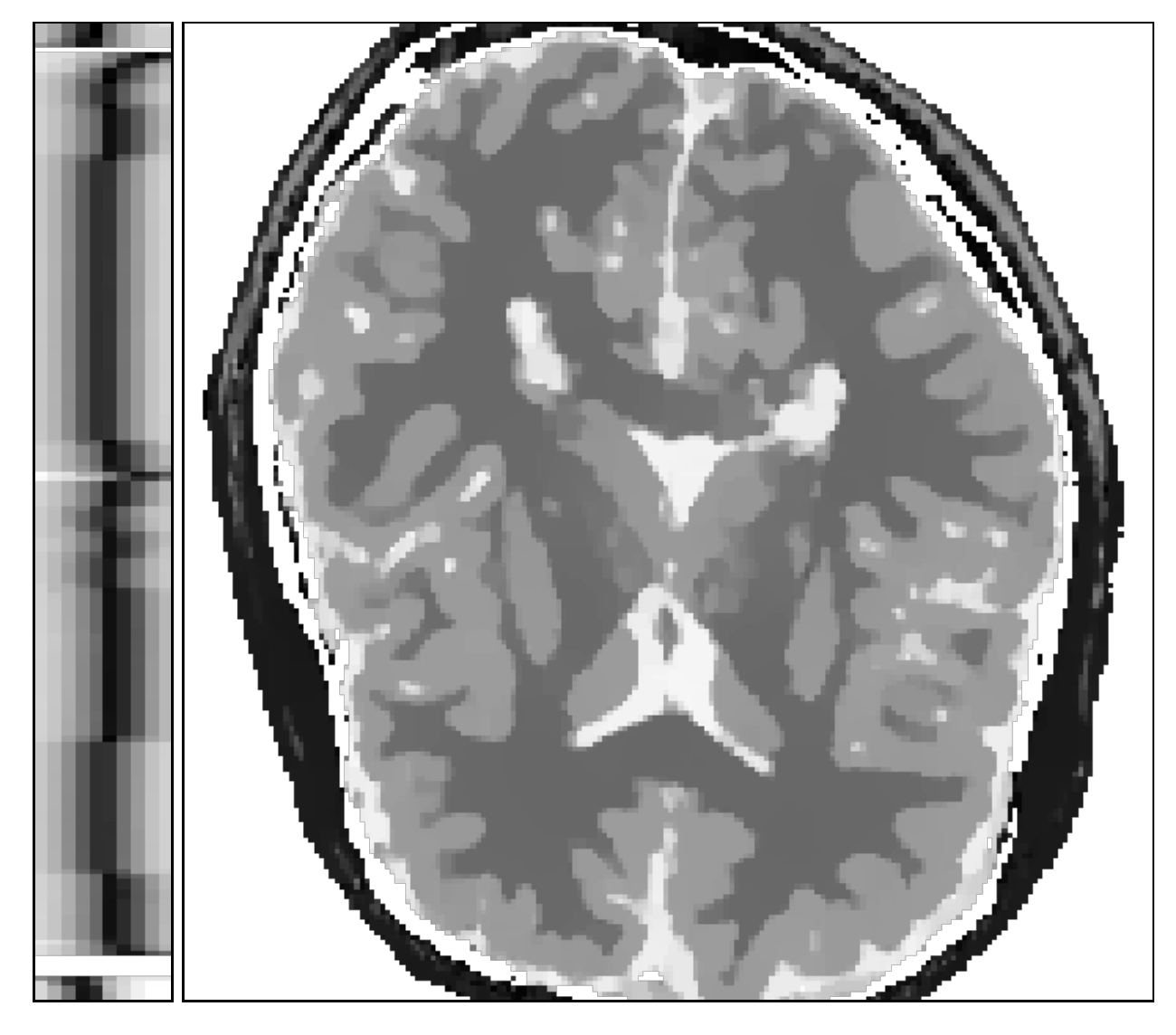}
		\includegraphics[width=0.15\textwidth]{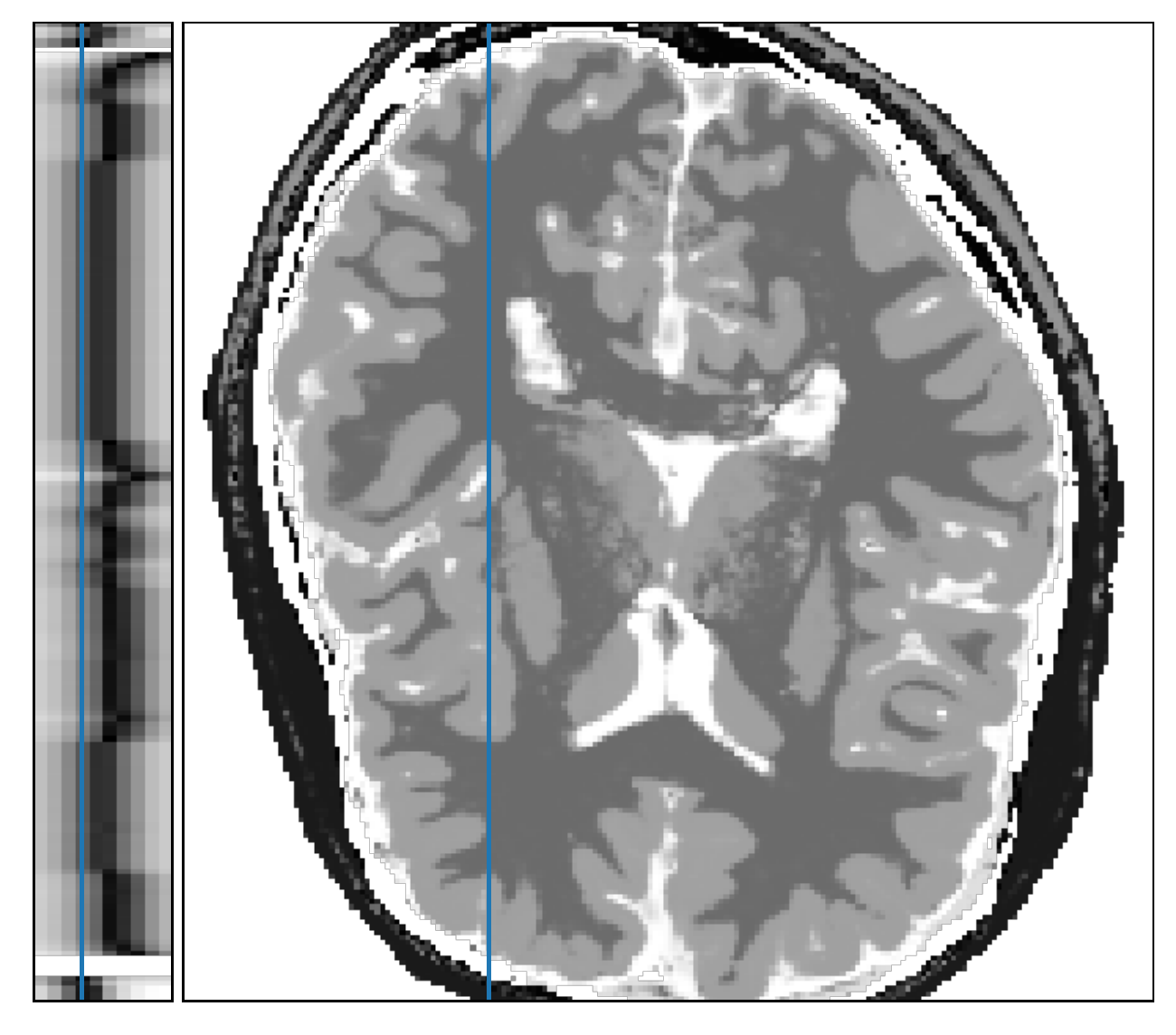}
		\includegraphics[width=0.15\textwidth]{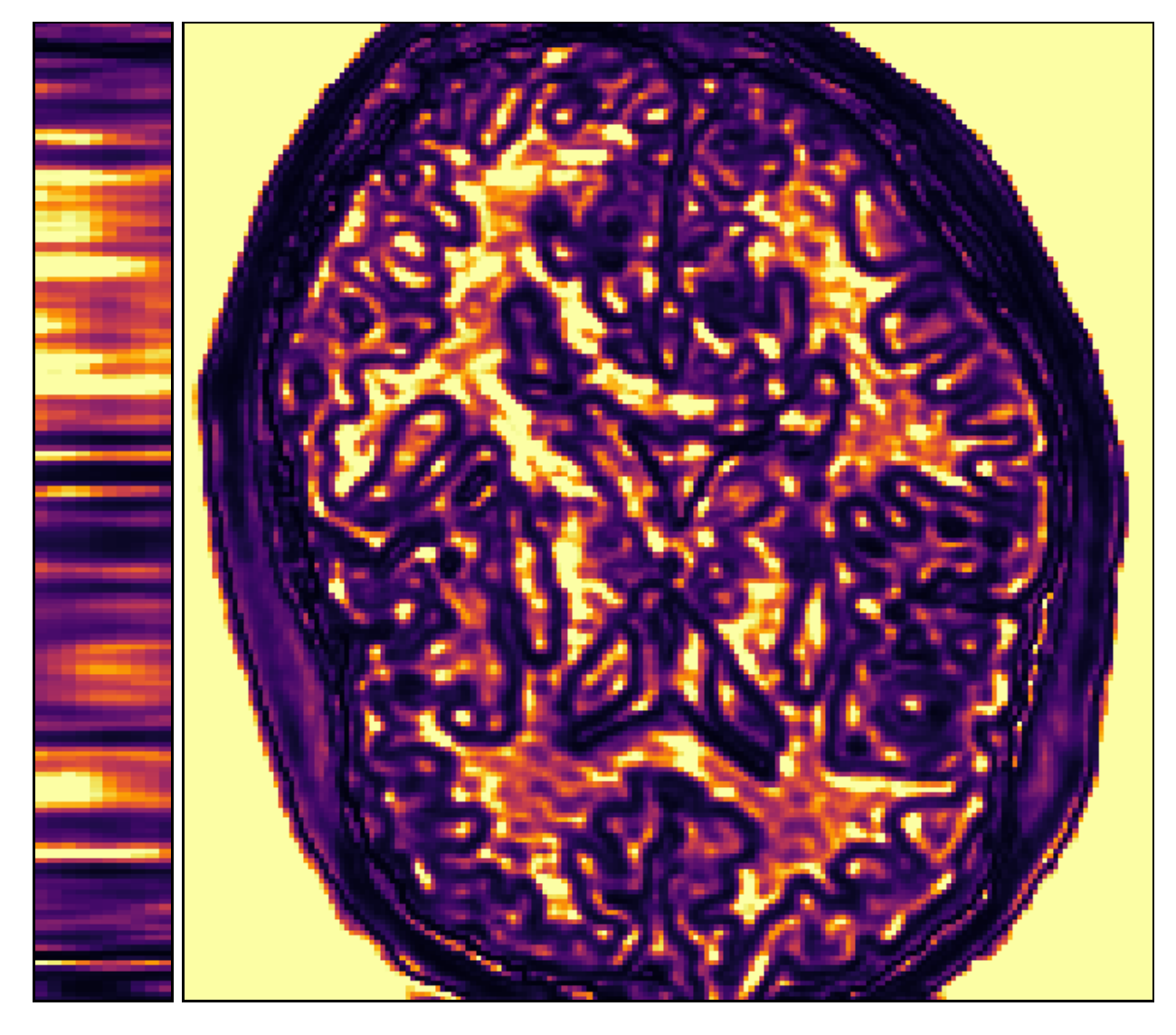}
		\vskip -3pt
		\includegraphics[width=0.15\textwidth]{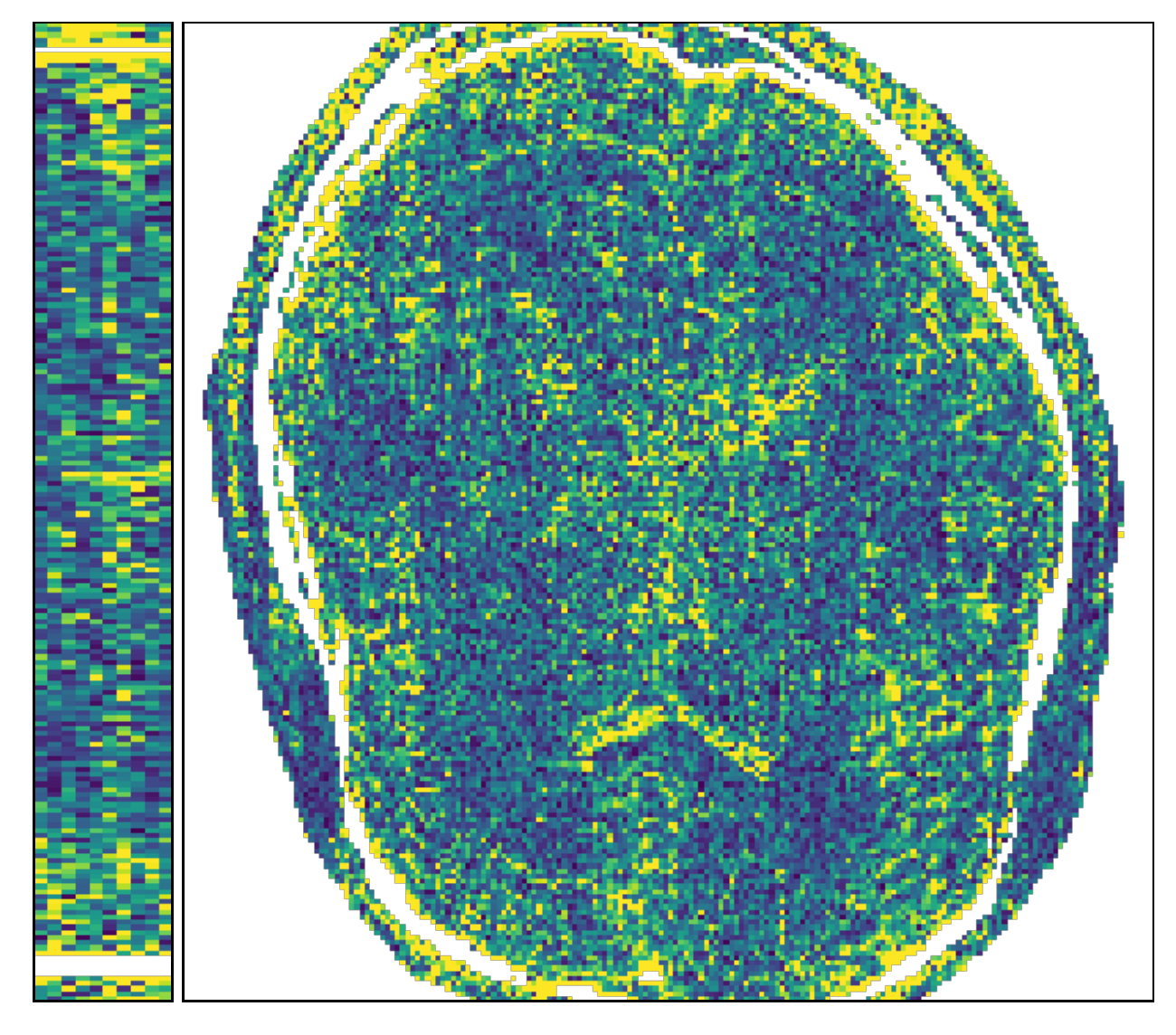}
		\includegraphics[width=0.15\textwidth]{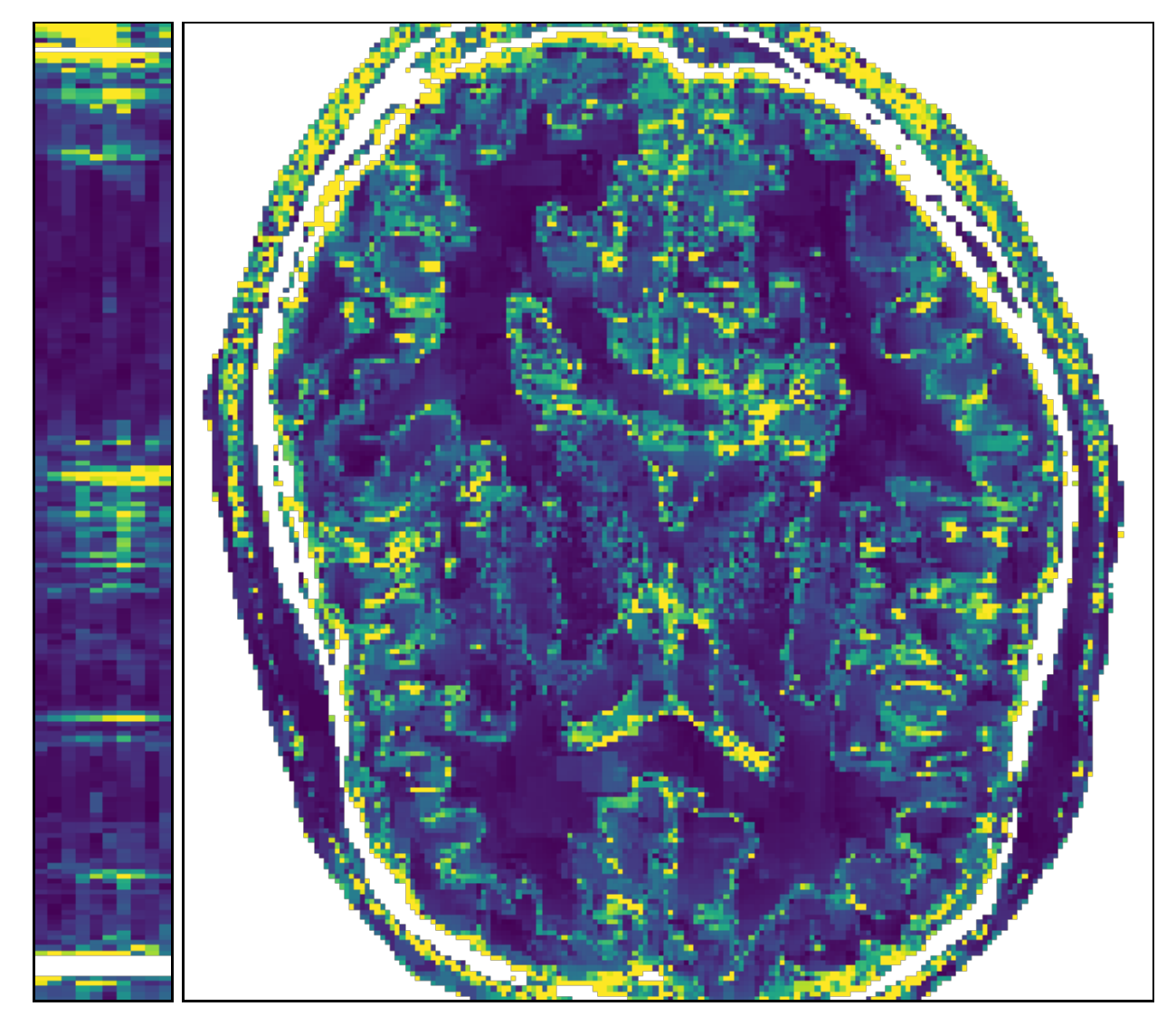}
		\includegraphics[width=0.15\textwidth]{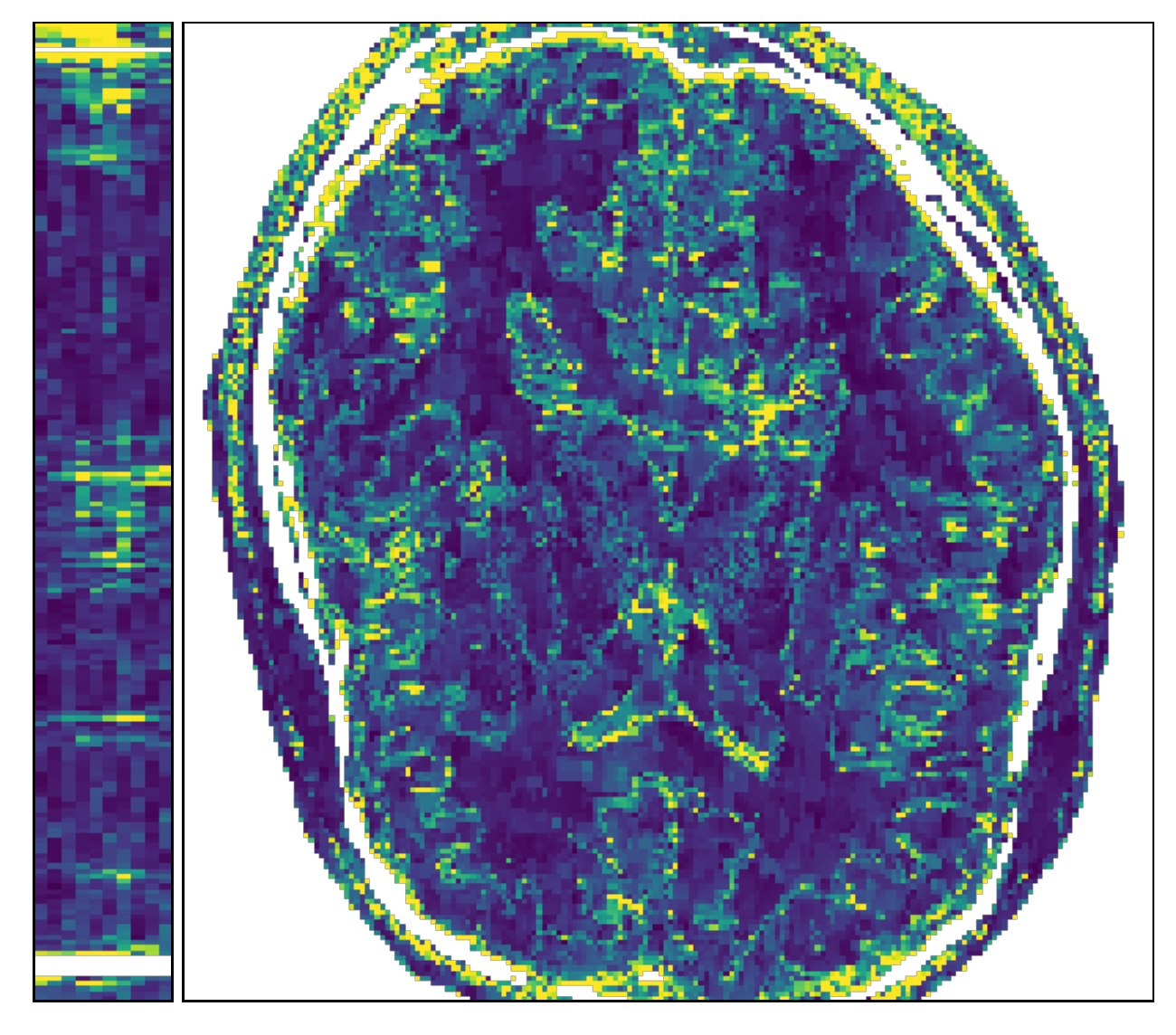}
		\includegraphics[width=0.15\textwidth]{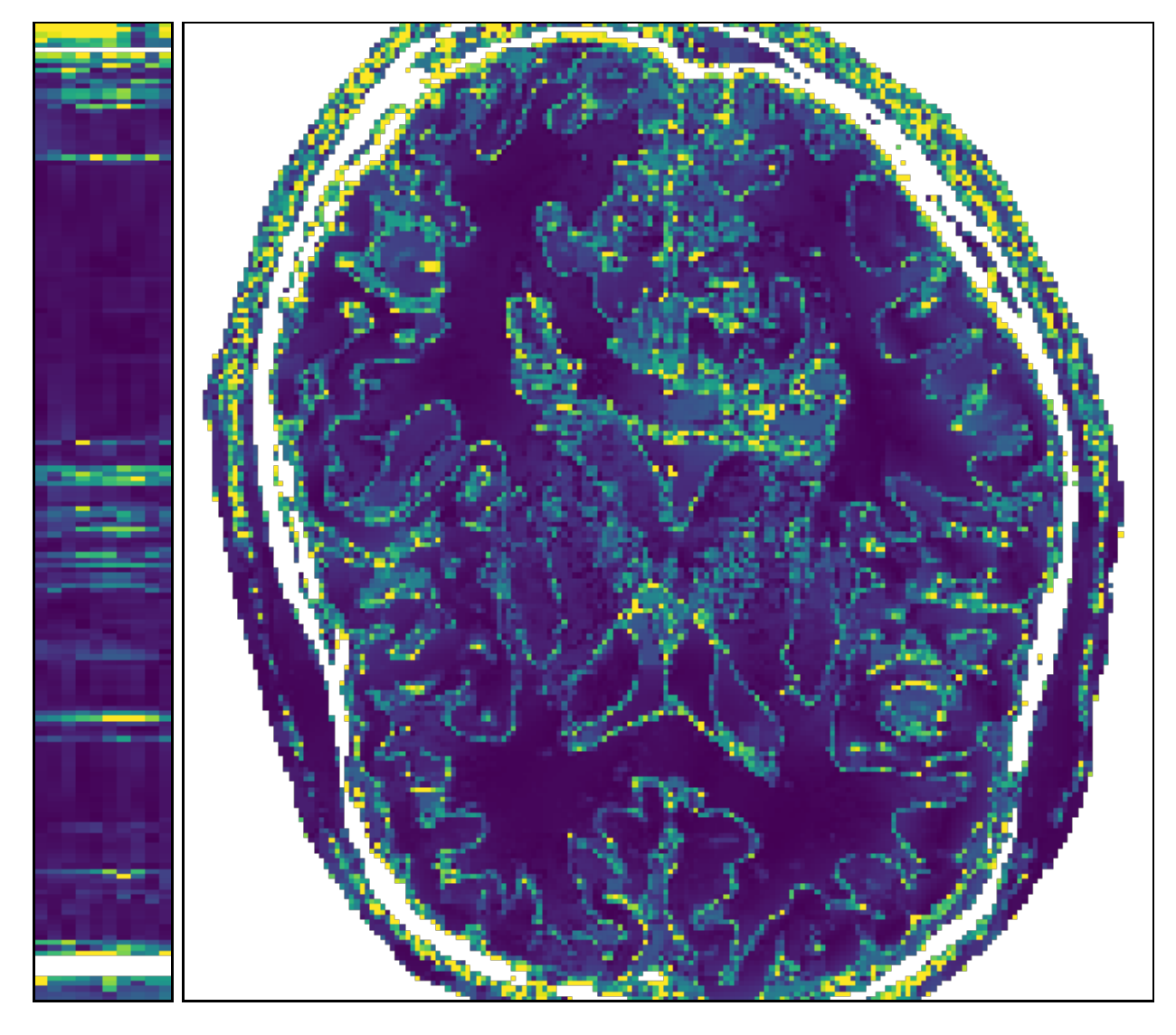}
		\includegraphics[width=0.15\textwidth]{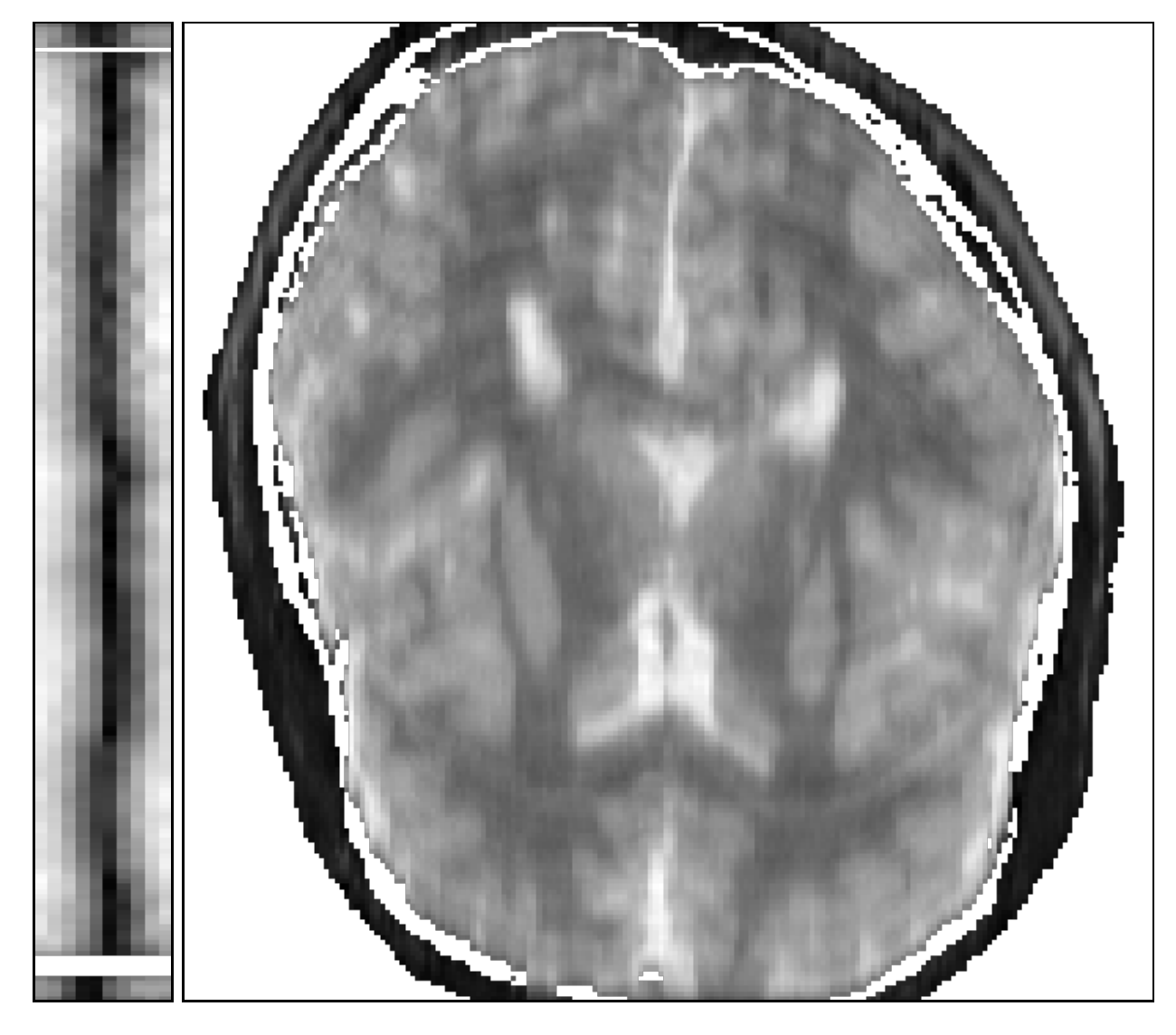}
		\includegraphics[width=0.15\textwidth]{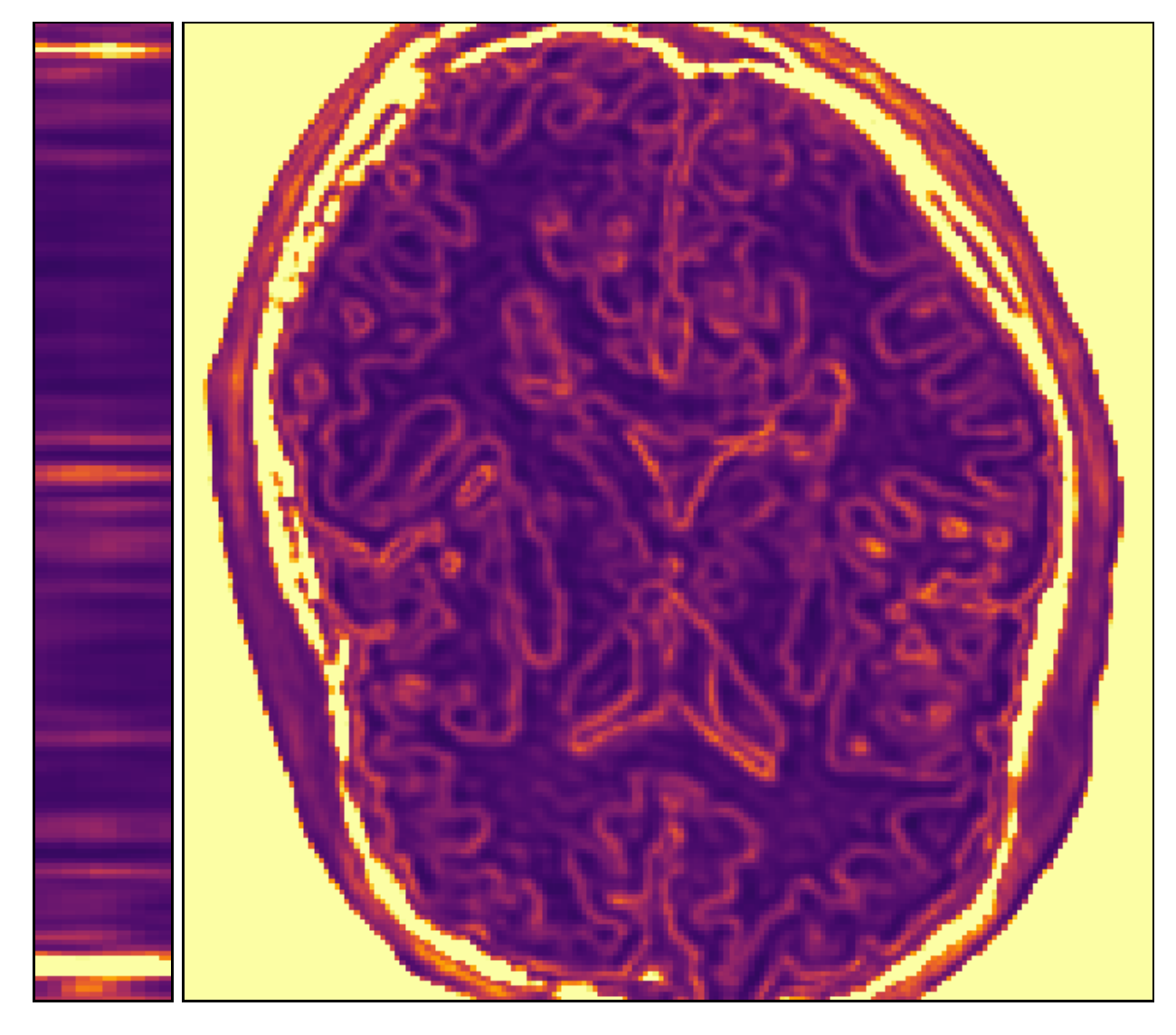}
	\end{subfigure}
	
	\vspace{0.01\textwidth}
	\rotatebox[origin=c]{90}{Example 2}\,\,\begin{subfigure}[p]{\textwidth}
    	\includegraphics[width=0.15\textwidth]{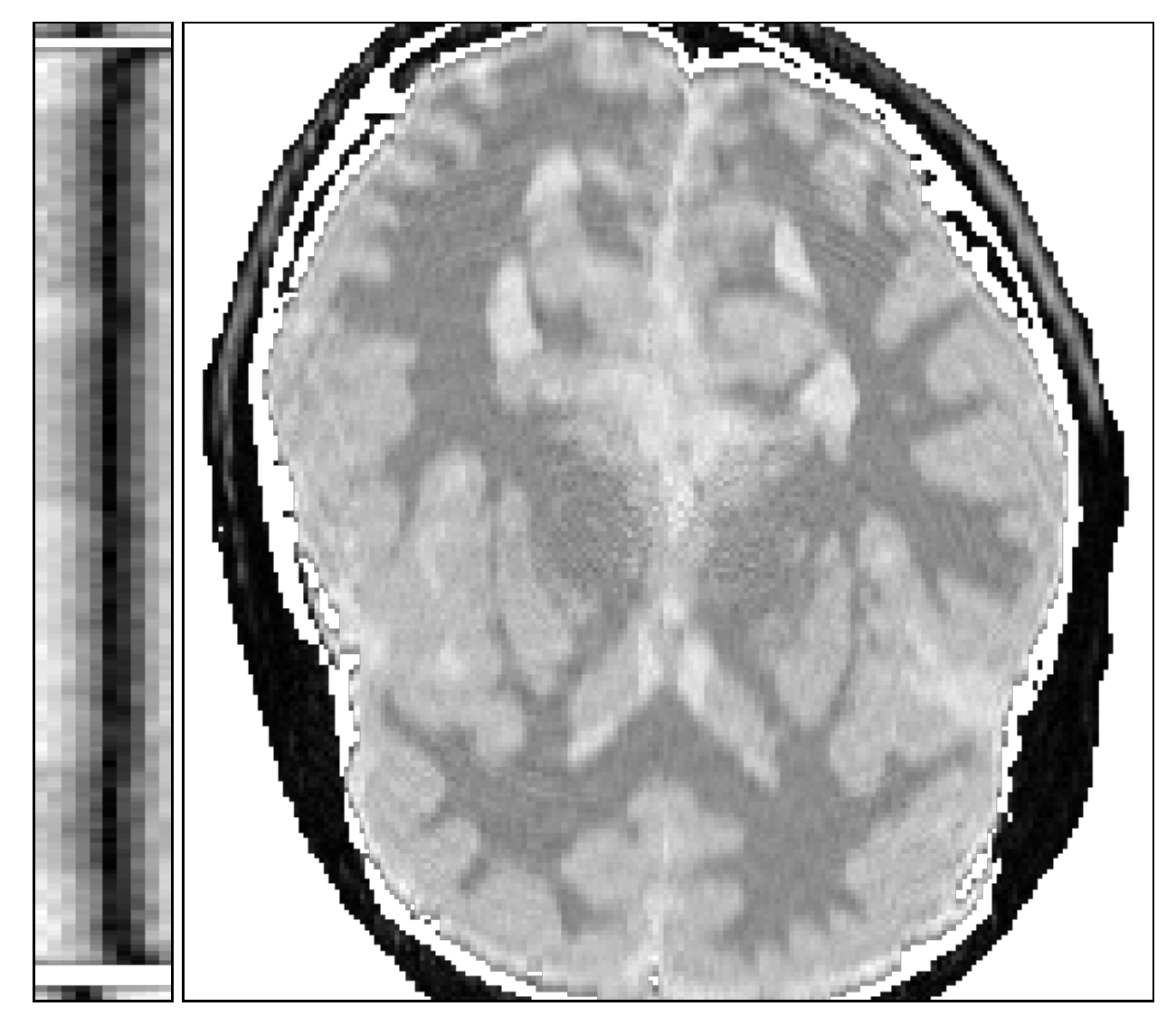} 
    	\includegraphics[width=0.15\textwidth]{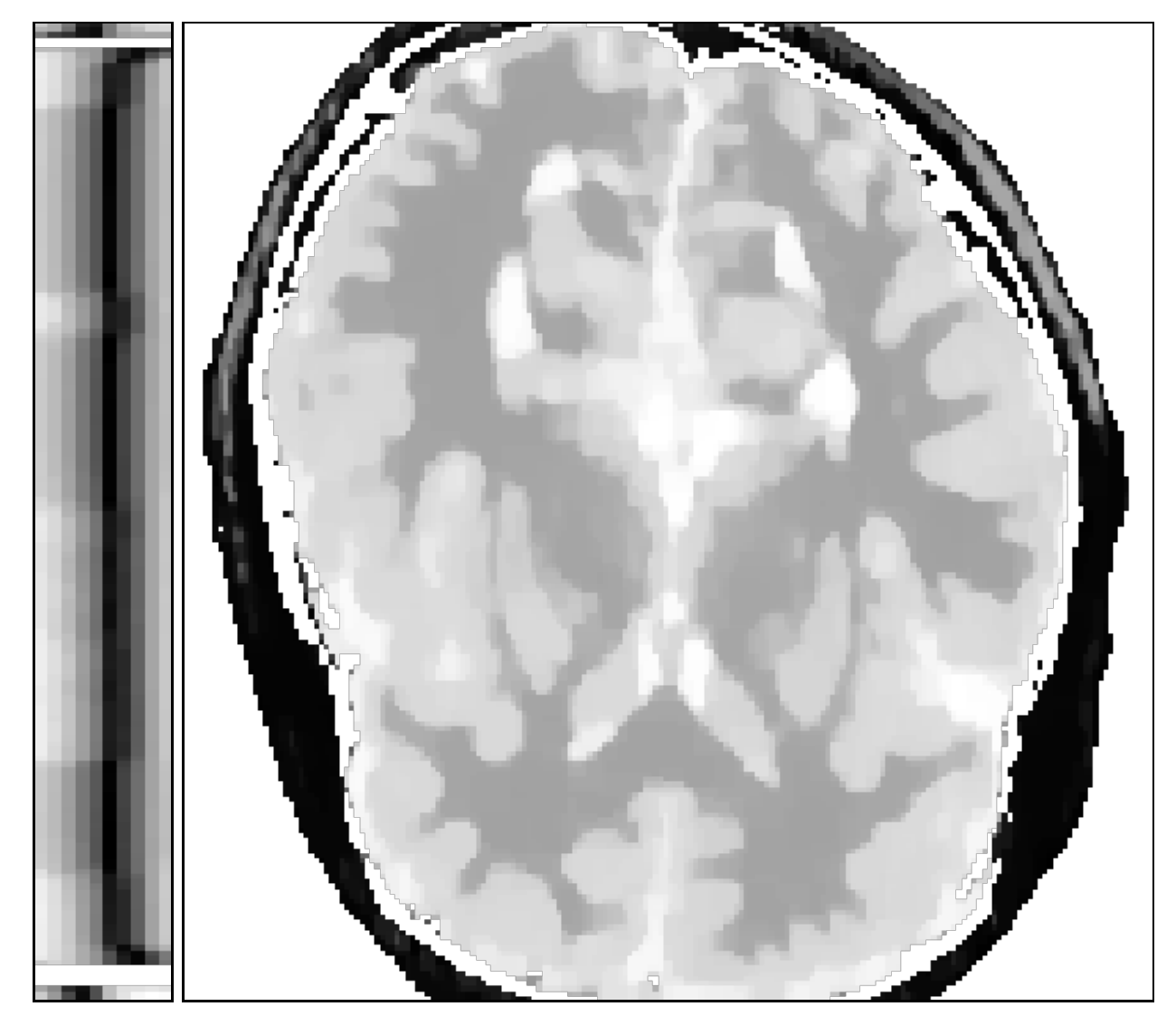} 
    	\includegraphics[width=0.15\textwidth]{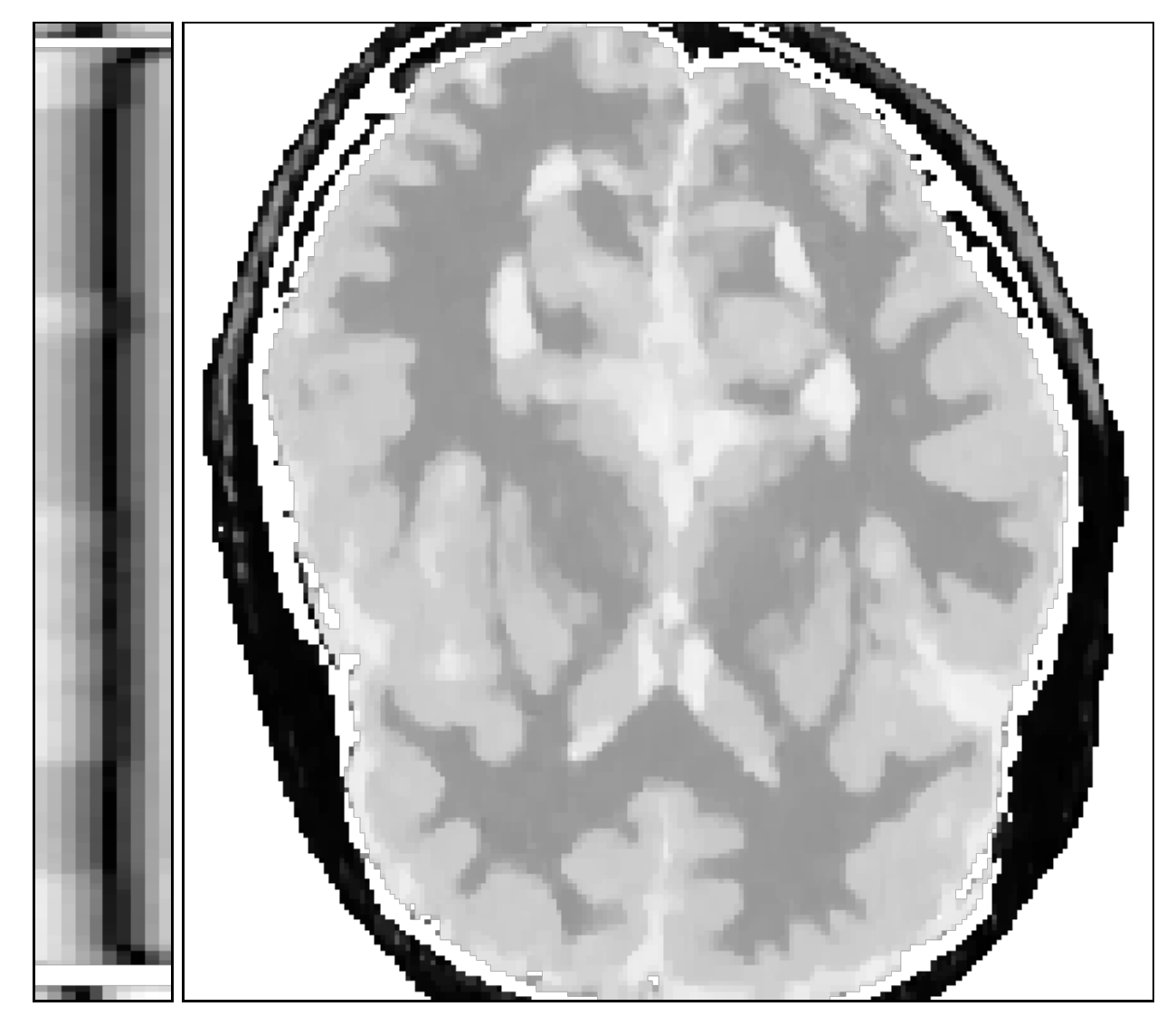} 
    	\includegraphics[width=0.15\textwidth]{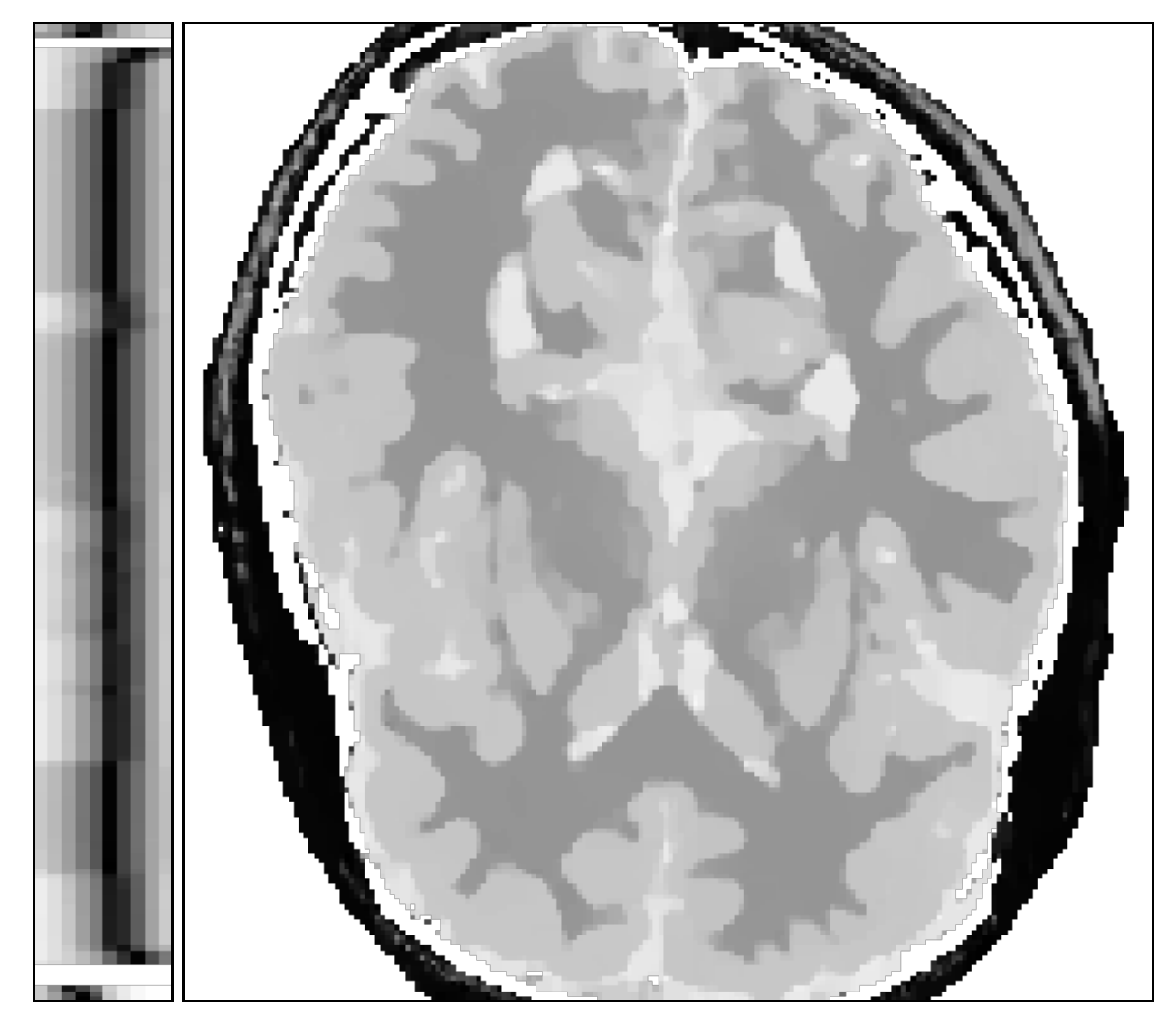} 
    	\includegraphics[width=0.15\textwidth]{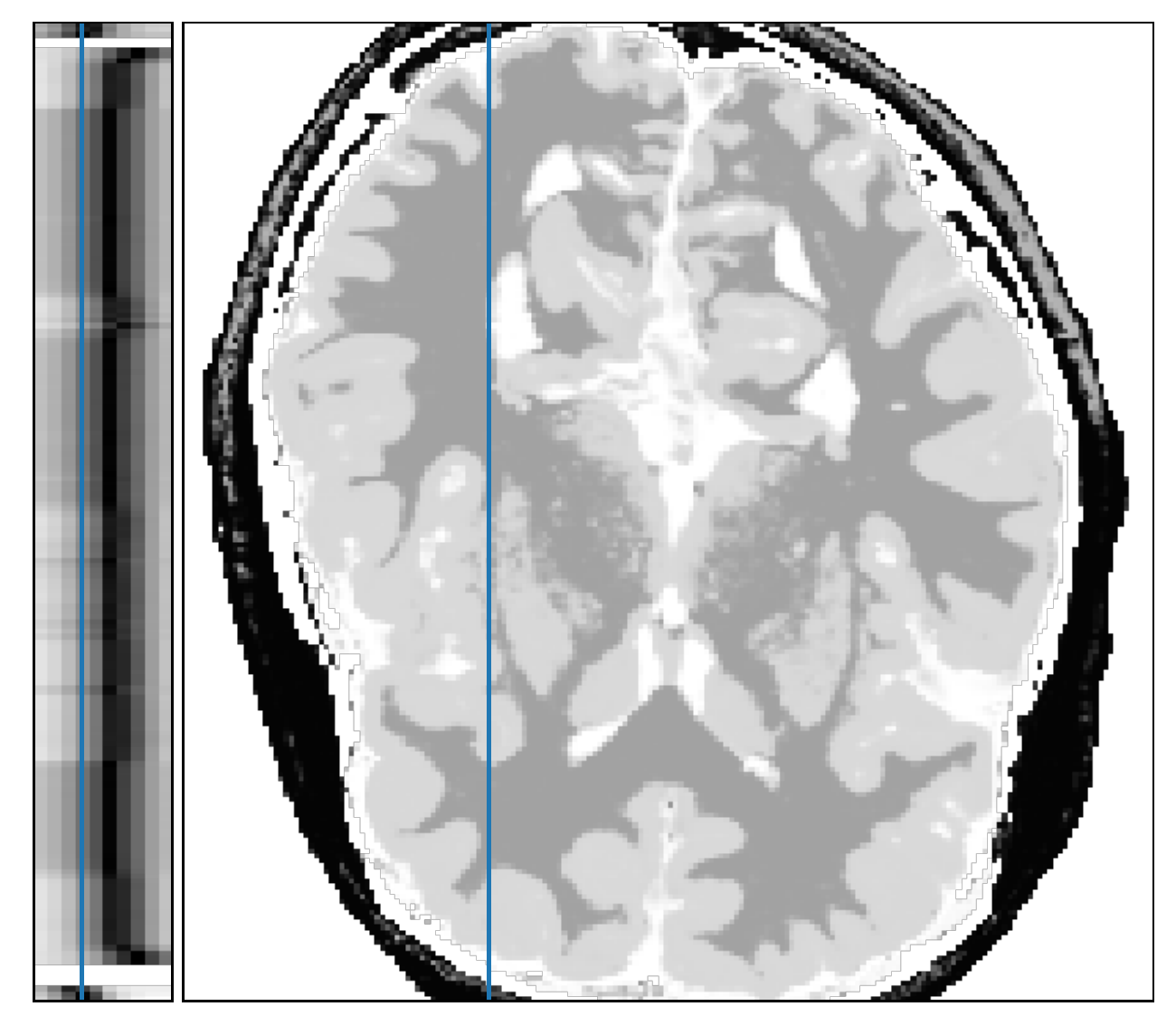} 
    	\includegraphics[width=0.15\textwidth]{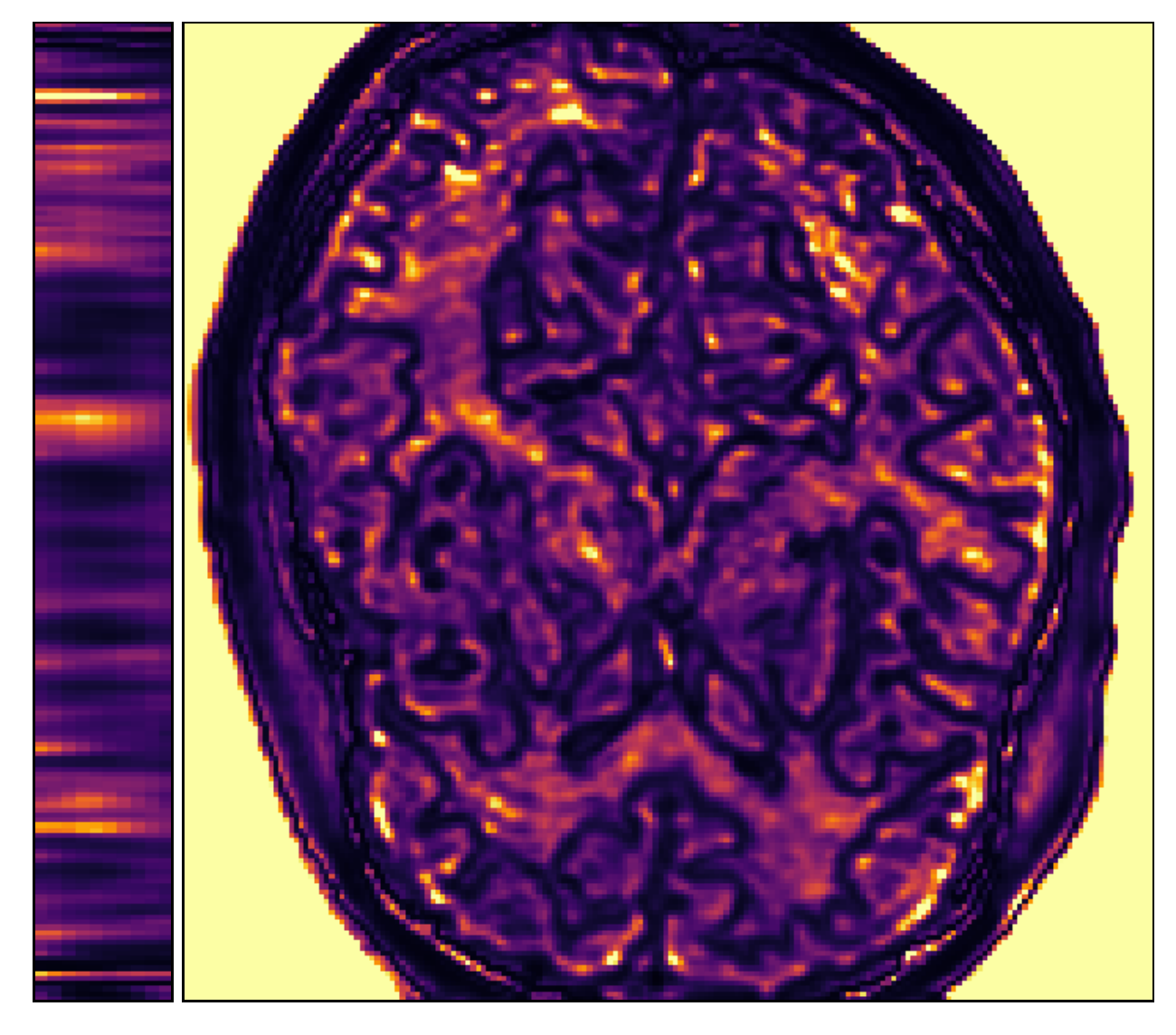}
    	\vskip -3pt
    	\includegraphics[width=0.15\textwidth]{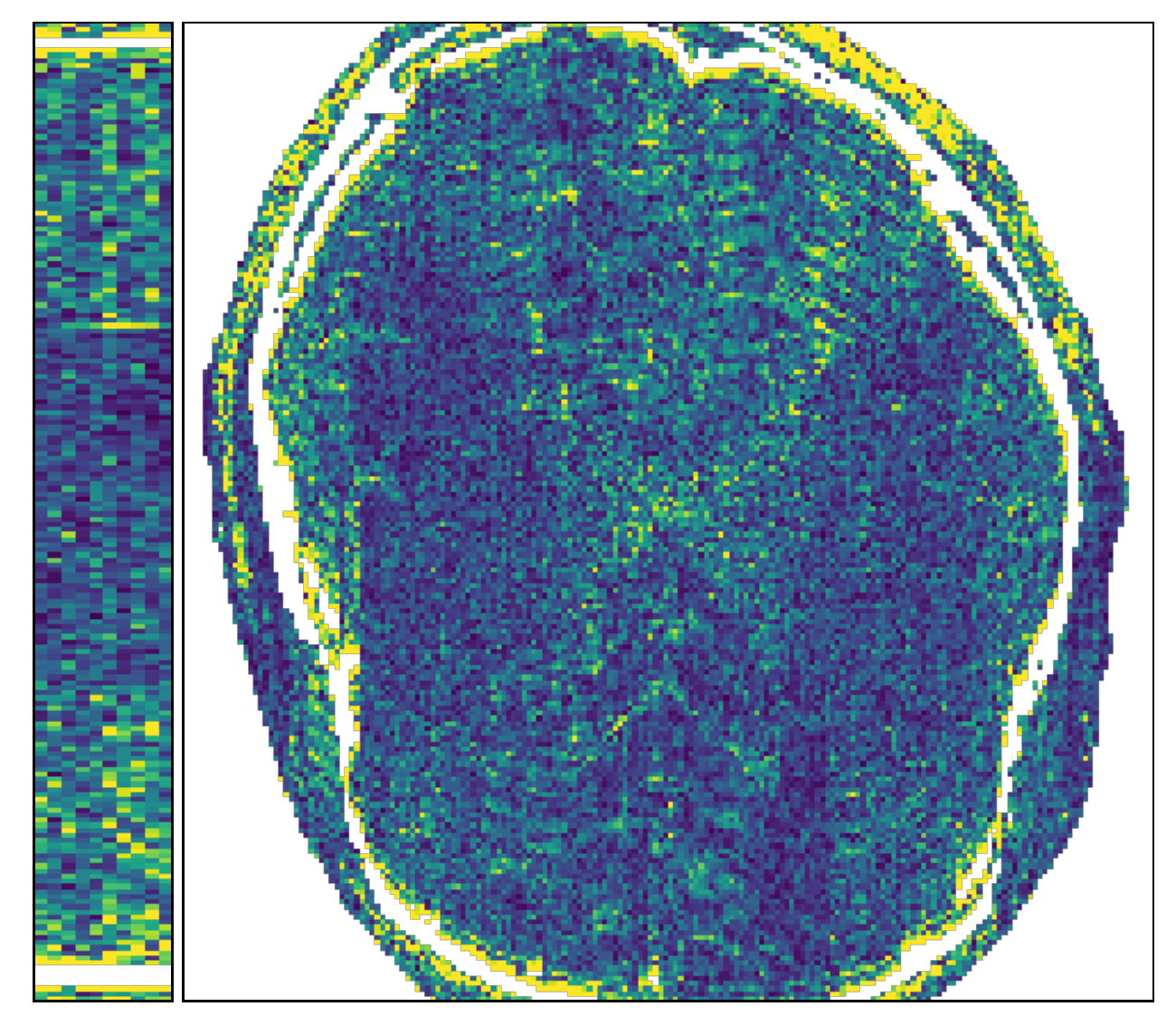} 
    	\includegraphics[width=0.15\textwidth]{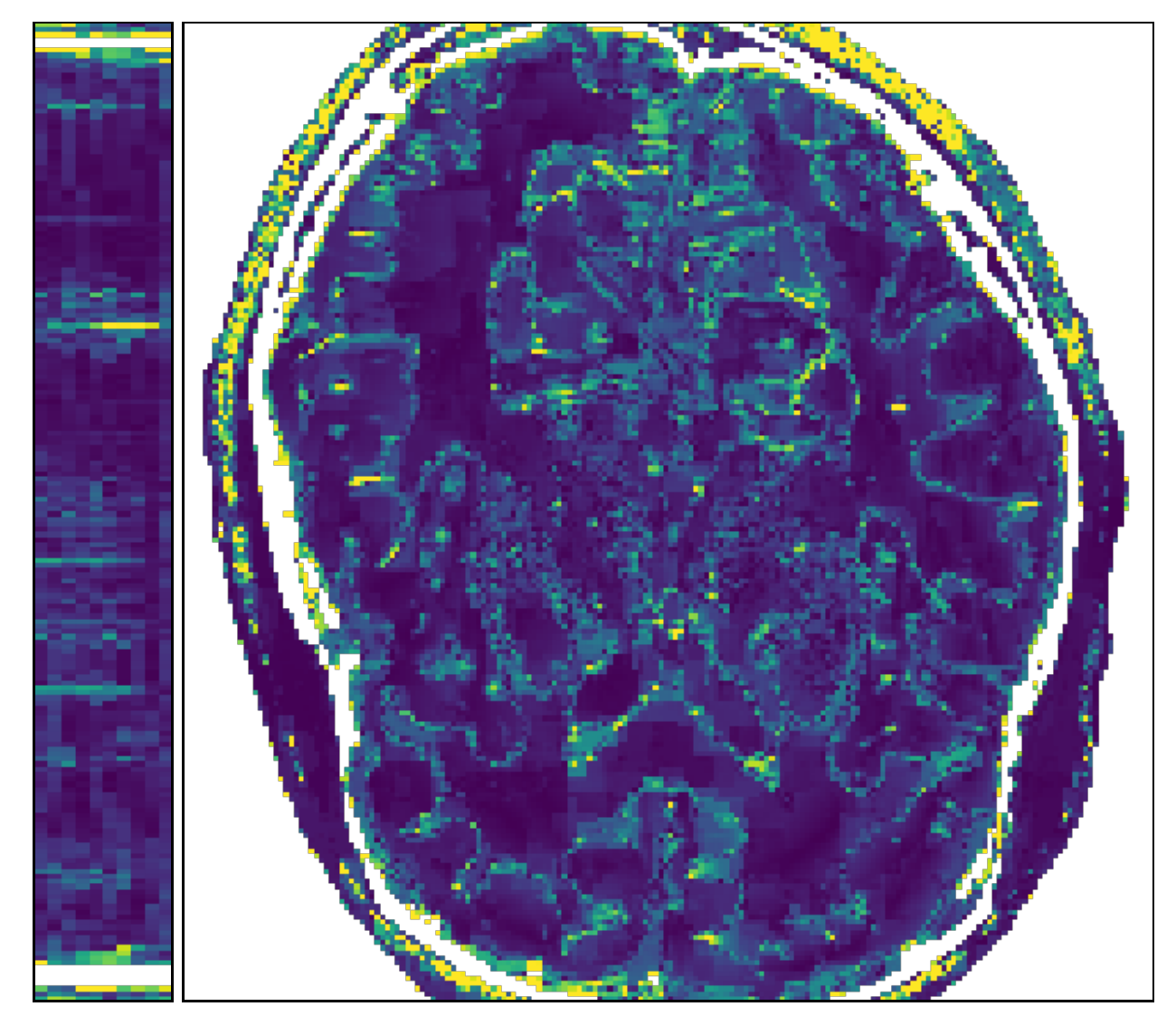} 
    	\includegraphics[width=0.15\textwidth]{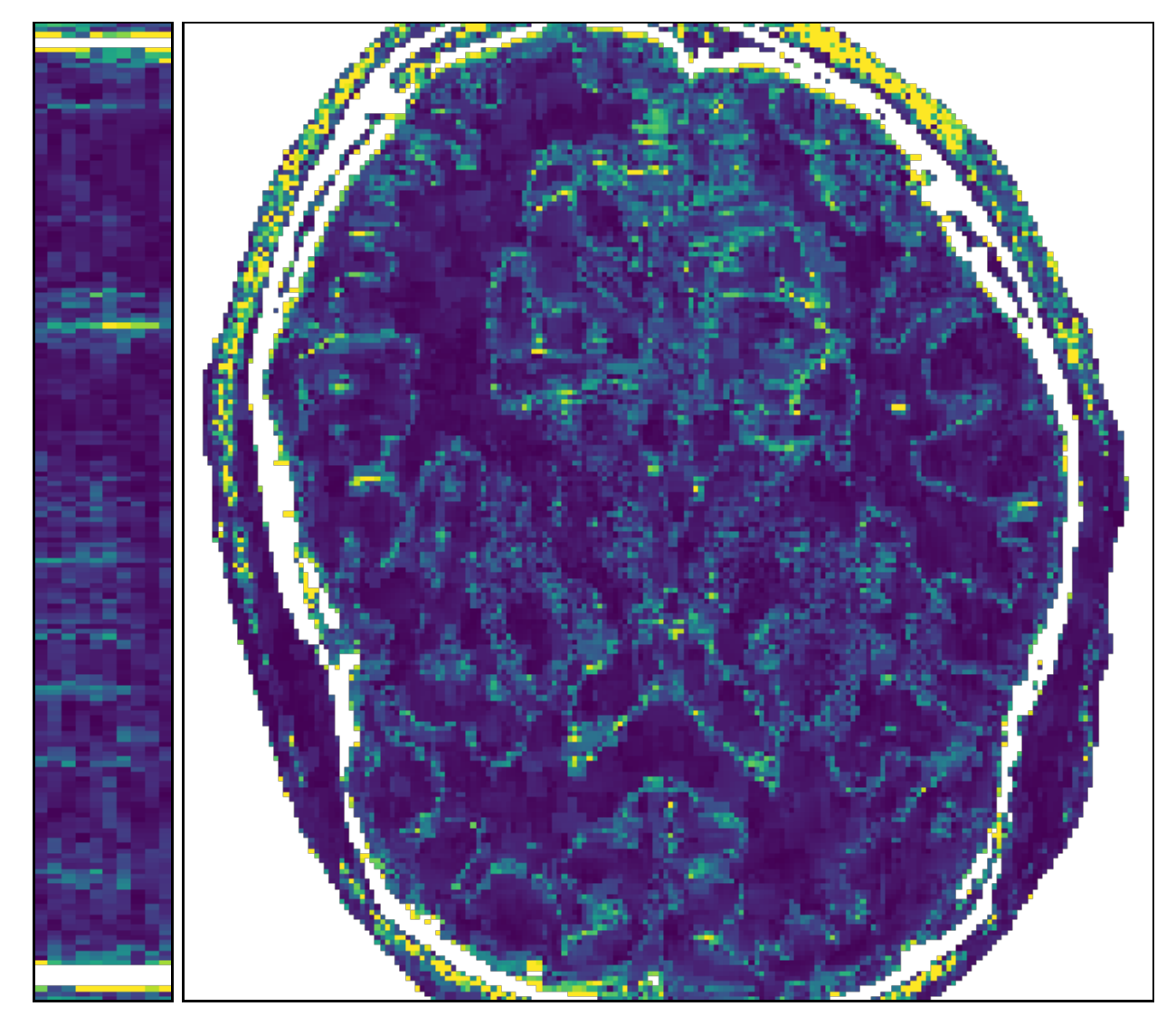} 
    	\includegraphics[width=0.15\textwidth]{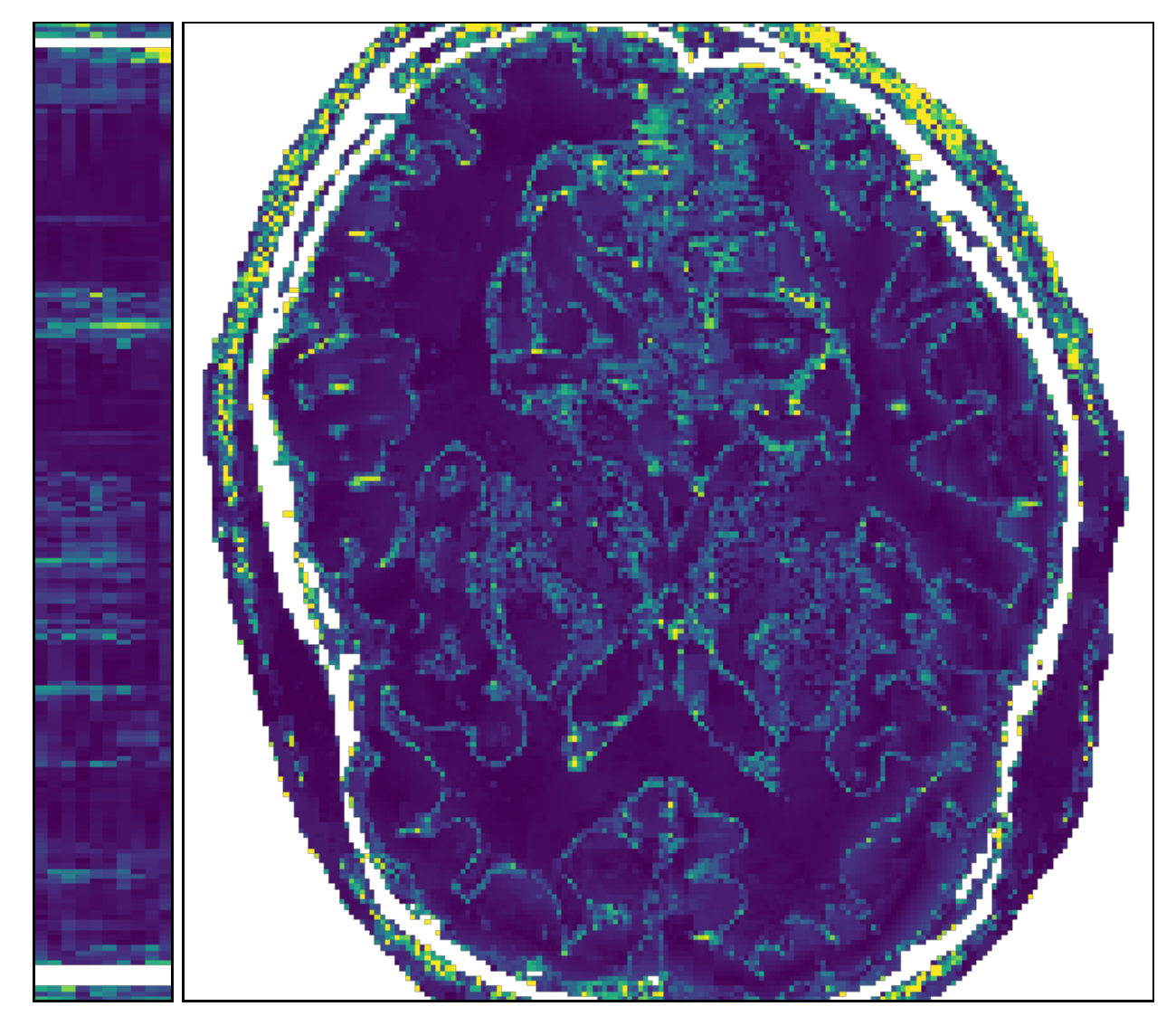}
    	\includegraphics[width=0.15\textwidth]{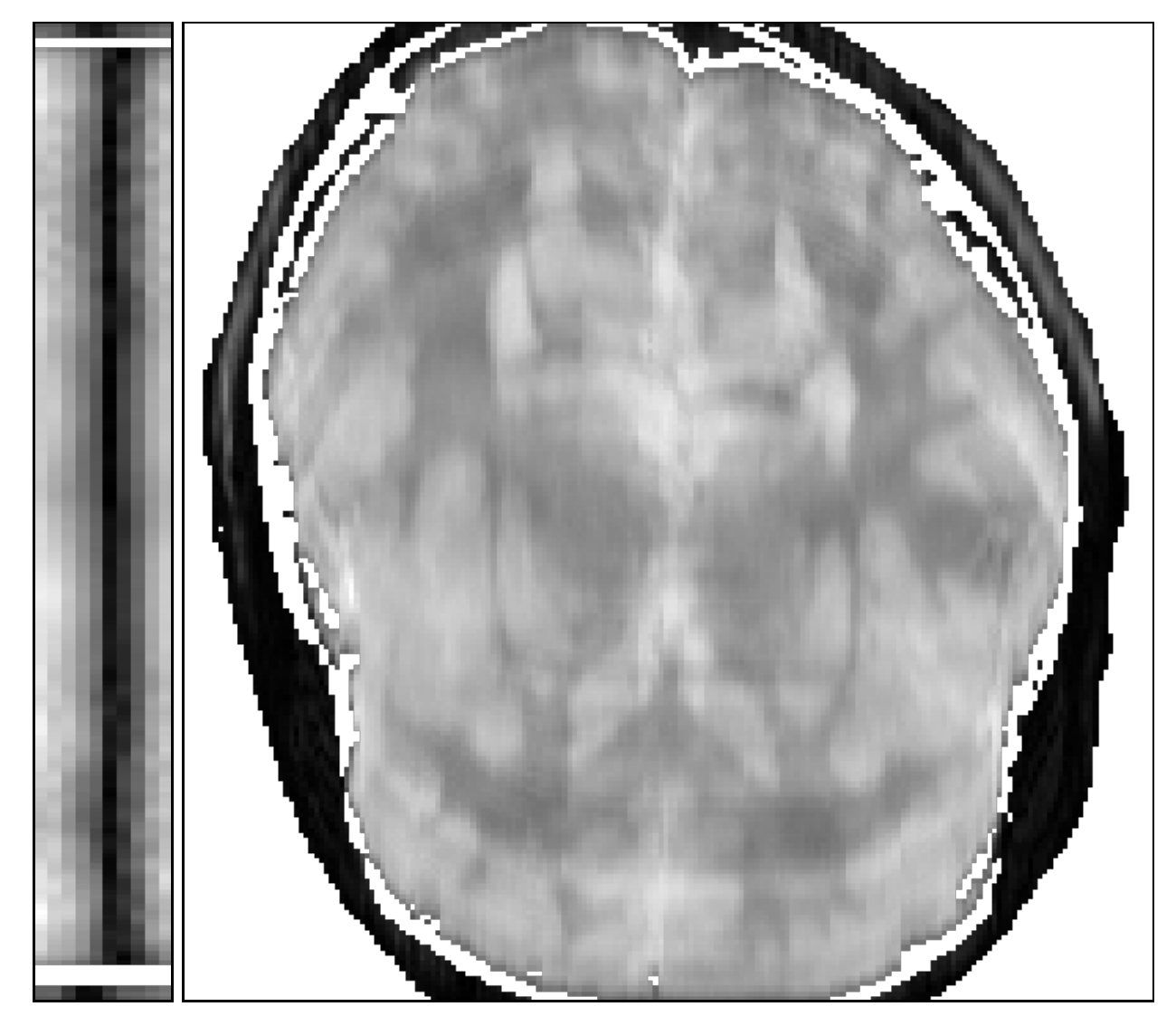} 
    	\includegraphics[width=0.15\textwidth]{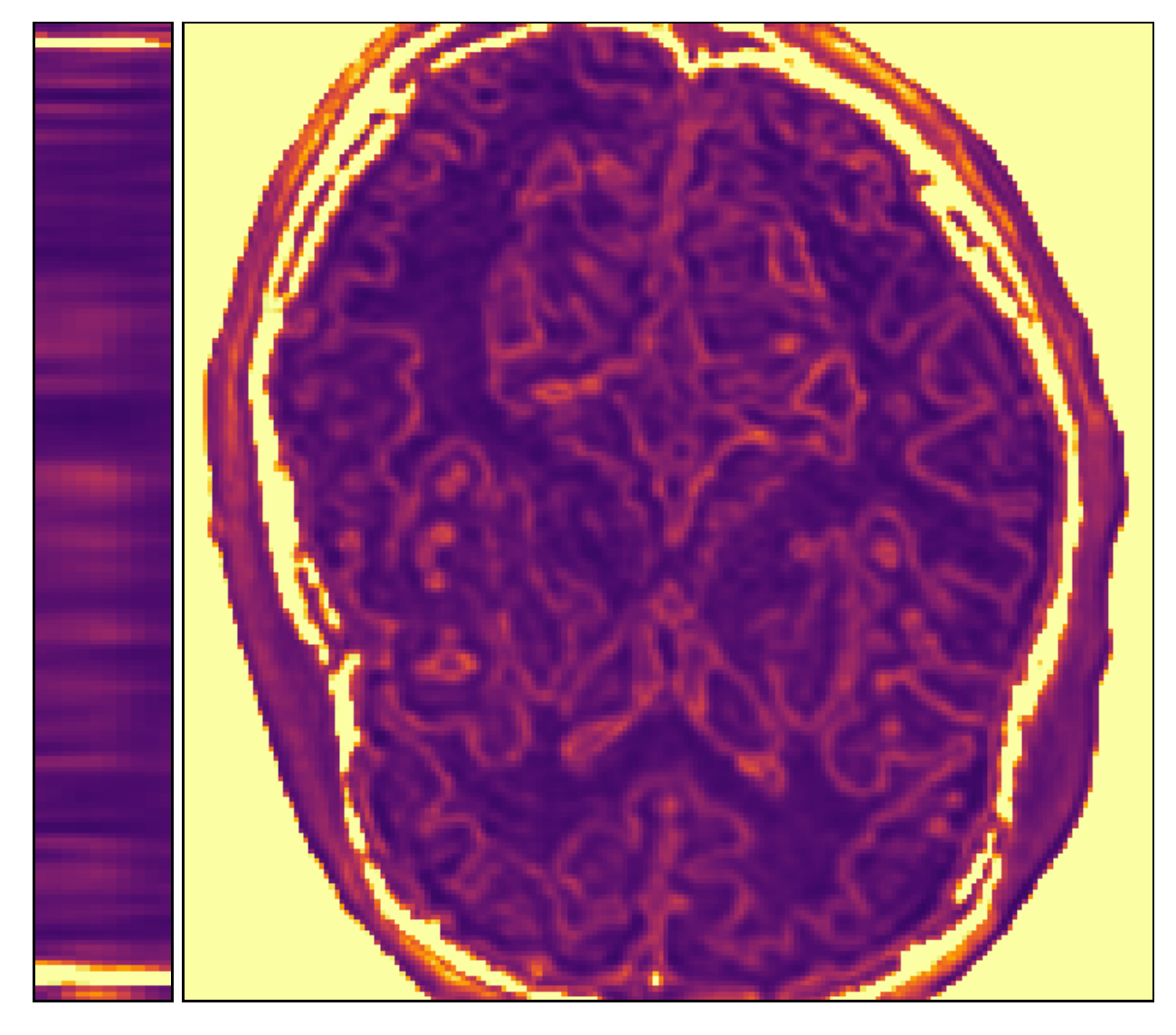}
    \end{subfigure}

    \vspace{0.01\textwidth}
    \rotatebox[origin=c]{90}{Example 3}\,\,\begin{subfigure}[p]{\textwidth}
        \includegraphics[width=0.15\textwidth]{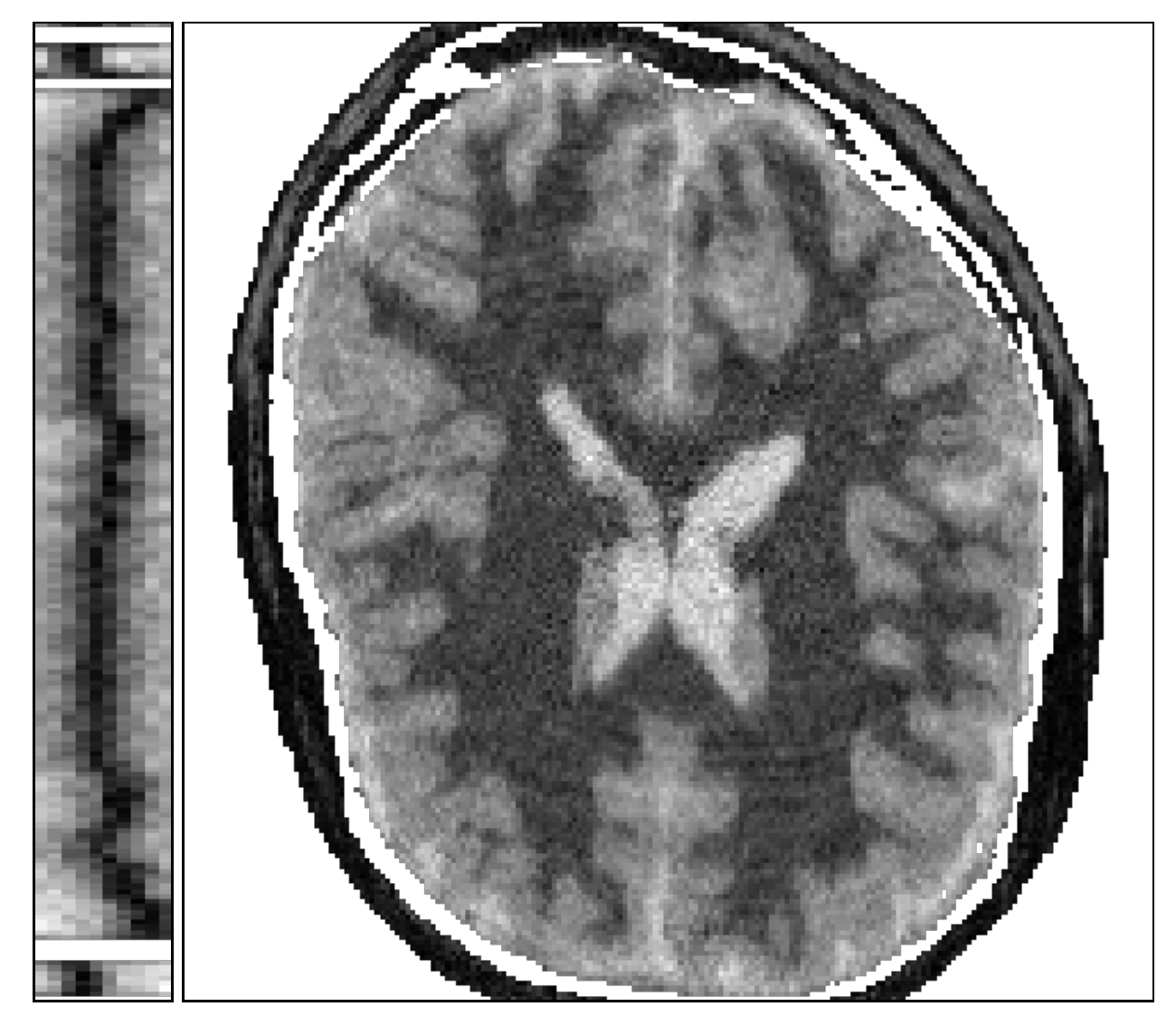} 
        \includegraphics[width=0.15\textwidth]{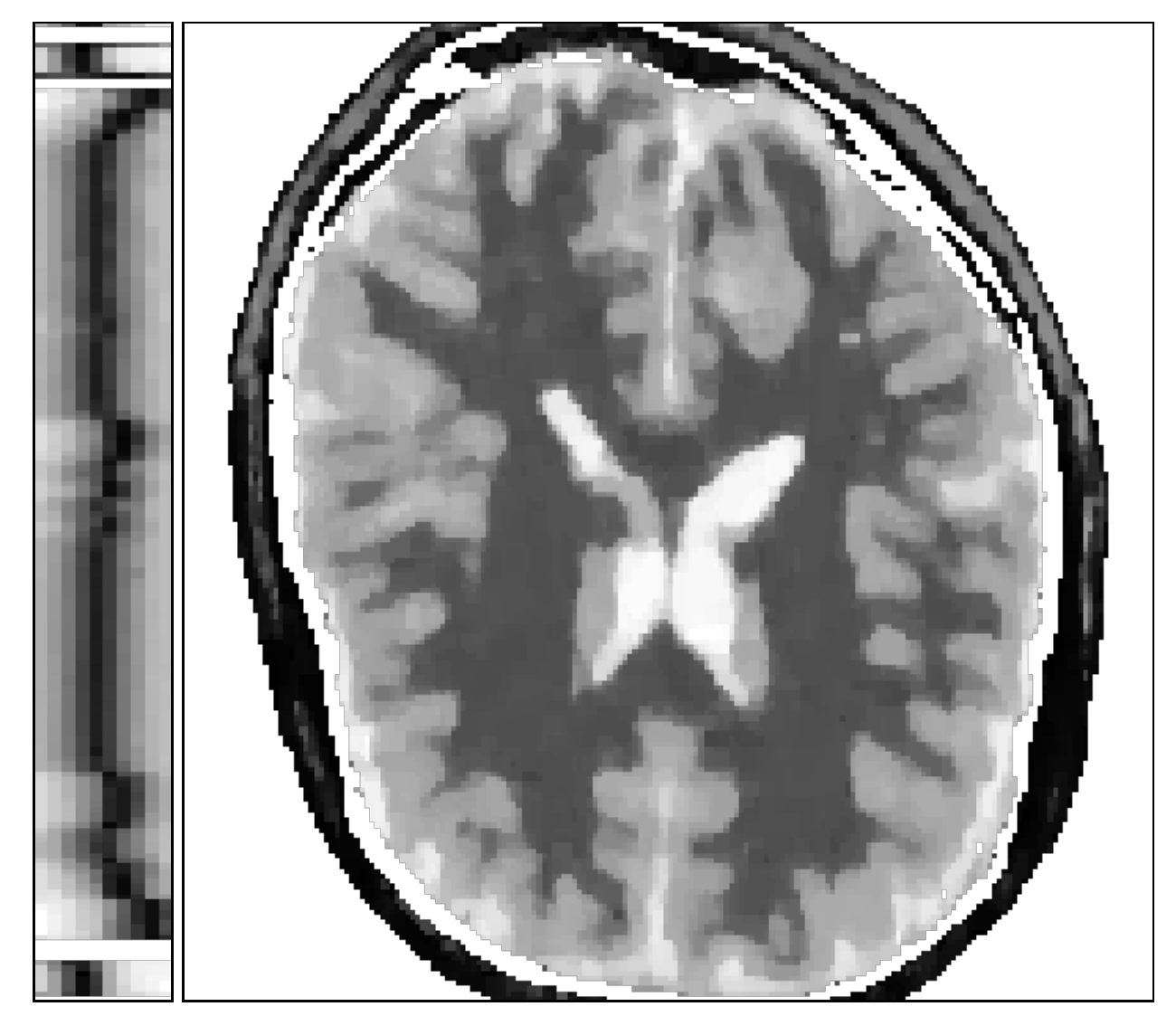} 
        \includegraphics[width=0.15\textwidth]{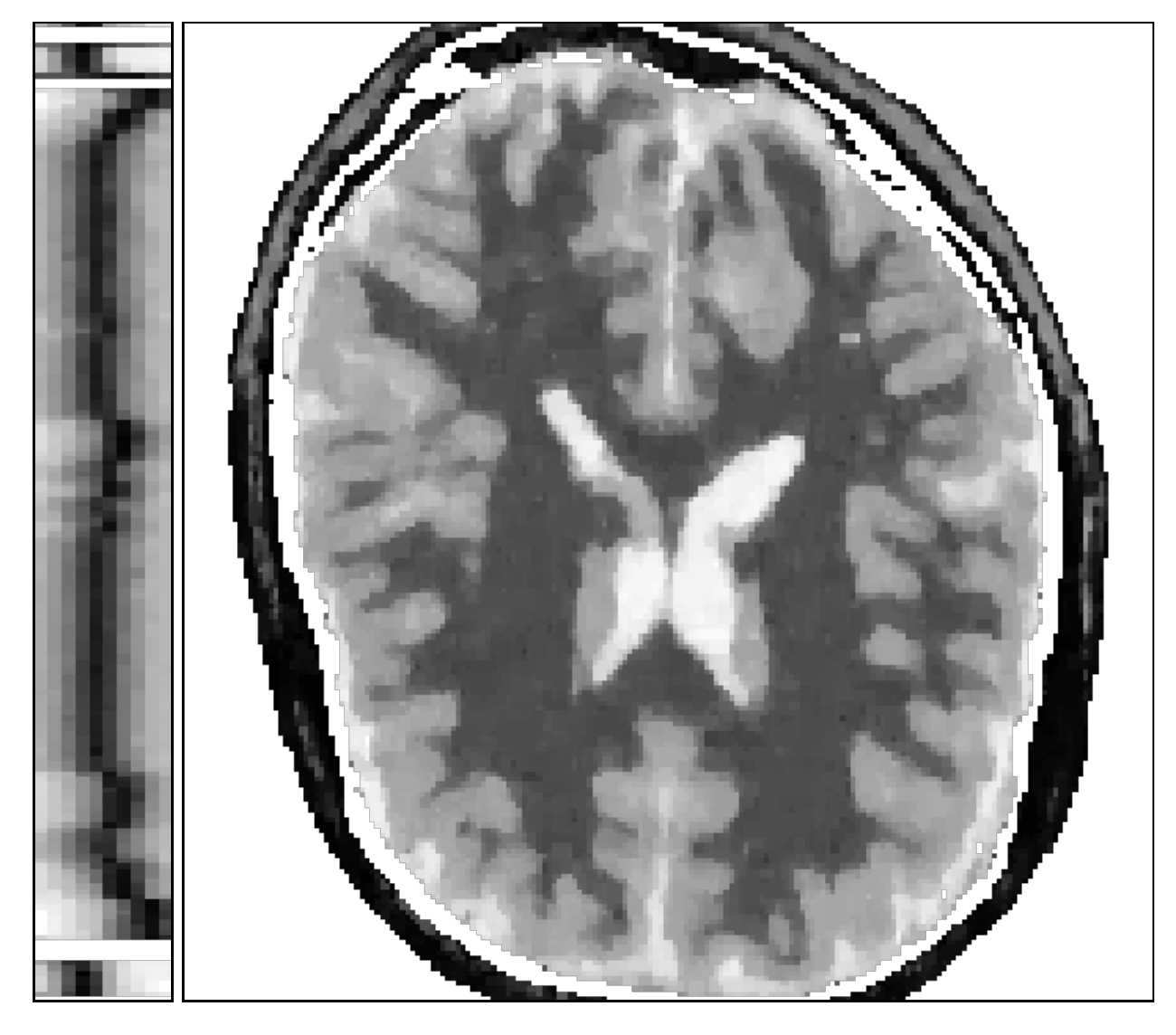} 
        \includegraphics[width=0.15\textwidth]{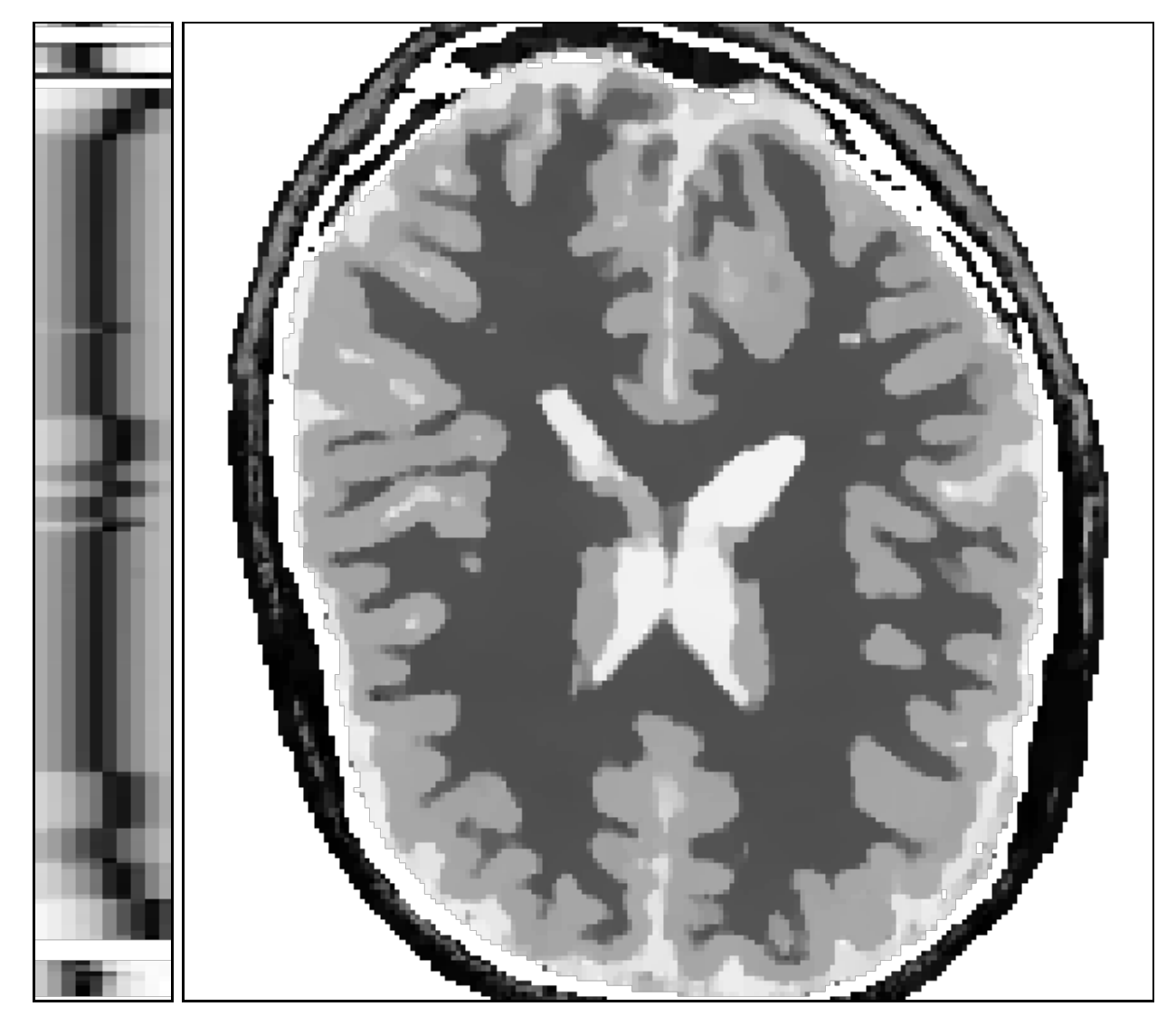} 
        \includegraphics[width=0.15\textwidth]{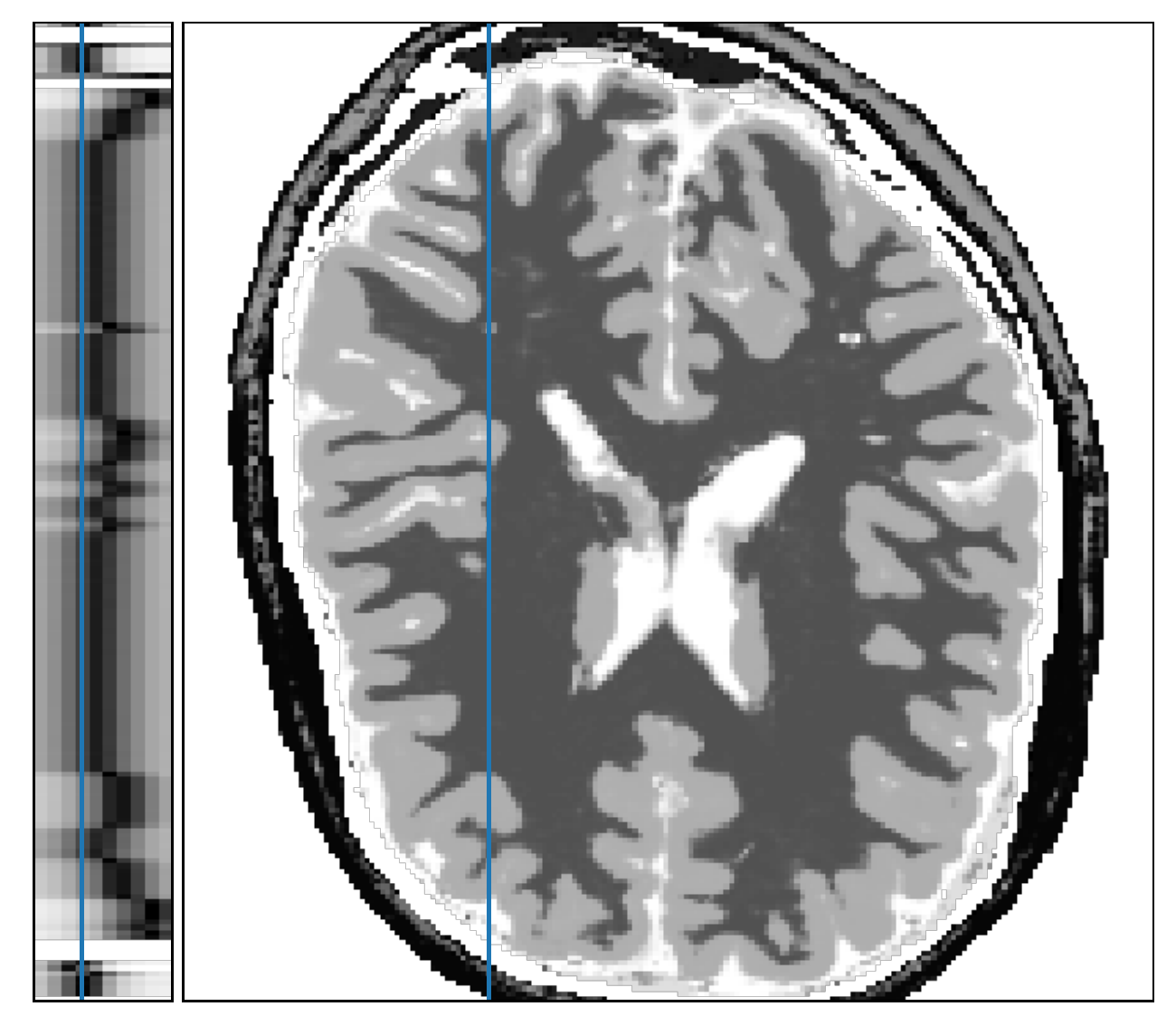} 
        \includegraphics[width=0.15\textwidth]{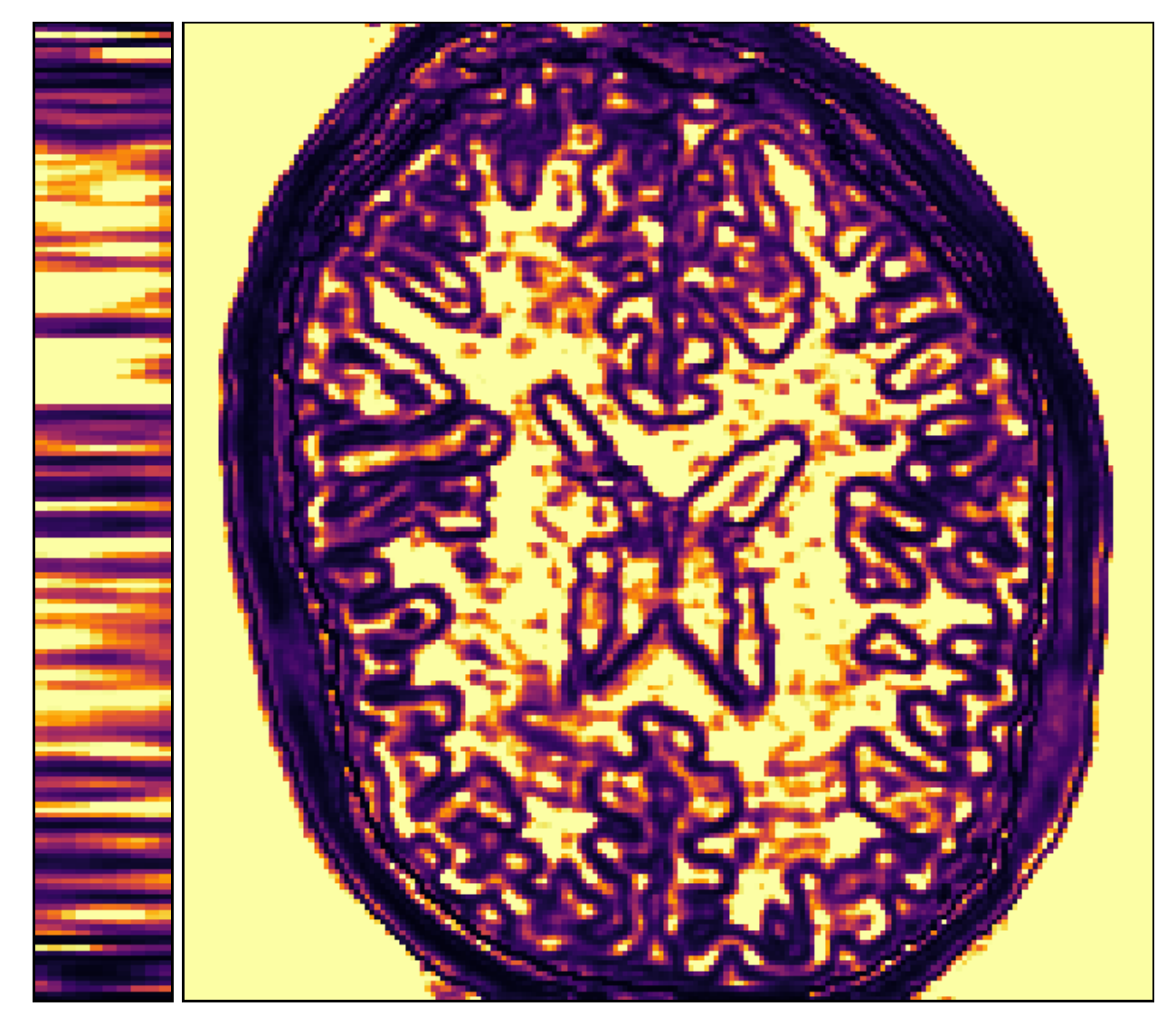}
        \vskip -3pt
        \includegraphics[width=0.15\textwidth]{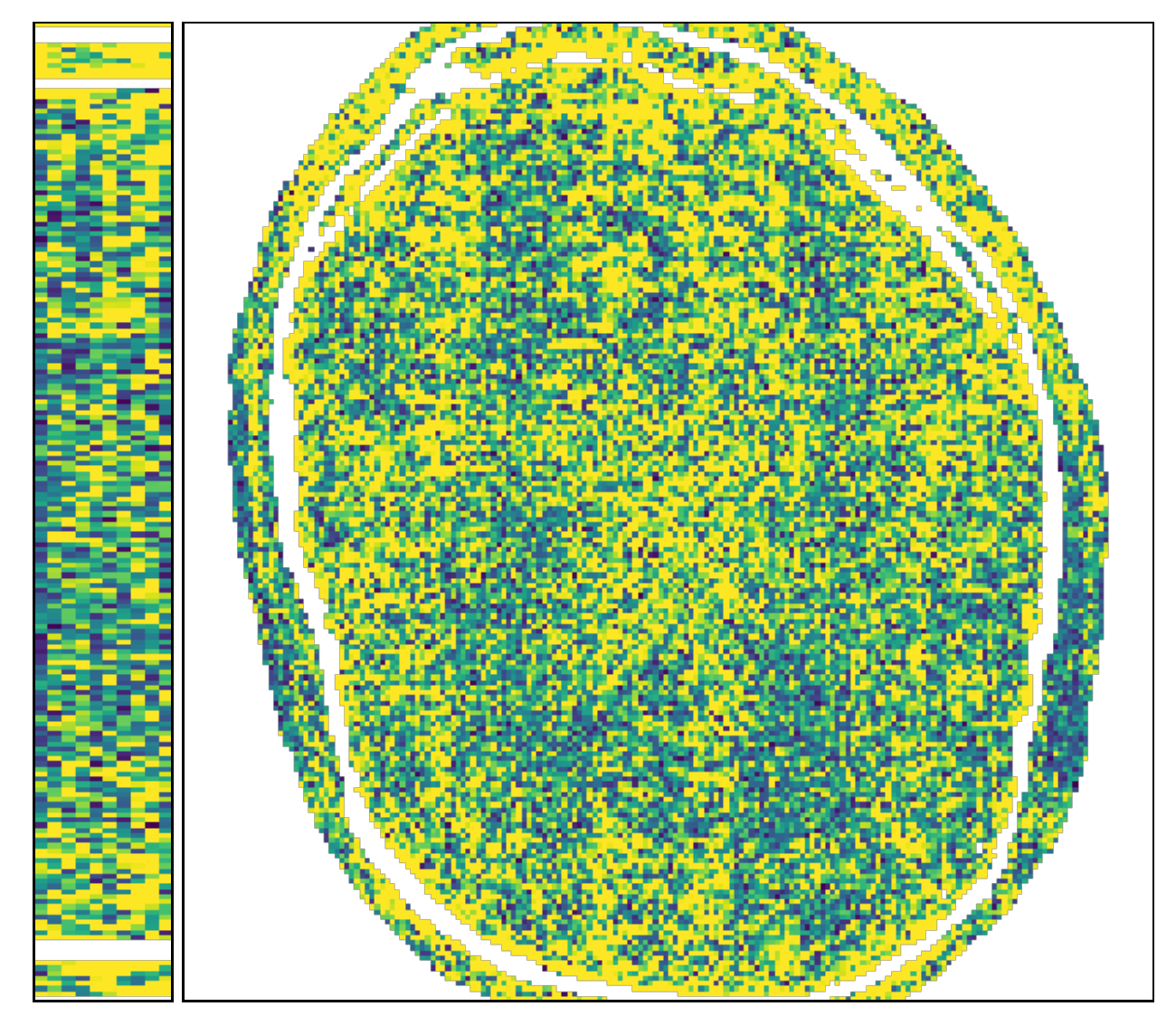} 
        \includegraphics[width=0.15\textwidth]{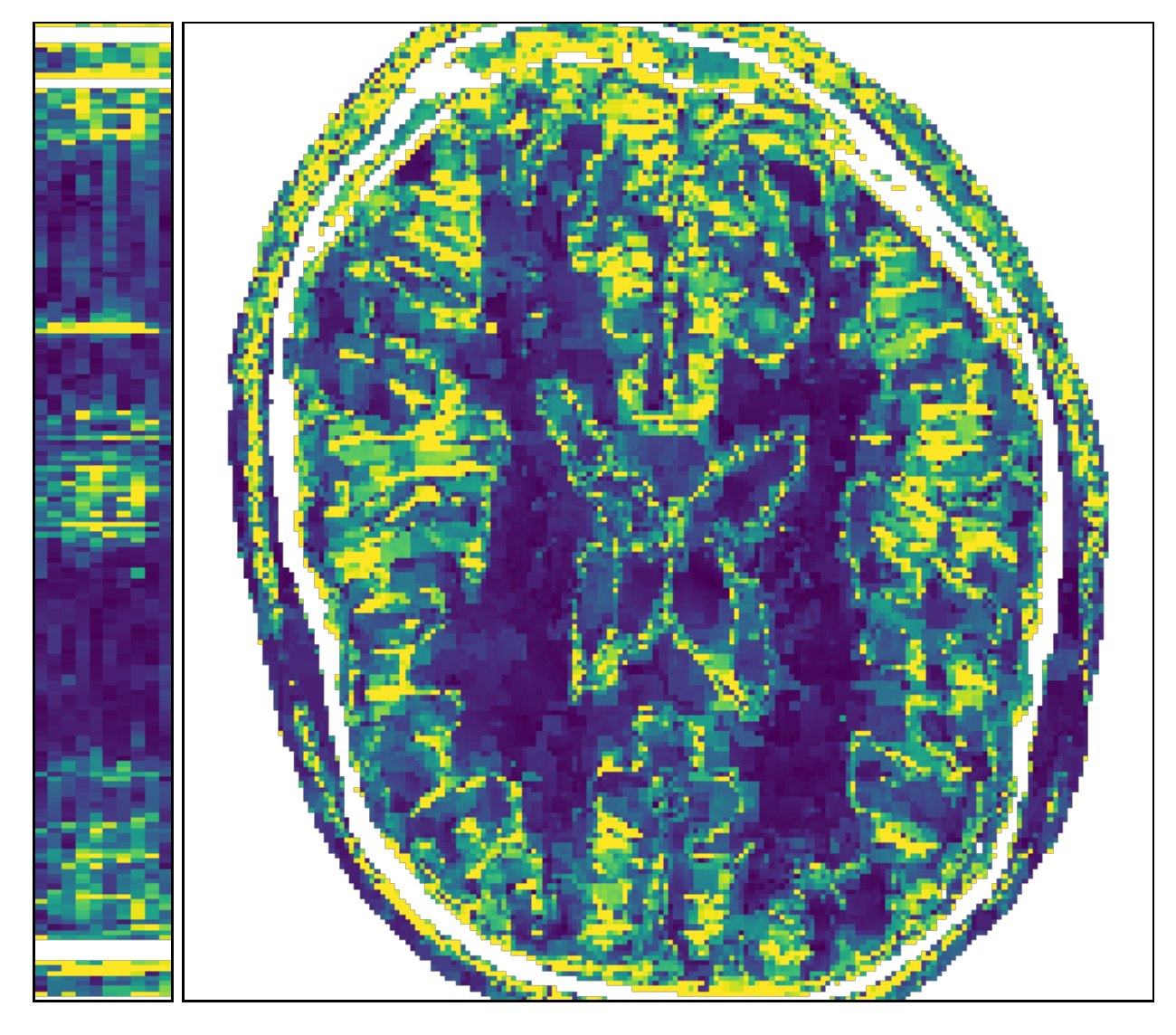} 
        \includegraphics[width=0.15\textwidth]{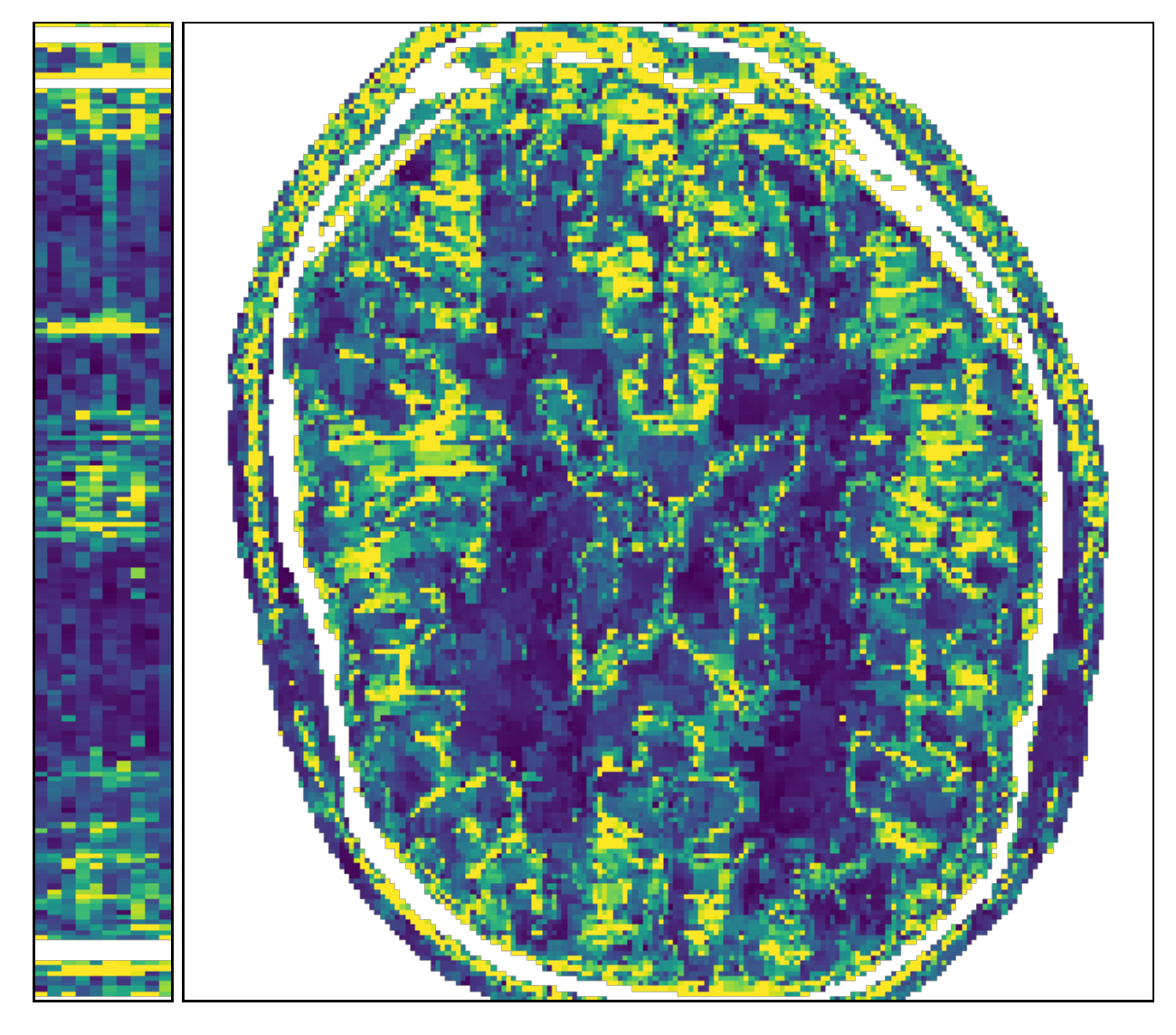} 
        \includegraphics[width=0.15\textwidth]{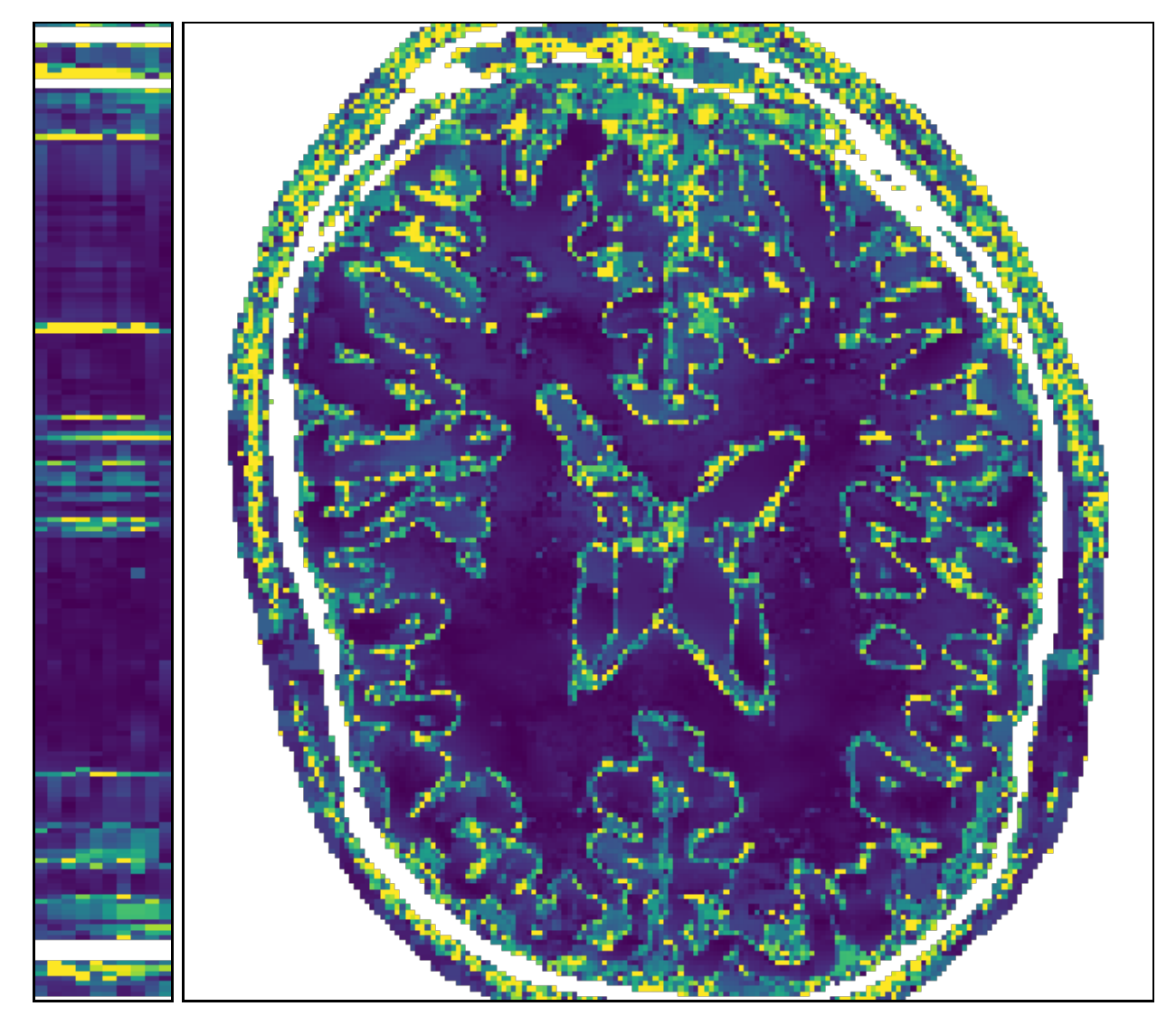}
        \includegraphics[width=0.15\textwidth]{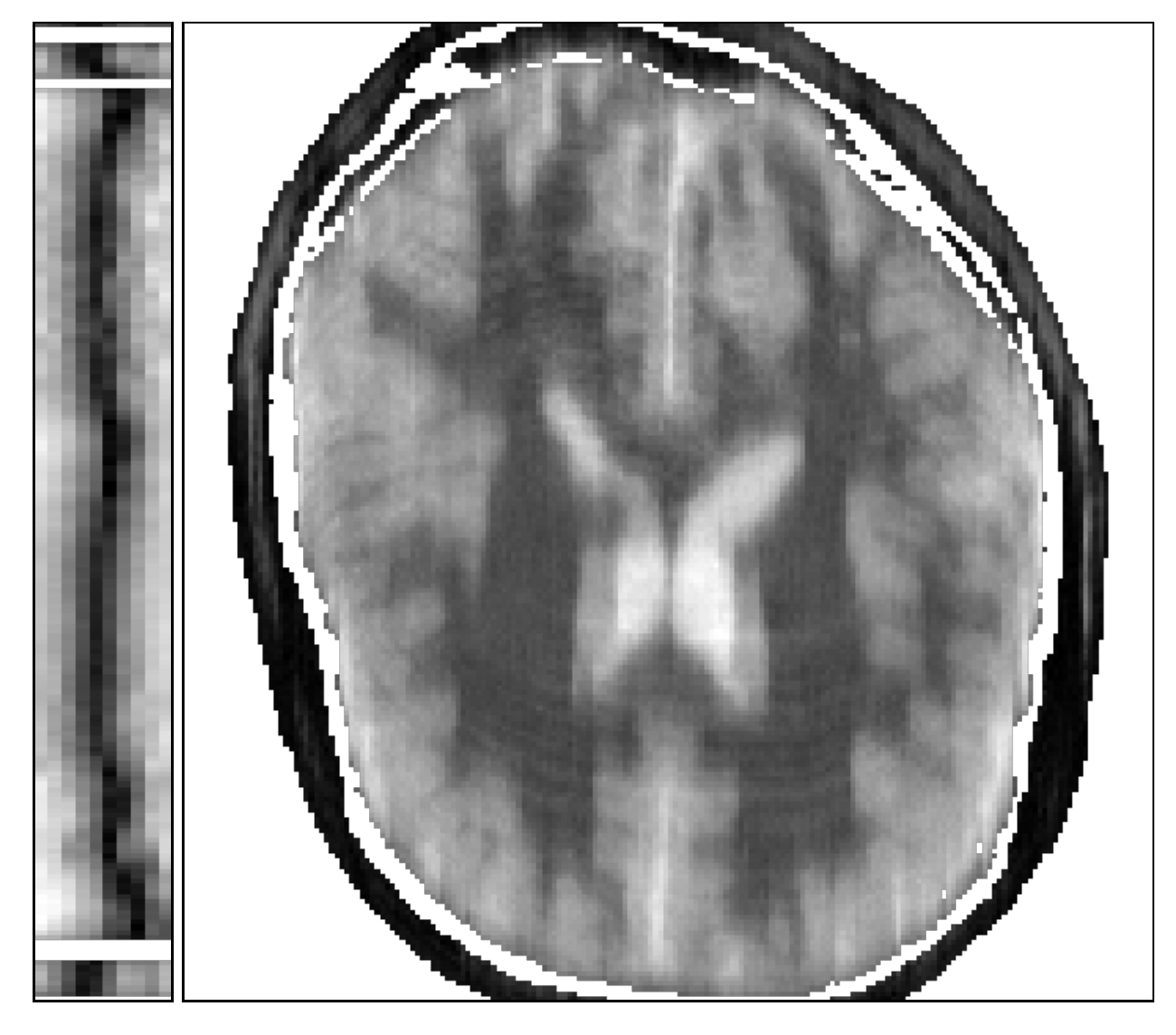} 
        \includegraphics[width=0.15\textwidth]{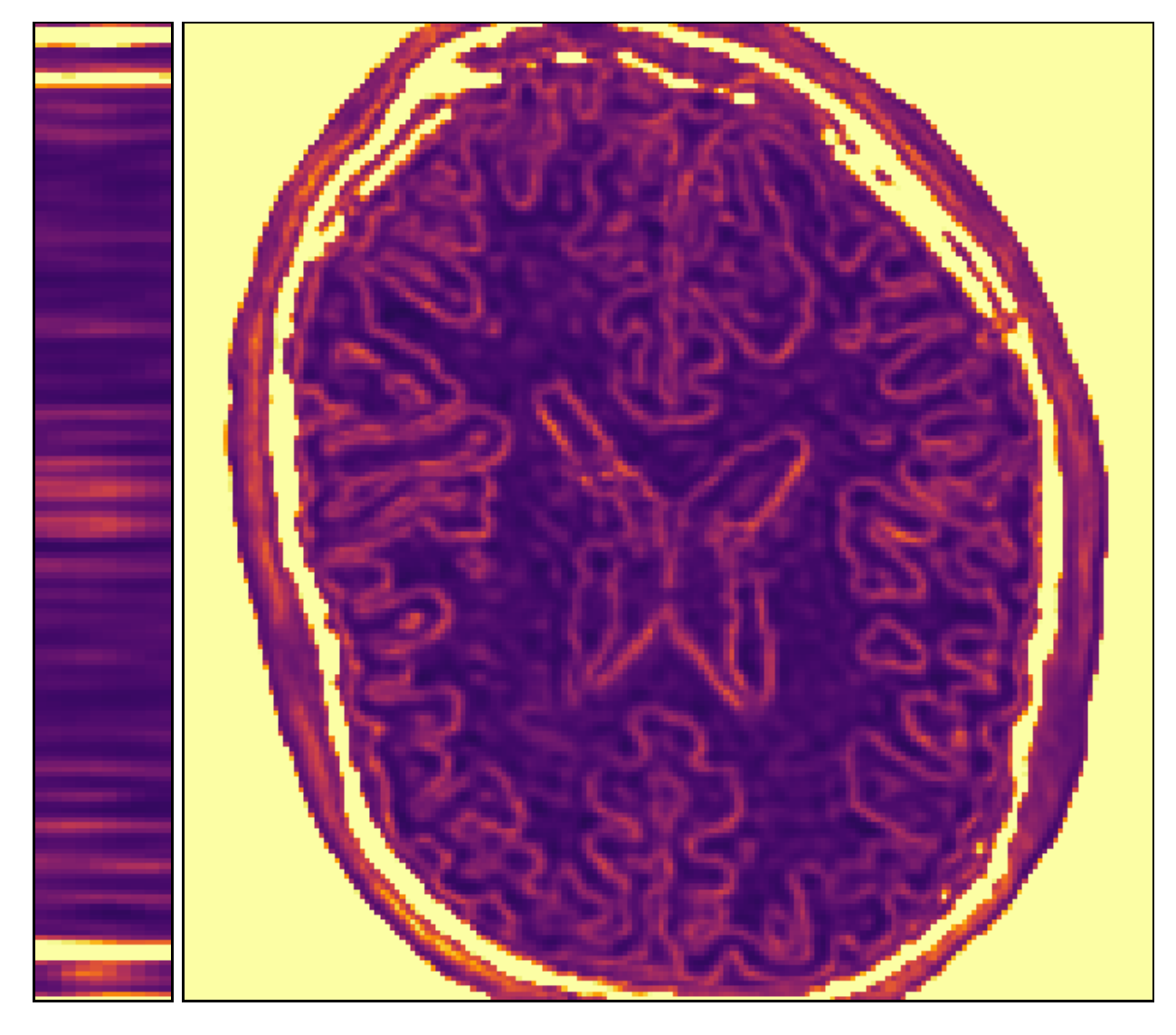}
        
        \includegraphics[width=0.46\textwidth]{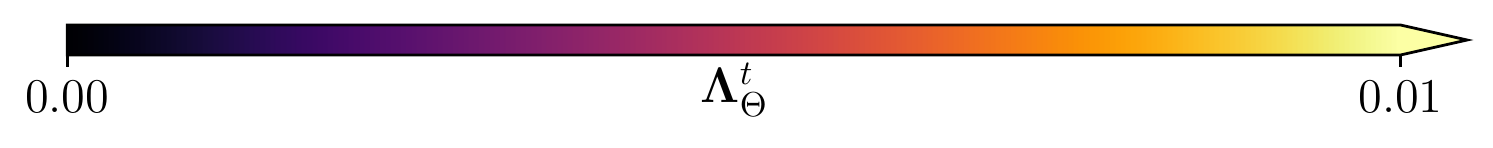}\includegraphics[width=0.46\textwidth]{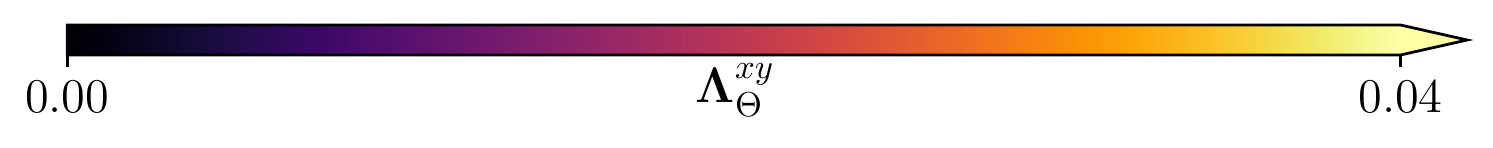}
    \end{subfigure} 
    	\caption{Three exemplary sets of reconstructed, qualitative magnitude images (all acceleration factor 6, $\sigma = $ 0.10 / 0.15 / 0.24 ) and absolute error compared to the ground truth (bottom rows). The last column shows spatial (top) and temporal (bottom) regularization strengths maps  for the PDHG reconstruction generated by the CNN, the second to last column shows the synthetic ground truth as well as the zero-filled reconstruction (bottom rows). Qualitatively, the reconstruction was reduced by the proposed method without resulting in noticeably increased blur. The regularization parameter-maps show strong spatial regularization within areas of homogeneous tissue and increased temporal regularization at the transitions between different tissues.}
	\label{fig:qmri_y}
\end{figure}

\begin{figure}
\hspace{0.16\textwidth}CG-SENSE\hspace{0.03\textwidth}PDHG $\lambda^{xy,t}_{\tilde{\mathrm{P}}}$\hspace{0.015\textwidth}PDHG $\lambda^{xy,t}_{\mathrm{P}}$\hspace{0.015\textwidth}PDHG $\LLambda^{xy,t}_\Theta$ \hspace{0.0\textwidth}Ground Truth 

\centering   
\def\which{20}

\def\acc{acc6}
\rotatebox[origin=c]{90}{$R=6$}\,\,\begin{subfigure}[p]{0.7\textwidth}
    \includegraphics[width=0.19\textwidth]{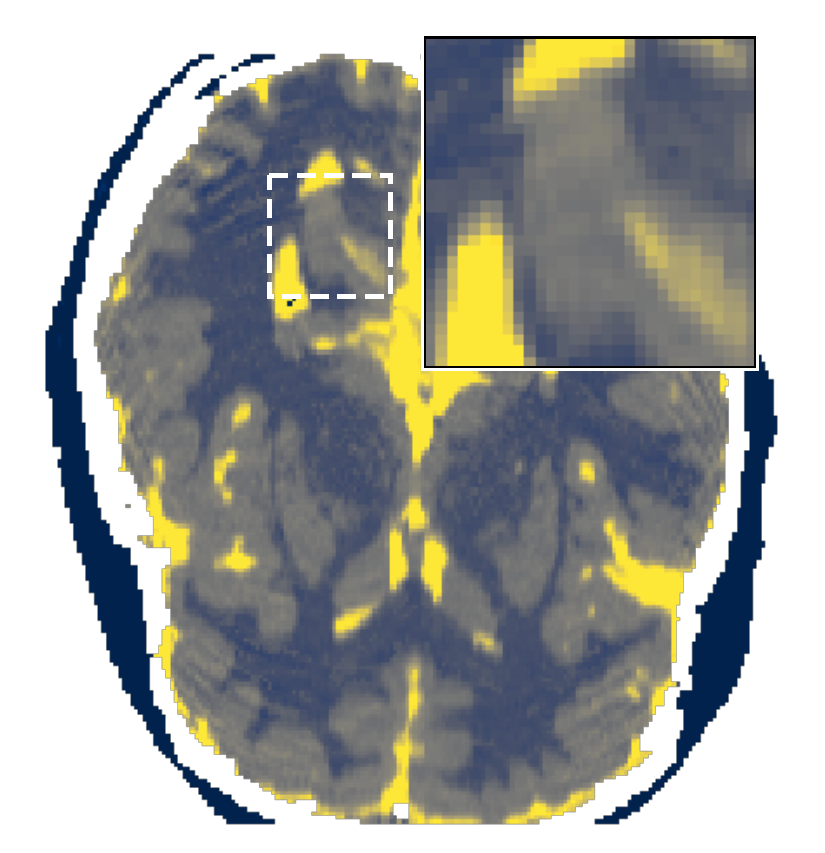} 
    \includegraphics[width=0.19\textwidth]{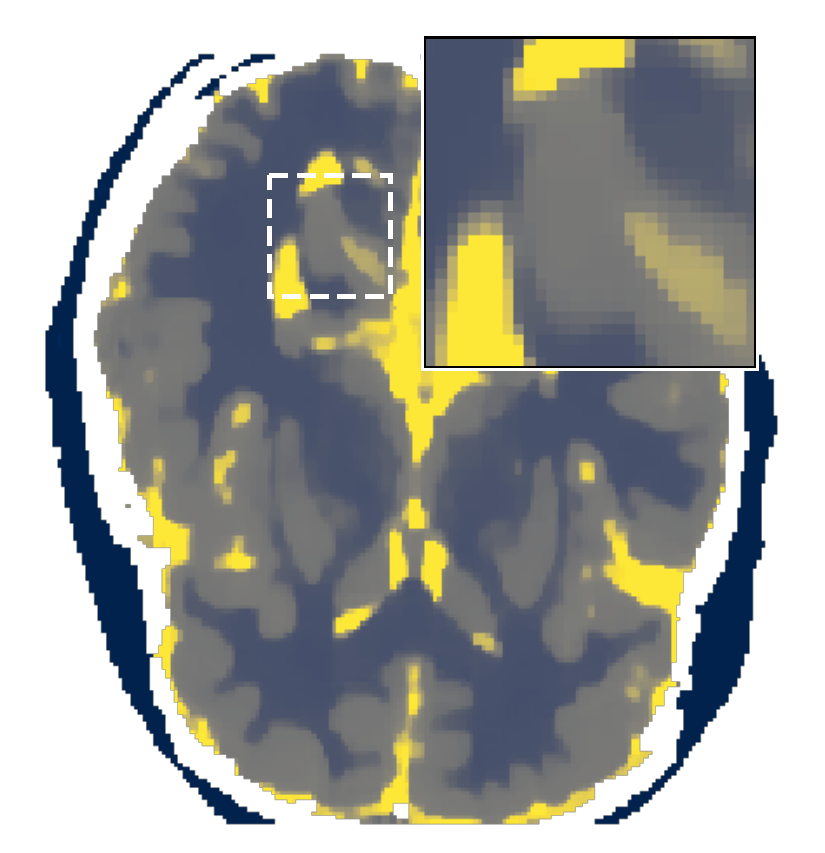} 
    \includegraphics[width=0.19\textwidth]{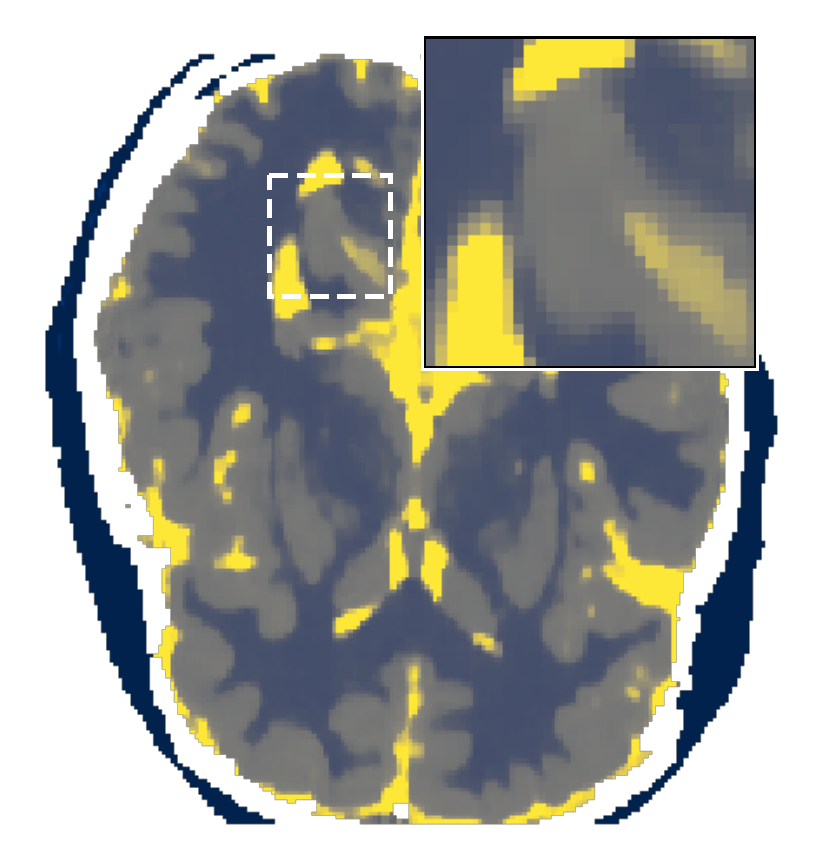} 
    \includegraphics[width=0.19\textwidth]{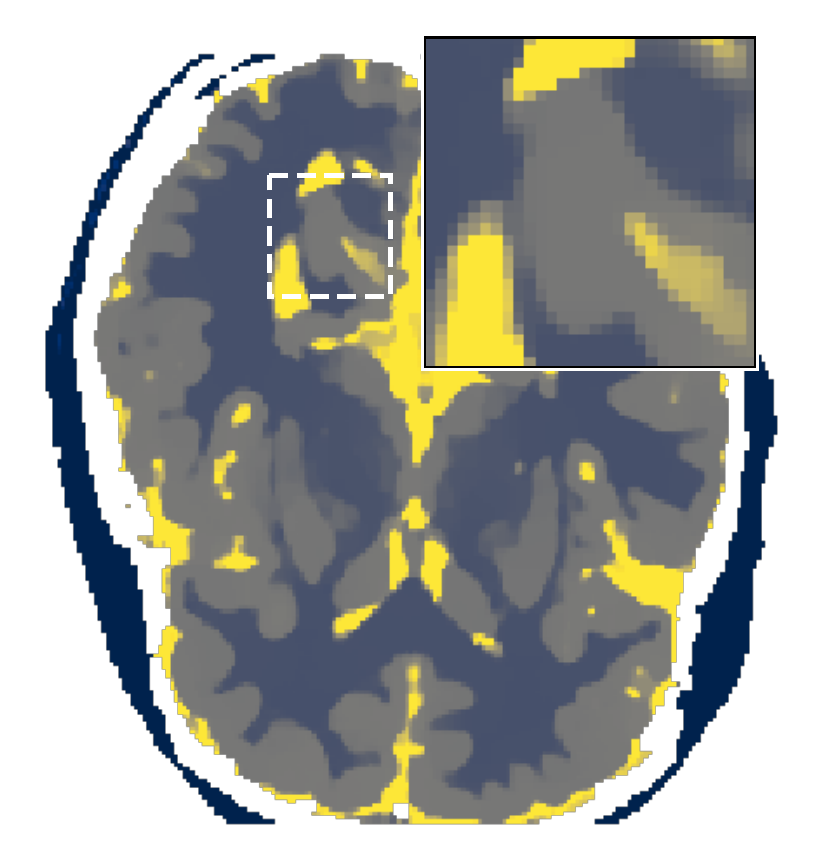} 
    \includegraphics[width=0.19\textwidth]{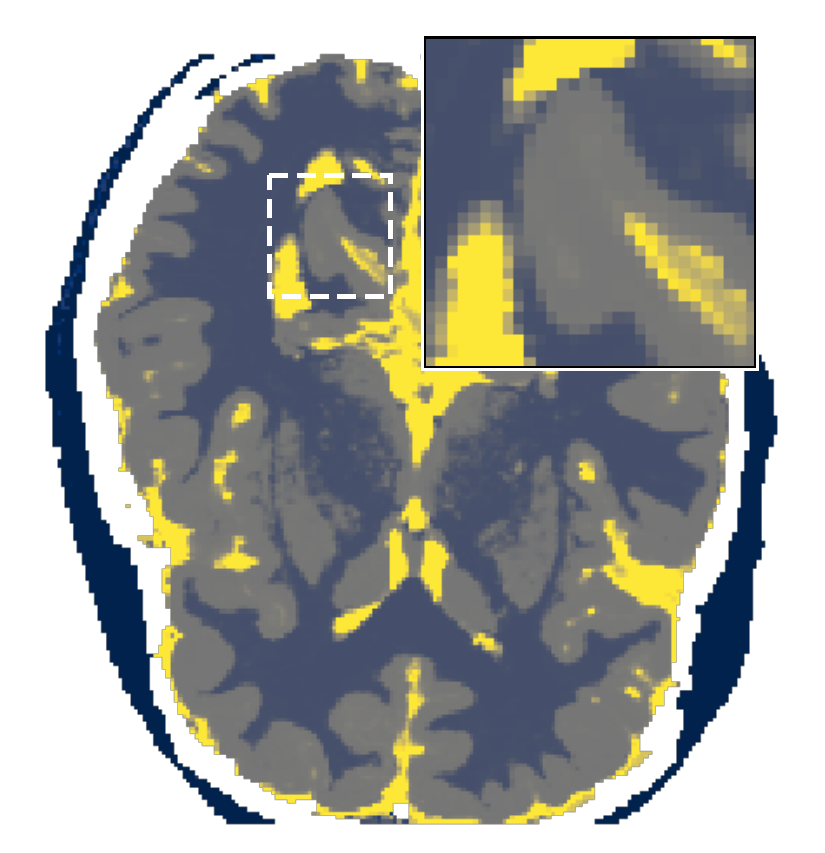}
 \vskip -5pt
    \includegraphics[width=0.19\textwidth]{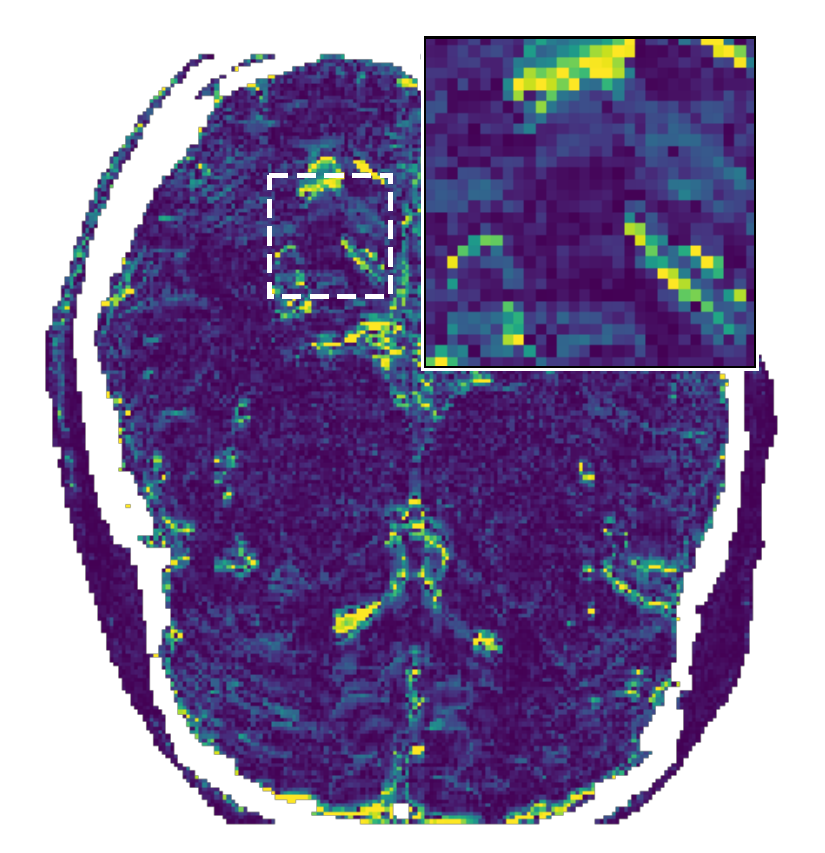} 
    \includegraphics[width=0.19\textwidth]{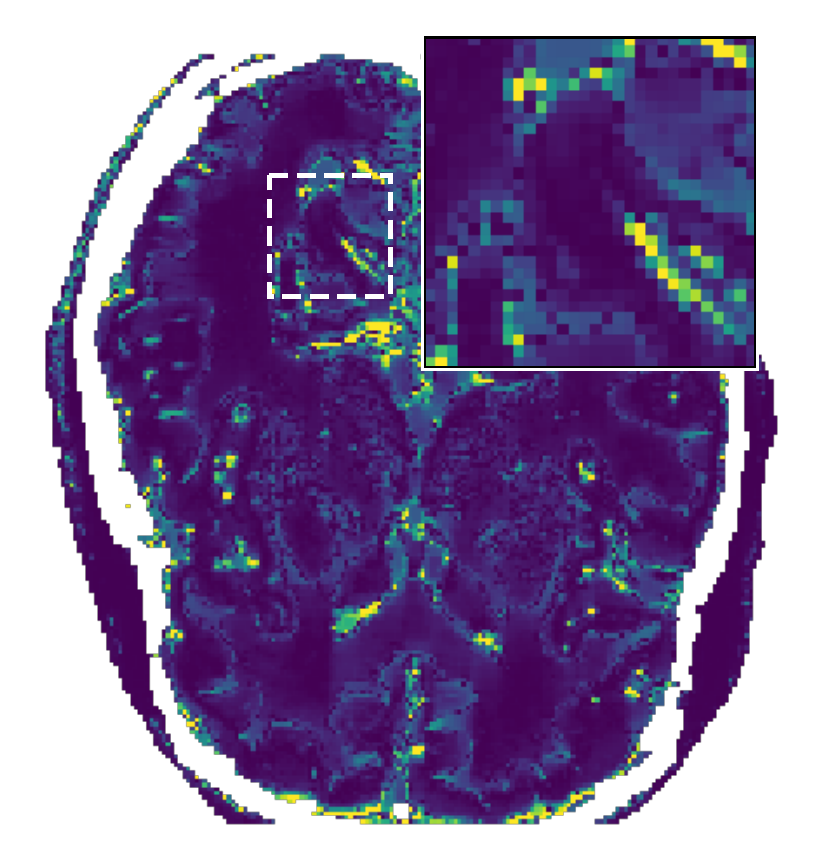} 
    \includegraphics[width=0.19\textwidth]{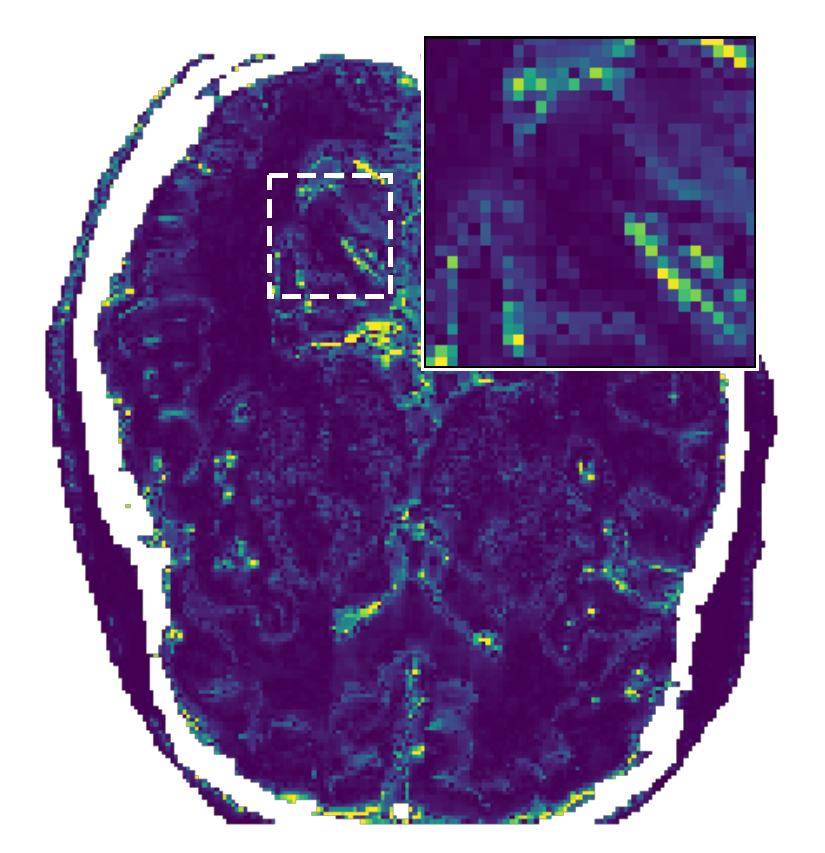} 
    \includegraphics[width=0.19\textwidth]{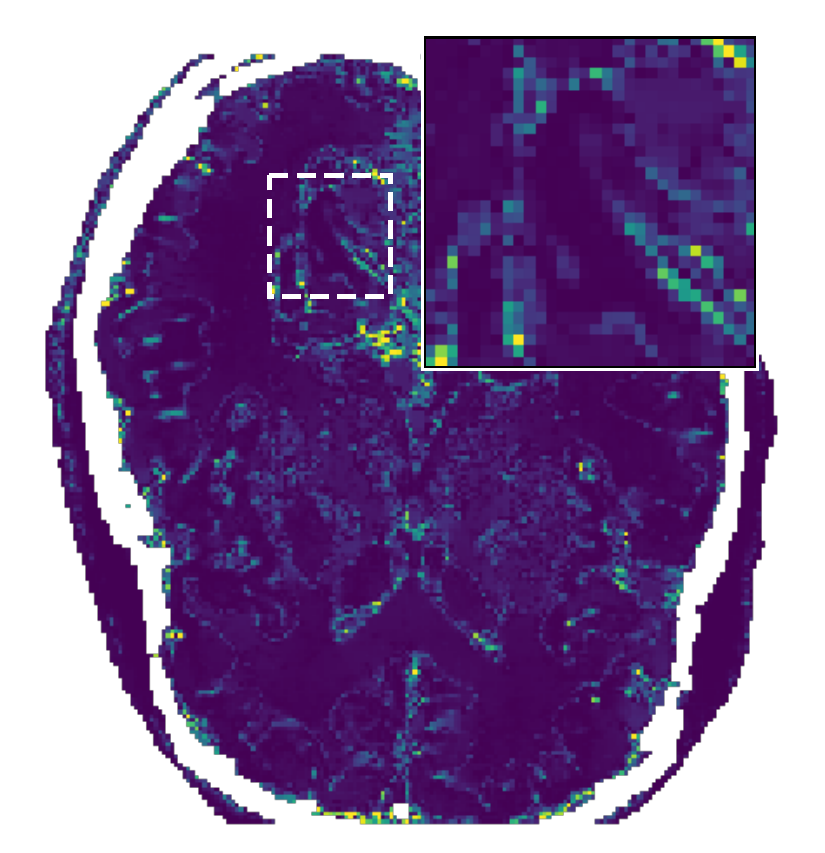}
    \hspace{0.19\textwidth} 
\end{subfigure}   
\vspace{0.005\textwidth}

\def\acc{acc8}
\rotatebox[origin=c]{90}{$R=8$}\,\,\begin{subfigure}[p]{0.7\textwidth}
    \includegraphics[width=0.19\textwidth]{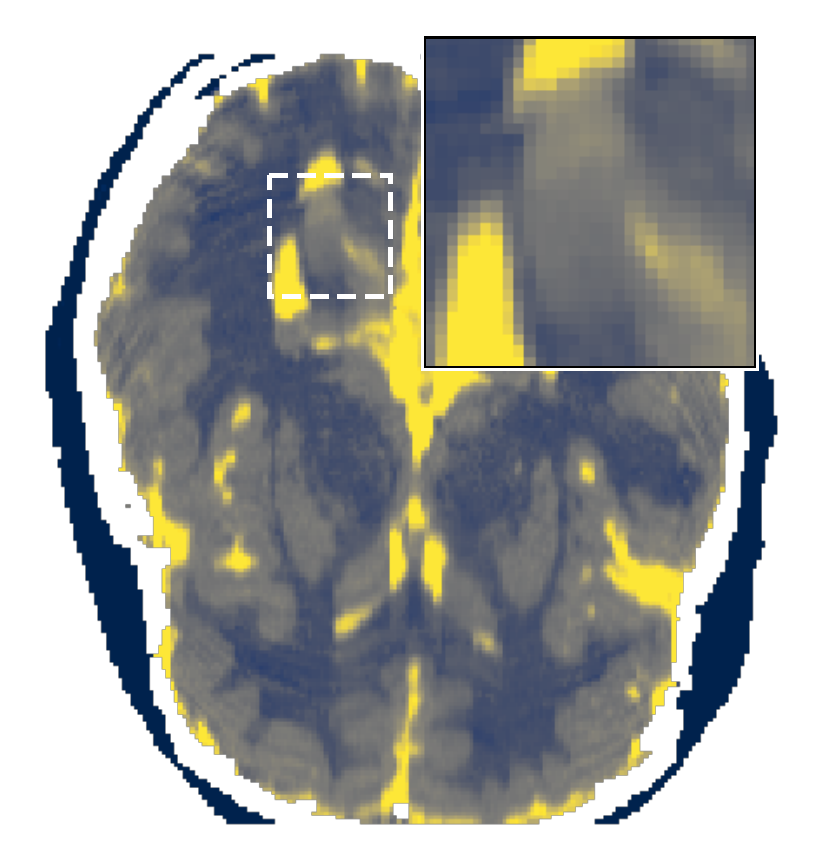} 
    \includegraphics[width=0.19\textwidth]{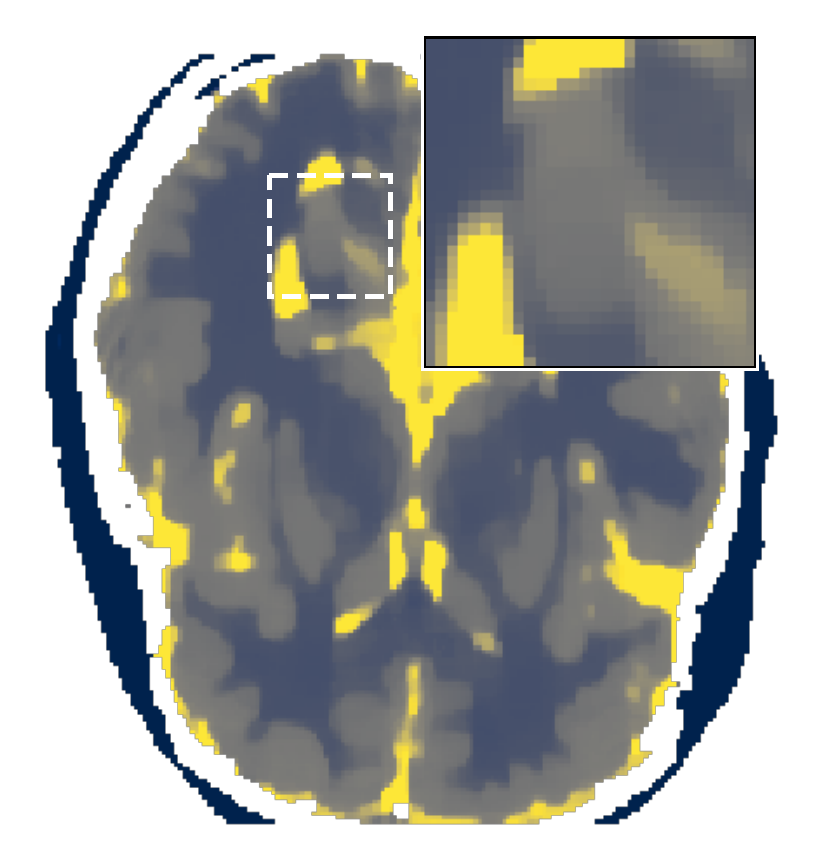} 
    \includegraphics[width=0.19\textwidth]{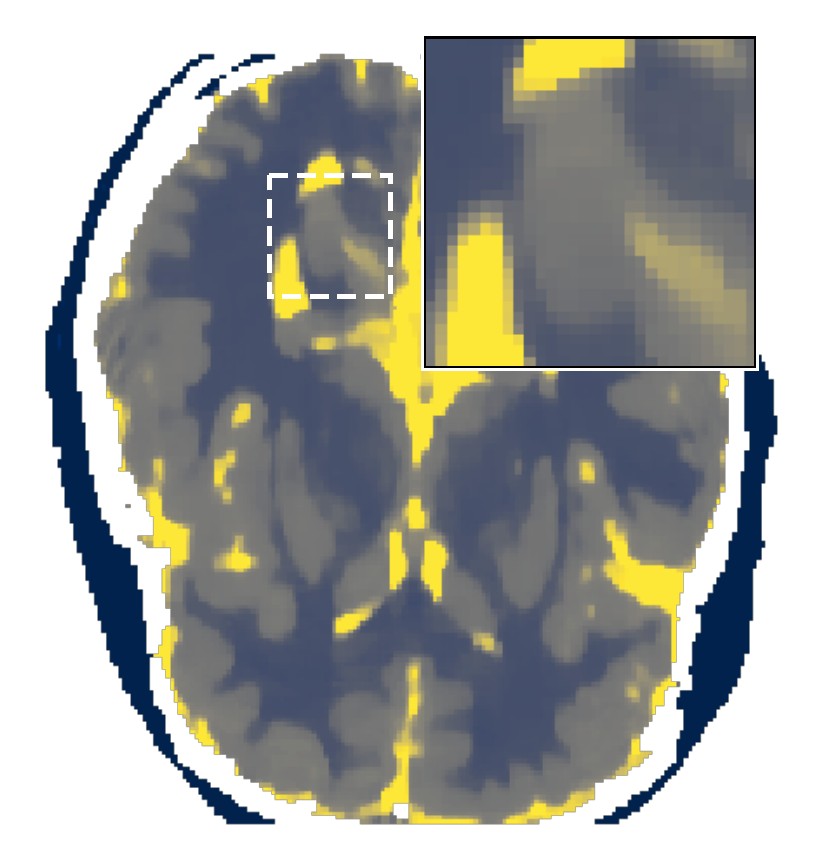} 
    \includegraphics[width=0.19\textwidth]{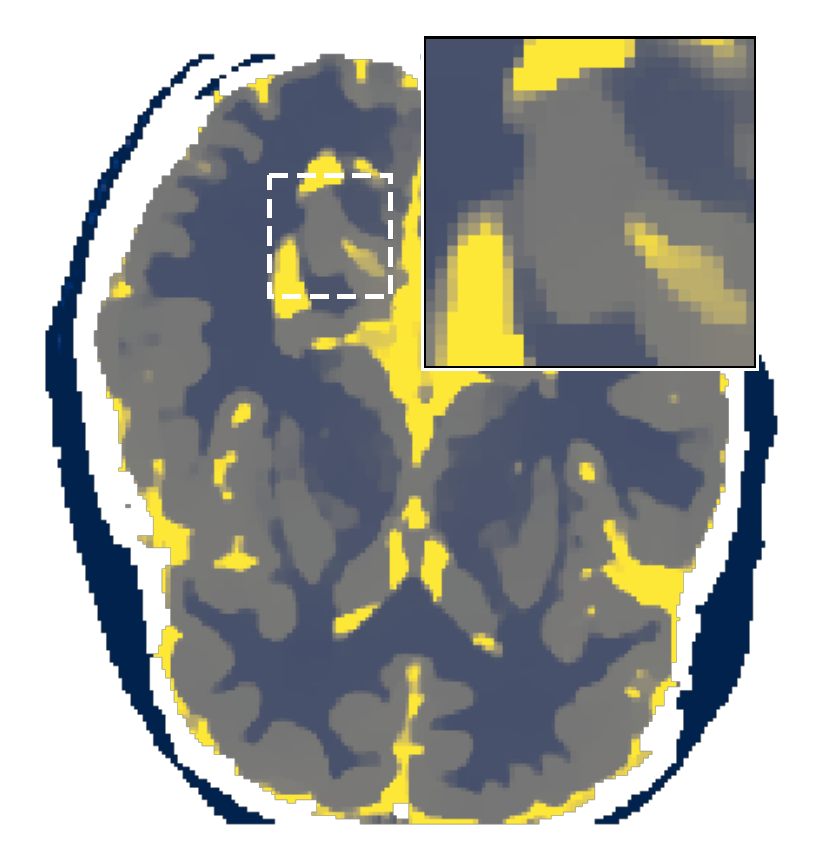} 
    \includegraphics[width=0.19\textwidth]{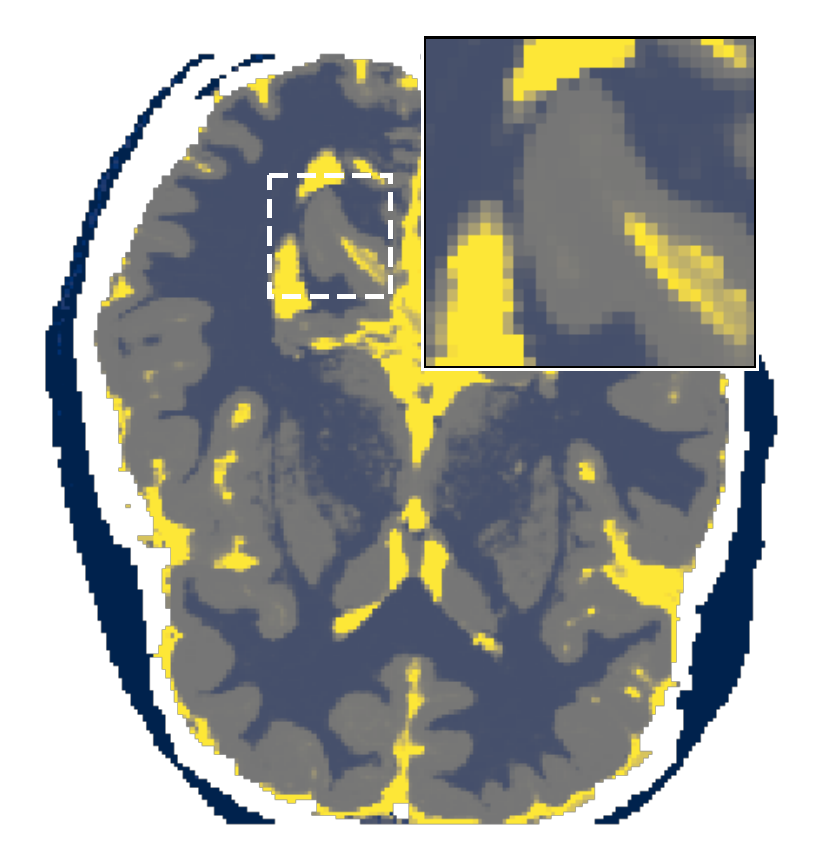}
\vskip -5pt
    \includegraphics[width=0.19\textwidth]{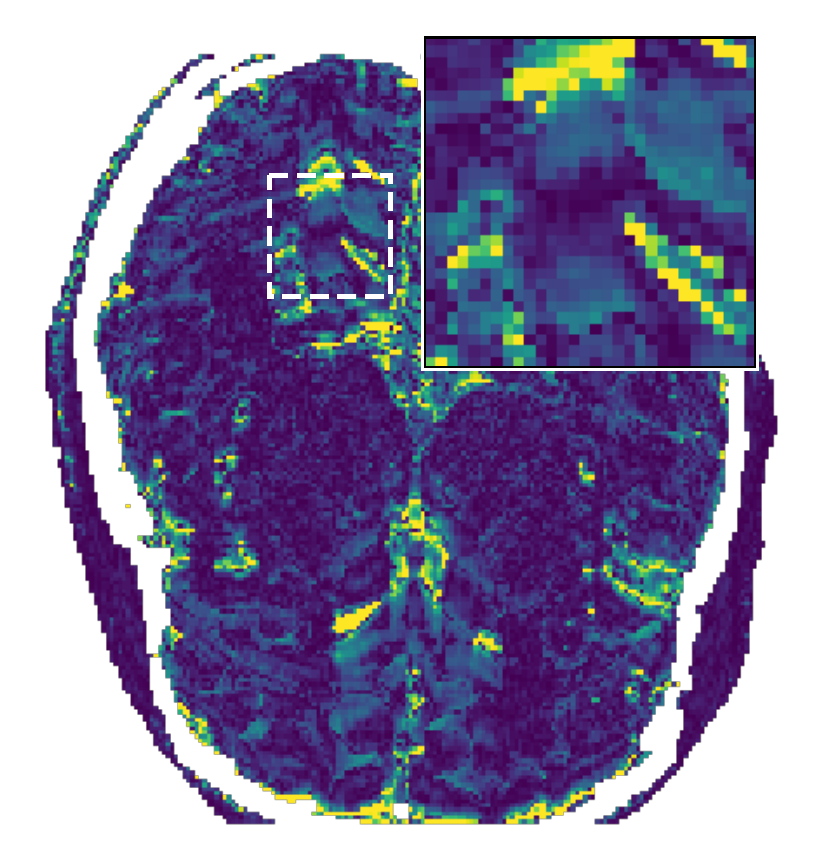} 
    \includegraphics[width=0.19\textwidth]{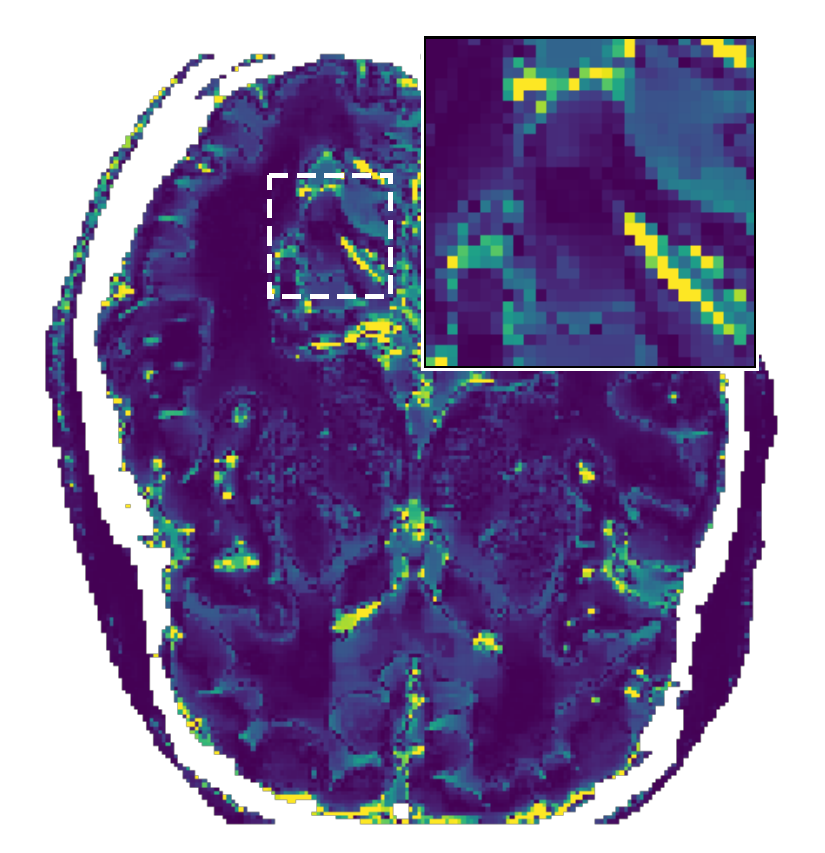} 
    \includegraphics[width=0.19\textwidth]{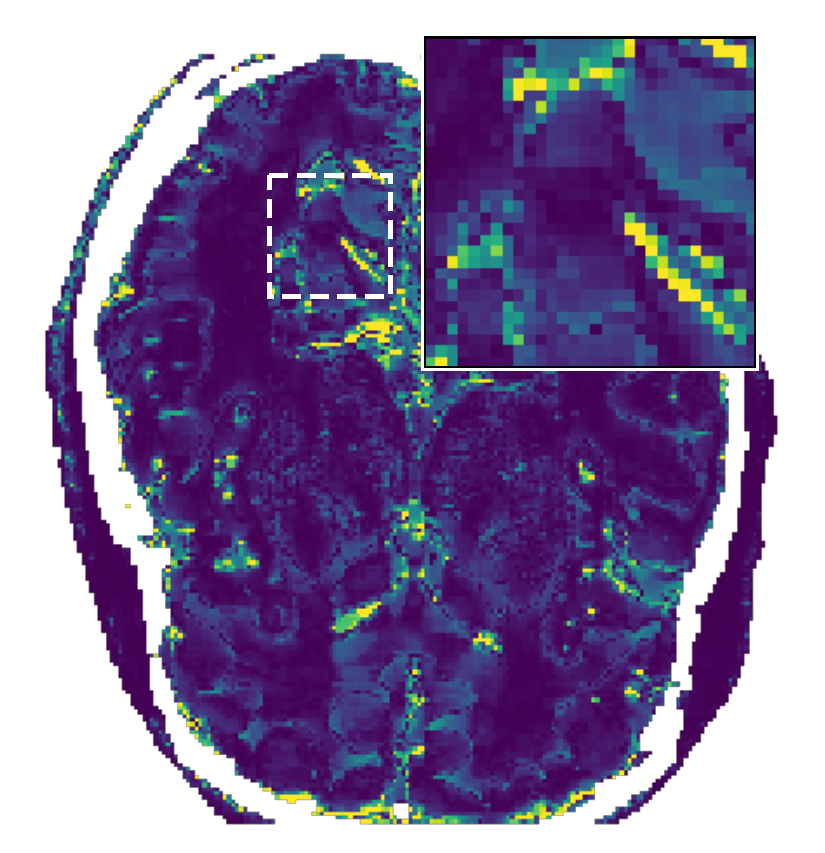} 
    \includegraphics[width=0.19\textwidth]{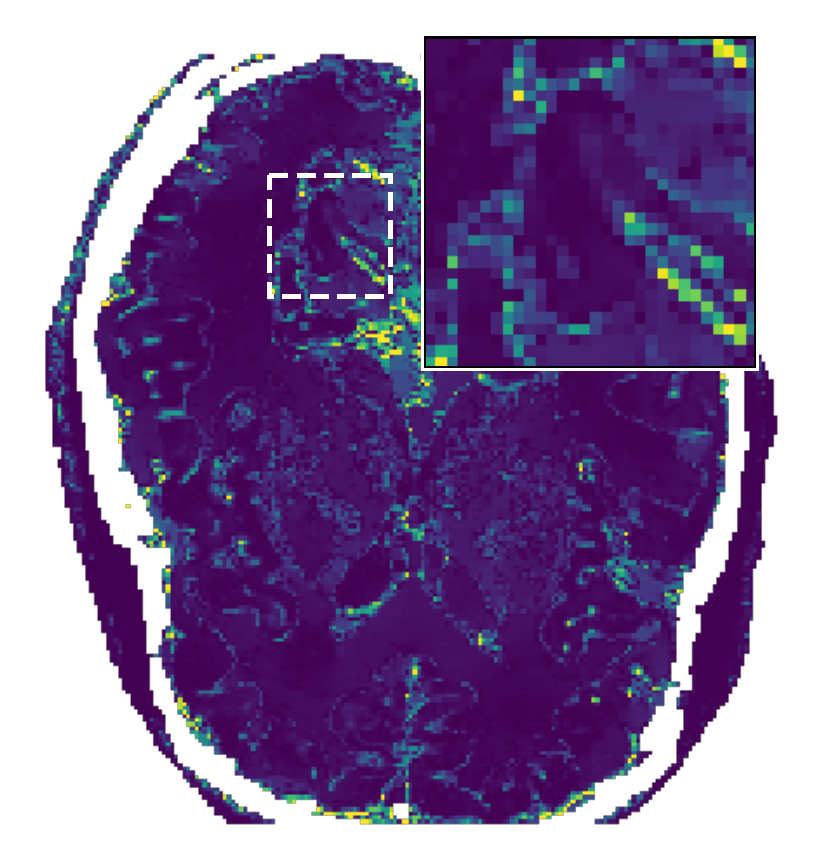}
    \includegraphics[width=0.07\textwidth]{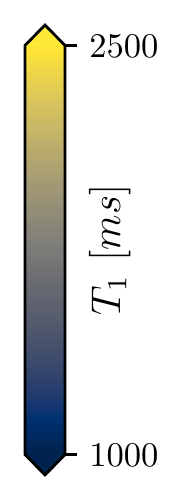}
    \includegraphics[width=0.06\textwidth]{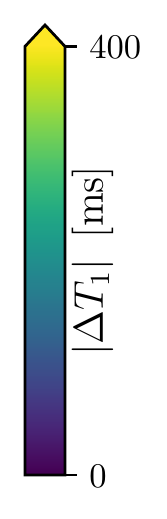} 
\end{subfigure} 

    \caption{Resulting quantitative $T_1$ parameter-maps after performing the regression on the reconstructed magnitude images and absolute errors compared to the ground truth (bottom row). Shown is \textit{Example 2} of Figure \ref{fig:qmri_y}, at acceleration 6 and 8.}
    \label{fig:qmri_t1}
\end{figure}

\begin{figure}
    \centering
\includegraphics{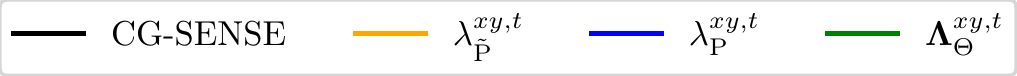}

    \includegraphics[width=0.3\textwidth]{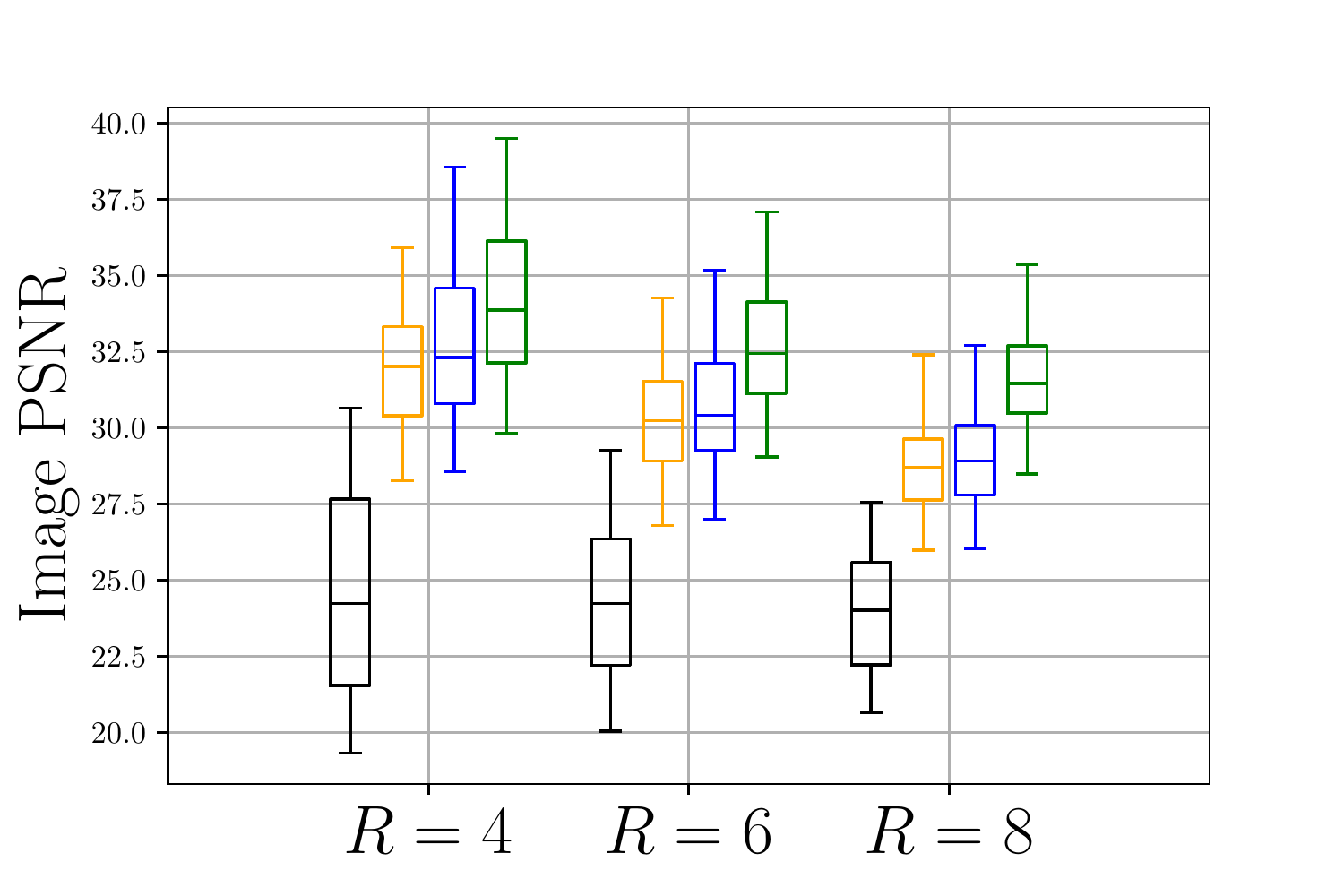}
    \includegraphics[width=0.3\textwidth]{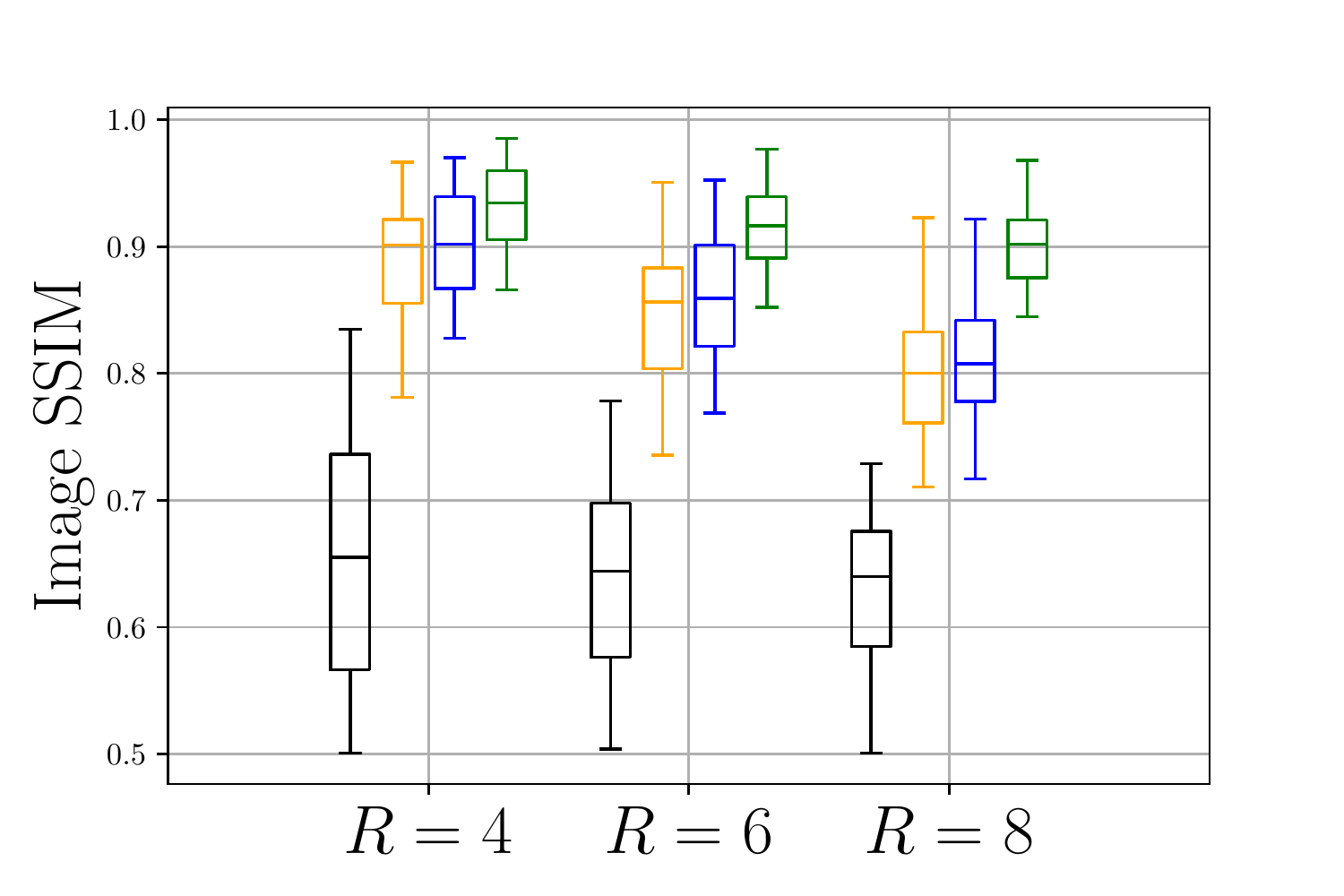}
  \includegraphics[width=0.3\textwidth]{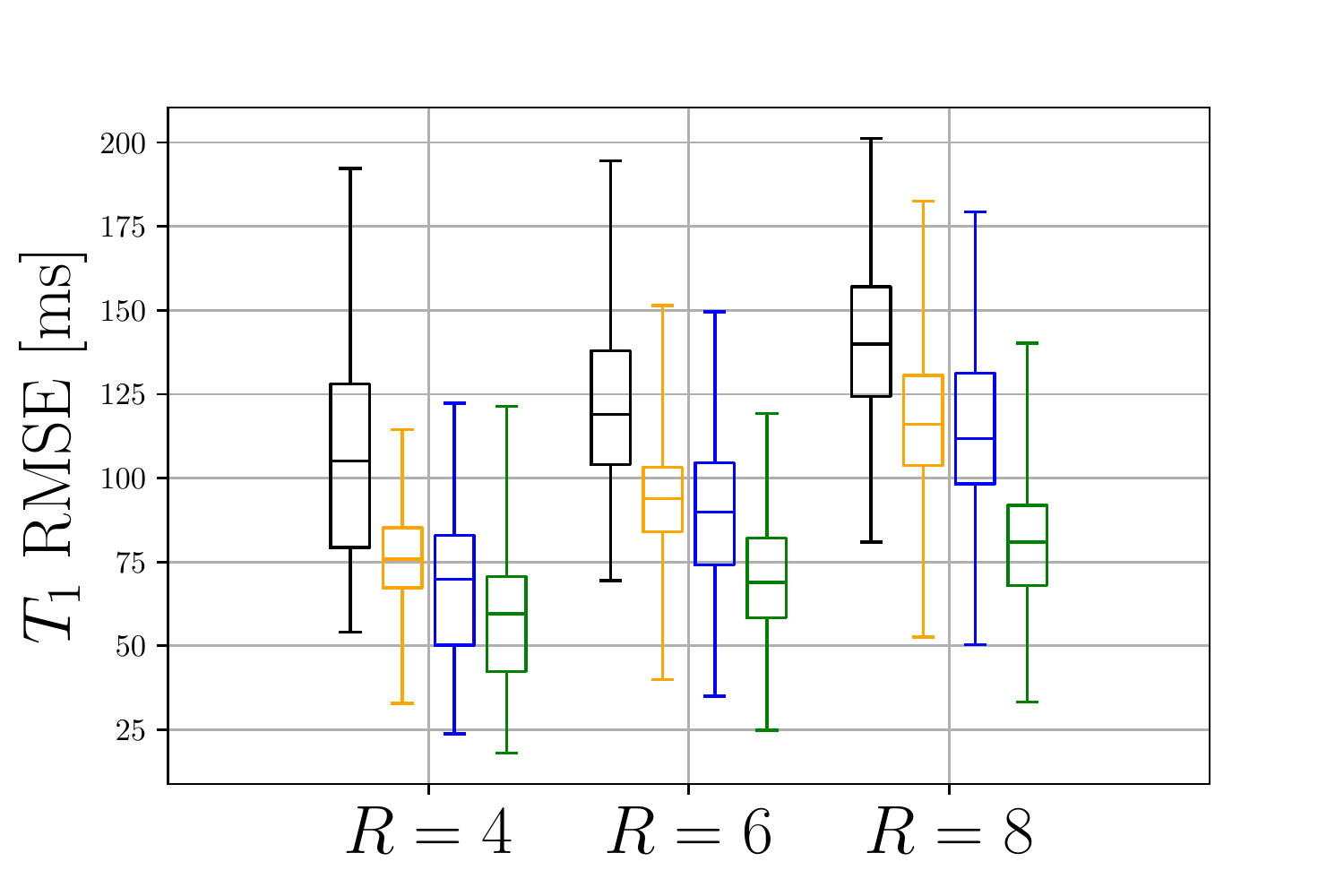}
    \caption{Quantification result of the quality of the reconstructed magnitude images in terms of PSNR and SSIM, and RMSE of the $T_1$ values after performing the signal regression. We compare our proposed CNN-based generation of optimal spatial and temporally varying regularization maps for the PDHG reconstruction with a CG-SENSE reconstruction\cite{pruessmann2001advances} (early stopping for solving the normal equations), two scalar regularization values chosen on the whole test dataset by grid-search ($\lambda_{\tilde{\mathrm{P}}}^{xy,t}$), as well as two scalar regularization values that were chosen for each image of the test dataset by grid-search with access to the ground truth images ($\lambda_{\mathrm{P}}^{xy,t}$).}
        \label{fig:box_plots_qmri}
\end{figure}

Figure \ref{fig:qmri_y} shows examples of the quantitative (magnitude) images $\UU$ of three of the 112 simulated inversion recovery measurements in the test dataset. We also show the  regularization parameter-maps for regularization along the spatial directions and along the inversion-time direction generated by the network. The mean PSNR and SSIM of our proposed method is consistently higher for all considered acceleration factors, even compared to PDHG with regularization strength along spatial and inversion-time direction chosen by grid-search with access to the ground truth images (shown in Figure \ref{fig:box_plots_qmri} and Table \ref{tab:qmri_results}).
The resulting $T_1$ parameter-maps after performing the regression on the reconstructed images are shown in Figure \ref{fig:qmri_t1}. Again, our proposed method results in the lowest RMS deviation from the ground truth images (Table \ref{tab:qmri_results}).

\begin{table}[h]
\centering
\footnotesize\rm
\begin{tabular}{r|l|rcl|rcl|rcl|rcl}

&    & \multicolumn{3}{c|}{\textbf{CG-SENSE}}  & \multicolumn{3}{c|}{\textbf{PDHG - $\lambda_{\tilde{\mathrm{P}}}^{xy,t}$}}  & \multicolumn{3}{c|}{\textbf{PDHG - $\lambda_{\mathrm{P}}^{xy,t}$}} & \multicolumn{3}{c}{\textbf{PDHG - $\boldsymbol{\Lambda}_{\Theta}^{xy,t}$}}  \\[5pt]
\hline
&  \textbf{PSNR}&  $24.62$&$\pm$&$3.45$&  $31.89$&$\pm$&$1.70$&  $32.85$&$\pm$&$2.48$&  $\mathbf{34.23}$&$\pm$&$2.50$\\ 
$R=4$&  \textbf{SSIM}&  $0.654$&$\pm$&$0.095$&  $0.884$&$\pm$&$0.044$&  $0.902$&$\pm$&$0.042$&  $\mathbf{0.930}$&$\pm$&$0.031$\\ 
&  \textbf{RMSE} [ms]&  $107$&$\pm$&$34$&  $76$&$\pm$&$15$&  $68$&$\pm$&$21$&  $\mathbf{58}$&$\pm$&$18$\\ 
 \hline 
&  \textbf{PSNR}&  $24.25$&$\pm$&$2.49$&  $30.19$&$\pm$&$1.58$&  $30.62$&$\pm$&$1.80$&  $\mathbf{32.62}$&$\pm$&$1.79$\\ 
$R=6$&  \textbf{SSIM}&  $0.639$&$\pm$&$0.077$&  $0.843$&$\pm$&$0.048$&  $0.859$&$\pm$&$0.044$&  $\mathbf{0.914}$&$\pm$&$0.027$\\ 
&  \textbf{RMSE} [ms]&  $122$&$\pm$&$26$&  $94$&$\pm$&$19$&  $91$&$\pm$&$22$&  $\mathbf{70}$&$\pm$&$16$\\ 
 \hline 
&  \textbf{PSNR}&  $23.87$&$\pm$&$1.90$&  $28.70$&$\pm$&$1.42$&  $28.93$&$\pm$&$1.53$&  $\mathbf{31.61}$&$\pm$&$1.48$\\ 
$R=8$&  \textbf{SSIM}&  $0.623$&$\pm$&$0.063$&  $0.799$&$\pm$&$0.043$&  $0.810$&$\pm$&$0.043$&  $\mathbf{0.897}$&$\pm$&$0.026$\\ 
&  \textbf{RMSE} [ms]&  $140$&$\pm$&$25$&  $117$&$\pm$&$24$&  $114$&$\pm$&$25$&  $\mathbf{82}$&$\pm$&$18$\\ 
 \hline 
\end{tabular}
\caption{Mean and standard deviation of the measures PSNR and SSIM of the qualitative images and RMSE of the $T_1$ parameter-maps over the test set.}\label{tab:qmri_results}
\end{table}

\newpage
\subsection{Dynamic Image Denoising}
Here, we apply the proposed method to estimate voxel-wise dependent regularization parameter-maps to be used in a dynamic image denoising problem. An important difference to the previously considered cardiac MRI example is that, while in the latter, a clear inherent distinction between the black background and the object of interest is possible, for the next videos to be considered, this is not the case. The samples might show scenes with static camera position and only moving objects or scenes in which also the camera-position changes over time.

\subsubsection{Problem Formulation}

The real-valued noisy video samples are denoted by $\XX \in \R^n$ with $n=n_x\cdot n_y\cdot n_t$. The forward operator for the dynamic denoising problem is simply given by an $n \times n$ identity operator, i.e.\ $\Ad = \Id_{n}$.

\subsubsection{Experimental Set-Up}

For training and testing, we used video samples from the benchmark dataset for multiple object tracking  \cite{MOT16}, containing both dynamic and static camera scenes. For training and validation, we scaled the resolution of the video samples by $0.5$ in each direction and extracted patches of size $n_x\times n_y\times n_t = 192 \times 192 \times 32$. During the training process, we used $1751$ patches for training and $1000$ patches from different video samples for validation. We tested the trained model on scaled resolution but the full spatial dimension with $100$ time points per test sample. Gaussian noise with a random standard deviation in the range of $\sigma = 0.1, 0.2, 0.3$ was added to the samples during training. For simplicity and increased training speed, we use a grey-scaled version of the video samples. Because the grey-scaled image data is real-valued, the CNN $u_{\Theta}$ was constructed as described in Section \ref{subsubsec:experiments:cardiac}, but with only one output-channel per output-dimension. For this example, we use the same CNN-block $u_{\Theta}$ as in Figure \ref{fig:u_Theta.pdf}. For comparison, we also trained $\lambda^{xyt}$ which holds a single value for both, the spatial and temporal dimension, and $\lambda^{xy, t}$ which holds two different values for the spatial and temporal dimension. During training we used $T=128$ network iterations, for testing we increased the number of iterations to $T=1024$. We minimized the mean squared error (MSE) between denoised and ground truth patches using the Adam optimizer \cite{kingma2014adam} with an initial learning rate of $10^{-4}$. All the training was performed for $100$ epochs, where validation was performed every second epoch.

\subsubsection{Results}

\begin{figure}[!h]

    \begin{minipage}{\linewidth}
        \begin{minipage}{\linewidth}
            \hspace{1.6cm} PDHG $\lambda^{xyt}_{\tilde{\mathrm{P}}}$ \hspace{1.7cm} PDHG $\lambda^{xy,t}_{\tilde{\mathrm{P}}}$ \hspace{1.7cm} PDHG $\LLambda^{xy,t}_\Theta$ \hspace{1.7cm} Target/Noisy
        \end{minipage} 
        
        \begin{minipage}{\linewidth}
            \rotatebox{90}{
                \begin{minipage}{0.4\linewidth}
                    \hspace{-4.2cm }Static camera \hspace{3.5cm} Moving camera
                \end{minipage}
            }
            \begin{minipage}{\linewidth}
                
                \centering
                \resizebox{\linewidth}{!}{
                    \includegraphics[height=3cm]{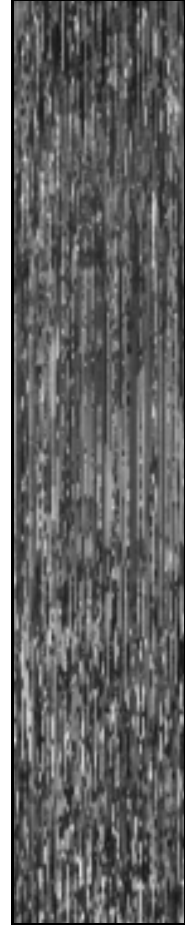}\hspace{-0.1cm}
                    \begin{tikzpicture}[spy using outlines={rectangle, white, magnification=2, size=1.cm, connect spies}]
                        \node[anchor=south west,inner sep=0]  at (0,0) {\includegraphics[height=3cm]{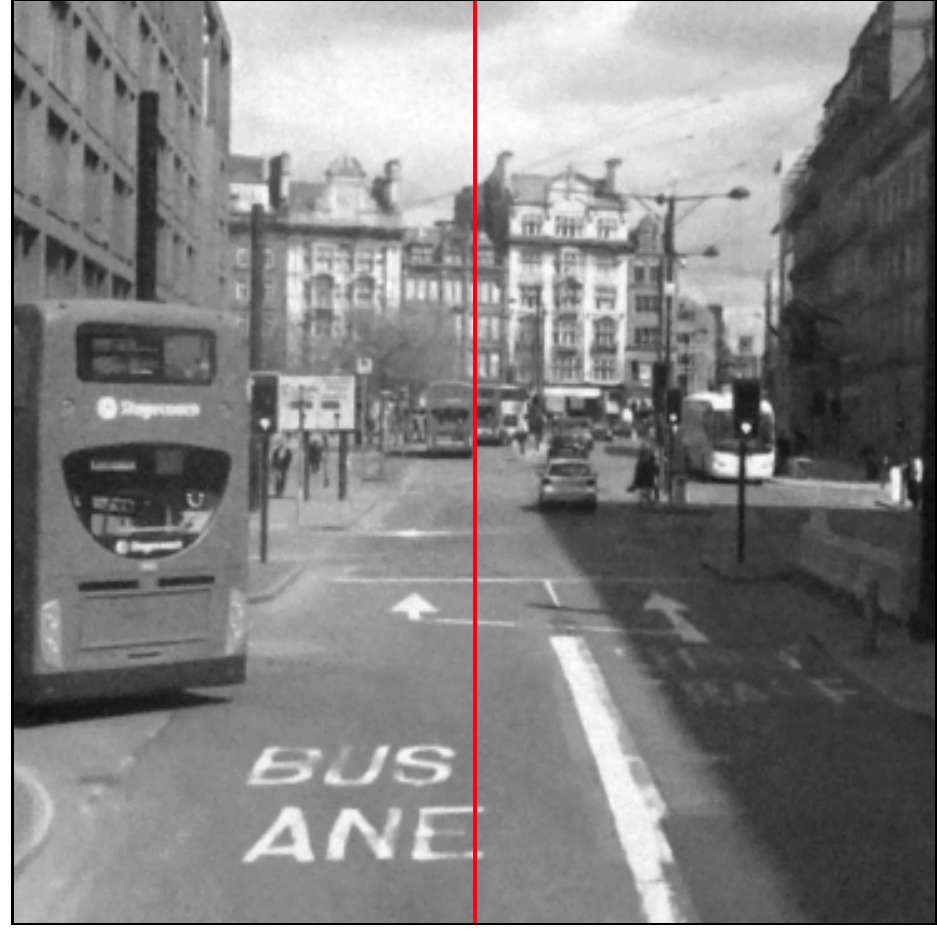}};
                        \spy on (0.8, 1.5) in node [left] at (3.0, 2.5);
                    \end{tikzpicture}
                    \includegraphics[height=3cm]{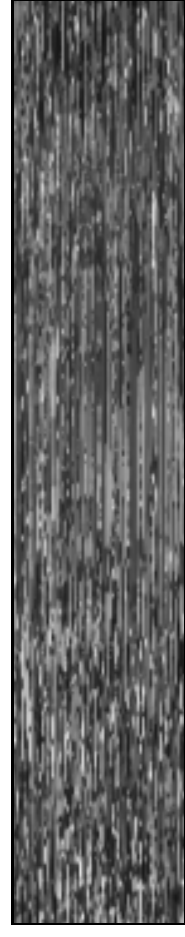}\hspace{-0.1cm}
                    \begin{tikzpicture}[spy using outlines={rectangle, white, magnification=2, size=1.cm, connect spies}]
                        \node[anchor=south west,inner sep=0]  at (0,0) {\includegraphics[height=3cm]{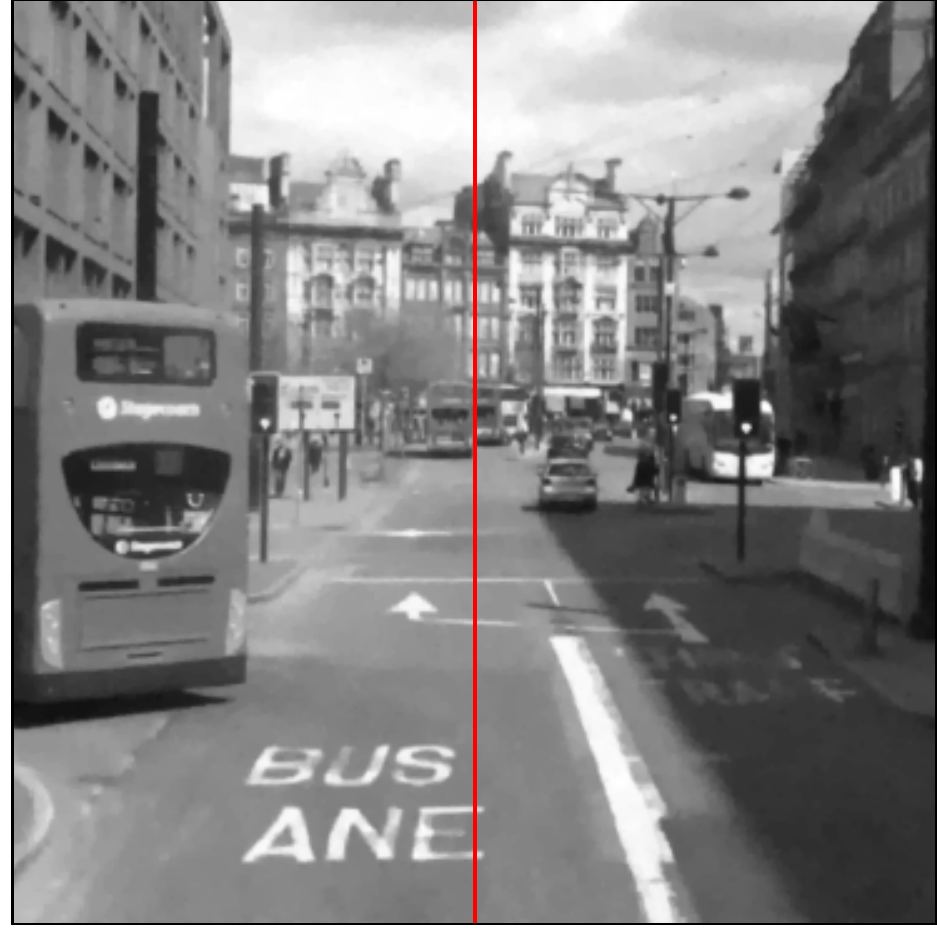}};
                        \spy on (0.8, 1.5) in node [left] at (3.0, 2.5);
                    \end{tikzpicture}
                    \includegraphics[height=3cm]{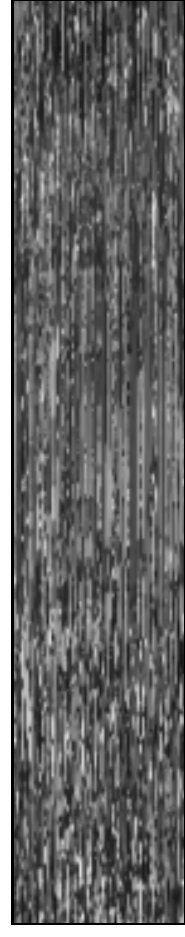}\hspace{-0.1cm}
                    \begin{tikzpicture}[spy using outlines={rectangle, white, magnification=2, size=1.cm, connect spies}]
                        \node[anchor=south west,inner sep=0]  at (0,0) {\includegraphics[height=3cm]{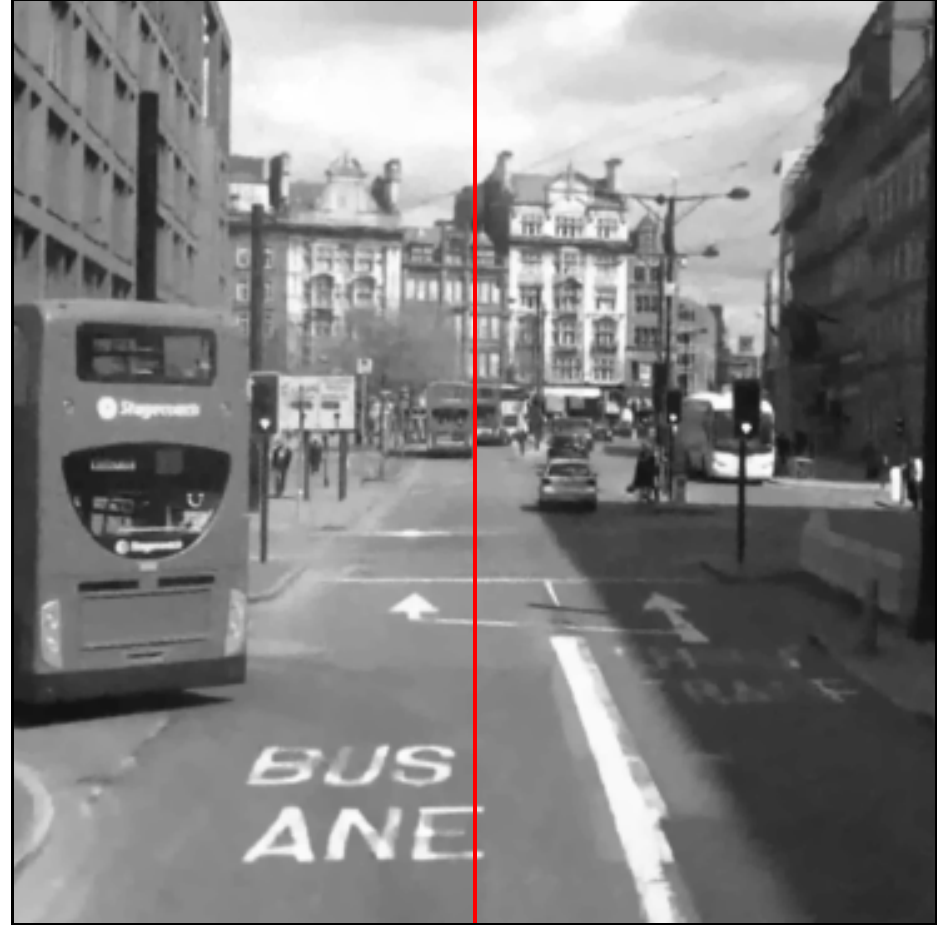}};
                        \spy on (0.8, 1.5) in node [left] at (3.0, 2.5);
                    \end{tikzpicture}
                    \includegraphics[height=3cm]{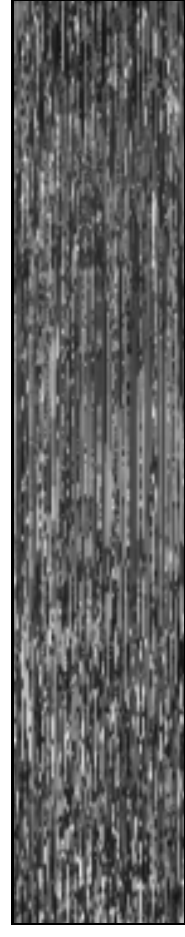}\hspace{-0.1cm}
                    \begin{tikzpicture}[spy using outlines={rectangle, white, magnification=2, size=1.cm, connect spies}]
                        \node[anchor=south west,inner sep=0]  at (0,0) {\includegraphics[height=3cm]{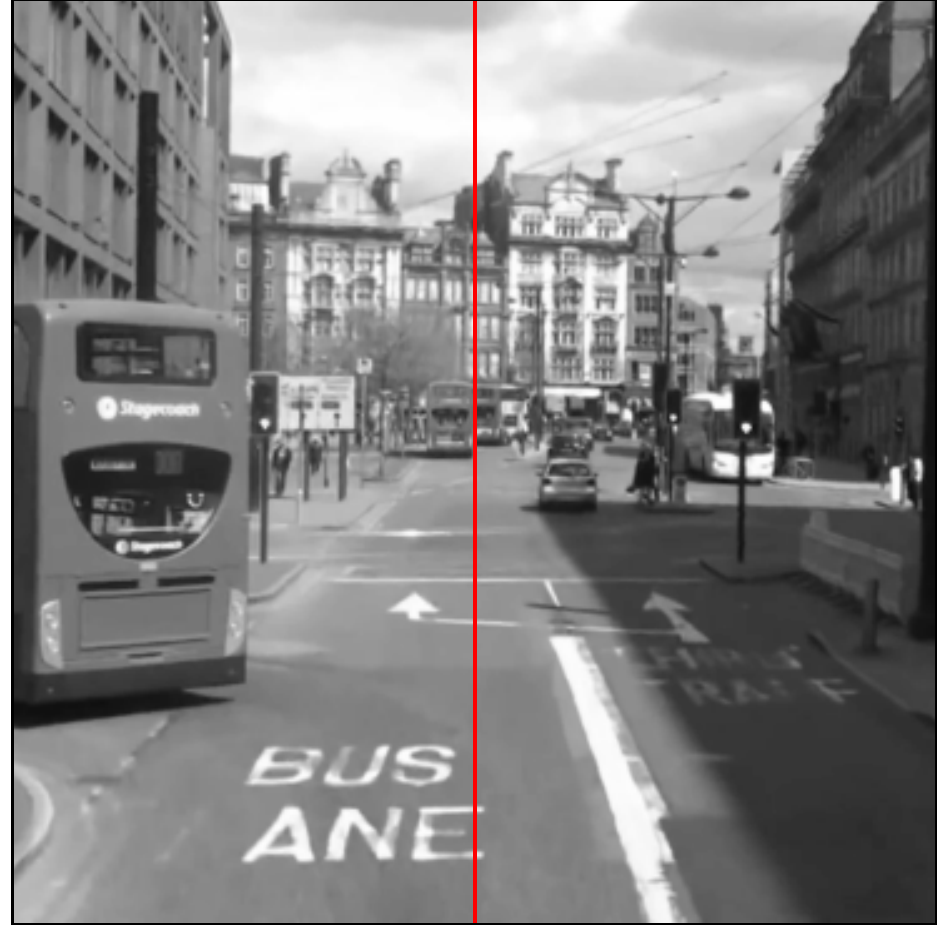}};
                        \spy on (0.8, 1.5) in node [left] at (3.0, 2.5);
                    \end{tikzpicture}
                }
                \resizebox{\linewidth}{!}{
                    \includegraphics[height=3cm]{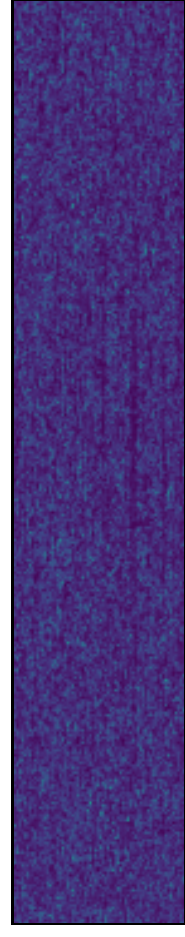}\hspace{-0.1cm}
                    \begin{tikzpicture}[spy using outlines={rectangle, white, magnification=2, size=1.cm, connect spies}]
                        \node[anchor=south west,inner sep=0]  at (0,0) {\includegraphics[height=3cm]{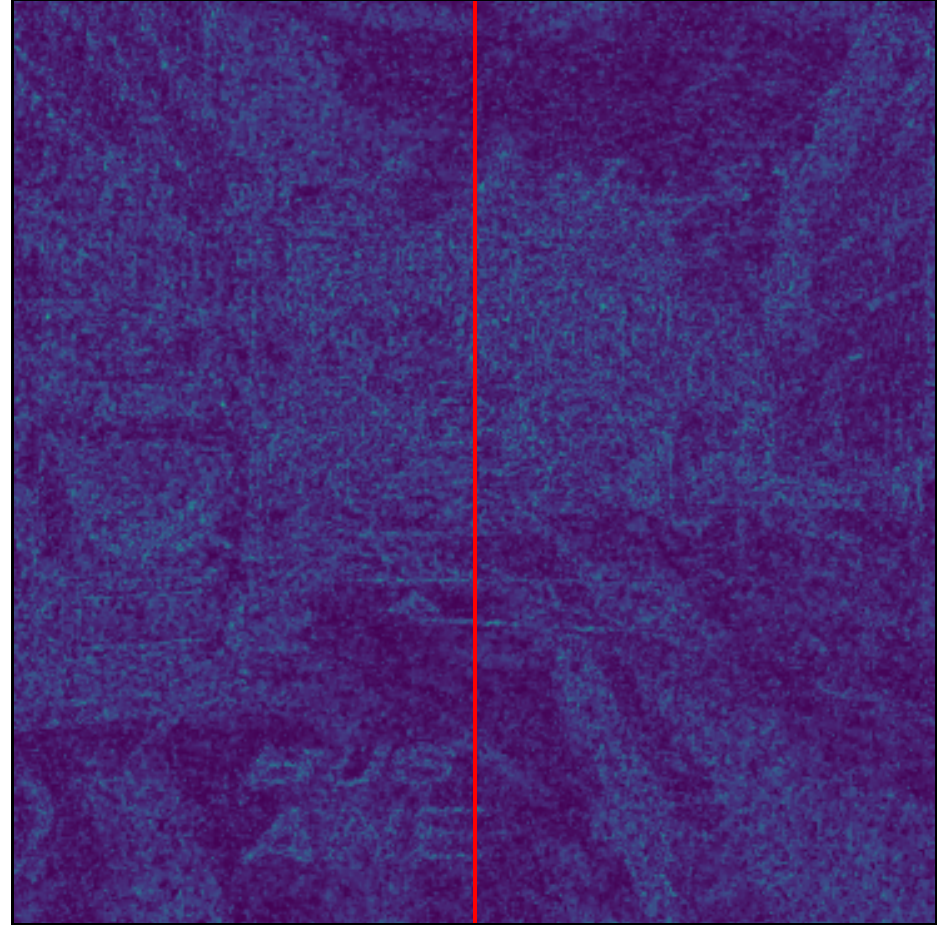}};
                        \spy on (0.8, 1.5) in node [left] at (3.0, 2.5);
                    \end{tikzpicture}
                    \includegraphics[height=3cm]{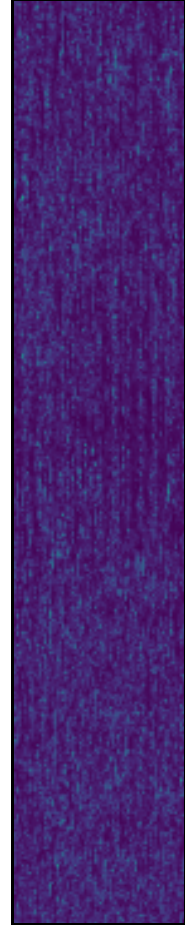}\hspace{-0.1cm}
                    \begin{tikzpicture}[spy using outlines={rectangle, white, magnification=2, size=1.cm, connect spies}]
                        \node[anchor=south west,inner sep=0]  at (0,0) {\includegraphics[height=3cm]{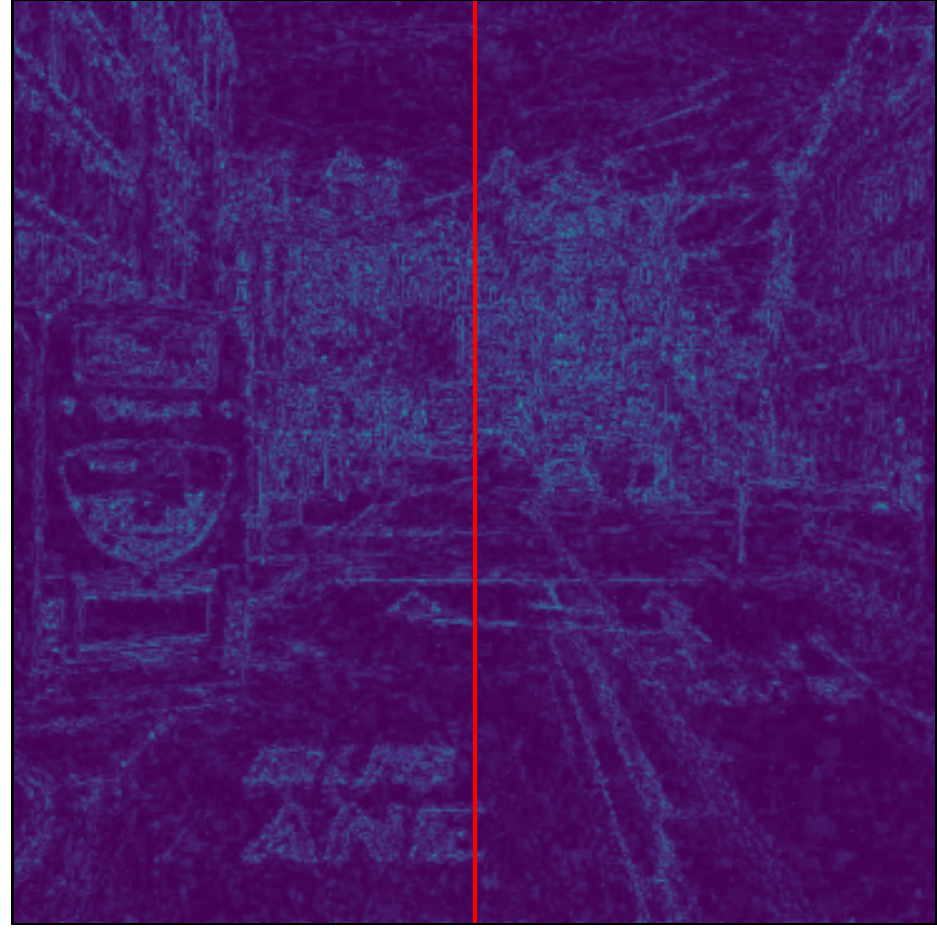}};
                        \spy on (0.8, 1.5) in node [left] at (3.0, 2.5);
                    \end{tikzpicture}
                    \includegraphics[height=3cm]{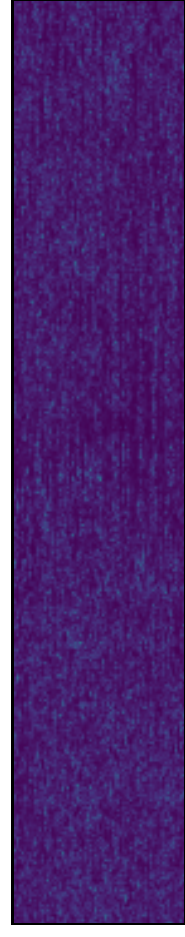}\hspace{-0.1cm}
                    \begin{tikzpicture}[spy using outlines={rectangle, white, magnification=2, size=1.cm, connect spies}]
                        \node[anchor=south west,inner sep=0]  at (0,0) {\includegraphics[height=3cm]{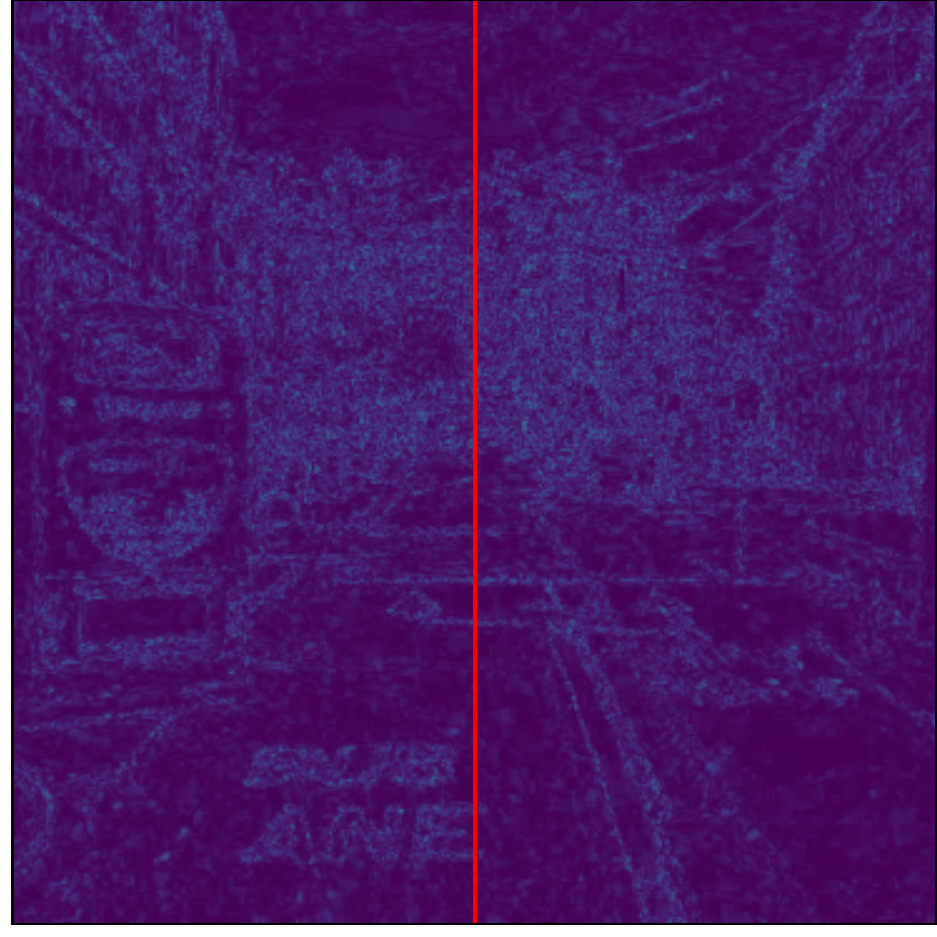}};
                        \spy on (0.8, 1.5) in node [left] at (3.0, 2.5);
                    \end{tikzpicture}
                    \includegraphics[height=3cm]{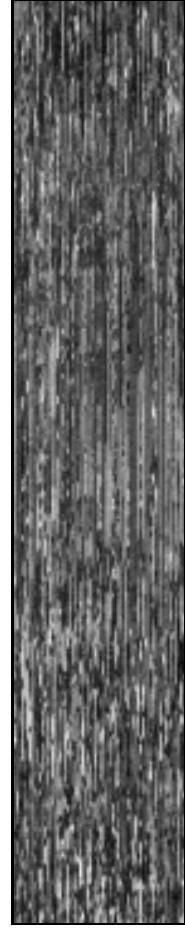}\hspace{-0.1cm}
                    \begin{tikzpicture}[spy using outlines={rectangle, white, magnification=2, size=1.cm, connect spies}]
                        \node[anchor=south west,inner sep=0]  at (0,0) {\includegraphics[height=3cm]{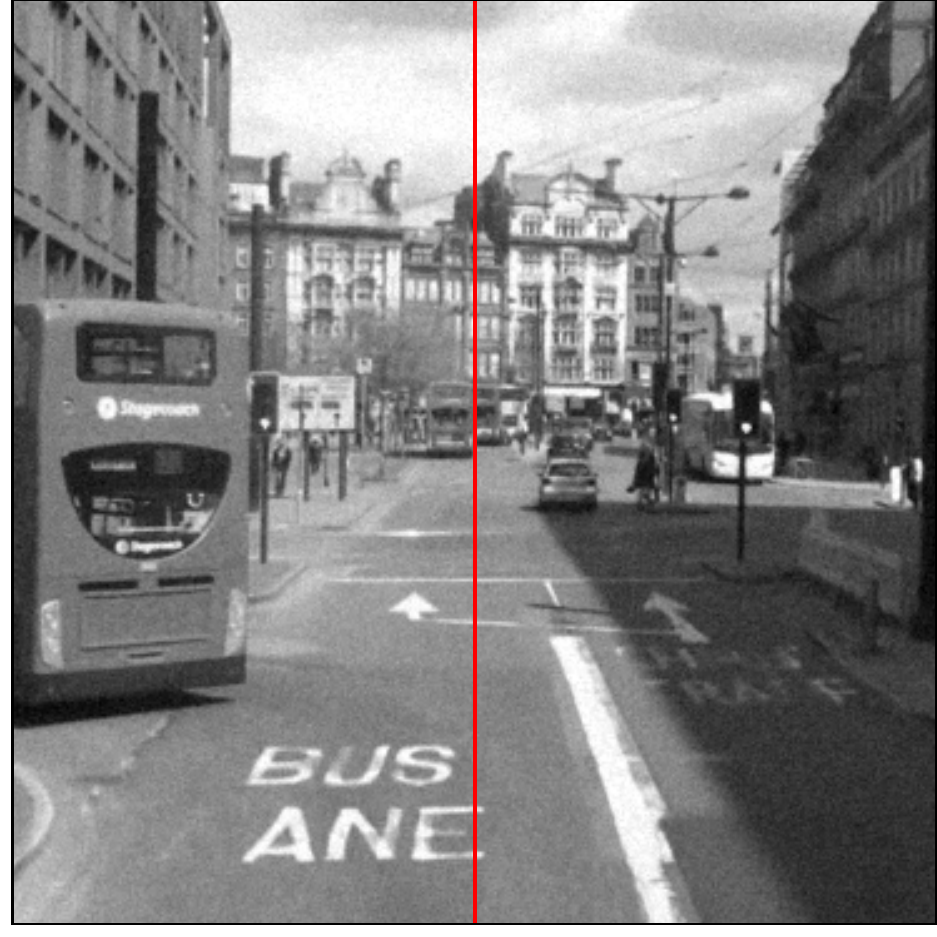}};
                        \spy on (0.8, 1.5) in node [left] at (3.0, 2.5);
                    \end{tikzpicture}
                }

                \resizebox{\linewidth}{!}{
                    \includegraphics[height=3cm]{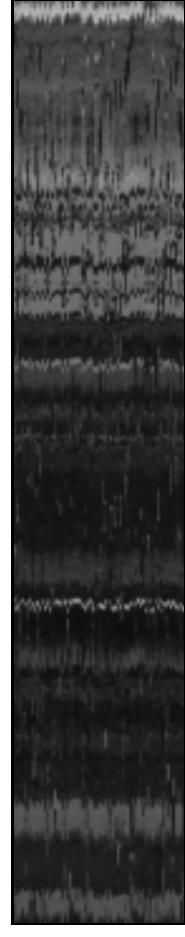}\hspace{-0.1cm}
                    \begin{tikzpicture}[spy using outlines={rectangle, white, magnification=2, size=1.cm, connect spies}]
                        \node[anchor=south west,inner sep=0]  at (0,0) {\includegraphics[height=3cm]{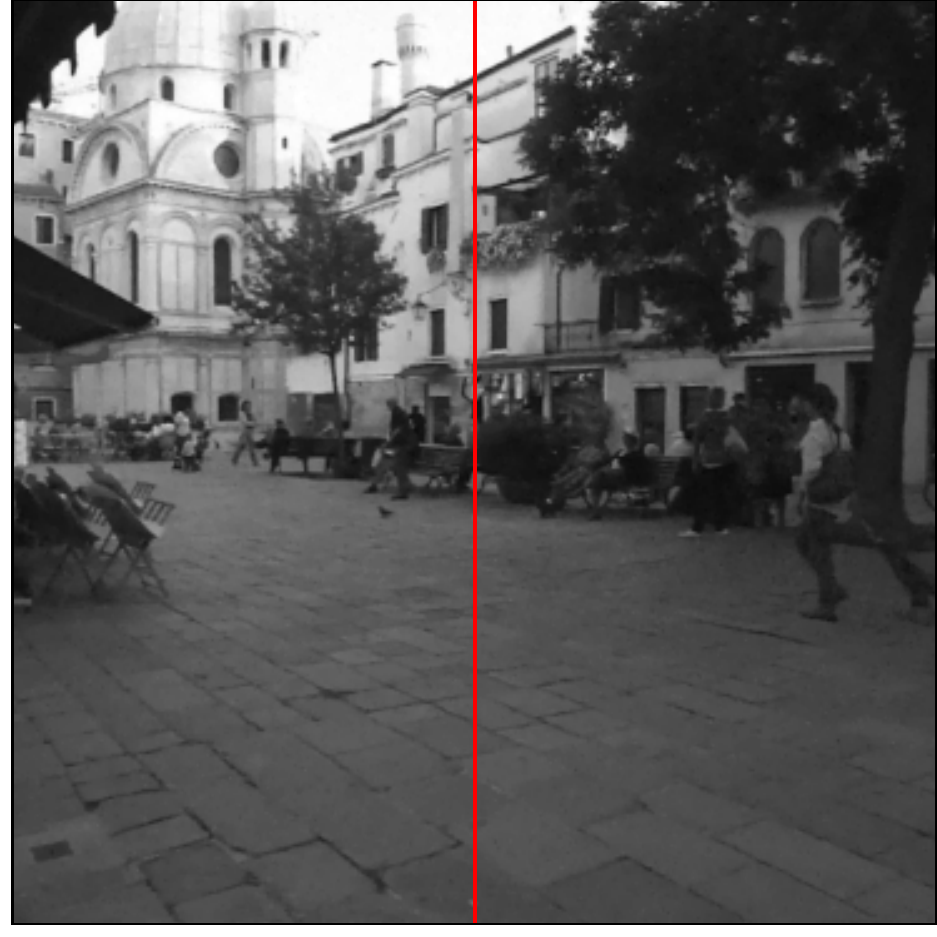}};
                        \spy on (0.8, 1.5) in node [left] at (3.0, 2.5);
                    \end{tikzpicture}
                    \includegraphics[height=3cm]{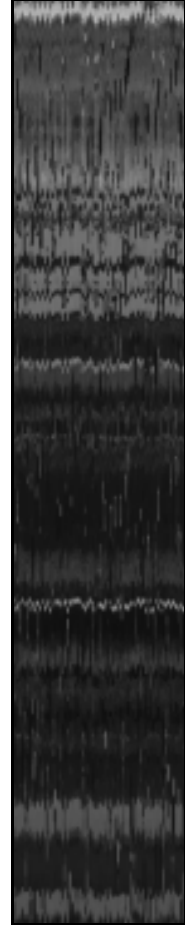}\hspace{-0.1cm}
                    \begin{tikzpicture}[spy using outlines={rectangle, white, magnification=2, size=1.cm, connect spies}]
                        \node[anchor=south west,inner sep=0]  at (0,0) {\includegraphics[height=3cm]{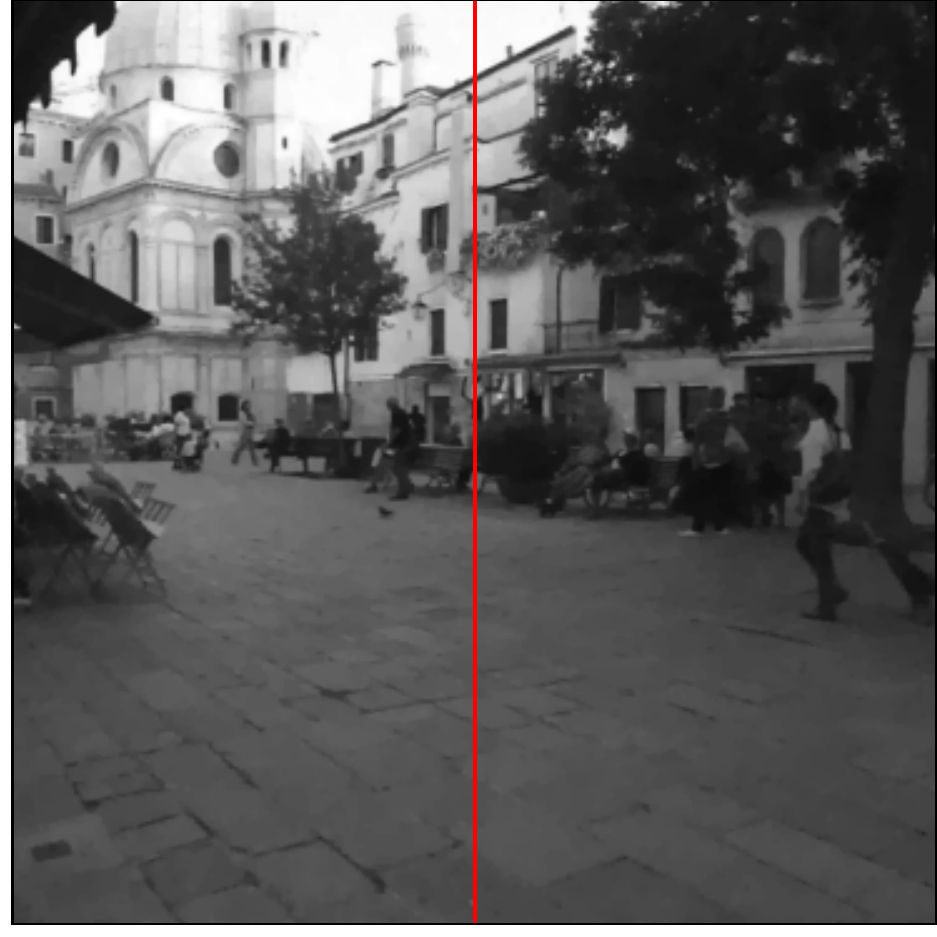}};
                        \spy on (0.8, 1.5) in node [left] at (3.0, 2.5);
                    \end{tikzpicture}
                    \includegraphics[height=3cm]{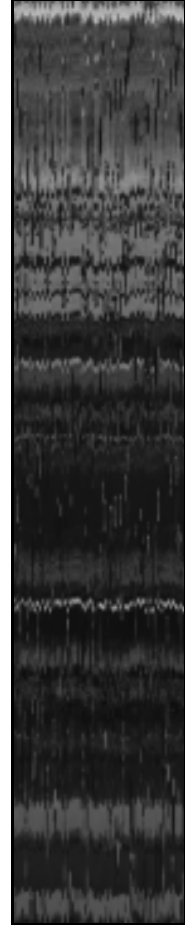}\hspace{-0.1cm}
                    \begin{tikzpicture}[spy using outlines={rectangle, white, magnification=2, size=1.cm, connect spies}]
                        \node[anchor=south west,inner sep=0]  at (0,0) {\includegraphics[height=3cm]{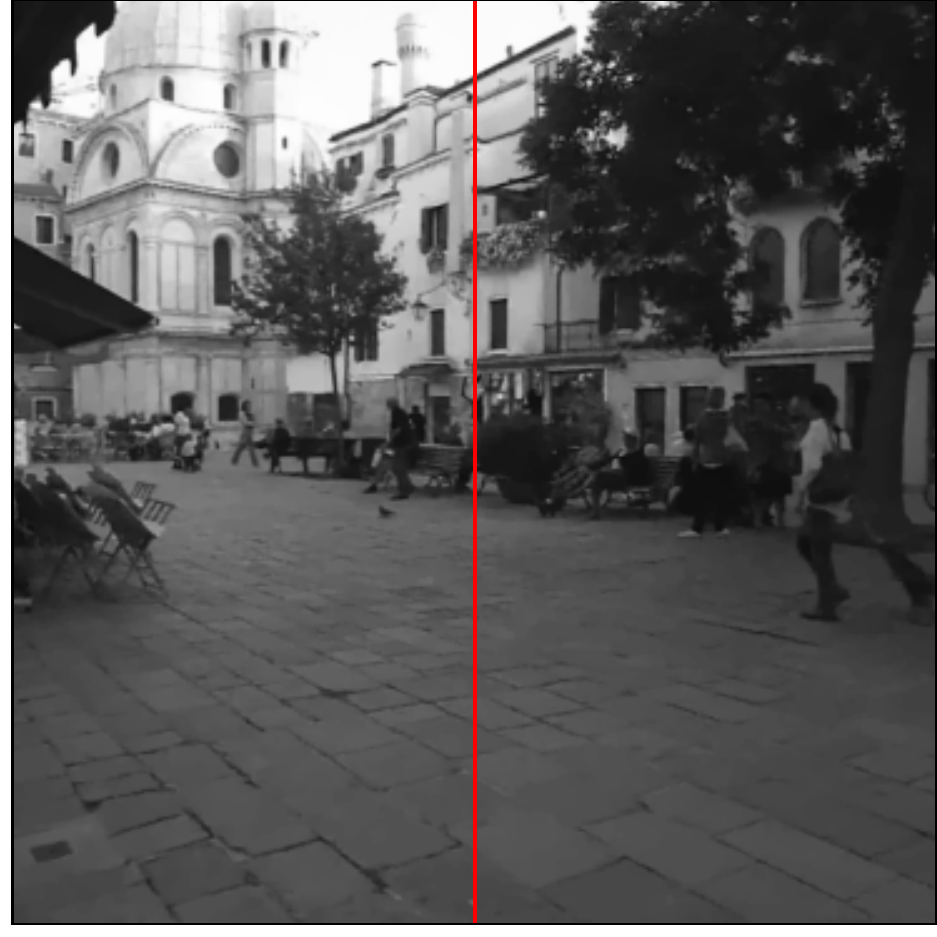}};
                        \spy on (0.8, 1.5) in node [left] at (3.0, 2.5);
                    \end{tikzpicture}
                    \includegraphics[height=3cm]{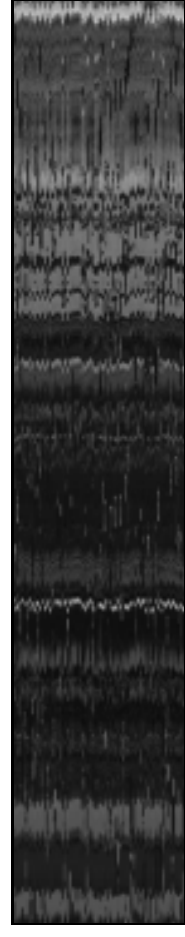}\hspace{-0.1cm}
                    \begin{tikzpicture}[spy using outlines={rectangle, white, magnification=2, size=1.cm, connect spies}]
                        \node[anchor=south west,inner sep=0]  at (0,0) {\includegraphics[height=3cm]{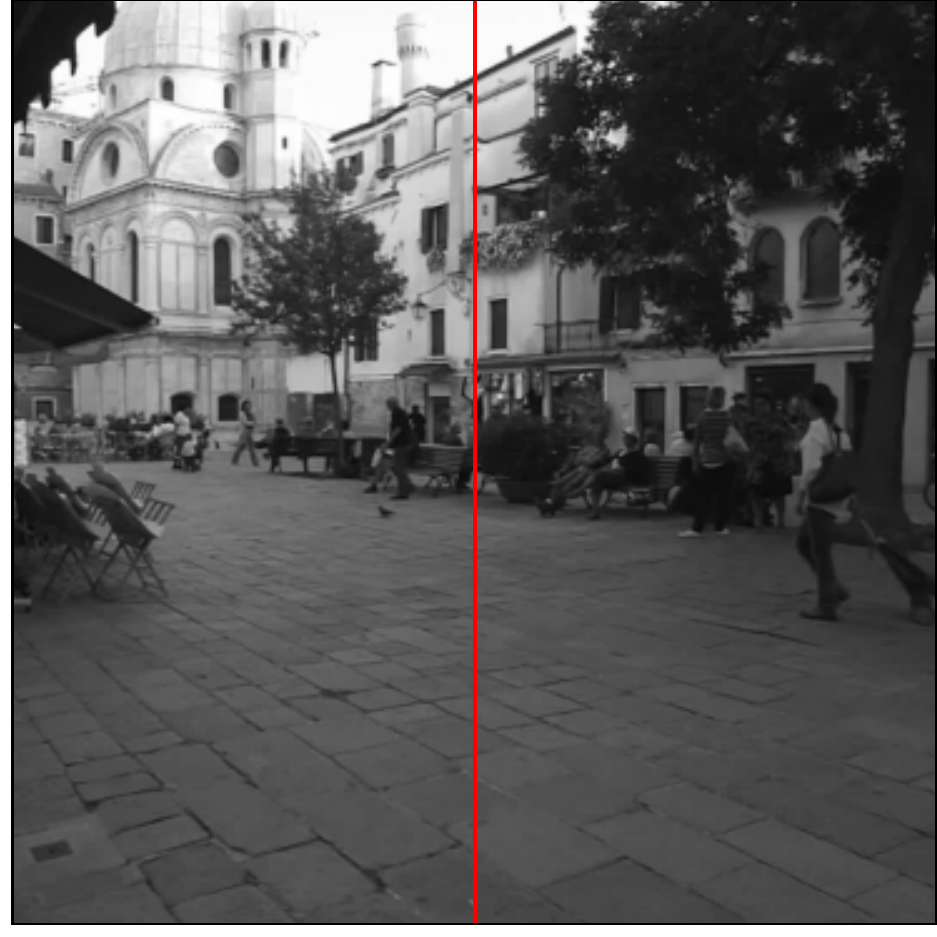}};
                        \spy on (0.8, 1.5) in node [left] at (3.0, 2.5);
                    \end{tikzpicture}
                }
                \resizebox{\linewidth}{!}{
                    \includegraphics[height=3cm]{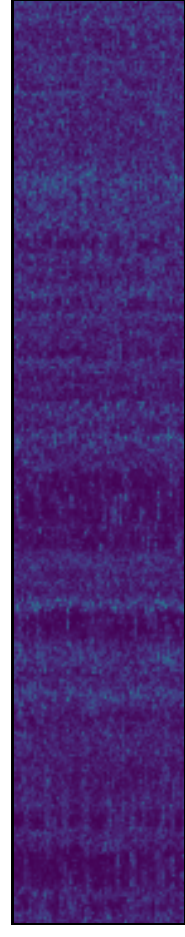}\hspace{-0.1cm}
                    \begin{tikzpicture}[spy using outlines={rectangle, white, magnification=2, size=1.cm, connect spies}]
                        \node[anchor=south west,inner sep=0]  at (0,0) {\includegraphics[height=3cm]{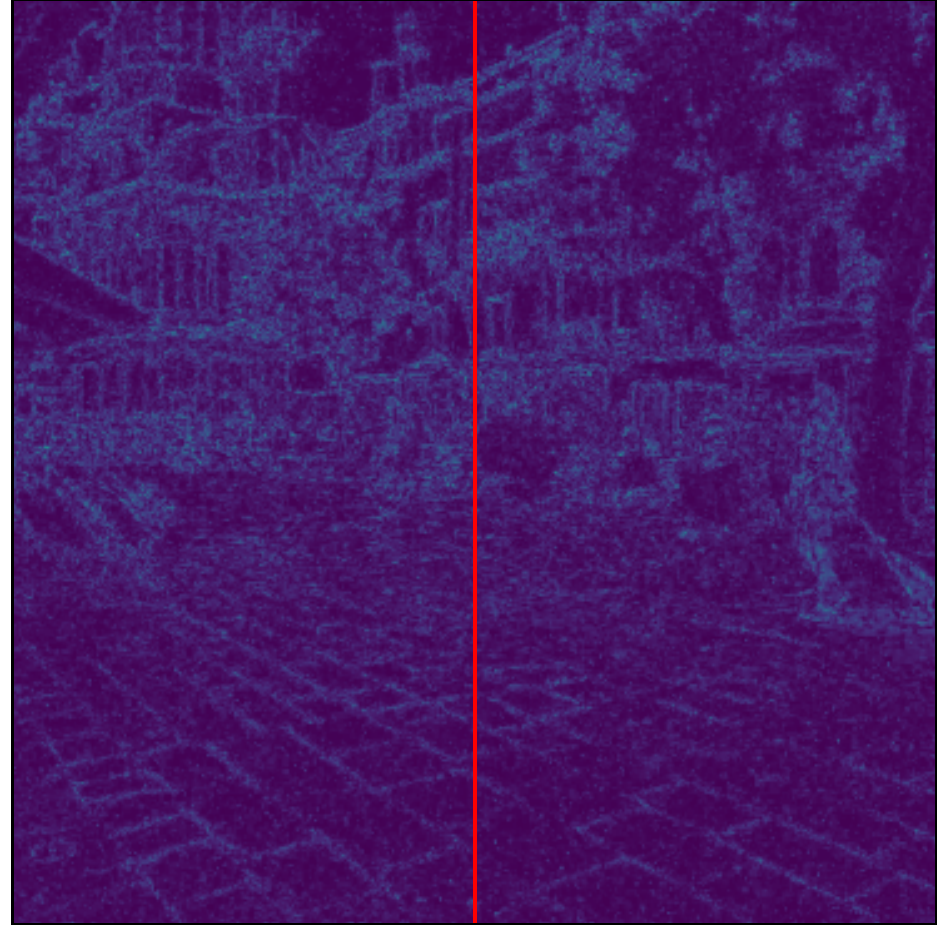}};
                        \spy on (0.8, 1.5) in node [left] at (3.0, 2.5);
                    \end{tikzpicture}
                    \includegraphics[height=3cm]{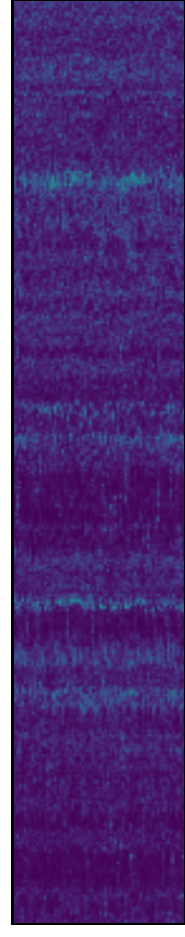}\hspace{-0.1cm}
                    \begin{tikzpicture}[spy using outlines={rectangle, white, magnification=2, size=1.cm, connect spies}]
                        \node[anchor=south west,inner sep=0]  at (0,0) {\includegraphics[height=3cm]{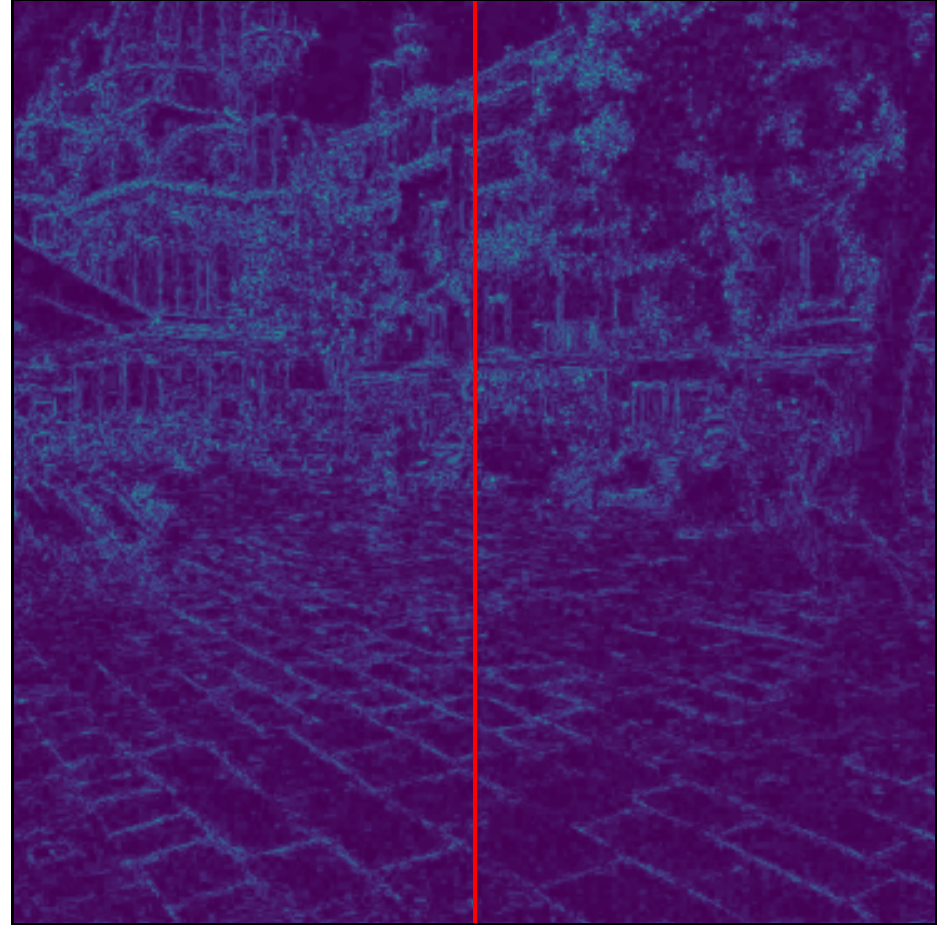}};
                        \spy on (0.8, 1.5) in node [left] at (3.0, 2.5);
                    \end{tikzpicture}
                    \includegraphics[height=3cm]{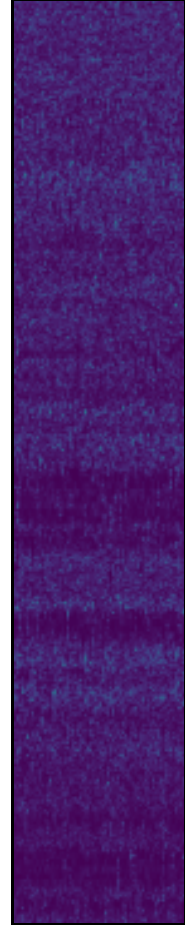}\hspace{-0.1cm}
                    \begin{tikzpicture}[spy using outlines={rectangle, white, magnification=2, size=1.cm, connect spies}]
                        \node[anchor=south west,inner sep=0]  at (0,0) {\includegraphics[height=3cm]{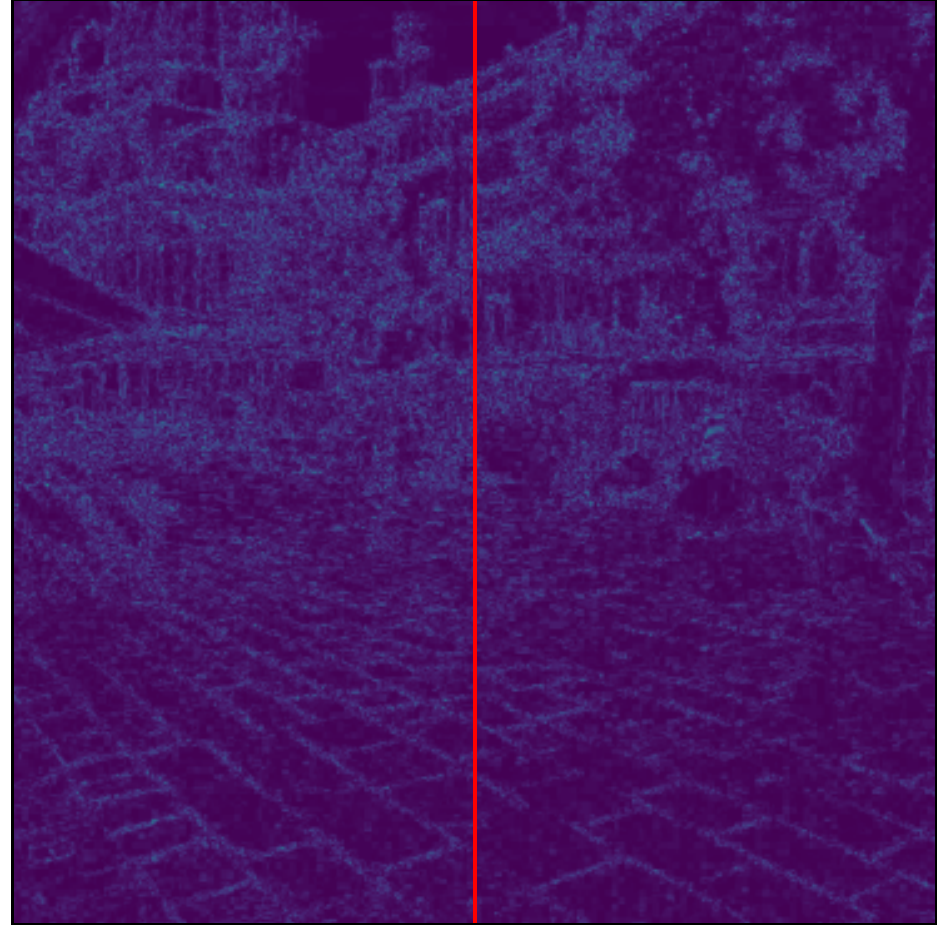}};
                        \spy on (0.8, 1.5) in node [left] at (3.0, 2.5);
                    \end{tikzpicture}
                    \includegraphics[height=3cm]{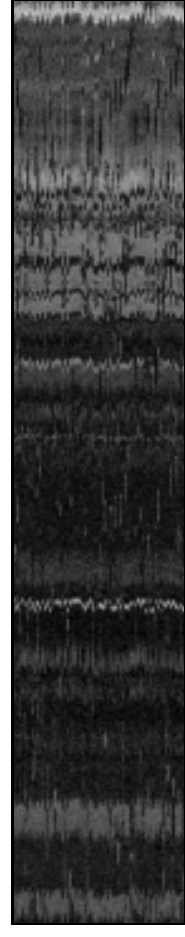}\hspace{-0.1cm}
                    \begin{tikzpicture}[spy using outlines={rectangle, white, magnification=2, size=1.cm, connect spies}]
                        \node[anchor=south west,inner sep=0]  at (0,0) {\includegraphics[height=3cm]{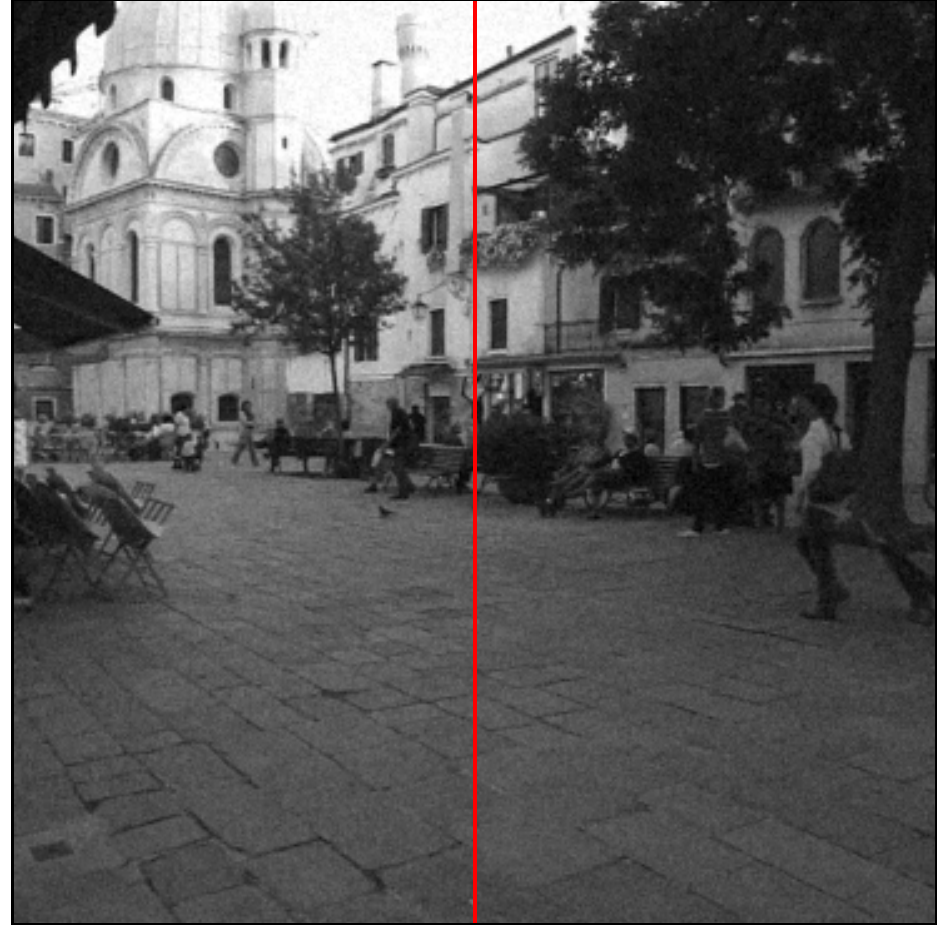}};
                        \spy on (0.8, 1.5) in node [left] at (3.0, 2.5);
                    \end{tikzpicture}
                }
            \end{minipage}
        \end{minipage}
    \end{minipage}
    
    \caption{An example of dynamic denoising with the PDHG from Algorithm \ref{algo:tv_reco_algo_mri} for different choices of regularization parameters at moving (\url{https://motchallenge.net/vis/MOT17-14}) and static (\url{https://motchallenge.net/vis/MOT17-01}) camera view. Single scalar regularization parameter $\lambda_{\mathrm{\tilde{P}}}^{xyt}$, two scalar regularization parameters for differently weighted spatial and temporal components $\lambda_{\mathrm{\tilde{P}}}^{xy,t}$, and the proposed spatially and temporally dependent parameter-map $\boldsymbol{\Lambda}_{\Theta}^{xy,t}$ obtained with the network $\mathcal{N}_{\Theta}^T$. The last column shows the target image and the noisy sample. The row underneath the denoised image shows the error map.}
    \label{fig:dynamic_denoising:method_comparison}
\end{figure}

\begin{figure}[!h]
    \centering
    \includegraphics[width=0.5\linewidth]{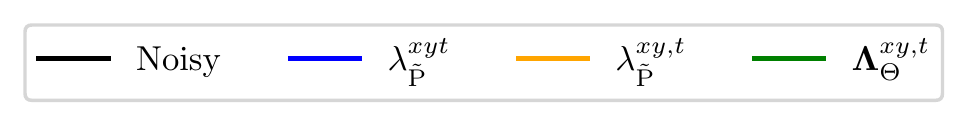}\\
    \includegraphics[width=0.325\linewidth]{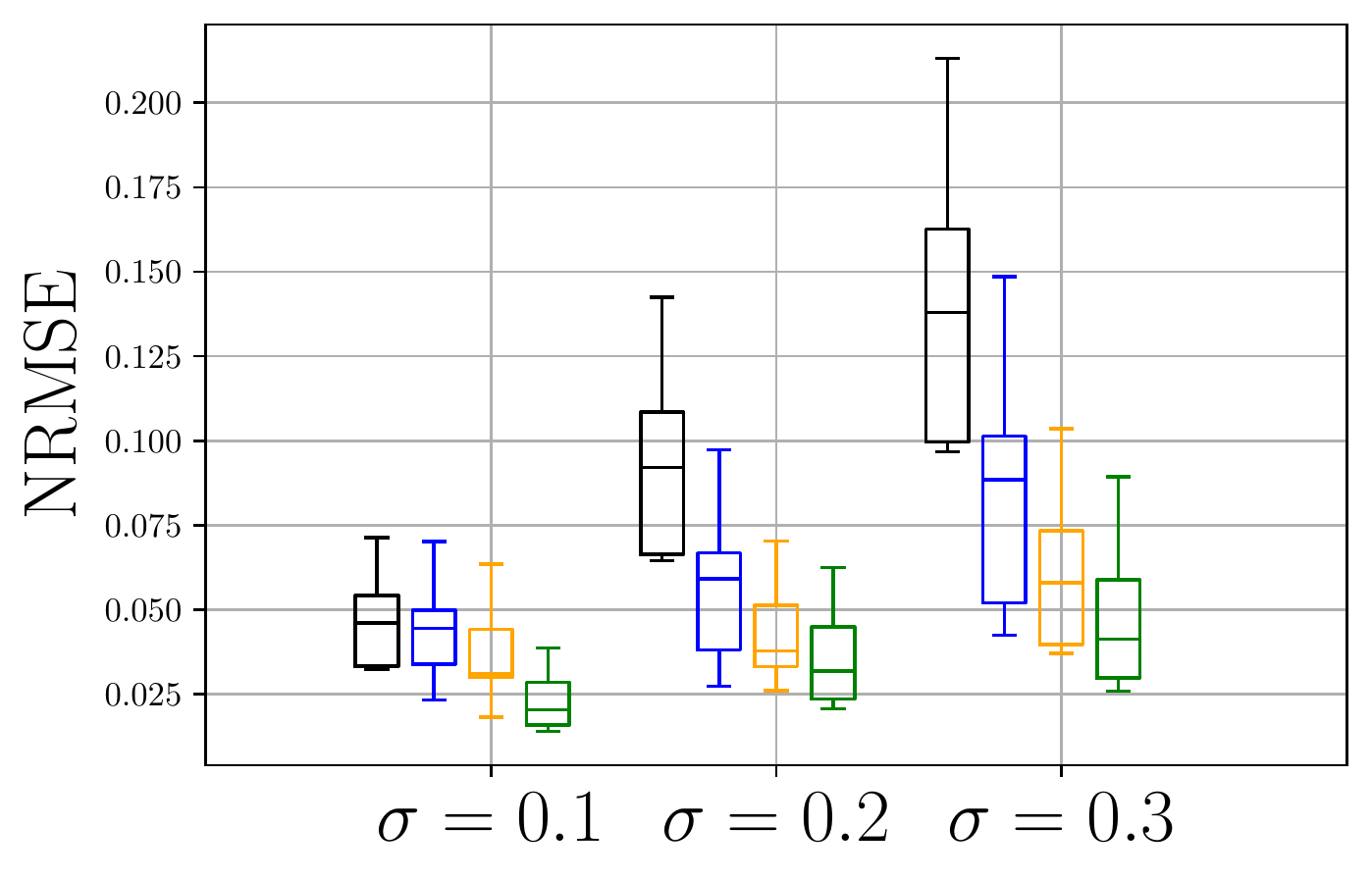}
    \includegraphics[width=0.325\linewidth]{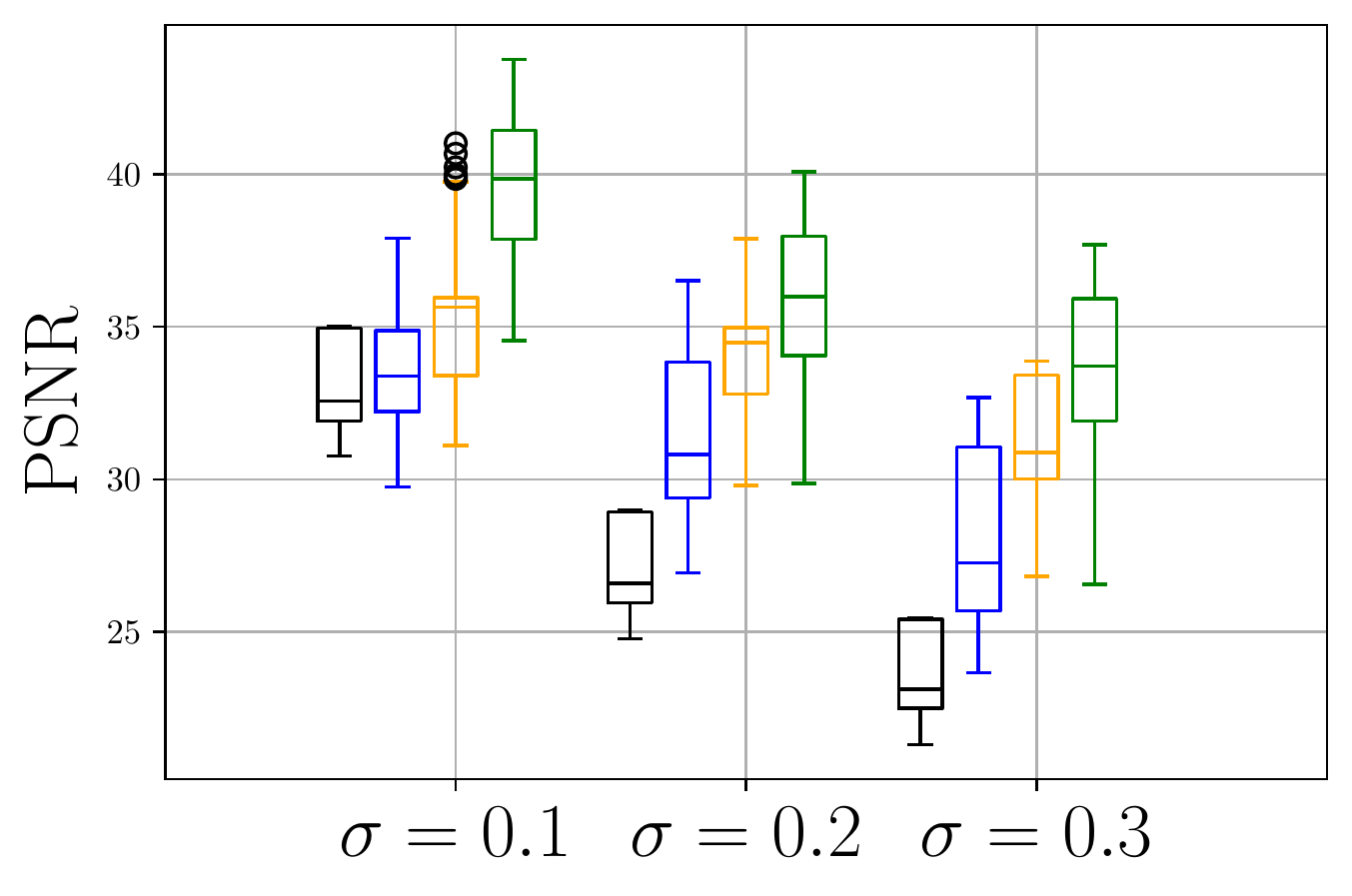}
    \includegraphics[width=0.325\linewidth]{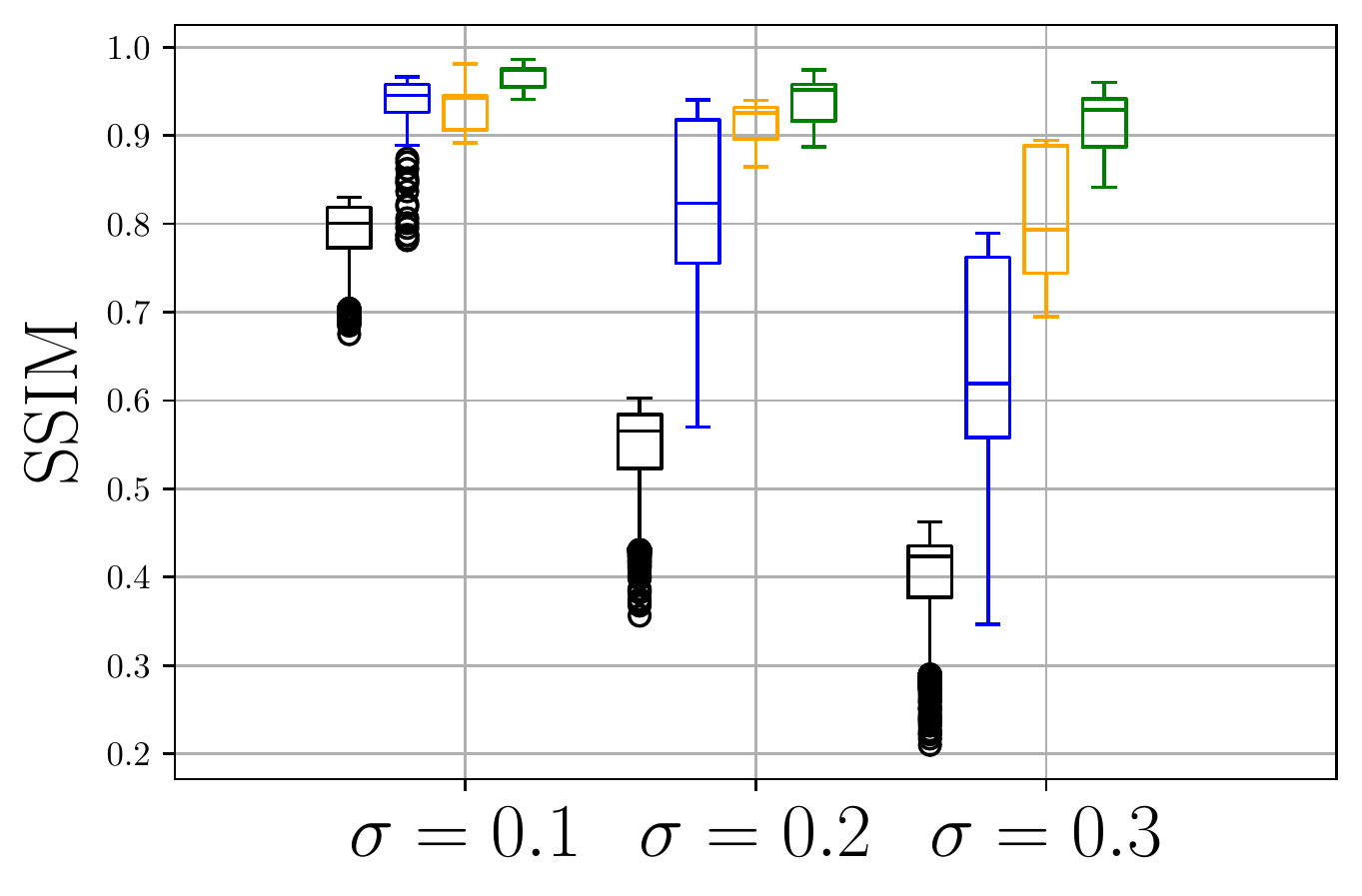}
    \caption{Box-plots summarizing the results in terms of PSNR, NRMSE and SSIM obtained with the PDHG algorithm for dynamic denoising. Single scalar regularization parameter  ($\lambda_{\tilde{\mathrm{P}}}^{xyt}$, blue), two scalar regularization parameters; one for the spatial $x$- and $y$-direction, one for the temporal $t$-direction ($\lambda_{\tilde{\mathrm{P}}}^{xy,t}$, orange) and the proposed spatially and temporally dependent parameter-map $\boldsymbol{\Lambda}_{\Theta}^{xy,t}$ obtained with a CNN ($\mathrm{NET}_{\Theta}$, green).}
    \label{fig:dynamic_denoising:box_plots}
\end{figure}

\begin{figure}[!h]

    \begin{minipage}{\linewidth}
        \begin{minipage}{\linewidth}
            \hspace{1.6cm} Noisy \hspace{1.4cm} PDHG $\LLambda^{xy,t}_\Theta$ \hspace{1.6cm} $\LLambda^{xy}_\Theta$ \hspace{2.2cm} $\LLambda^{t}_\Theta$ \hspace{1.8cm} Target
        \end{minipage} 
        
        \begin{minipage}{\linewidth}
            \rotatebox{90}{
                \begin{minipage}{0.16\linewidth}
                    \hspace{-2.4cm }Static camera \hspace{0.1cm} Moving camera
                \end{minipage}
            }
            \begin{minipage}{\linewidth}

                \centering
                \resizebox{\linewidth}{!}{
                    \includegraphics[height=3cm]{figures/DynDenoising/sample_14_sigma_02/noisy_yt.pdf}\hspace{-0.1cm}
                    \begin{tikzpicture}[spy using outlines={rectangle, white, magnification=2, size=1.cm, connect spies}]
                        \node[anchor=south west,inner sep=0]  at (0,0) {\includegraphics[height=3cm]{figures/DynDenoising/sample_14_sigma_02/noisy_xy.pdf}};
                    \end{tikzpicture}
                    \includegraphics[height=3cm]{figures/DynDenoising/sample_14_sigma_02/cnn_yt.pdf}\hspace{-0.1cm}
                    \begin{tikzpicture}[spy using outlines={rectangle, white, magnification=2, size=1.cm, connect spies}]
                        \node[anchor=south west,inner sep=0]  at (0,0) {\includegraphics[height=3cm]{figures/DynDenoising/sample_14_sigma_02/cnn_xy.pdf}};
                    \end{tikzpicture}
                    \includegraphics[height=3cm]{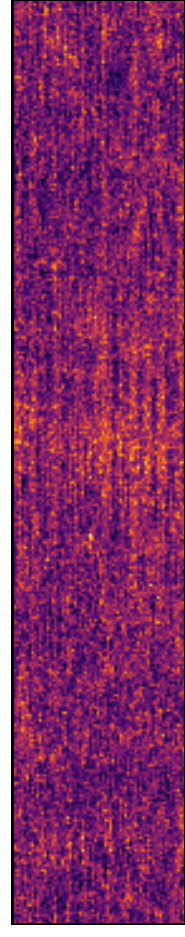}\hspace{-0.1cm}
                    \begin{tikzpicture}[spy using outlines={rectangle, white, magnification=2, size=1.cm, connect spies}]
                        \node[anchor=south west,inner sep=0]  at (0,0) {\includegraphics[height=3cm]{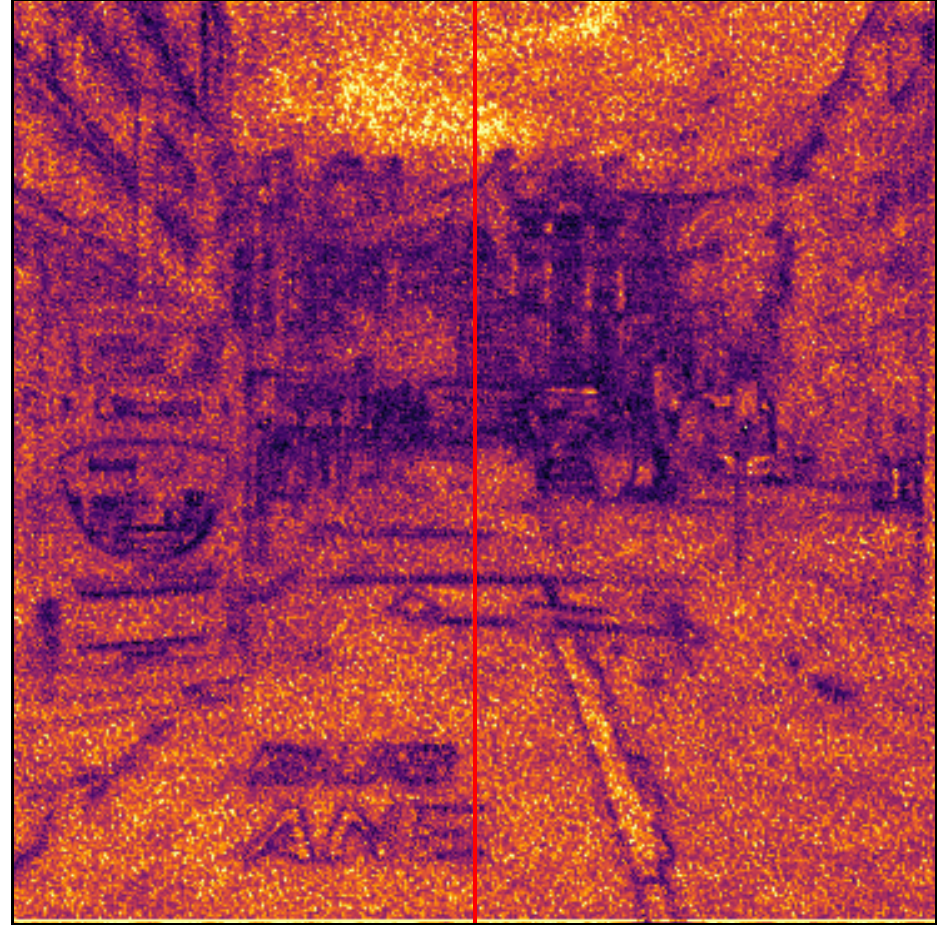}};
                    \end{tikzpicture}
                    \includegraphics[height=3cm]{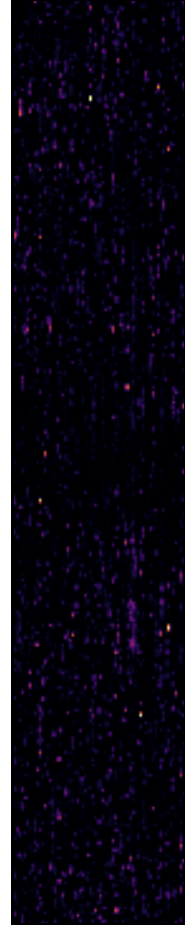}\hspace{-0.1cm}
                    \begin{tikzpicture}[spy using outlines={rectangle, white, magnification=2, size=1.cm, connect spies}]
                        \node[anchor=south west,inner sep=0]  at (0,0) {\includegraphics[height=3cm]{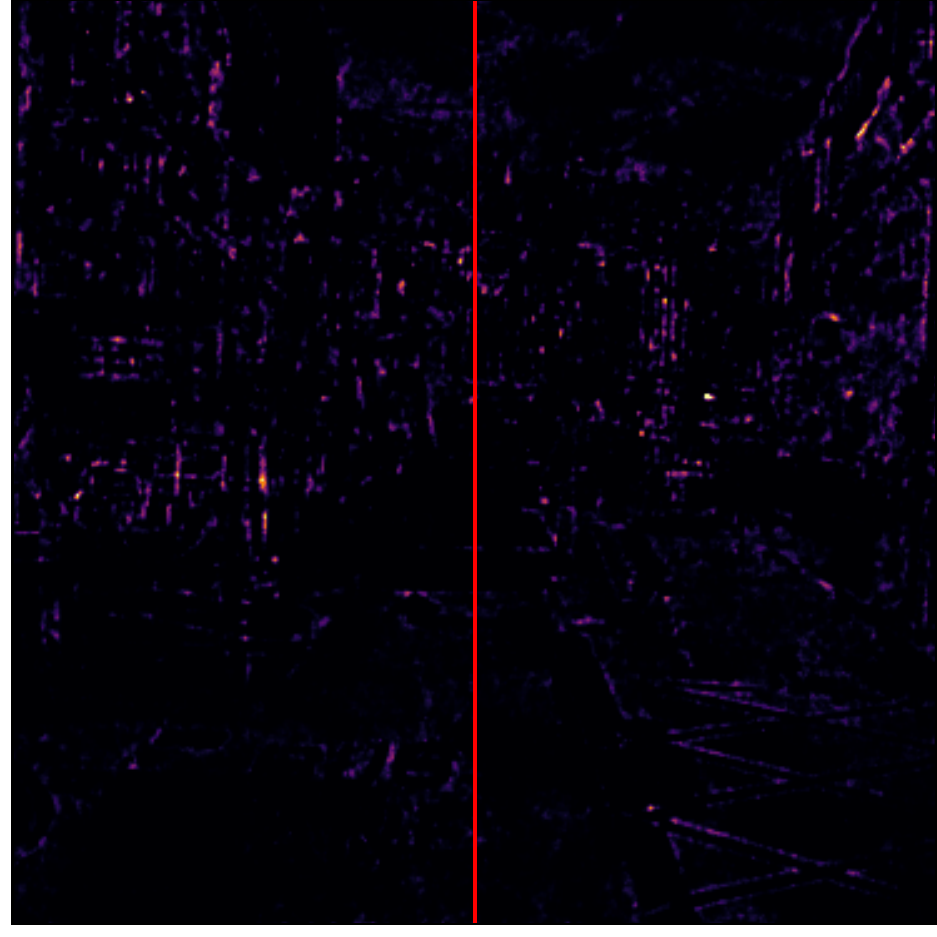}};
                    \end{tikzpicture}
                    \includegraphics[height=3cm]{figures/DynDenoising/sample_14_sigma_02/ground_truth_yt.pdf}\hspace{-0.1cm}
                    \begin{tikzpicture}[spy using outlines={rectangle, white, magnification=2, size=1.cm, connect spies}]
                        \node[anchor=south west,inner sep=0]  at (0,0) {\includegraphics[height=3cm]{figures/DynDenoising/sample_14_sigma_02/ground_truth_xy.pdf}};
                    \end{tikzpicture}
                }
                \resizebox{\linewidth}{!}{
                    \includegraphics[height=3cm]{figures/DynDenoising/sample_01_sigma_02/noisy_yt.pdf}\hspace{-0.1cm}
                    \begin{tikzpicture}[spy using outlines={rectangle, white, magnification=2, size=1.cm, connect spies}]
                        \node[anchor=south west,inner sep=0]  at (0,0) {\includegraphics[height=3cm]{figures/DynDenoising/sample_01_sigma_02/noisy_xy.pdf}};
                    \end{tikzpicture}
                    \includegraphics[height=3cm]{figures/DynDenoising/sample_01_sigma_02/cnn_yt.pdf}\hspace{-0.1cm}
                    \begin{tikzpicture}[spy using outlines={rectangle, white, magnification=2, size=1.cm, connect spies}]
                        \node[anchor=south west,inner sep=0]  at (0,0) {\includegraphics[height=3cm]{figures/DynDenoising/sample_01_sigma_02/cnn_xy.pdf}};
                    \end{tikzpicture}
                    \includegraphics[height=3cm]{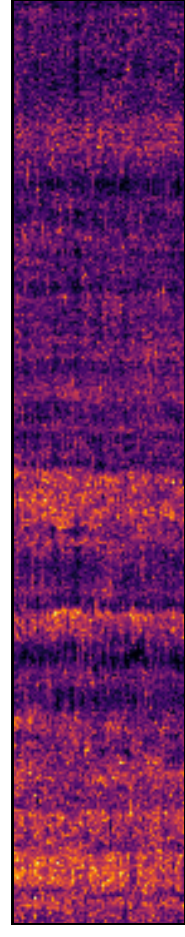}\hspace{-0.1cm}
                    \begin{tikzpicture}[spy using outlines={rectangle, white, magnification=2, size=1.cm, connect spies}]
                        \node[anchor=south west,inner sep=0]  at (0,0) {\includegraphics[height=3cm]{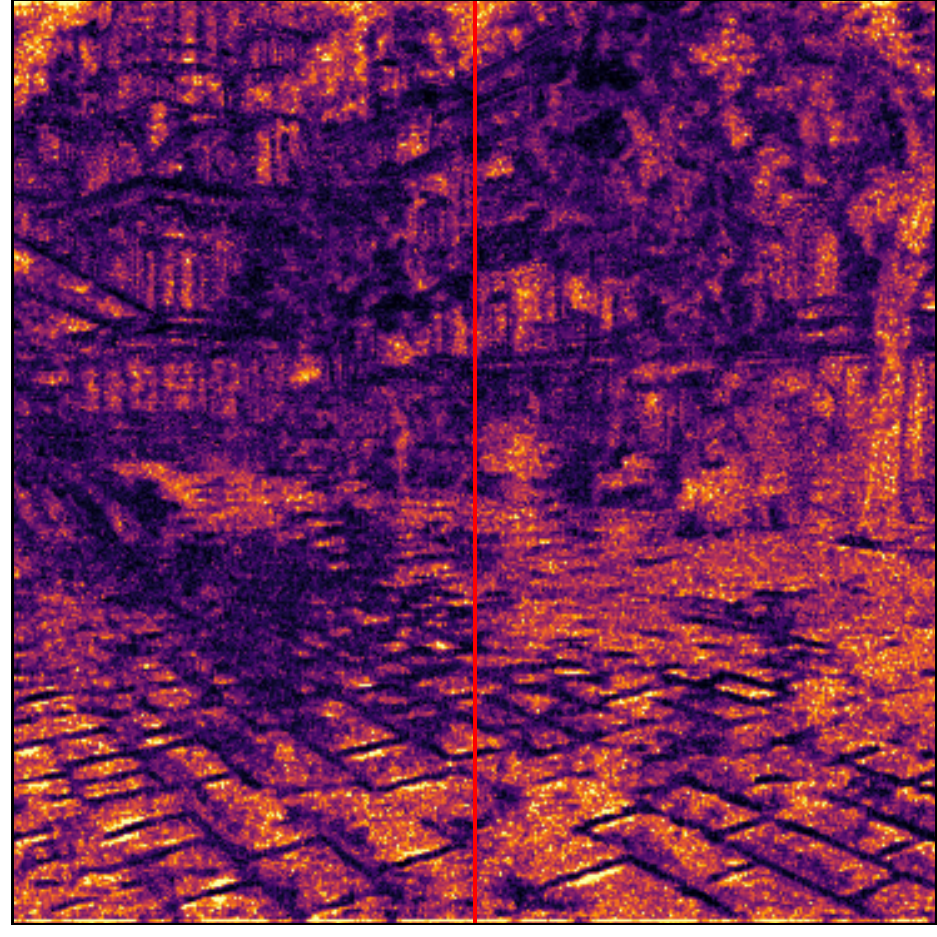}};
                    \end{tikzpicture}
                    \includegraphics[height=3cm]{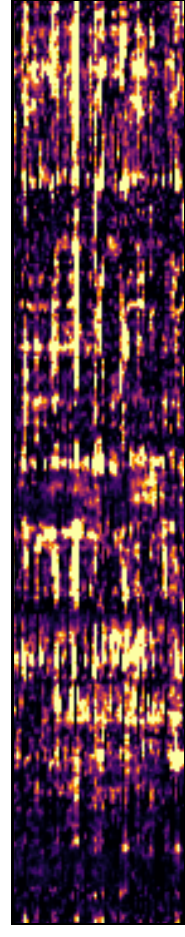}\hspace{-0.1cm}
                    \begin{tikzpicture}[spy using outlines={rectangle, white, magnification=2, size=1.cm, connect spies}]
                        \node[anchor=south west,inner sep=0]  at (0,0) {\includegraphics[height=3cm]{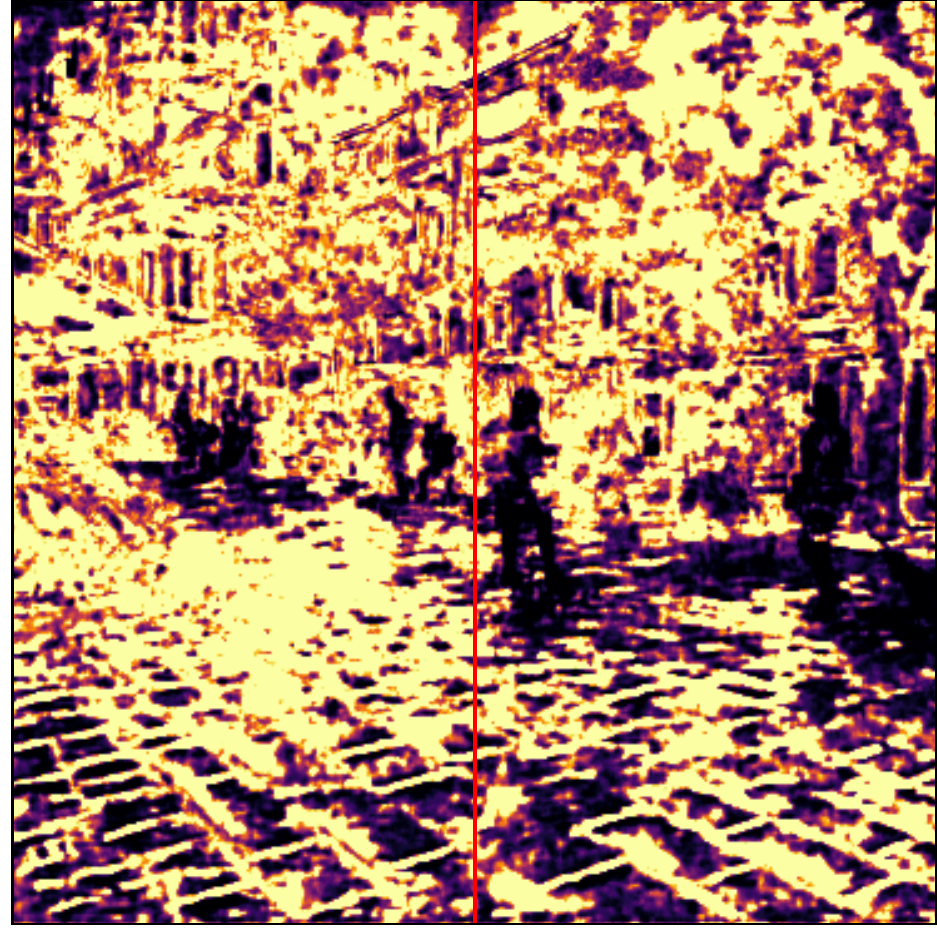}};
                    \end{tikzpicture}
                    \includegraphics[height=3cm]{figures/DynDenoising/sample_01_sigma_02/ground_truth_yt.pdf}\hspace{-0.1cm}
                    \begin{tikzpicture}[spy using outlines={rectangle, white, magnification=2, size=1.cm, connect spies}]
                        \node[anchor=south west,inner sep=0]  at (0,0) {\includegraphics[height=3cm]{figures/DynDenoising/sample_01_sigma_02/ground_truth_xy.pdf}};
                    \end{tikzpicture}
                }
            \end{minipage}
        \end{minipage}

        \begin{minipage}{\linewidth}
            \vspace{0.2cm}
            \begin{center}
                $\LLambda^{xy}_\Theta$\\[0.1cm]
                \includegraphics[width=0.8\linewidth]{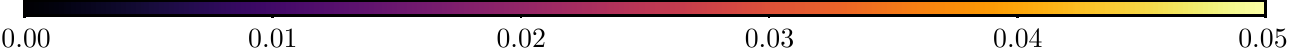}
            \end{center}
        \end{minipage}

        \begin{minipage}{\linewidth}
            \vspace{0.2cm}
            \begin{center}
                $\LLambda^{t}_\Theta$\\[0.1cm]
                \includegraphics[width=0.8\linewidth]{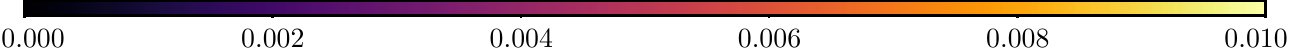}
            \end{center}
        \end{minipage}
        
    \end{minipage}
    
    \caption{Two examples with moving (top row) and static (bottom row) camera view. The columns show the noisy sample, the denoised sample, spatial and temporal dependent parameter-maps $\boldsymbol{\Lambda}_{\Theta}^{xy}$ and $\boldsymbol{\Lambda}_{\Theta}^{t}$ and the ground truth sample. By comparing the top and the bottom row, we can see that the CNN-block $\mathrm{NET}_{\Theta}$ is able to differentiate between scenes with static and dynamic camera positions as it exploits the higher temporal correlation in the first by assigning higher temporal regularization values and at the same time provides lower temporal regularization values when there is less temporal correlation to exploit due to the moving camera position. }
    \label{fig:dynamic_denoising:parameter_maps}    
\end{figure}

We compare the 2D time frames of the video samples from the test dataset to the denoised frames, regularized by $\lambda_{\tilde{\mathrm{P}}}^{xyt}$, $\lambda_{\tilde{\mathrm{P}}}^{xy,t}$ and by the spatio-temporal parameter-map $\boldsymbol{\Lambda}_{\Theta}^{xy,t}$. The metrics were calculated frame-wise for all samples at three different noise levels, characterized by the standard deviation of the Gaussian distribution. From the box-plots in Figure \ref{fig:dynamic_denoising:box_plots} we see that the PDHG reconstructions using the proposed spatio-temporal regularization parameter-map yield superior reconstructions compared to $\lambda_{\tilde{\mathrm{P}}}^{xyt}$ and $\lambda_{\tilde{\mathrm{P}}}^{xy,t}$ with respect to all measures.Table \ref{tab:dyn_denoising_results} quantitatively summarizes the results . In Figure \ref{fig:dynamic_denoising:method_comparison}, we compare two samples from the test dataset with a static camera view in the first row and a dynamic camera view in the third row. The vertical red lines in Figure \ref{fig:dynamic_denoising:method_comparison} indicate the $x$-location of the $yt$-excerpt shown to the left of each image. The second and the fourth row show the pointwise absolute errors of the respective images. For both samples, the lowest error is achieved by the $\boldsymbol{\Lambda}_{\Theta}^{xy,t}$ parameter-map. The spatial and the temporal components of the obtained regularization parameter-maps $\boldsymbol{\Lambda}_{\Theta}^{xy,t}$  are visualized in Figure \ref{fig:dynamic_denoising:parameter_maps}. Here, the noisy samples, the results obtained with PDGH using  $\boldsymbol{\Lambda}_{\Theta}$, the spatially and temporally dependent $\boldsymbol{\Lambda}_{\Theta}^{xy,t}$ parameter-maps and the ground truth-images are depicted. By comparing the static and dynamic case, we see that the trained CNN is able  to differentiate between the two inherently different cases. Thereby, for the video sample with the static camera position, where the background remains constant over time and only objects are changing position, the CNN imposes an overall higher temporal regularization. For the video sample where the camera position also changes over time, the CNN is able to predict the less prominent potential to exploit the temporal gradient-sparsity and thus assigns relatively low

\begin{table}[h]
\centering
\footnotesize\rm
\begin{tabular}{r|r|rcl|rcl|rcl}

&    & \multicolumn{3}{c}{\textbf{PDHG - $\lambda_{\tilde{\mathrm{P}}}^{xyt}$}}  & \multicolumn{3}{c}{\textbf{PDHG - $\lambda_{\tilde{\mathrm{P}}}^{xy,t}$}} & \multicolumn{3}{c}{\textbf{PDHG - $\boldsymbol{\Lambda}_{\Theta}^{xy,t}$}}  \\[5pt]
\hline
                &  \textbf{SSIM}    & $0.941$ & $\pm$ & $0.017$ &           $0.934$ & $\pm$ & $0.021$ &     $\textbf{0.968}$ & $\pm$ & $0.011$    \\
$\sigma=0.1$    &  \textbf{PSNR}    & $33.26$ & $\pm$ & $1.53$ &            $34.65$ & $\pm$ & $1.75$ &      $\textbf{39.30}$ & $\pm$ & $2.06$     \\
                &  \textbf{NRMSE}   & $0.044$ & $\pm$ & $0.009$ &           $0.038$ & $\pm$ & $0.011$ &     $\textbf{0.022}$ & $\pm$ & $0.007$    \\
\hline
                &  \textbf{SSIM}    & $0.837$ & $\pm$ & $0.075$ &           $0.914$ & $\pm$ & $0.021$ &     $\textbf{0.940}$ & $\pm$ & $0.021$  \\
$\sigma=0.2$    &  \textbf{PSNR}    & $31.32$ & $\pm$ & $2.25$ &            $33.52$ & $\pm$ & $1.63$ &      $\textbf{35.537}$ & $\pm$ & $2.24$  \\
                &  \textbf{NRMSE}   & $0.056$ & $\pm$ & $0.014$ &           $0.043$ & $\pm$ & $0.011$ &     $\textbf{0.035}$ & $\pm$ & $0.011$  \\
\hline
                &  \textbf{SSIM}    & $0.649$ & $\pm$ & $0.105$ &           $0.814$ & $\pm$ & $0.067$ &     $\textbf{0.915}$ & $\pm$ & $0.028$  \\
$\sigma=0.3$    &  \textbf{PSNR}    & $28.09$ & $\pm$ & $2.54$ &            $31.04$ & $\pm$ & $2.06$ &      $\textbf{33.36}$ & $\pm$ & $2.33$   \\
                &  \textbf{NRMSE}   & $0.082$ & $\pm$ & $0.024$ &           $0.058$ & $\pm$ & $0.016$ &     $\textbf{0.045}$ & $\pm$ & $0.014$  \\
\hline
\end{tabular}
\caption{Mean and standard deviation of the measures SSIM, PSNR and NRMSE and Blur obtained over the test set for the dynamic image denoising example. The TV-reconstruction using the proposed spatio-temporal parameter-maps $\boldsymbol{\Lambda}_{\Theta}^{xy,t}$ improves the results especially in terms of SSIM and PSNR.}\label{tab:dyn_denoising_results}
\end{table}

\subsection{Low-Dose Computerized Tomography}\label{sec:CT}
In the last section, we show an application of our proposed method to a static 2D low-dose CT reconstruction problem. Because of the different noise statistics and as a result, a different fidelity term, see \eqref{eq_CTdatafidelity} below, the problem requires the use of a reconstruction algorithm different than PDHG which shows that our proposed method can be used in conjunction with any iterative scheme. We mention however that since this fidelity term is not strongly convex, the consistency results of Section
\ref{sec:consistency_ana_scheme} cannot be applied in this case. We leave the corresponding  extension of these results to the CT case for future work.

\subsubsection{Problem Formulation}

We consider the proposed NN for the low-dose Computerized Tomography (CT) setting. Here a two-dimensional parallel beam geometry is chosen and the corresponding ray transform is given by the Radon transform \cite{Radon86}. As forward operator, we then consider the discretized Radon transformation, which is a finite-dimensional linear map $\Ad\colon \R^n \to \R^m$, where $n$ is the dimension of the image space and $m$ is the product between the number of angles of the measurement and the number of the equidistant detector bins. Then we can formulate the inverse problem as 
\begin{align*}
\ZZ= \Ad \XX + \mathbf{e},~\text{where} ~ \mathbf{e} = -\Ad \XX - \log (\mathbf{\tilde{N}_1} / N_0) ~ \text{and} ~ \mathbf{\tilde{N}_1} \sim \text{Pois}(N_0 \exp(-\Ad \XX \mu)),
\end{align*}
where $\mu$ is a normalization constant and $N_0$ denotes the mean photon count per detector bin without attenuation. Note that here we do not have Gaussian noise, but some noise which follows the negative log-transformation of a Poisson distribution. Therefore, the data-discrepancy in \eqref{eq:tv_min_problem} is not the L2-error, and the correct term can be derived from a Bayesian viewpoint, where the data-discrepancy corresponds to the negative log-likelihood $p_{Y|X=x}$. Using that the negative log-likelihood of a Poisson distributed random variable is given by the Kullback-Leibler divergence, the resulting data-discrepancy can be written as
\begin{align} \label{eq_CTdatafidelity}
d(\Ad \XX ,\ZZ) = \sum_{i=1}^m e^{-(\Ad \XX )_i \mu} N_0 - e^{-\ZZ_i \mu} N_0 \big(-(\Ad \XX)_i \mu + \log(N_0) \big),
\end{align}
see e.g. \cite{ADHHMS2022,Leuschner21} for more details. Consequently, we cannot use Algorithm~\ref{algo:tv_reco_algo_mri} for reconstruction and a reformulation of Algorithm~\ref{algo:pdhg_algo} for this data-discrepancy does not lead to a closed form. This can be seen by the proximal operator of $f^*$ in line 3 of Algorithm~\ref{algo:pdhg_algo}. In our case the convex functional $f$ would be given by $f(\XX) = d(\XX,\ZZ)$, i.e.,
\begin{align} \label{eq:f_data_fidelity}
f(\XX) = \sum_{i=1}^m e^{-\XX_i \mu} N_0 - \tilde{\ZZ} N_0 \big(-\XX_i \mu + \log(N_0) \big),
\end{align}
where we set $\tilde{\ZZ} = e^{-\ZZ_i \mu}$ for simplicity.
Then the convex conjugate is given by
\begin{align} \label{eq:dual_general}
f^*(\PP) = \max_\XX \sum_{i=1}^m \XX_i \PP_i - f(\XX) = \max_\XX \sum_{i=1}^m \XX_i \PP_i - e^{-\XX_i \mu} N_0 - \tilde{\ZZ} N_0 \XX_i \mu,
\end{align}
where we used in the second equality that $\tilde{\ZZ} \log(N_0)$ is independent of $\XX$. Differentiation with respect to $\XX$ shows that for a maximizer $\hat{\XX}$ it holds
\begin{align*}
\hat{\XX}_i = -\frac{1}{\mu} \log \big(\ZZ- \frac{\ZZ_i}{\mu N_0} \big).
\end{align*} 
Inserting this in \eqref{eq:dual_general} yields the convex conjugate of $f$
\begin{align*}
f^*(\PP) &= \sum_{i=1}^m -\frac{1}{\mu} \log \big(\ZZ- \frac{\PP_i}{\mu N_0} \big) \PP_i - \big(\ZZ- \frac{\PP_i}{\mu N_0} \big) + \ZZ\log \big(\ZZ- \frac{\PP_i}{\mu N_0} \big) \\ 
&=\sum_{i=1}^m \big(\ZZ- \frac{\PP_i}{\mu} \big) \log \big(\ZZ- \frac{\PP_i}{\mu N_0} \big)  -\big(\ZZ- \frac{\PP_i}{\mu N_0} \big).
\end{align*} 
Then for this $f^*$ the proximal operator does not have a simple closed form.

As a remedy we consider the primal-dual algorithm PD3O \cite{Yan2018}, which is a generalization of the PDHG method. The PD3O aims to minimize the sum of proper, lower semi-continuous and convex functions
\begin{align*}
    \min_\XX f(\mathbf{K}\XX) + g(\XX) + h(\XX),
\end{align*}
where $\mathbf{K}\colon V^n \to V^{\tilde{m}}$ is a bounded linear operator, $h$ is differentiable with a Lipschitz continuous gradient and for both $g$ and $f^*$ the proximal operator has a analytical solution. The general scheme of PD3O is described in Algorithm~\ref{algo:general_pd3o_algo}

\begin{algorithm}[h]
\caption{Unrolled PD3O algorithm (adapted from \cite{Yan2018})}\label{algo:general_pd3o_algo}
  \begin{algorithmic}[1]
  \INPUT $L = \text{Lip}(\nabla h)$,\; $\tau = 2/L$, \; $\sigma = 1/(\tau \Vert \mathbf{K} \mathbf{K}^\trans\| )$, \; $\text{initial guess}~ \bar{\XX}_0$
  \OUTPUT reconstructed image $\XX_{\mathrm{TV}}$
  \STATE $\PP_0 = \mathbf{\bar{x}_0}$
  \STATE $\QQ_0 = \mathbf{0}$
    \FOR {$k < T$ }
    \STATE $\QQ_{k+1} = \mathrm{prox}_{\sigma f^{*}}(\QQ_k + \sigma \mathbf{K}\bar{x}_k)$
    \STATE $\PP_{k+1} = \mathrm{prox}_{\tau g}(\PP_k - \tau \nabla h(\PP_k) - \tau \mathbf{K}^\trans \QQ_{k+1})$
    \STATE $\bar{\XX}_{k+1} = 2\PP_{k+1} - \PP_k + \tau\nabla h(\PP_k) - \tau \nabla h(\PP_{k+1})$ 
    \ENDFOR
    \STATE $\XX_{\mathrm{TV}} = \XX_T$
  \end{algorithmic}
\end{algorithm}

Note that we recover the PDHG algorithm if we set $h = 0$. For application of PD3O to our CT case we define
\begin{align*}
    f(\QQ) &= \Vert \mathbf{\Lambda} \, \QQ \Vert_1, \\
    g(\PP) &= \iota_{\{\PP \ge 0 \}} (\PP) = \begin{cases}
    0\quad &\text{if} ~\PP \ge 0, \\
    + \infty & \text{else},
    \end{cases} \\
    h(\PP) &= \sum_{i=1}^m e^{-\PP_i \mu} N_0 - e^{-\ZZ_i \mu} N_0 \big(-\PP_i \mu + \log(N_0) \big), \\
    \mathbf{K} &= \nabla.
\end{align*}
The proximal operator of $f^*$ is already given in \eqref{eq:clipping_op}, the proximal operator of $g$ is given by
\begin{align*}
    \text{prox}_{\tau g}(\ZZ) = \text{ReLU}(\ZZ) = \begin{cases}
    \ZZ \quad &\text{if}~ \ZZ \ge 0, \\
    0 &\text{else},
    \end{cases}
\end{align*}
and $\nabla h$ is given by
\begin{align*}
    \nabla h(\PP) = \mu N_0 \Ad^\trans \big( - e^{-\Ad \PP \mu} + e^{- \ZZ\mu} \big).
\end{align*}
Note that $\nabla h$ is not globally Lipschitz continuous, but due to the non-negativity constraint $g$ we only have to consider $\nabla h$ for $\PP$ with non-negative entries. Consequently, we can find an upper bound of the Lipschitz constant of $\nabla h$ by $\text{Lip}(\nabla h) \le \Vert \Ad \Vert^2 \mu^2 N_0$.
The resulting scheme we use for CT reconstruction is summarized in Algorithm~\ref{algo:tv_pd3o_algo}.
\begin{algorithm}[h]
        \caption{Unrolled PD3O algorithm  for general linear operator $\Ad$ with $h(\, \cdot \,) = d(\Ad \cdot \, , \ZZ)$, $d$ as in \eqref{eq_CTdatafidelity} and \textit{fixed} regularization parameter-map $\boldsymbol{\Lambda}$ (adapted from \cite{Yan2018})}\label{algo:tv_pd3o_algo}
  \begin{algorithmic}[1]
  \INPUT $L = \text{Lip}(\nabla h)$, \; $\tau = 2/L$, \; $\sigma = 1/(\tau\Vert \nabla \Vert )$, \; $\text{initial guess}~ \mathbf{\bar{x}_0}$
   \OUTPUT reconstructed image $\XX_{\mathrm{TV}}$ 
  \STATE $\PP_0 = \mathbf{\bar{x}_0}$
  \STATE $\QQ_0 = \mathbf{0}$
    \FOR {$k < T$ }
   \STATE $\QQ_{k+1}  = \mathrm{clip}_{\boldsymbol{\Lambda}}(\QQ_k + \sigma \nabla \bar{\XX}_k)$  
    \STATE $\PP_{k+1} = \text{ReLU} ( \PP_k - \tau \mu N_0 \Ad^\trans ( e^{-\ZZ\mu} - e^{-\Ad\PP_k \mu}) - \tau \nabla^\trans \QQ_{k+1})$ 
     \STATE $\bar{\XX}_{k+1} = 2 \PP_{k+1} - \PP_k + \tau \mu N_0 \Ad^\trans(e^{-\Ad\PP_{k+1} \mu} - e^{-\Ad\PP_k \mu})$ 
    \ENDFOR
   \STATE $\XX_{\mathrm{TV}} = \bar{\XX}_T$
  \end{algorithmic}
\end{algorithm}

\subsubsection{Experimental Set-Up}

We use the LoDoPaB dataset \cite{LoDoPaB21}\footnote{available at \url{https://zenodo.org/record/3384092##.Ylglz3VBwgM}} for low-dose CT imaging. It is based on scans of the Lung Image Database Consortium and Image Database Resource Initiative \cite{Armato11} which serve as ground truth images, while the measurements are simulated. The dataset contains 35820 training images, 3522 validation images and 3553 test images. Here the ground truth images have a resolution of $362\times 362$ on a  domain of $26\text{cm}\times26\text{cm}$. We only use the first 300 training images and the first 10 validation images.
For the forward operator we consider a normalization constant $\mu = 81.35858$, the mean photon count per detector bin $N_0 = 4096$ as well as 513 equidistant detector bins and 1000 equidistant angles between 0 and $\pi$. A detailed description of the data generation process is given in \cite{LoDoPaB21}.
Following the naming convention of Figure~\ref{fig:u_Theta.pdf}, the network $u_\theta$ is a 2D U-Net\footnote{available at \url{https://jleuschn.github.io/docs.dival/_modules/dival/reconstructors/networks/unet.html}}, where the number of channels at different scales are 32, 32, 64, 64 and 128 resulting in 610673 trainable parameters. For training we use Adam optimizer \cite{kingma2014adam} with a learning rate of $10^{-4}$, a batch size of 1 and train for 50 epochs. Then we used the model configuration for which the MSE on the validation set was lowest. The number of iterations of PD3O is set to $T=512$ resulting in a training time of around 24 hours on a single NVIDIA GeForce RTX 2080 super GPU with 8 GB GPU memory. At test time, we use $T=1024$ iterations for reconstruction. The forward and the adjoint operator as well as the FBP were implemented using the publicly available library \texttt{ODL}\cite{adler2017operator}.

\subsubsection{Results}

\begin{figure}
\centering
\begin{subfigure}[t]{.24\textwidth}
\caption*{PD3O  $\lambda_{\tilde{\mathrm{P}}}^{xy}$}
\begin{tikzpicture}[spy using outlines={rectangle,white,magnification=3,size=2.3cm, connect spies}]
\node[anchor=south west,inner sep=0]  at (0,0) {\includegraphics[width=\linewidth]{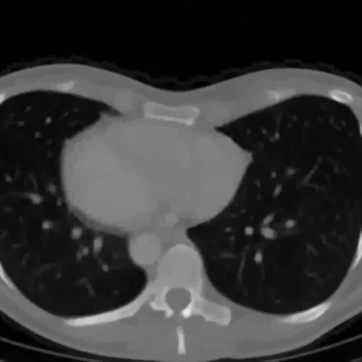}};
\spy on (1.8,.8) in node [left] at (2.33,2.49);
\end{tikzpicture}
\end{subfigure}%
\hfill
\begin{subfigure}[t]{.24\textwidth}
  \caption*{PD3O  $\lambda_{\mathrm{P}}^{xy}$} 
  \begin{tikzpicture}[spy using outlines={rectangle,white,magnification=3,size=2.3cm, connect spies}]
\node[anchor=south west,inner sep=0]  at (0,0) {  \includegraphics[width=\linewidth]{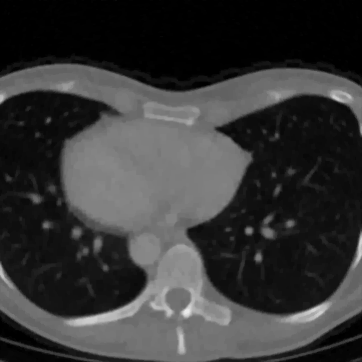}};
\spy on (1.8,.8) in node [left] at (2.33,2.49);
\end{tikzpicture}
\end{subfigure}%
\hfill
\begin{subfigure}[t]{.24\textwidth}
  \caption*{PD3O  $\boldsymbol{\Lambda}_{\Theta}^{xy}$}  
  \begin{tikzpicture}[spy using outlines={rectangle,white,magnification=3,size=2.3cm, connect spies}]
\node[anchor=south west,inner sep=0]  at (0,0) {  \includegraphics[width=\linewidth]{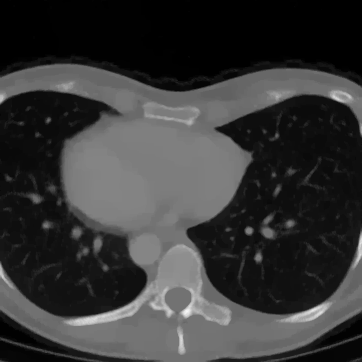}};
\spy on (1.8,.8) in node [left] at (2.33,2.49);
\end{tikzpicture}
\end{subfigure}%
\hfill
\begin{subfigure}[t]{.24\textwidth}
  \caption*{Target / FBP}  
  \begin{tikzpicture}[spy using outlines={rectangle,white,magnification=3,size=2.3cm, connect spies}]
\node[anchor=south west,inner sep=0]  at (0,0) {  \includegraphics[width=\linewidth]{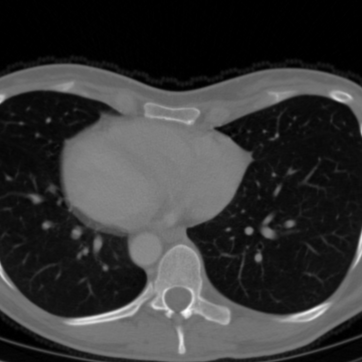}};
\spy on (1.8,.8) in node [left] at (2.33,2.49);
\end{tikzpicture}
\end{subfigure}%

\begin{subfigure}[t]{.24\textwidth}
  \begin{tikzpicture}[spy using outlines={rectangle,white,magnification=3,size=2.3cm, connect spies}]
\node[anchor=south west,inner sep=0]  at (0,0) {  \includegraphics[width=\linewidth]{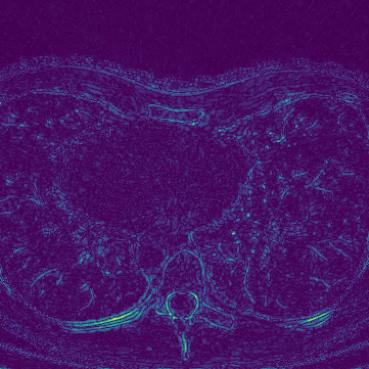}};
\spy on (1.8,.8) in node [left] at (2.33,2.49);
\end{tikzpicture}
\end{subfigure}%
\hfill
\begin{subfigure}[t]{.24\textwidth}
  \begin{tikzpicture}[spy using outlines={rectangle,white,magnification=3,size=2.3cm, connect spies}]
\node[anchor=south west,inner sep=0]  at (0,0) {  \includegraphics[width=\linewidth]{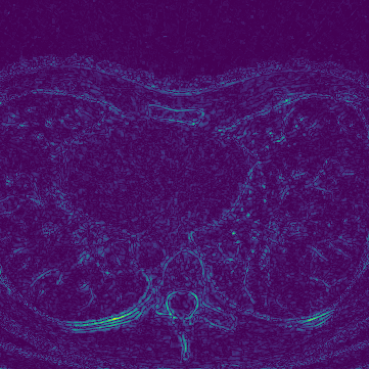}};
\spy on (1.8,.8) in node [left] at (2.33,2.49);
\end{tikzpicture}
\end{subfigure}%
\hfill
\begin{subfigure}[t]{.24\textwidth}
  \begin{tikzpicture}[spy using outlines={rectangle,white,magnification=3,size=2.3cm, connect spies}]
\node[anchor=south west,inner sep=0]  at (0,0) {  \includegraphics[width=\linewidth]{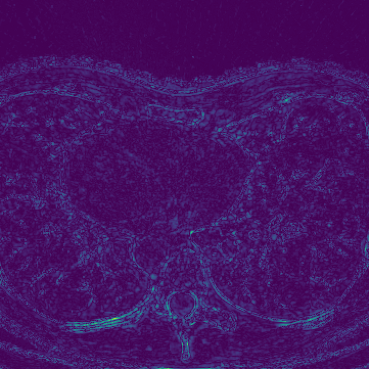}};
\spy on (1.8,.8) in node [left] at (2.33,2.49);
\end{tikzpicture}
\end{subfigure}%
\hfill
\begin{subfigure}[t]{.24\textwidth}
  \begin{tikzpicture}[spy using outlines={rectangle,white,magnification=3,size=2.3cm, connect spies}]
\node[anchor=south west,inner sep=0]  at (0,0) {  \includegraphics[width=\linewidth]{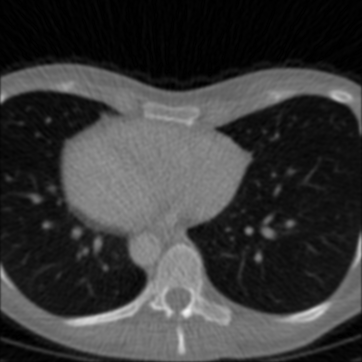}};
\spy on (1.8,.8) in node [left] at (2.33,2.49);
\end{tikzpicture}

\end{subfigure}%
\caption{Reconstruction of the ground truth CT image using different choices of regularization parameters. The last column shows the ground truth and the FBP reconstruction.
\textit{Top}: full image. \textit{Bottom}: difference to the ground truth.} \label{fig:CT_img}
\end{figure}

\begin{figure}[!h]
    \centering
    \includegraphics[width=0.325\linewidth]{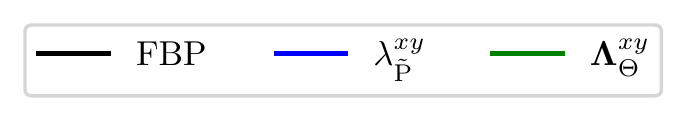}\\
    \includegraphics[width=0.24\linewidth]{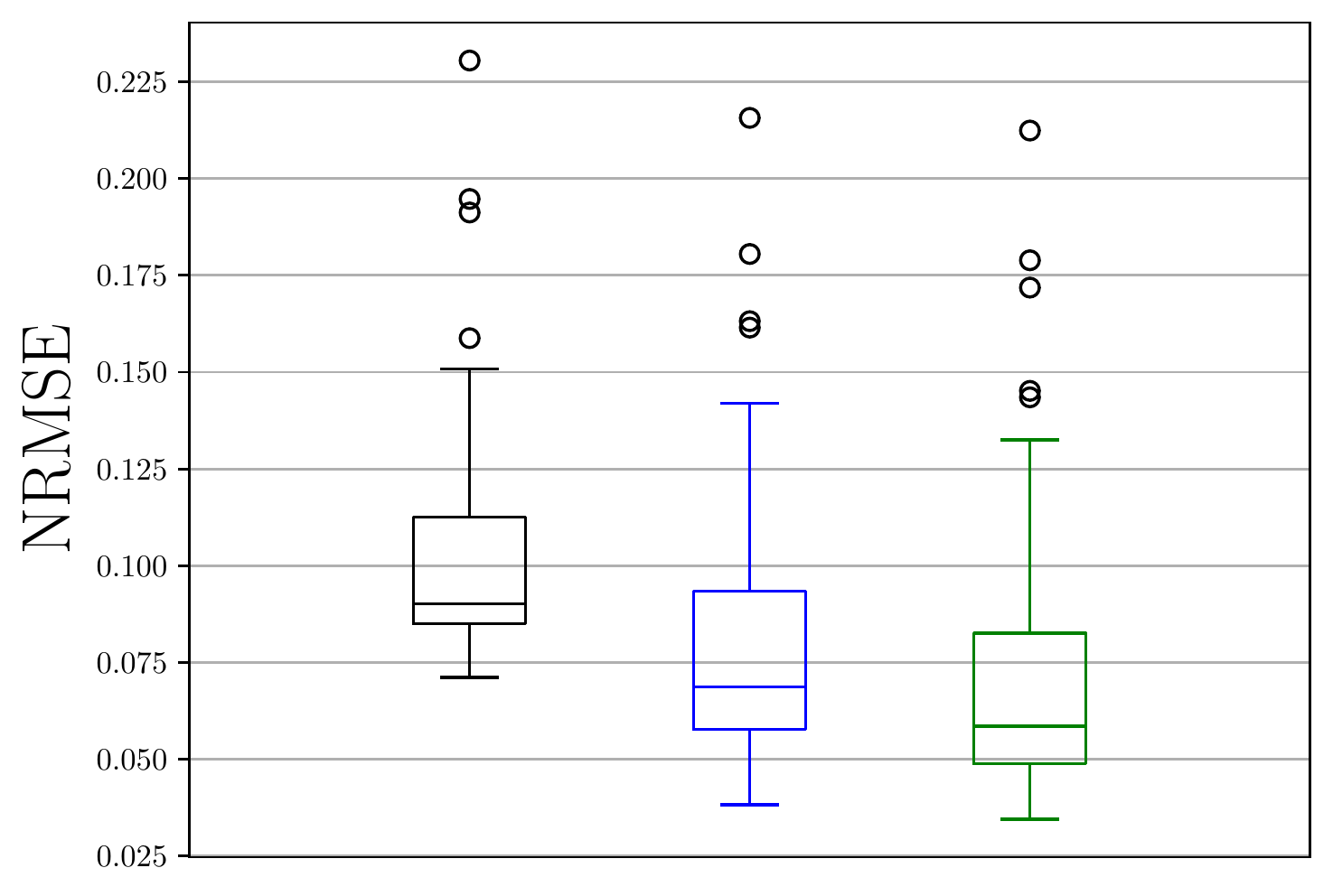}
    \includegraphics[width=0.24\linewidth]{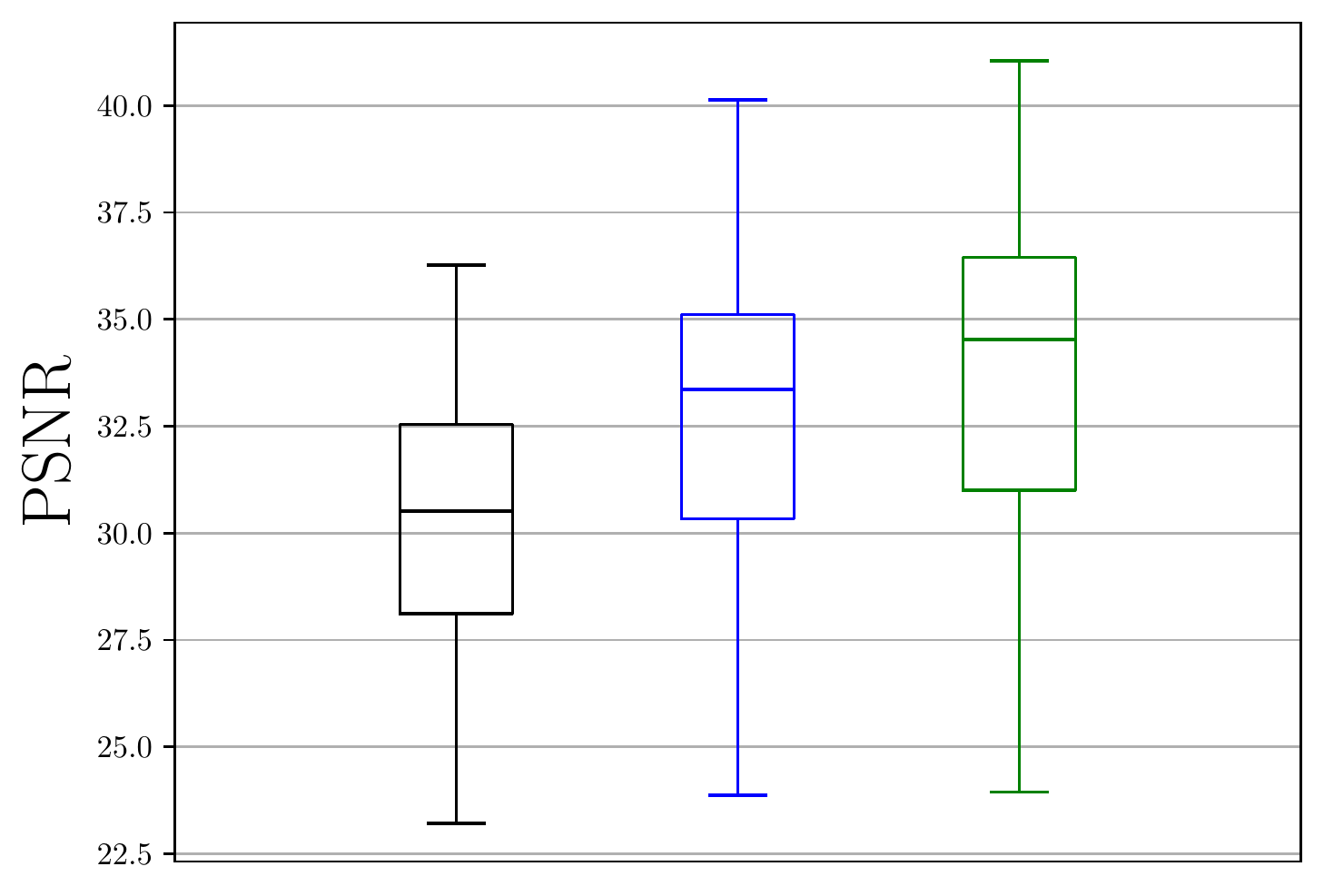}
    \includegraphics[width=0.24\linewidth]{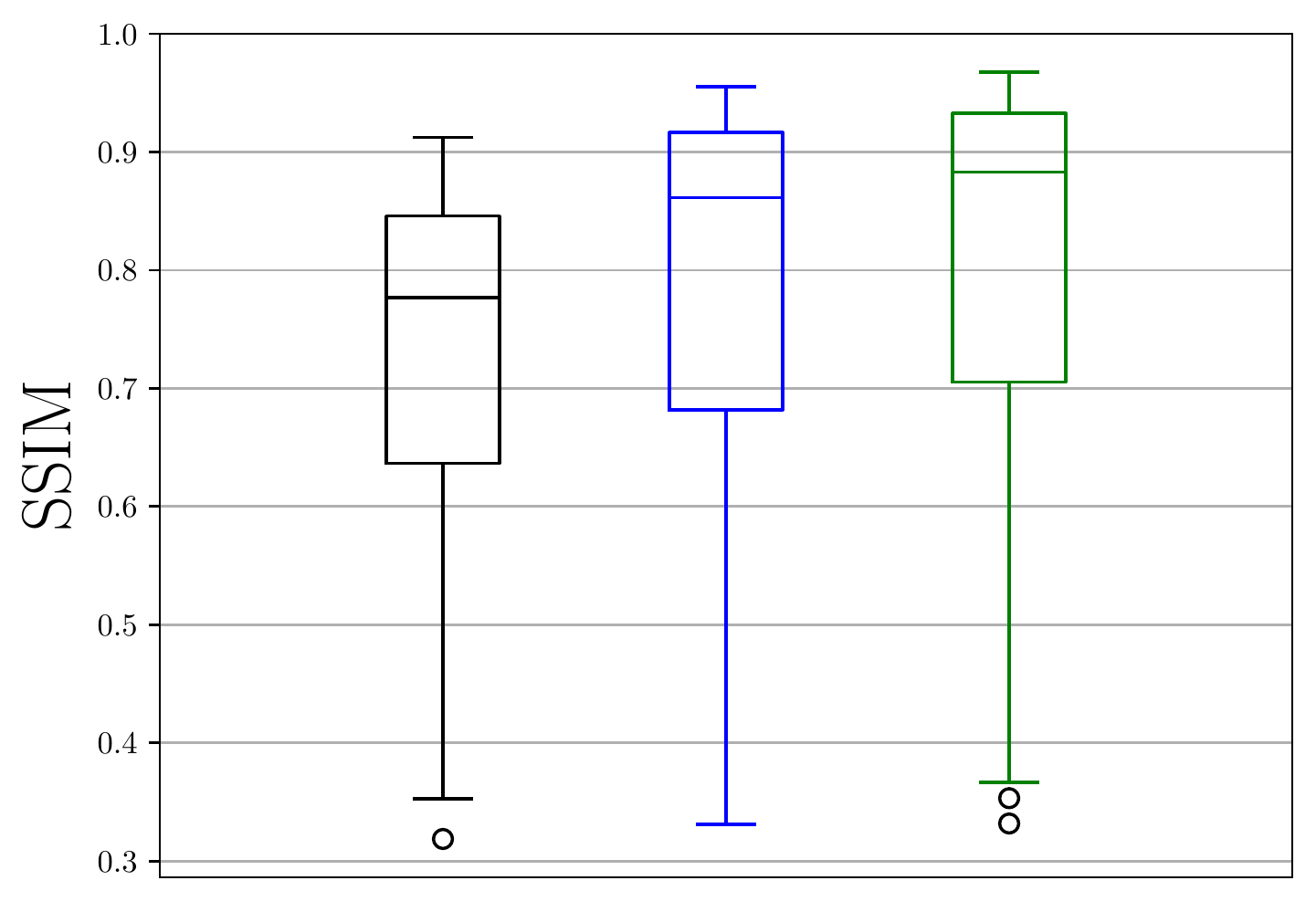}
    \includegraphics[width=0.24\linewidth]{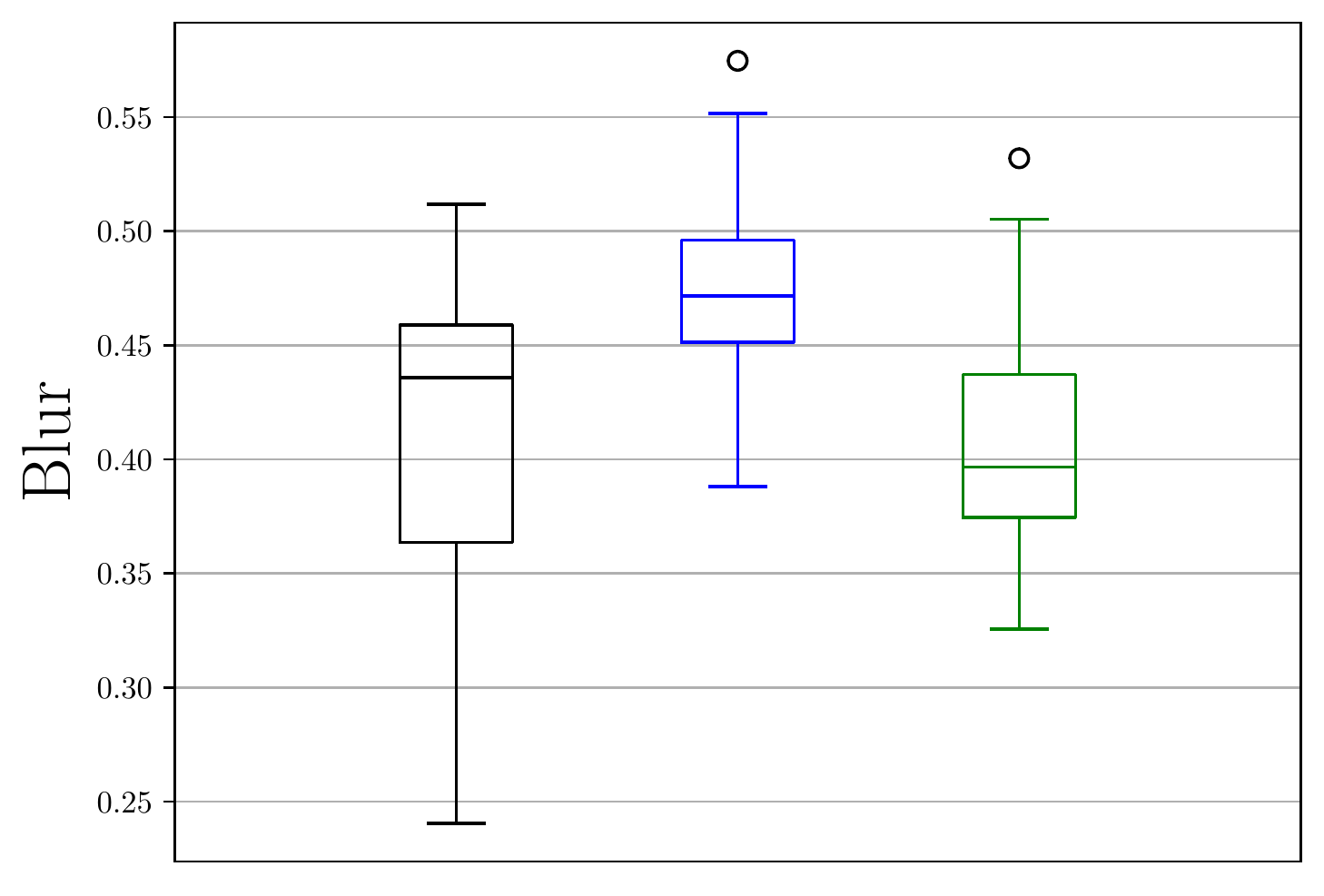}
    \caption{Box-plots summarizing the reconstruction results in terms of PSNR, NRMSE, SSIM and blur effect obtained with the PD3O algorithm for a CT reconstruction problem for different choices of the regularization parameter. Filtered back-projection (black), single scalar regularization parameter  ($\lambda_{\tilde{\mathrm{P}}}^{xy}$, blue) and the proposed parameter-map $\boldsymbol{\Lambda}_{\Theta}$ obtained with a CNN ($\mathrm{NET}_{\Theta}$, green).}
    \label{fig:box_plots_CT}
\end{figure}

In Figure~\ref{fig:CT_img} we compare the PD3O reconstructions (top) and their corresponding errors with respect to the ground truth (bottom) using different regularization parameter choices $\lambda_{\tilde{P}}^{xy}$, $\lambda_{P}^{xy}$ and $\boldsymbol\Lambda_{\Theta}$ for PD3O. Obviously, using the estimated parameter-map $\boldsymbol\Lambda_{\Theta}$ leads to a significant improvement of the reconstruction. In particular, sharp edges are retained, while using a constant regularizing parameter results in a significant blur. This can be also seen in Table~\ref{tab:ct_results}, where we compare the NRMSE, PSNR, SSIM and blur and  evaluated on the first 100 test images of the LoDoBaP dataset. These results are visualized in Figure~\ref{fig:box_plots_CT} using box-plots. Note that the FBP seems to better than PD3O-$\lambda_{\tilde{P}}^{xy}$ in terms of the blur effect, but this can be explained by the fact that FBP reconstructions admit a lot of high-frequency artefacts leading to a small blur effect.

\begin{figure}[!h]
\centering
\begin{subfigure}[t]{.24\textwidth}
\caption*{FBP}
  \begin{tikzpicture}[spy using outlines={rectangle,white,magnification=2,size=1.5cm, connect spies}]
\node[anchor=south west,inner sep=0]  at (0,0) { \includegraphics[width=\linewidth]{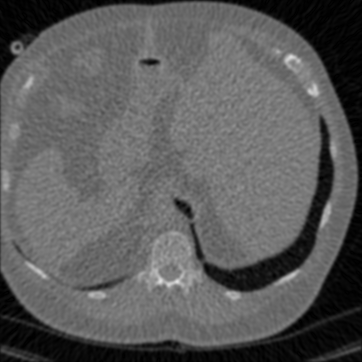}};
\spy on (1.75,1.) in node [left] at (1.52,2.9);
\end{tikzpicture}
\end{subfigure}%
\hfill
\begin{subfigure}[t]{.24\textwidth}
  \caption*{PD3O  $\boldsymbol{\Lambda}_{\Theta}^{xy}$} 
  \begin{tikzpicture}[spy using outlines={rectangle,white,magnification=2,size=1.5cm, connect spies}]
\node[anchor=south west,inner sep=0]  at (0,0) { \includegraphics[width=\linewidth]{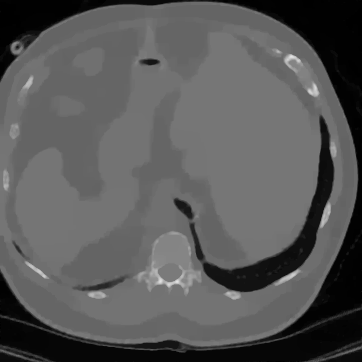}};
\spy on (1.75,1.) in node [left] at (1.52,2.9);
\end{tikzpicture}
\end{subfigure}%
\hfill
\begin{subfigure}[t]{.24\textwidth}
  \caption*{$\boldsymbol{\Lambda_{\Theta}}^{xy}$}  
  \begin{tikzpicture}[spy using outlines={rectangle,white,magnification=2,size=1.5cm, connect spies}]
\node[anchor=south west,inner sep=0]  at (0,0) {\includegraphics[width=\linewidth]{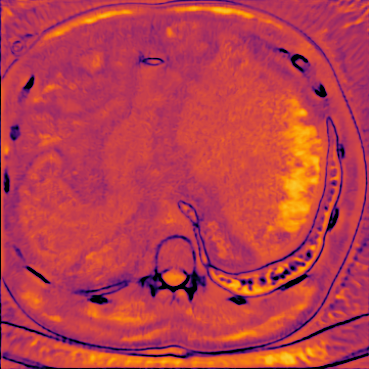}};
\spy on (1.75,1.) in node [left] at (1.52,2.9);
\end{tikzpicture}
\end{subfigure}%
\hfill
\begin{subfigure}[t]{.24\textwidth}
  \caption*{Target}  
  \begin{tikzpicture}[spy using outlines={rectangle,white,magnification=2,size=1.5cm, connect spies}]
\node[anchor=south west,inner sep=0]  at (0,0) {\includegraphics[width=\linewidth]{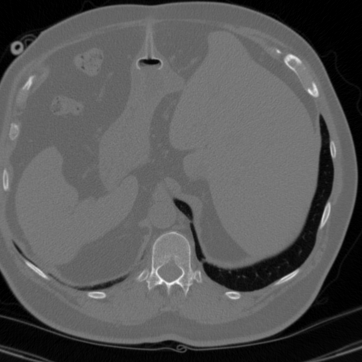}};
\spy on (1.75,1.) in node [left] at (1.52,2.9);
\end{tikzpicture}
\end{subfigure}%

\begin{subfigure}[t]{.24\textwidth}
\begin{tikzpicture}[spy using outlines={rectangle,white,magnification=2,size=1.5cm, connect spies}]
\node[anchor=south west,inner sep=0]  at (0,0) { \includegraphics[width=\linewidth]{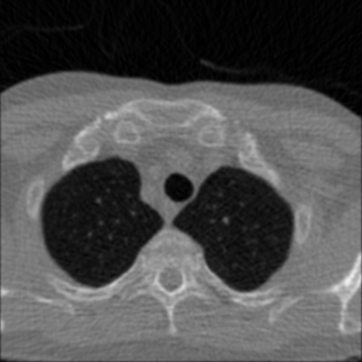}};
\spy on (1.75,.65) in node [left] at (1.52,2.9);
\end{tikzpicture}
\end{subfigure}%
\hfill
\begin{subfigure}[t]{.24\textwidth}
  \begin{tikzpicture}[spy using outlines={rectangle,white,magnification=2,size=1.5cm, connect spies}]
\node[anchor=south west,inner sep=0]  at (0,0) { \includegraphics[width=\linewidth]{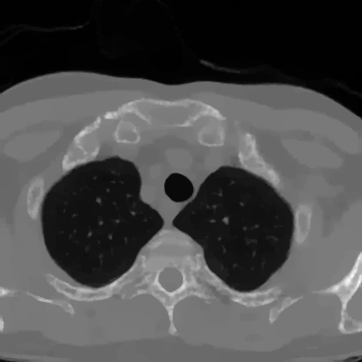}};
\spy on (1.75,.65) in node [left] at (1.52,2.9);
\end{tikzpicture}
\end{subfigure}%
\hfill
\begin{subfigure}[t]{.24\textwidth} 
  \begin{tikzpicture}[spy using outlines={rectangle,white,magnification=2,size=1.5cm, connect spies}]
\node[anchor=south west,inner sep=0]  at (0,0) {\includegraphics[width=\linewidth]{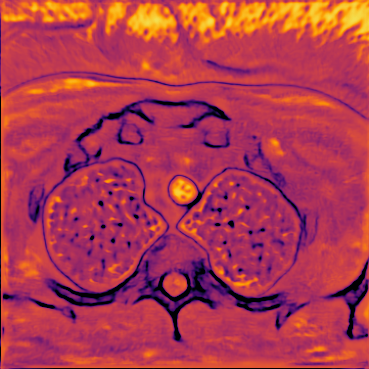}};
\spy on (1.75,.65) in node [left] at (1.52,2.9);
\end{tikzpicture}
\end{subfigure}%
\hfill
\begin{subfigure}[t]{.24\textwidth}
  \begin{tikzpicture}[spy using outlines={rectangle,white,magnification=2,size=1.5cm, connect spies}]
\node[anchor=south west,inner sep=0]  at (0,0) {\includegraphics[width=\linewidth]{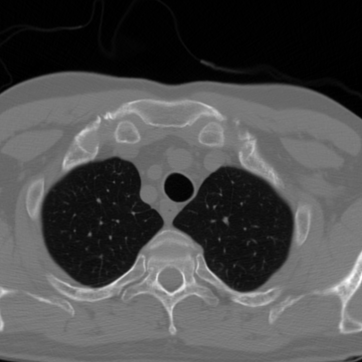}};
\spy on (1.75,.65) in node [left] at (1.52,2.9);
\end{tikzpicture}
\end{subfigure}%

\begin{subfigure}[t]{.24\textwidth}
\begin{tikzpicture}[spy using outlines={rectangle,white,magnification=1.8,size=1.5cm, connect spies}]
\node[anchor=south west,inner sep=0]  at (0,0) { \includegraphics[width=\linewidth]{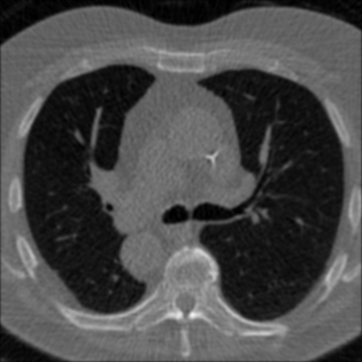}};
\spy on (1.9,.8) in node [left] at (1.52,2.9);
\end{tikzpicture}
\end{subfigure}%
\hfill
\begin{subfigure}[t]{.24\textwidth}
  \begin{tikzpicture}[spy using outlines={rectangle,white,magnification=1.8,size=1.5cm, connect spies}]
\node[anchor=south west,inner sep=0]  at (0,0) { \includegraphics[width=\linewidth]{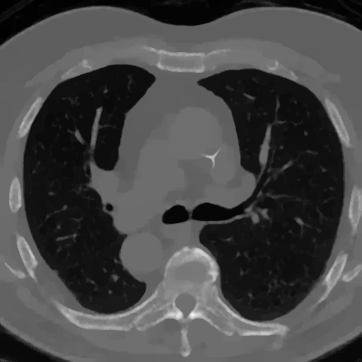}};
\spy on (1.9,.8) in node [left] at (1.52,2.9);
\end{tikzpicture}
\end{subfigure}%
\hfill
\begin{subfigure}[t]{.24\textwidth}
  \begin{tikzpicture}[spy using outlines={rectangle,white,magnification=1.8,size=1.5cm, connect spies}]
\node[anchor=south west,inner sep=0]  at (0,0) {\includegraphics[width=\linewidth]{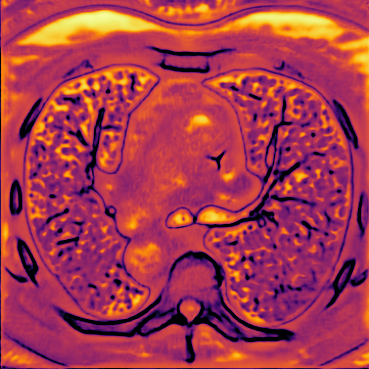}};
\spy on (1.9,.8) in node [left] at (1.52,2.9);
\end{tikzpicture}
\end{subfigure}%
\hfill
\begin{subfigure}[t]{.24\textwidth} 
  \begin{tikzpicture}[spy using outlines={rectangle,white,magnification=1.8,size=1.5cm, connect spies}]
\node[anchor=south west,inner sep=0]  at (0,0) {\includegraphics[width=\linewidth]{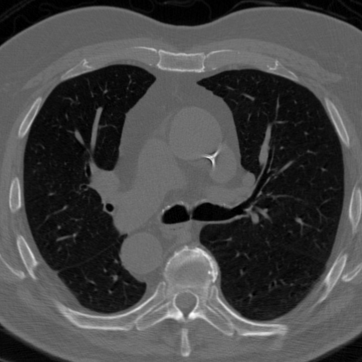}};
\spy on (1.9,.8) in node [left] at (1.52,2.9);
\end{tikzpicture}
\end{subfigure}%

\includegraphics[width=0.9\linewidth]{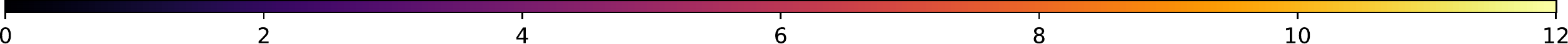}
\caption{Different reconstructions obtained with PD3O employing the regularization parameter-maps obtained with the  proposed CNN. From left to right: initial FBP-reconstruction, PD3O-reconstruction, spatial regularization parameter-map and ground truth image. As can be seen, the network attributes higher regularization parameters to image content with smooth structures while it yields lower regularization parameters at the edges to prevent smoothing.} \label{fig:CT_img_lambdamaps}
\end{figure}

Further PD3O-$\boldsymbol\Lambda_{\Theta}$ reconstructions with their corresponding estimated parameter-maps $\boldsymbol\Lambda_{\Theta}$ are shown in Figure~\ref{fig:CT_img_lambdamaps}. Note that the parameter-maps are given in a logarithmic scale. As expected, the regularization is strong in constant areas and less strong on edges or finer details in order to reduce a smoothing in these regions.

\begin{table}[h]
\centering
\footnotesize\rm
\begin{tabular}[t]{c|rcl|rcl|rcl} 
            & \multicolumn{3}{c}{\textbf{FBP}}                & \multicolumn{3}{c}{\textbf{PD3O - $\lambda_{\tilde{\mathrm{P}}}^{xy}$}}          & \multicolumn{3}{c}{\textbf{PD3O - $\LLambda_{\Theta}^{xy}$}} \\[5pt]
\hline
\textbf{PSNR}        & 30.37 &$\pm$& 2.95   & 32.87 &$\pm$& 3.59  & \textbf{33.90} &$\pm$& 3.94   \\ 
\textbf{SSIM}        & 0.739 &$\pm$& 0.141  & 0.796 &$\pm$& 0.152 & \textbf{0.809} &$\pm$& 0.157  \\
\textbf{NRMSE}        & 0.101 &$\pm$& 0.028  & 0.079 &$\pm$& 0.032 & \textbf{0.071} &$\pm$& 0.033  \\
\textbf{Blur Effect}    & 0.412 &$\pm$& 0.067  & 0.472 &$\pm$& 0.038 & \textbf{0.407} &$\pm$& 0.043  \\
\hline
\textbf{Training time} &  & - &  & & 5 h &         &  & 24 h      \\
\textbf{Runtime}     &  & 0.03 s   &           & &  5.08 s  &    &   & 5.08 s    
\end{tabular}
\caption{Mean and standard deviation of the measures SSIM, PSNR, NRMSE and Blur effect obtained over the CT test set. The best value ist marked in bold. The TV-reconstruction using the proposed parameter-maps $\boldsymbol{\Lambda}_{\Theta}^{xy}$ improve the results in every quality measure.}\label{tab:ct_results}
\end{table}

\newpage
\subsection{Choosing the Number of Iterations $T$}\label{subsec:choosing_T}
Since our proposed method to obtain the regularization parameter-map $\boldsymbol{\Lambda}_{\Theta}$ is based on unrolling an iterative algorithm as PDHG or PD3O using a fixed number of iterations $T$, questions about how to choose $T$ during training as well as at inference time are relevant. 
Recall from Section \ref{sec:consistency_ana_scheme} that $\XX_T := S^T(\ZZ, \LLambda)$ and $\XX^* := S^*(\ZZ, \LLambda)$ denote the $T$-th iterate and the exact solution of problem \eqref{eq:tv_min_problem_spatial}, respectively.
For addressing questions about the choice of $T$ at training and testing time, here, we emphasize the dependence of the solutions $\XX_T$ and $\XX^\ast$ on the set of parameters $\Theta$, by writing $\XX_T(\Theta) := S^T(\ZZ, \LLambda_{\Theta})$ and $\XX^{\ast}(\Theta) := S^\ast(\ZZ, \LLambda_{\Theta})$ and by denoting $\Theta_T$ as the set of trainable parameters which is obtained by training the  network which unrolls using $T$ iterations of PDHG or PD3O.

We illustrate the following considerations relying on results obtained for the dynamic cardiac MRI application shown in Subsection \ref{subsec:dyn_mri} but point out that these could be derived from the other applications examples as well.

\subsubsection{Choosing $T$ at Training Time}

\begin{figure}[h]
    \centering
    \includegraphics[width=0.8\linewidth]{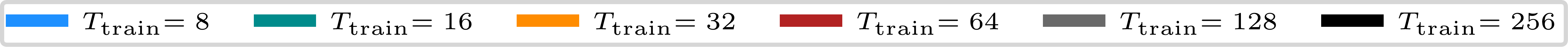}
    \includegraphics[width=0.6\linewidth]{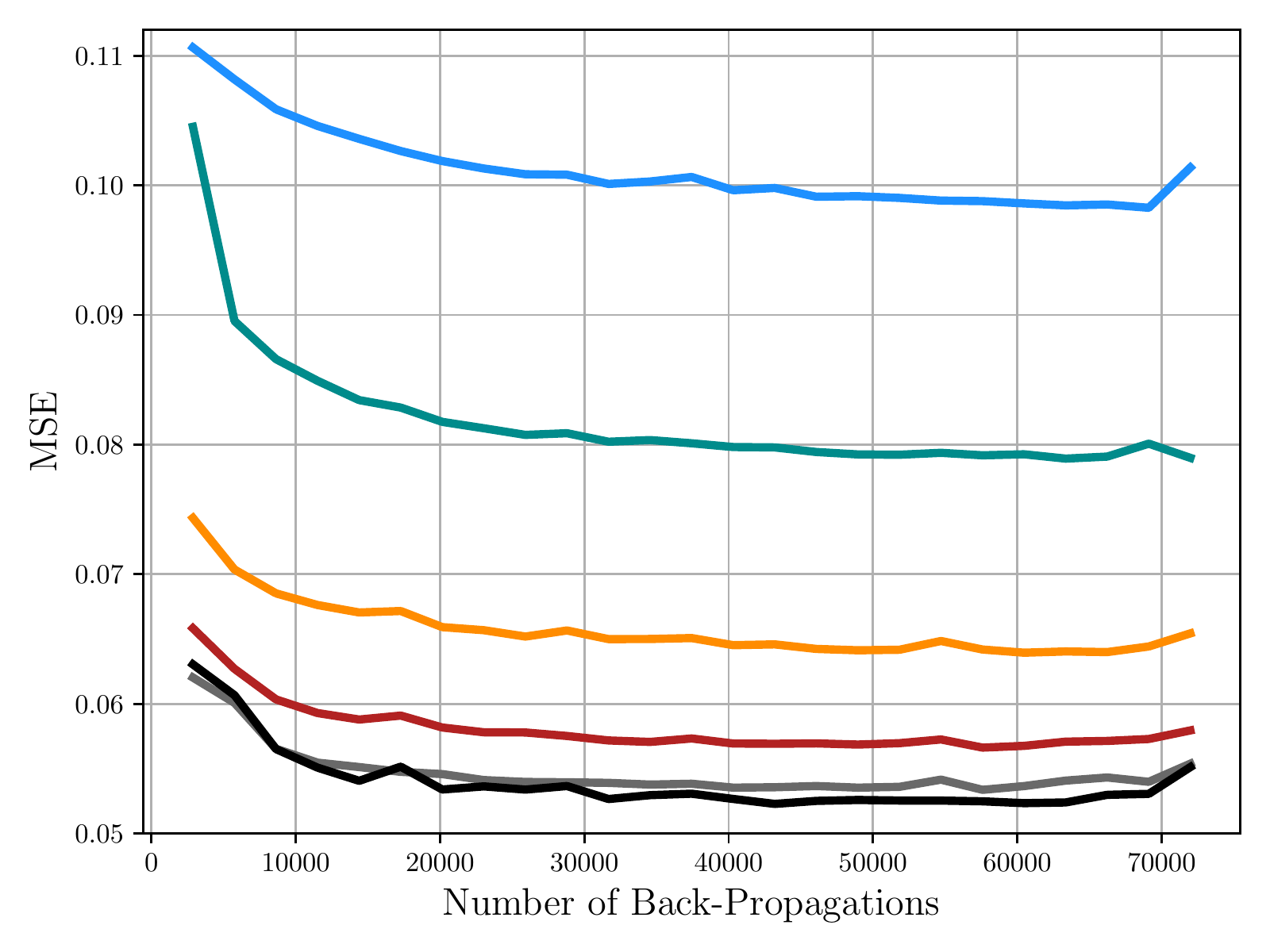}
    \caption{The validation error during training of the proposed method for the previously shown dynamic cardiac MRI example for different number of iterations $T$ which were used for PDHG. Using higher $T$ results in more accurate reconstructions.}
    \label{fig:T_training_curves}
\end{figure}

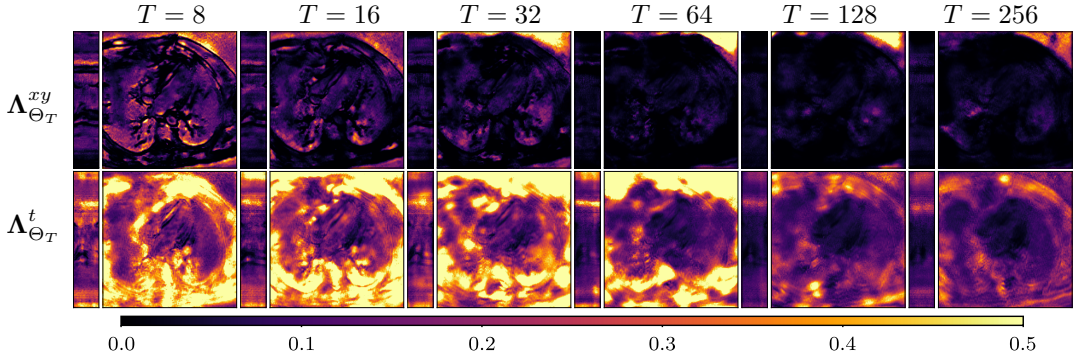
\begin{figure}[!h]
    \centering
    \begin{overpic}[width=0.87\linewidth]{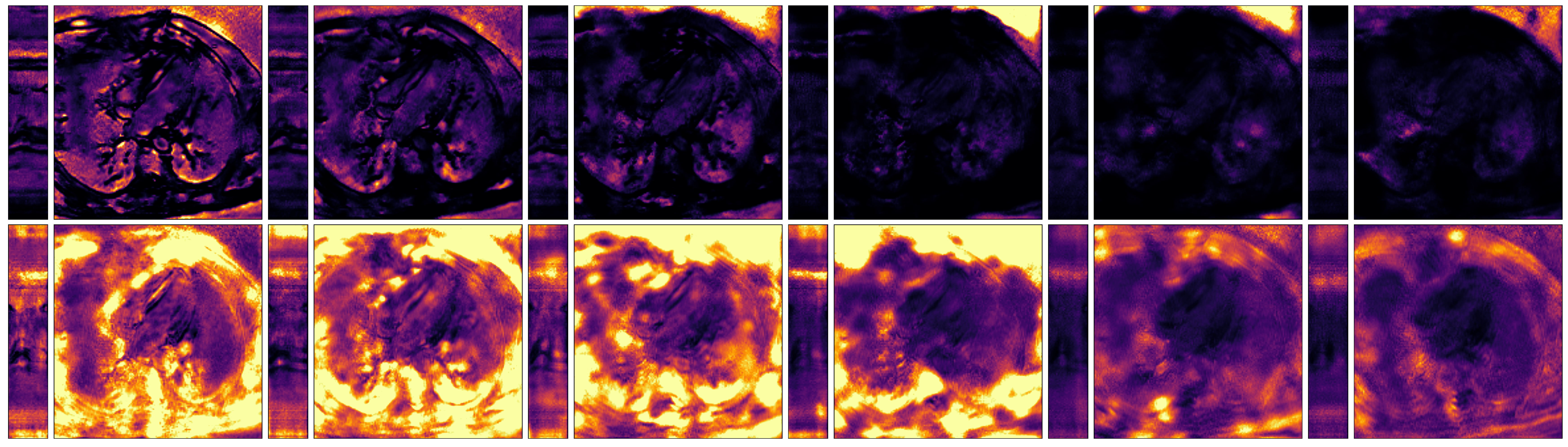}
    \put(-7,20){ {\textcolor{black}{\normalsize{\bf{$\boldsymbol{\Lambda}_{\Theta_T}^{xy}$ }}}}} 
    \put(-7,8){ {\textcolor{black}{\normalsize{\bf{$\boldsymbol{\Lambda}_{\Theta_T}^t$ }}}}} 
    \put(6,29){ {\textcolor{black}{\normalsize{\bf{$T=8$}}}}} 
    \put(22,29){ {\textcolor{black}{\normalsize{\bf{$T=16$}}}}}
    \put(38,29){ {\textcolor{black}{\normalsize{\bf{$T=32$}}}}}
    \put(55,29){ {\textcolor{black}{\normalsize{\bf{$T=64$}}}}}
    \put(70,29){ {\textcolor{black}{\normalsize{\bf{$T=128$}}}}}
    \put(86,29){ {\textcolor{black}{\normalsize{\bf{$T=256$}}}}}
    \end{overpic}
    \includegraphics[width=0.8\linewidth]{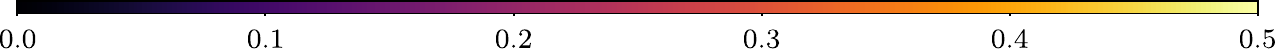}
    \caption{An example of different spatio-temporal regularization parameter-maps $\boldsymbol{\Lambda}_{\Theta_T}^{xy,t}$ for the dynamic cardiac MRI example obtained by training the unrolled network $\mathcal{N}_{\Theta}^T$ using different numbers of iterations $T$. All images are shown on the scale $[0,0.5]$. Although for different $T$ the  regularization parameter-maps have a similar structure - for example, the cardiac region is consistently regularized less strongly along time - for lower $T$, the values tend to be higher in general. Intuitively speaking, the network estimates higher regularization values in order to compensate for the lower $T$.}
    \label{fig:lambda_maps_Theta_T}
\end{figure}

Clearly, the obvious choice is to set $T$ as high as possible during training. The reason is to hope to be able to have $\XX_T(\Theta)\approx \XX^\ast(\Theta) $ and therefore to optimize $\Theta$ such that its optimal when the reconstruction algorithm given by $\mathcal{N}_{\Theta}^T$ is run until convergence. However, choosing a high $T$ increases training times as well as hardware requirements. Thus, one could on purpose choose or be forced to choose to use a lower $T$ for training and hope that the sub-network $\mathrm{NET}_{\Theta}$ is flexible enough to compensate for that.\\
Figure \ref{fig:T_training_curves} shows the validation error during training of an unrolled PDHG for the dynamic MRI example for different $T$. As can be seen, for smaller $T$, the NNs' ability to accurately reconstruct the images is clearly reduced. Further, Figure \ref{fig:lambda_maps_Theta_T} shows an example of different regularization parameter-maps which were obtained by training using a different number of iterations $T$. It shows that indeed, the obtained regularization parameter-maps vary depending on the number of iterations chosen for training, although they seem to share local features. Clearly, when $T$ is set too low, the network tends to yield higher regularization parameter-maps to try to compensate for the low number of iterations. However, from Figure \ref{fig:T_training_curves} one cannot infer whether  the limited reconstruction accuracy is attributable to the too low number of iterations, a resulting sub-optimal $\Theta_T$ or combinations thereof.

\subsubsection{Choosing $T$ at Test Time}

\begin{figure}
    \centering
    \resizebox{\linewidth}{!}{
    \begin{overpic}[width=0.95\linewidth]{figures/T_test_curves/legend_T8_T256.pdf}
    \put(14,-4){ {\textcolor{black}{\large{\bf{$R=4$}}}}} 
    \put(46,-4){ {\textcolor{black}{\large{\bf{$R=6$}}}}} 
    \put(80,-4){ {\textcolor{black}{\large{\bf{$R=8$}}}}} 
    
    \end{overpic} 
    }\\
    \vspace{0.7cm}
    \includegraphics[width=0.315\linewidth]{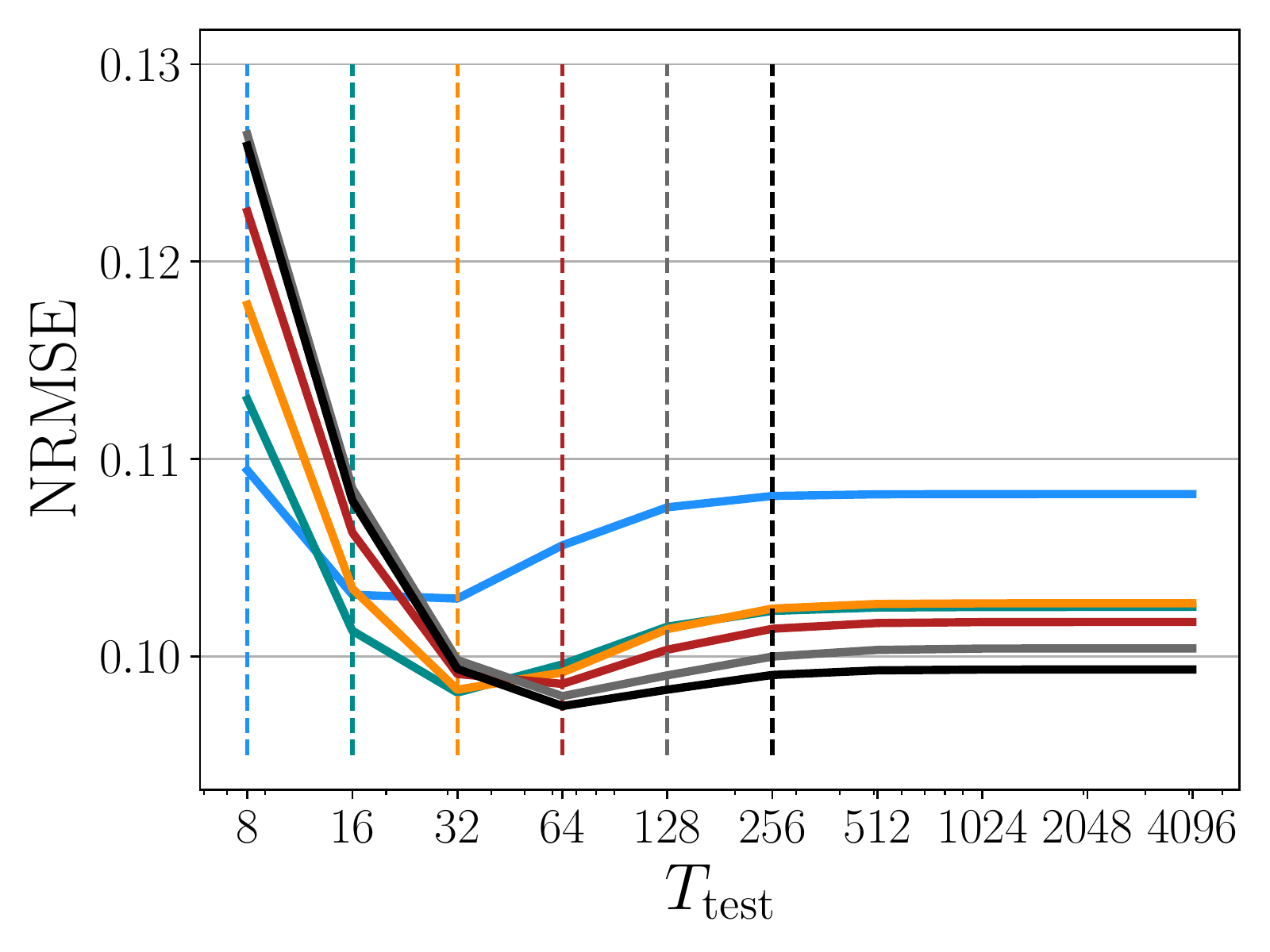}
    \includegraphics[width=0.315\linewidth]{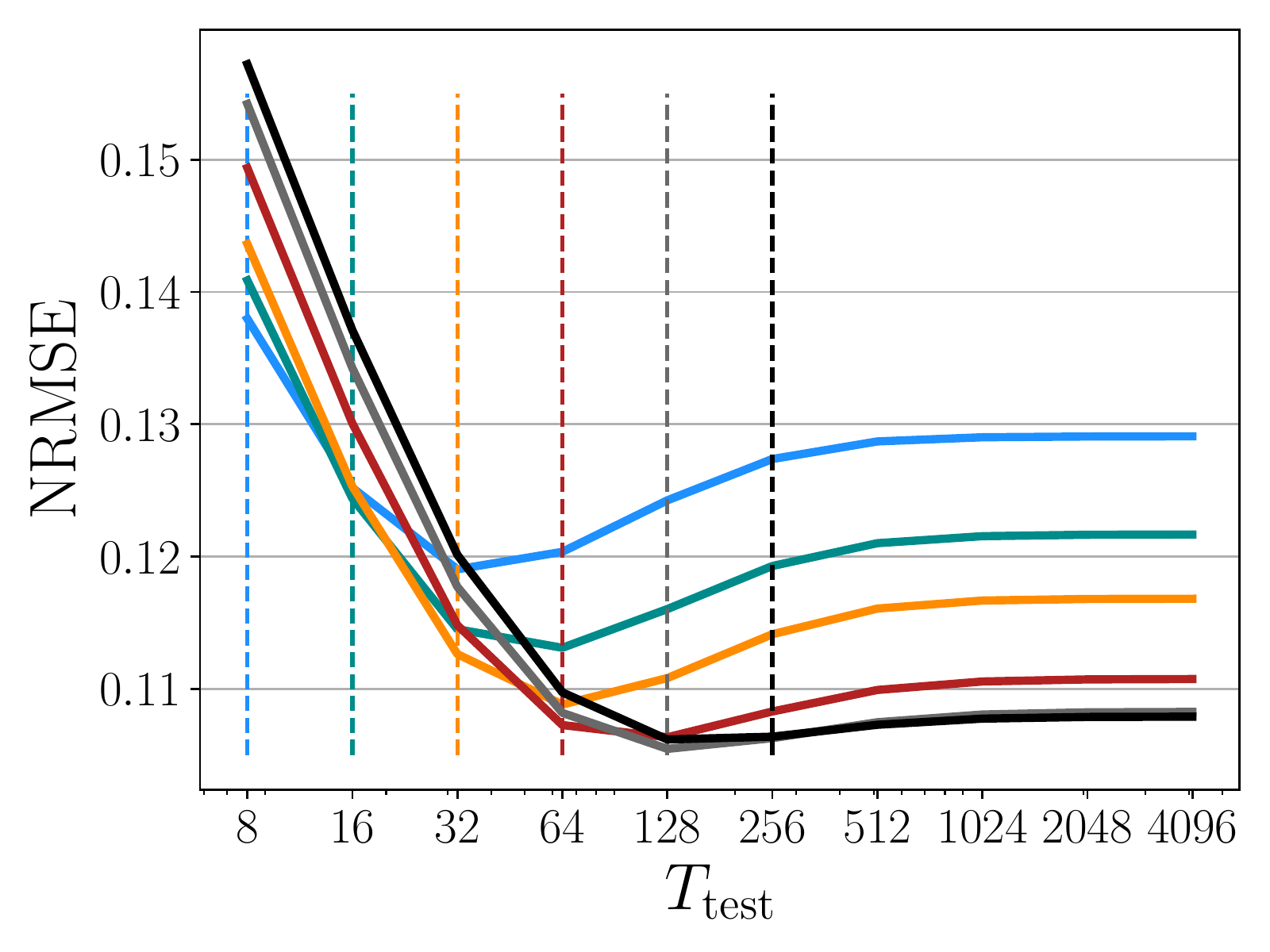}
    \includegraphics[width=0.315\linewidth]{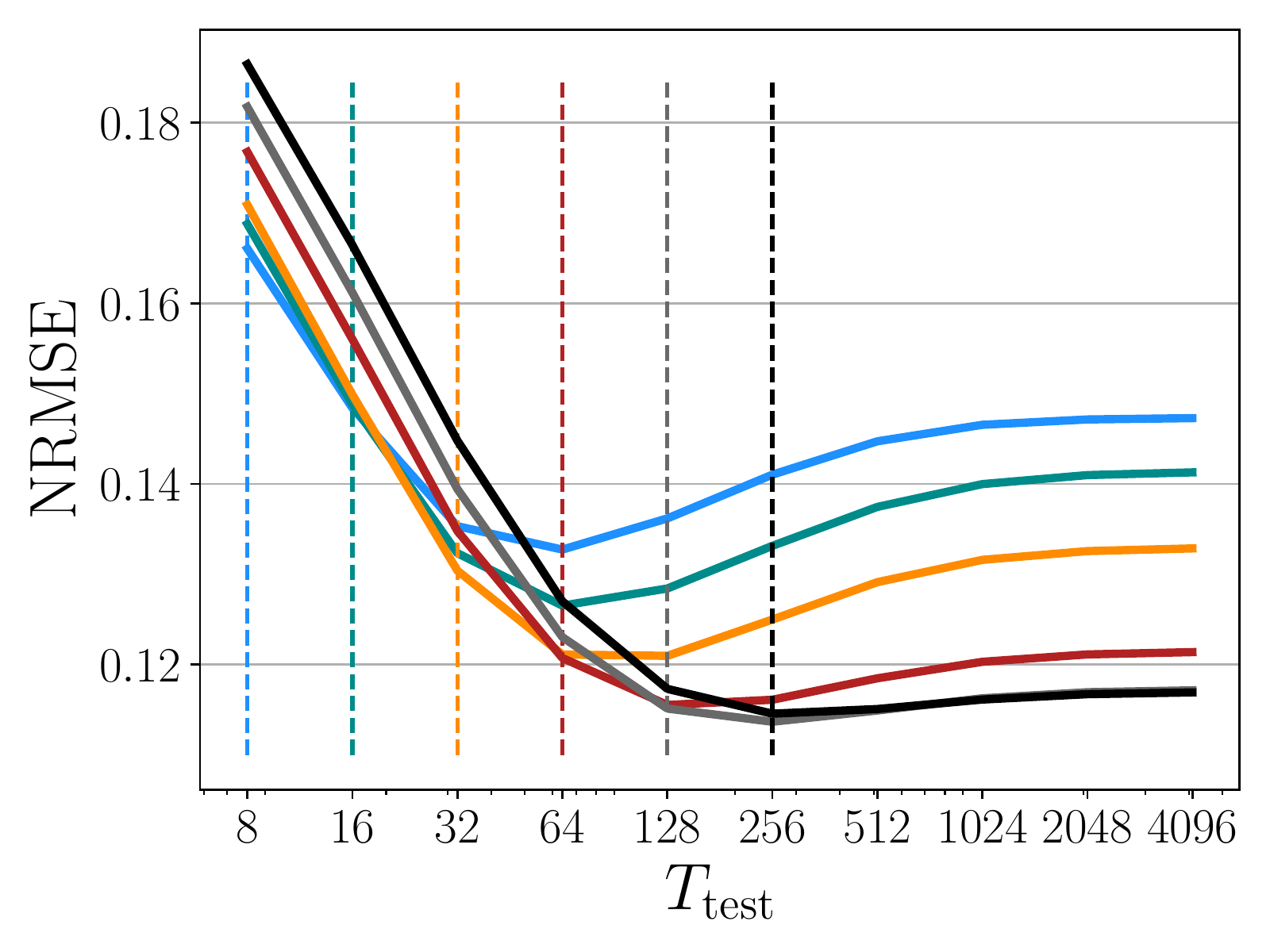}\\
    \includegraphics[width=0.315\linewidth]{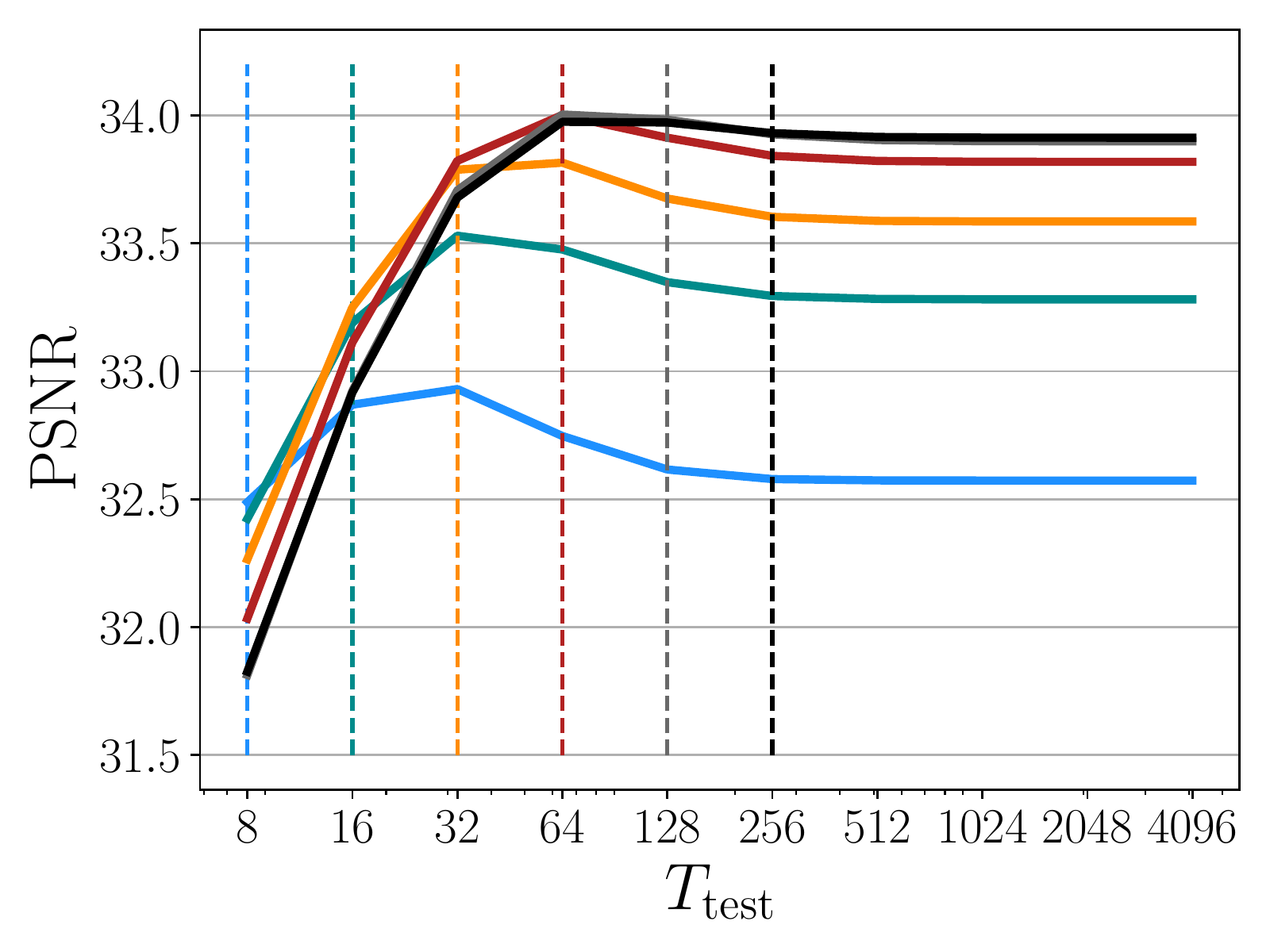}
    \includegraphics[width=0.315\linewidth]{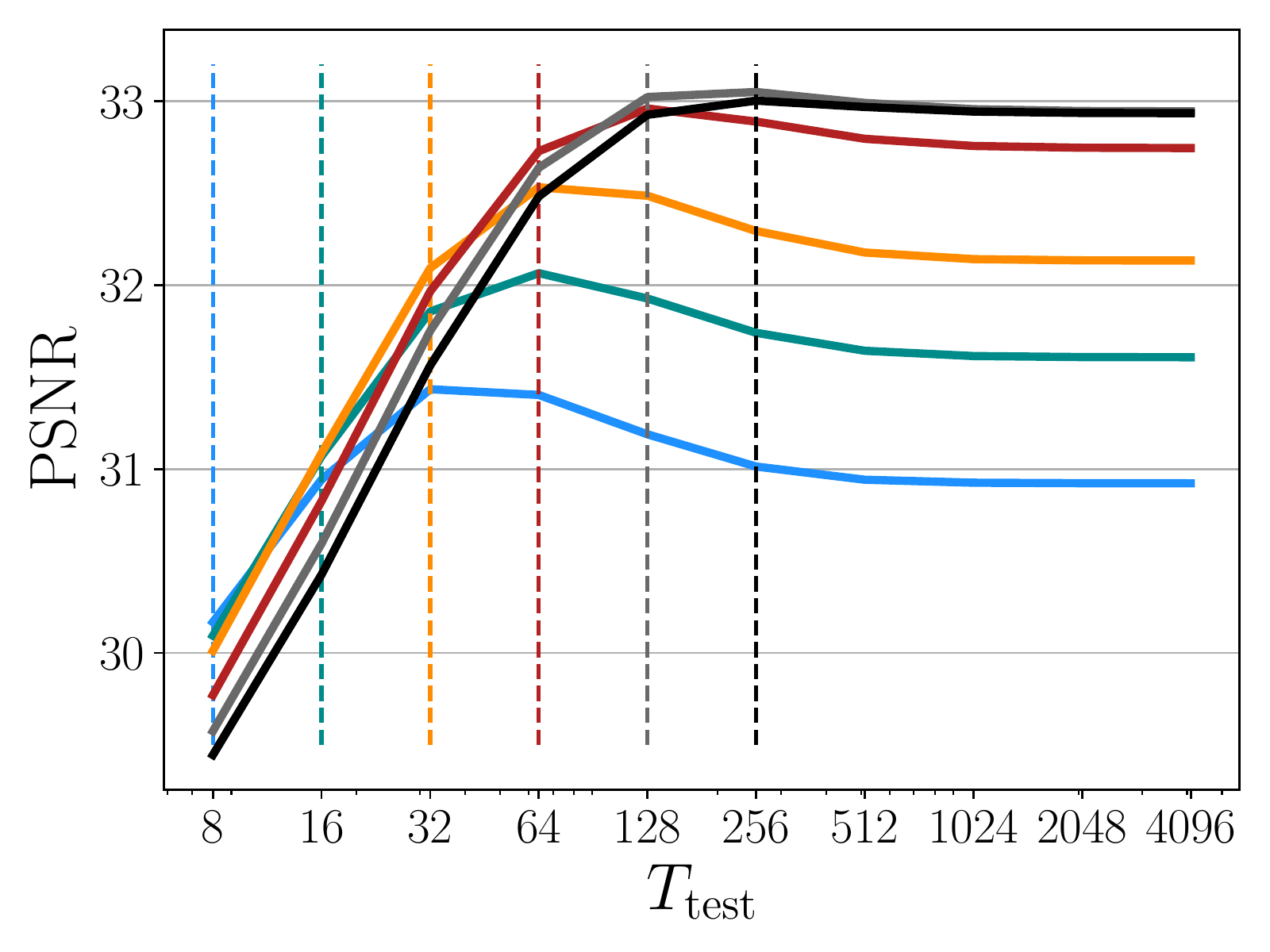}
    \includegraphics[width=0.315\linewidth]{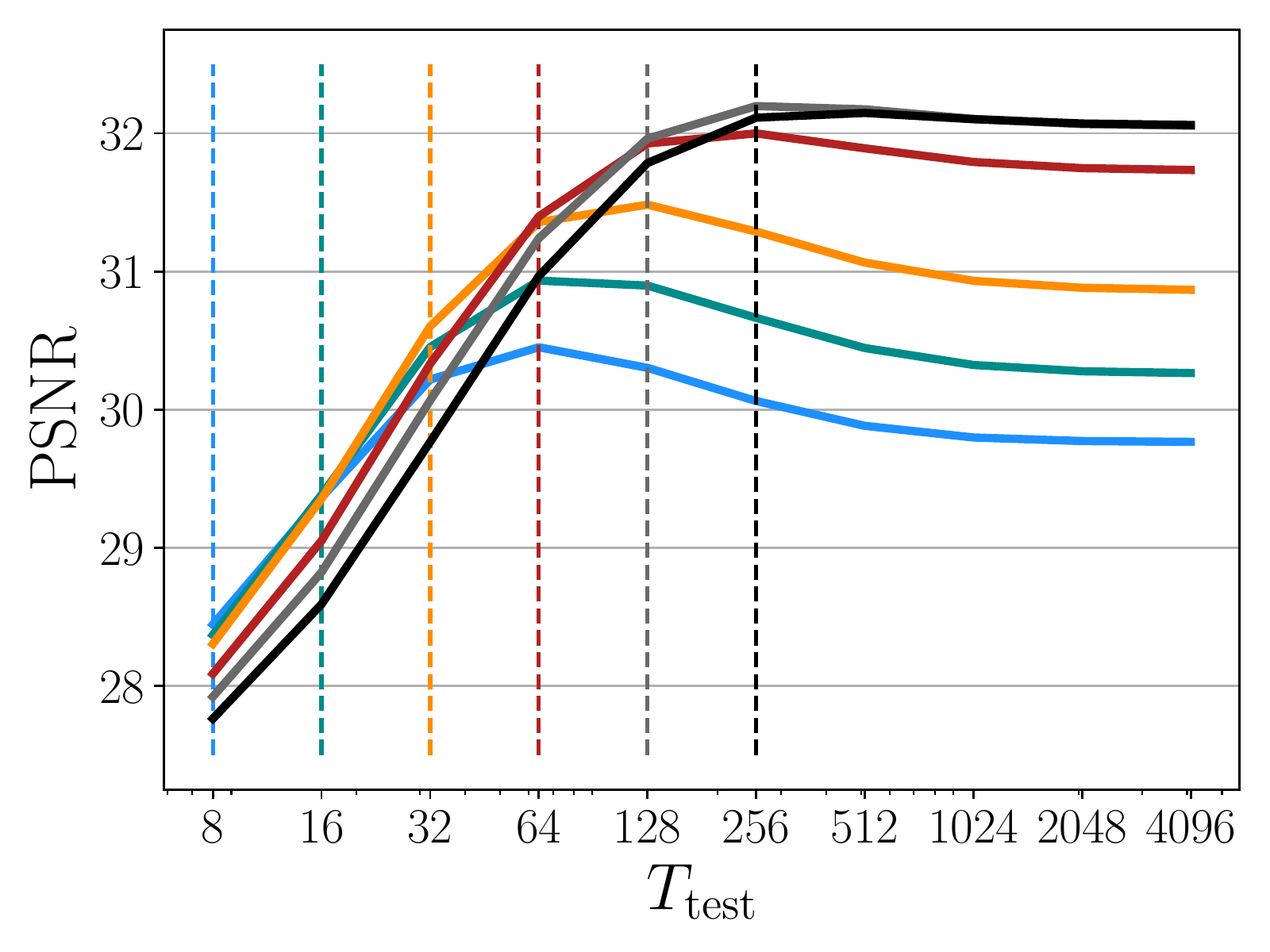}\\
    \includegraphics[width=0.315\linewidth]{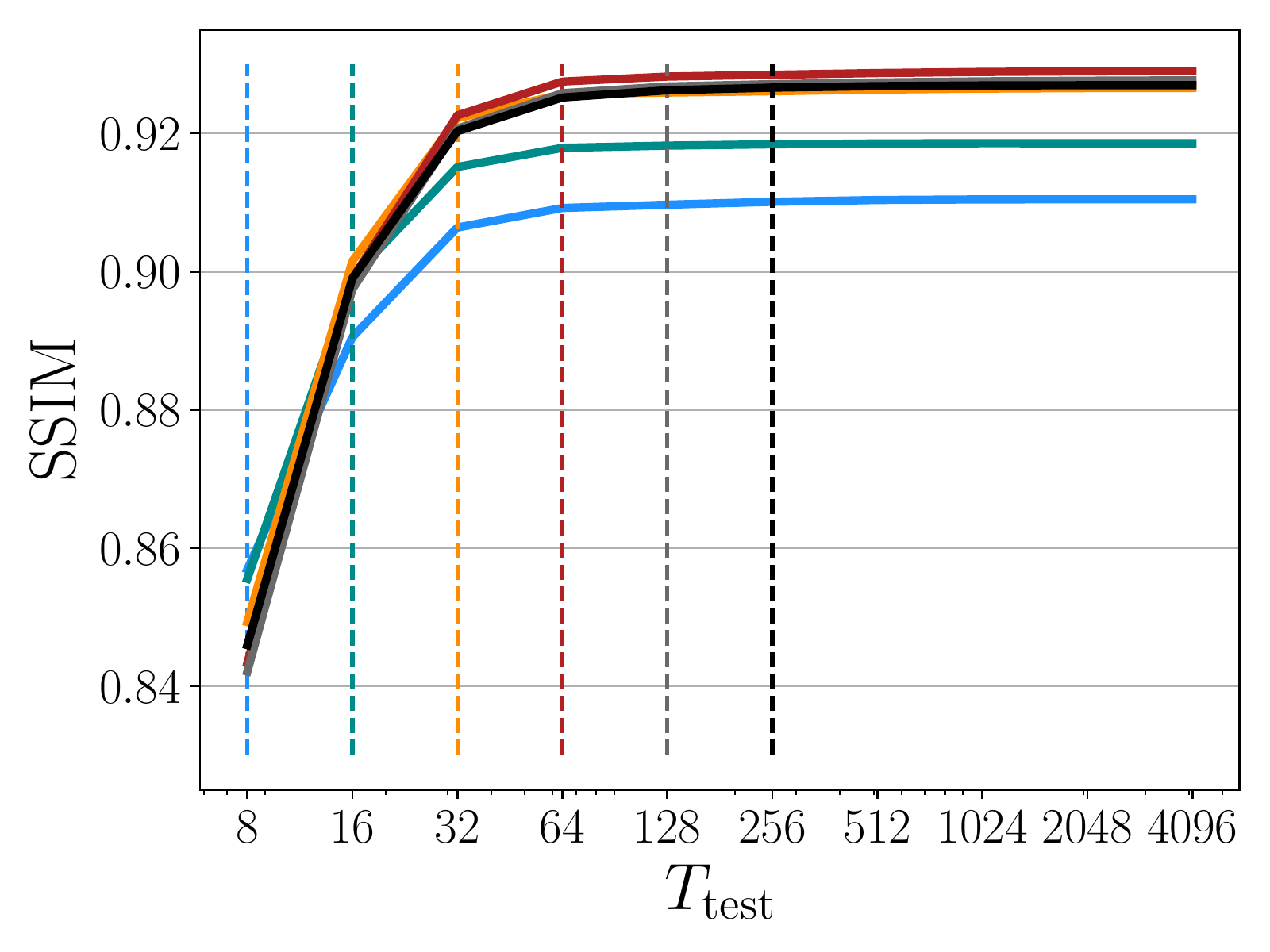}
    \includegraphics[width=0.315\linewidth]{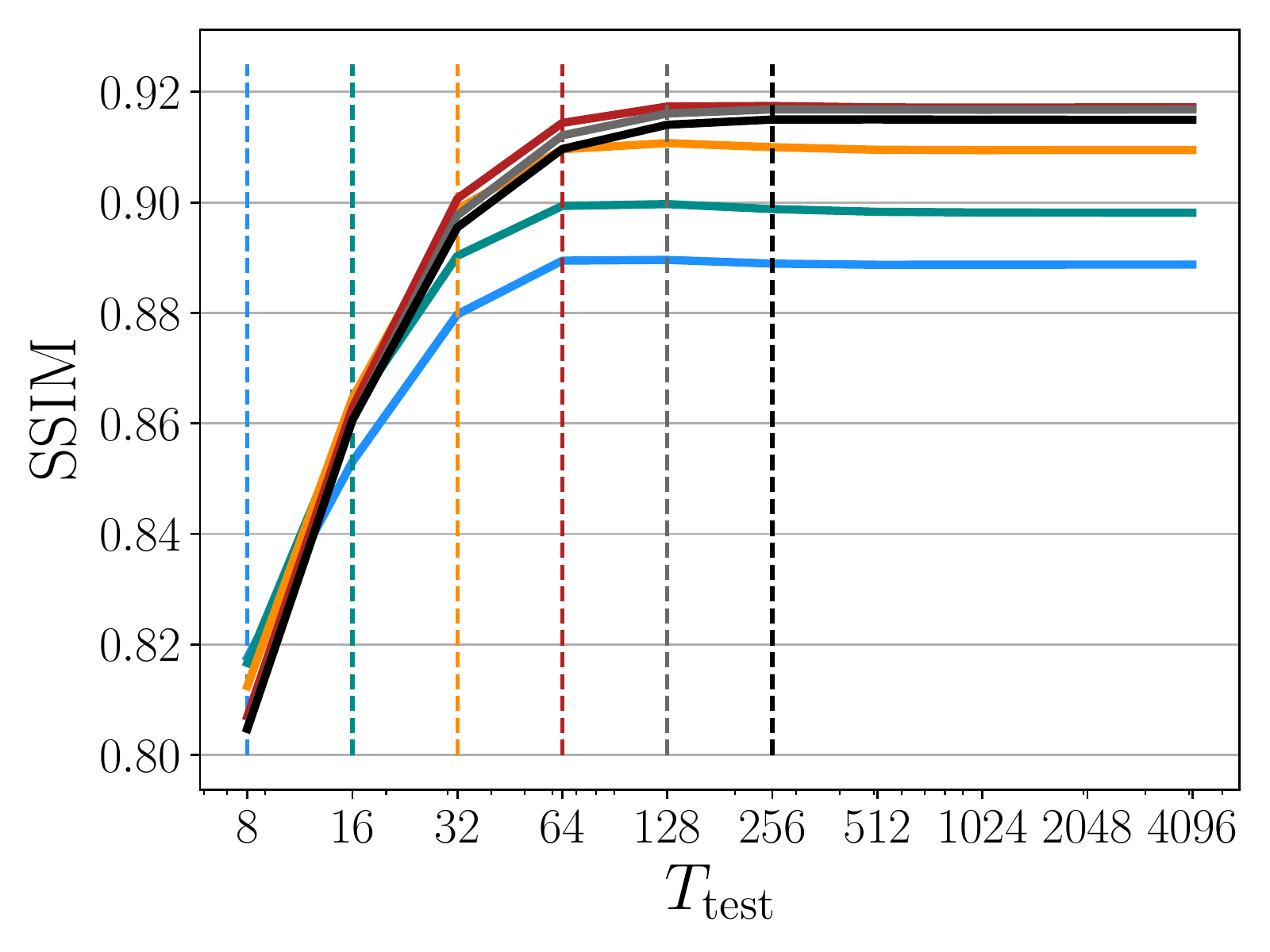}
    \includegraphics[width=0.315\linewidth]{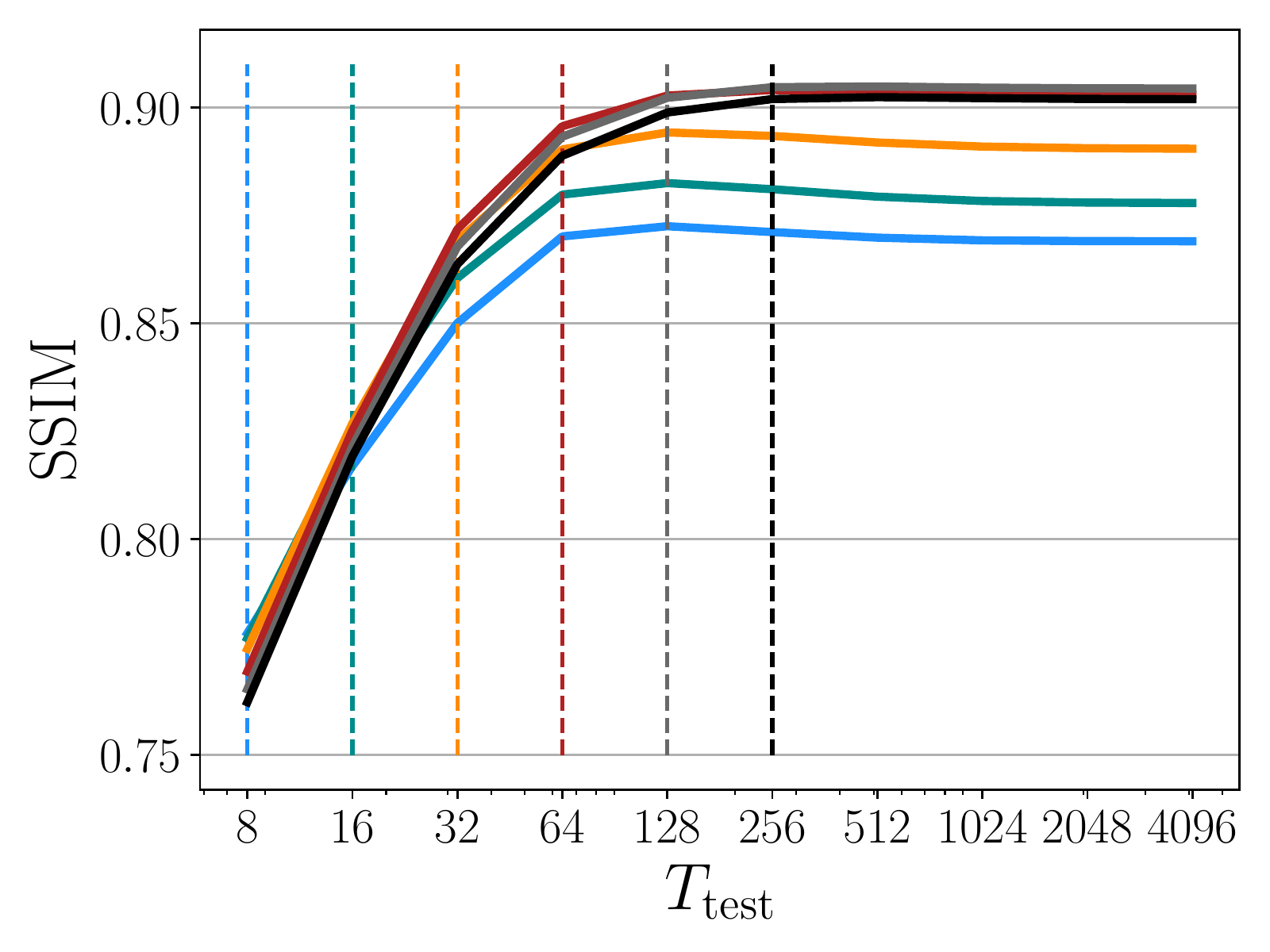}\\
   
    \caption{Average NRMSE, PSNR and SSIM for different PDHG-networks $\mathcal{N}_{{\Theta}_{T_{\mathrm{train}}}}^{T_{\mathrm{test}}}$ trained and tested with different combinations of  $T_{\mathrm{train}}$ and $T_{\mathrm{test}}$ shown for the accelerated dynamic cardiac MRI example for acceleration factors $R=4,6,8$. The color of the dashed lines encodes the number of iterations which were used for training, and consequently, the number of iterations used at test time for which one could expect the corresponding network to yield the best measure. We see that, however, this does not consistently hold, especially for lower $T_{\mathrm{train}}$, where choosing $T_{\mathrm{test}}>T_{\mathrm{train}}$ quite consistently improves the results up to some point with respect to all measures.} 
    \label{fig:T_train_vs_test_curves}
\end{figure}

Recall that at test time, the proposed method generates an input-dependent regularization parameter-map $\boldsymbol{\Lambda}_{\Theta_T}$ which is inherently dependent on the number of iterations $T$ the network was trained with. With $\boldsymbol{\Lambda}_{\Theta_T}$, we can then formulate the reconstruction problem \eqref{eq:tv_min_problem_parameter_map}. Conceptually, it might be desirable to exactly solve problem \eqref{eq:tv_min_problem_parameter_map}, i.e., to run the network $\mathcal{N}_{\Theta_T}^{T^{\prime}}$ until convergence by setting $T^{\prime}$ high enough, i.e. conceptually let\ $T^{\prime}\rightarrow \infty$. By the triangle inequality, we have 
\begin{equation}\label{eq:approx_accuracy_test}
\|\XX^{\ast}(\Theta_T) - \XX_{\mathrm{true}}\|_2 \leq  \| \XX^{\ast}(\Theta_T) - \XX_{T^{\prime}}(\Theta_T) \|_2 + \|\XX_{T^{\prime}}(\Theta_T)- \XX_{\mathrm{true}} \|_2,
\end{equation}
where $\| \XX^{\ast}(\Theta_T) - \XX_{T^{\prime}}(\Theta_T) \|_2 \rightarrow 0$ as $T^{\prime}\rightarrow \infty$ due to the convergence of the algorithm that the network $\mathcal{N}_{\Theta_T}^{T^{\prime}}$ implicitly defines. Note however, that the contribution of the second term is not necessarily the smallest for $T^{\prime}=T$ (see also Figure \ref{fig:T_train_vs_test_curves}, especially for lower $T_{\mathrm{test}}$), since $\Theta_T$ was obtained by training only $\Theta$ (see \eqref{eq:loss_fct}) and not $\Theta$ \textit{and} $T^{\prime}$ jointly. This means that, in general, there could exist a configuration which  further improves the results, i.e.\ $\|\XX_{T}(\Theta_T)- \XX_{\mathrm{true}} \|_2 \geq \|\XX_{T^{\prime}}(\Theta_T)- \XX_{\mathrm{true}} \|_2$ for some $T^{\prime}\neq T$.  This is clearly visible from Figure \ref{fig:T_train_vs_test_curves}, which shows the average NRMSE, PSNR and SSIM for different combinations of $T_{\mathrm{train}}$ and $T_{\mathrm{test}}$. Interestingly, it reveals that setting $T_{\mathrm{train}}=T_{\mathrm{test}}$ is indeed not consistently the best choice. Especially for lower $T_{\mathrm{train}}$,  setting $T_{\mathrm{test}}>T_{\mathrm{train}}$  introduces further regularization and yields more accurate reconstructions. For higher $T_{\mathrm{train}}$, in contrast, we see that setting $T_{\mathrm{test}}>T_{\mathrm{train}}$ possibly introduces reconstruction errors which, however, can be entirely attributed to the regularization inherently imposed by the TV, i.e., coming from the term $\|\XX^{\ast}(\Theta_T)- \XX_{\mathrm{true}} \|_2$.\\
This means that, in general, one can view our proposed method in two different flavours. From one point of view, with the proposed method, we can generate a regularization parameter-map $\boldsymbol{\Lambda}_{\Theta_T}$ which we then use to formulate the reconstruction problem \eqref{eq:tv_min_problem_parameter_map} and  which we then aim to subsequently solve \textit{exactly}, i.e., the corresponding algorithm defined by  $\mathcal{N}_{\Theta_T}^{T^{\prime}}$ is run until convergence by letting $T^{\prime} \rightarrow \infty$. Thereby, at test time, we implicitly accept the inherent model-error which is made by choosing the TV-minimization as the underlying regularization method. All the results shown in the paper were generated following this strategy.
From a second, more applied perspective, we can view the proposed approach which yields regularization parameter-maps which are also tailored to the specific number of iterations the network was trained with. Therefore, assuming one was able to use a high-enough $T_{\mathrm{train}}$ for training, at test time, one simply uses $T_{\mathrm{test}}=T_{\mathrm{train}}$ or, if $T_{\mathrm{train}}$ had to be strongly compromised during training (for example, due to limited GPU-memory) one can manually adjust an appropriate $T_{\mathrm{test}}$ on a validation set to compensate for the effects seen in Figure \ref{fig:T_train_vs_test_curves}.

\section{Conclusion}\label{sec:conclusion}

We have presented a simple yet efficient data-driven approach to automatically select data/patient-adaptive spatial/spatio-temporal dependent regularization parameter-maps for the variational regularization approach focusing on TV-minimization. This constitutes  a simple yet efficient and elegant way to combine variational methods with the versatility of deep learning-based approaches, yielding an entirely interpretable reconstruction algorithm which inherits all theoretical properties of the scheme the network implicitly defines. We showed consistency results of the proposed unrolled scheme and we  applied the proposed method to a dynamic MRI reconstruction problem, a quantitative MRI reconstruction problem, a dynamic image denoising problem and a low-dose CT reconstruction problem. In the following, we discuss possible future research directions and we also comment on the limitations of our approach.

We can immediately identify several different components worth further investigations. 
First of all, for a fixed problem formulation and choice of regularization method (i.e.\ the TV-minimization considered in this work) there exist several different reconstruction algorithms, all with their theoretical and practical advantages and limitations, see e.g.\ \cite{wang2008new, goldstein2009split, chambolle2011first, stadler}. It might be interesting to investigate whether our approach yields similar regularization maps regardless of the chosen reconstruction method and if not, to what extent they differ in. 
Second, in this work, we have considered the TV-minimization as the regularization method of choice. However, also TV minimization-based methods are known to have limitations, e.g.\ in producing staircasing effects. We hypothesize that the proposed method could as well be expanded to TGV-based methods \cite{TGV} to overcome these limitations. In addition, the parameter-map learning can be applied when a combination of regularizers is considered. For example, similar to the dynamic MRI and denoising case studies, the proposed method can be used for Hyperspectral X-ray CT, where the spatial and spectral domains are regularized differently, see e.g., \cite{Warr2021,Warr2022}. Further, other regularization methods as for example Wavelet-based methods \cite{Donoho_1994, Chang_2000} could be considered as well, where instead of employing the finite differences operator $\nabla$, a Wavelet-operator $\mathbf{\Psi}$ would be the sparsity-transform of choice. Thereby, the multi-scale decomposition of the U-Net which we have used in our work also naturally fits the problem and could be utilized to estimate different parameter-maps for each different level of the Wavelet-decomposition.
Third, although we have used a plain U-Net \cite{Ronneberger2015} for the estimation of the regularization parameter-maps, there exist nowadays more sophisticated network architectures, e.g.\ transformers \cite{lu2022transformer, liang2021swinir}, which could be potentially adopted as well. Lastly, from the theoretical prospective, future work can include extension of the consistency results to stationary points instead for minimizers only as well as extension to the non-strongly convex fidelity terms in order to cover the CT case as well. It would be also interesting to investigate theoretically in what degree CNN-produced artefacts in the parameter-maps can affect or create artefacts to the corresponding reconstructions.

The main limitations of the proposed approach are the ones which are common for every unrolled NN: the large GPU-memory consumption to store intermediate results and their corresponding gradients while training on the one hand and the possibly long training times which are attributed to the need to repeatedly apply the forward and the adjoint operator during training on the other hand.
As we have seen from Figure \ref{fig:T_training_curves}, to be able to learn the regularization parameter-map with a CNN as proposed, one must be able to use a certain number of iterations $T$ for the unrolled NN to ensure that the output image of the reconstruction network has sufficiently converged to the solution of problem \eqref{eq:tv_min_problem}. How large this number needs to be depends on the considered application as well as the convergence rate of the unrolled algorithm which is used for the reconstruction.

\section*{Acknowledgments}
The authors acknowledge the support of the German Research Foundation (DFG) under Germany's Excellence Strategy – The Berlin Mathematics
Research Center MATH+ (EXC-2046/1, project ID: 390685689) as this work was initiated during the Hackathon event “Maths Meets Image”, Berlin, March 2022, which was part of the MATH+ Thematic Einstein Semester on ``Mathematics of Imaging in Real-World Challenges". This work is further supported  via the MATH+ project EF3-7.

This work was funded by the UK EPSRC grants the ``Computational Collaborative Project in Synergistic Reconstruction for Biomedical Imaging" (CCP SyneRBI) EP/T026693/1; ``A Reconstruction Toolkit for Multichannel CT" EP/P02226X/1 and ``Collaborative Computational Project in tomographic imaging’' (CCPi) EP/M022498/1 and EP/T026677/1. This work made use of computational support by CoSeC, the Computational Science Centre for Research Communities, through CCP SyneRBI and CCPi.

This work is part of the Metrology for Artificial Intelligence for Medicine (M4AIM) project that is funded by the German Federal Ministry for Economic Affairs and Energy (BMWi) in the framework of the QI-Digital initiative.

\bibliographystyle{plain}
\bibliography{refs.bib}

\end{document}